\documentclass[letter,11pt]{book}
\usepackage[T1]{fontenc}
\usepackage{lmodern}
\usepackage{etex} 
\usepackage{hologo}
\usepackage{amsmath}
\usepackage{amsfonts}
\usepackage{amssymb}
\usepackage{mathtools}
\usepackage{graphicx}
\usepackage[misc]{ifsym}
\usepackage[utf8]{inputenc}
\usepackage{multicol}
\usepackage{float}
\usepackage[titletoc]{appendix}
\usepackage{makeidx}
\PassOptionsToPackage{hyphens}{url}\usepackage[unicode,linktocpage=true]{hyperref}
\usepackage[spanish,es-noquoting]{babel} 
\usepackage{tocloft}
\usepackage[explicit]{titlesec}
\usepackage{enumerate}
\usepackage{perpage}
\PassOptionsToPackage{dvipsnames,svgnames,table}{xcolor}
\usepackage{colortbl}
\usepackage{systeme}
\usepackage{pst-plot} 
\usepackage{pstricks}
\usepackage{pgfplots}       
\usepackage{tikz}
\usepackage[american,cuteinductors,smartlabels]{circuitikz}
\usetikzlibrary{matrix,arrows,positioning,decorations.pathmorphing}
\usetikzlibrary{shapes,shadows,calc}
\usepackage{tkz-euclide}
\usepackage{tikz-network}
\usetkzobj{all}
\usepackage{kpfonts}
\usepackage{fancyhdr}
\usepackage{ifthen}\newboolean{color}
\usepackage[breakable,skins]{tcolorbox}
\definecolor{Border}{HTML}{6699FF}
\definecolor{Background}{HTML}{F8F8F8}
\newtcolorbox{maxcode}{
    breakable,
    enhanced,
    colback=Background,
    colframe=Border,
    top=4mm,
    enlarge top by=\baselineskip/2+1mm,
    pad at break=\baselineskip,
    fontupper=\normalsize,}

\pgfplotsset{compat=1.14}

\definecolor{labelcolor}{RGB}{100,0,0}
\definecolor{LightCyan}{rgb}{0.88,1,1}

\newcommand{\degree}{\ensuremath{^\circ}}
\def\R{\mathbb{R}}
\def\rr{0.1}

\urldef\freqletras\url{https://es.wikipedia.org/wiki/Frecuencia_de_aparici%C3%B3n_de_letras#Frecuencia_de_aparici%C3%B3n_de_letras_en_espa%C3%B1ol}

\newsavebox{\picturebox}
\newlength{\pictureboxwidth}
\newlength{\pictureboxheight}
\newcommand{\includeimage}[1]{
    \savebox{\picturebox}{\includegraphics{#1}}
    \settoheight{\pictureboxheight}{\usebox{\picturebox}}
    \settowidth{\pictureboxwidth}{\usebox{\picturebox}}
    \ifthenelse{\lengthtest{\pictureboxwidth > .95\linewidth}}
    {
        \includegraphics[scale=0.6]{#1}
    }
    {
        \ifthenelse{\lengthtest{\pictureboxheight>.80\textheight}}
        {
            \includegraphics[scale=0.6]{#1}
            
        }
        {
            \includegraphics{#1}
        }
    }
}

\setlength{\oddsidemargin}{0.0in}
\setlength{\evensidemargin}{0.0in}
\setlength{\textwidth}{6.5in}
\setlength{\topmargin}{-0.5in}
\setlength{\textheight}{8.5in}
\setlength{\columnsep}{-12pt}

\setlength\cftparskip{-8pt}
\setlength\cftbeforesecskip{10pt}
\cftsetindents{section}{0.185in}{24pt}

\makeatletter
\newcommand\testcolorednode[3]{
           \pgfutil@ifundefined{pgf@sh@ns@#1}{#2}{#3}}
\makeatother

\tikzset{>=stealth'}
\tikzset{
	MyPersp/.style={scale=1.8,x={(-0.8cm,-0.4cm)},y={(0.8cm,-0.4cm)},
    z={(0cm,1cm)}},
	MyPoints/.style={fill=white,draw=black,thick}
		}

\newcommand{\itemequ}[1]{%
        \begingroup%
        \setlength{\abovedisplayskip}{0pt}%
        \setlength{\belowdisplayskip}{0pt}%
        \parbox[c]{\linewidth}{\begin{flalign*}#1&&\end{flalign*}}%
        \endgroup}

\let\origdoublepage\cleardoublepage
\newcommand{\clearemptydoublepage}{%
  \clearpage
  {\pagestyle{empty}\origdoublepage}%
}
\let\cleardoublepage\clearemptydoublepage

\pagestyle{fancy}
\fancyhf{}
\fancyhead[LE,RO]{\thepage}
\fancyhead[LO]{\itshape\nouppercase{\rightmark}}
\fancyhead[RE]{\itshape\nouppercase{\leftmark}}

\MakePerPage{footnote}

\newcommand*\chapterlabel{}
\titleformat{\chapter}
  {\gdef\chapterlabel{}
   \normalfont\sffamily\Huge\bfseries\scshape}
  {\gdef\chapterlabel{\thechapter\ }}{0pt}
  {\begin{tikzpicture}[remember picture,overlay]
    \node[yshift=-3cm] at (current page.north west)
      {\begin{tikzpicture}[remember picture, overlay]
        \draw[fill=LightSlateGray] (0,0) rectangle
          (\paperwidth,3cm);
        \node[anchor=east,xshift=.9\paperwidth,rectangle,
              rounded corners=20pt,inner sep=11pt,
              fill=MidnightBlue]
              {\color{white}\chapterlabel #1};
       \end{tikzpicture}
      };
   \end{tikzpicture}
  }
\titlespacing*{\chapter}{0pt}{20pt}{-60pt}

\newenvironment{dedication}
{
   \cleardoublepage
   \thispagestyle{empty}
   \vspace*{\stretch{1}}
   \hfill\begin{minipage}[t]{0.33\textwidth}
   \raggedright
}
{
   \end{minipage}
   \vspace*{\stretch{3}}
   \clearpage
}

\makeatletter
\renewcommand{\@chapapp}{}

\makeatother

\makeatletter
\let\sv@endpart\@endpart
\def\@endpart{\thispagestyle{empty}\sv@endpart}
\makeatother

\makeindex

\begin{document}

\begin{tikzpicture}[remember picture,overlay]
    \node[yshift=-8cm] at (current page.north west)
      {\begin{tikzpicture}[remember picture, overlay]
        \draw[fill=LightSlateGray] (0,0) rectangle
          (\paperwidth,10cm);
        \node[anchor=east,xshift=.9\paperwidth,rectangle,
              rounded corners=20pt,inner sep=11pt,
              fill=MidnightBlue]
              {\color{white} \normalfont\sffamily\Huge\bfseries\scshape 
              Matem\'aticas B\'asicas (con Software Libre)};
       \end{tikzpicture}
      };
   \end{tikzpicture}

\vspace{6cm}   
\begin{flushright}
{\Large\sffamily\bfseries\scshape Jos\'e A Vallejo}\\
{\Large\sffamily\bfseries\scshape Antonio Morante}\\[8pt]
{\normalfont\sffamily\large\bfseries\scshape Facultad de Ciencias\\
Universidad Aut\'onoma de San Luis Potos\'i}\\[8pt]
Versi\'on del \today
\end{flushright}
 
  

\newpage\thispagestyle{empty}

\vspace{3cm}
\noindent Jos\'e Antonio Vallejo Rodr\'iguez\\
\Letter\, \href{mailto:jvallejo@fc.uaslp.mx}{\tt jvallejo@fc.uaslp.mx} \\
URL: \url{http://galia.fc.uaslp.mx/~jvallejo}\\

\noindent Antonio Morante Lezama\\
\Letter\, \href{mailto:amorante@fciencias.uaslp.mx}{\tt amorante@fciencias.uaslp.mx}\\

\noindent Facultad de Ciencias\\
Universidad Aut\'onoma de San Luis Potos\'i\\
Lateral Av Salvador Nava s/n Colonia Lomas\\
CP 78290 San Luis Potos\'i (SLP)\\
M\'exico

\newpage

\begin{dedication}
A mis hijas, Alexandra y Gabriela, por si alg\'un d\'ia pudiera serles \'util (JAV)
\end{dedication}
\newpage\null\thispagestyle{empty}\newpage

\setcounter{tocdepth}{1}
\tableofcontents
\thispagestyle{empty}
\addtocontents{toc}{\protect\thispagestyle{empty}}
\mainmatter

\addcontentsline{toc}{chapter}{Prefacio}
\chapter*{Prefacio}
\markboth{Prefacio}{} 
\vspace{2.5cm}

Estas notas surgieron de la necesidad de disponer de unos materiales adecuados a los
alumnos que ingresan a las carreras impartidas en la Facultad de Ciencias de la
Universidad Aut\'onoma de San Luis Potos\'i (aunque seguramente esta
necesidad sea la misma en cualquier otra universidad), alumnos que por lo general
tienen una formaci\'on deficiente tanto en la parte cient\'ifica como en la human\'istica.
El desconocimiento de las t\'ecnicas para trabajar con la parte m\'as elemental de las matem\'aticas
les impide cursar con \'exito las asignaturas especializadas como \'algebra o c\'alculo, donde
esos conocimientos se dan por supuestos, pero es a\'un peor la falta de comprensi\'on lectora.
Incluso aquellos estudiantes que podr\'ian resolver la parte de c\'alculo de un problema
(digamos, despejar una ecuaci\'on y hallar el valor de una inc\'ognita), fallan porque no saben traducir el enunciado en la ecuaci\'on que saben resolver mec\'anicamente.
Por eso, los ejercicios que se presentan son de dos tipos: unos puramente t\'ecnicos 
y otros basados en un enunciado del que hay que extraer las ecuaciones a resolver. Es important\'isimo que el alumno practique los ejercicios del segundo tipo (hay m\'as de
$165$ ejercicios propuestos, de variada dificultad,
de modo que material sobra para ello). Una parte de estos ejercicios se debe a \'Alvaro
P\'erez Raposo, de la Universidad Polit\'ecnica de Madrid.

Todos los cap\'itulos tienen la misma estructura: se comienza con una secci\'on de
nociones b\'asicas, en la que se introducen unos conceptos matem\'aticos con el m\'inimo 
te\'orico indispensable para poder comprenderlos y utilizarlos.
Tratando de que el alumno aprecie desde un primer momento la importancia que tienen las
matem\'aticas en la vida cotidiana, se han inclu\'ido unas secciones de aplicaciones,
donde se ha buscado presentar situaciones que puedan resultar atractivas
y en las cuales se haga uso de las ideas expuestas, ofreciendo un
panorama de las matem\'aticas en acci\'on (en situaciones simplificadas, claro est\'a).
Los t\'opicos tratados van desde la estructura del c\'odigo ISBN a los fundamentos de
la teor\'ia de la informaci\'on y la entrop\'ia, pasando por el algoritmo criptogr\'afico
RSA, la navegaci\'on a\'erea por deriva, los c\'odigos correctores de
Reed--Solomon o el procesamiento b\'asico de im\'agenes. 
En realidad, esperamos que despu\'es de trabajar este texto el estudiante
pueda dar cumplida respuesta a preguntas como `?`para qu\'e se usa la trigonometr\'ia
en la vida real?, ?`por qu\'e hay que estudiar los polinomios?' y otras semejantes. 
Una  tercera secci\'on de ejercicios sirve para comprobar si se han asimilado
los conceptos y t\'ecnicas adecuadamente. 
Por \'ultimo, cada cap\'itulo acaba con otra secci\'on en la que se utiliza  
tecnolog\'ia inform\'atica. Se sugiere que las actividades de esta parte se
realicen en una sala de c\'omputo, por parejas o grupos a lo sumo de tres 
integrantes, incluso si los alumnos disponen de su propia laptop.
El uso de estas tecnolog\'ias requiere de alguna aclaraci\'on. Frecuentemente,
se emplean \'unicamente para mejorar las representaciones gr\'aficas, sustituyendo los
antiguos dibujos sobre el pizarr\'on con una gr\'afica de Excel\textregistered\, o la escritura
sobre ese pizarr\'on por una presentaci\'on en PowerPoint\textregistered\,. Decir que esto es
una `innovaci\'on tecnol\'ogica' es tanto como decir que el cambio del gis (la tiza) al
marcador (el rotulador) tambi\'en fue un paso innovador. Esto es simplemente una actualizaci\'on
de los medios, pero no supone ning\'un cambio en la actitud docente, ni en la orientaci\'on
de los contenidos. Lo que se pretende en este texto es entrelazar la ense\~nanza de
conceptos te\'oricos con su uso pr\'actico a trav\'es de un equipo de c\'omputo. Esto
supone que el alumno debe conocer el fundamento te\'cnico de un determinado c\'alculo
y tambi\'en debe saber c\'omo llevarlo a cabo eficientemente. Es en este \'ultimo paso
donde intervienen de manera fundamental las computadoras, ante la imposibilidad de
realizar c\'alculos elementales pero muy largos en un tiempo razonable. Todo esto se ha
discutido en innumerables foros y no es \'este el lugar adecuado para debatir sobre
el tema, pero una buena muestra del contexto real en el que se dan estas situaciones
se presenta en los cap\'itulos \ref{cap2}, \ref{cap3} y \ref{finalcap}. 
Aqu\'i haremos uso de dos programas de software
libre: el sistema de
\'algebra computacional Maxima y el de geometr\'ia din\'amica GeoGebra (descritos en el ap\'endice). Ambos son sobradamente conocidos en el entorno acad\'emico por su potencia y versatilidad, sustituyendo exitosamente
a programas comerciales de elevado coste como Mathematica\textregistered\ o Geometer 
Sketchpad\textregistered{}. La curva
de aprendizaje en ambos casos es m\'inima y los resultados que se obtienen justifican
en cualquier caso la inversi\'on de tiempo que se haga en ellos. Como ventaja adicional,
ambos programas se pueden instalar en cualquier dispositivo con alguno de los sistemas operativos Windows\textregistered{},
Android\textregistered\ o Apple OS\textregistered{}, lo cual los convierte
en candidatos id\'oneos para usarse en un entorno de aprendizaje m\'ovil (en ingl\'es,
\emph{mobile learning} o \emph{m--learning}).

Con el fin de reforzar la conexi\'on de las matem\'aticas con el mundo que nos rodea, se ha
hecho un uso intensivo de recursos de internet. Cuando ha sido posible, se han utilizado datos
reales extra\'idos de organismos oficiales, como la ONU, NASA, etc. Se proporciona la
URL (\textbf{U}niform \textbf{R}esource \textbf{L}ocator, la direcci\'on de internet) de todas las 
p\'aginas a las que se hace
referencia, y se ha comprobado que todas estuvieran vigentes a la fecha de edici\'on, \today.
Tambi\'en, con la finalidad de hacer m\'as agradable el texto, se han inclu\'ido numerosos
elementos gr\'aficos, cuya procedencia se cita expresamente (todos son o bien de dominio
p\'ublico o bien disponen de una licencia de uso que permite su inclusi\'on en esta obra).
Para que el libro fuera efectivamente utilizable en un curso, los autores nos impusimos un 
l\'imite de unas $200$ p\'aginas de extensi\'on, de manera que los contenidos se puedan cubrir
sin dificultad en un semestre.
En cuanto al nivel acad\'emico, la idea ha sido partir de unos conocimientos totalmente 
elementales (a nivel incluso de
secundaria o preparatoria) y dotar al estudiante con las herramientas que le permitan superar una asignatura de c\'alculo diferencial o de \'algebra lineal. Por eso la dificultad de los ejercicios y los temas tratados, as\'i como de las t\'ecnicas de c\'omputo,
va aumentando en sucesivos cap\'itulos. El curso se ha dise\~nado
en la l\'inea de los `precalculus' de algunos centros en USA o la conocida serie de textos
`cap\'itulo 0'.\\

En todo el libro se usa la notaci\'on del \emph{punto decimal}, como en $1{.}4142.$\\

Finalmente, cabe destacar que el dise\~no (que usa ejemplos de Stefan Kottwitz,
\url{http://texblog.net/}) y maquetaci\'on se han realizado enteramente con 
software libre. En particular, los gr\'aficos usan la librer\'ias TikZ y el formato se ha
hecho con \hologo{TeX} Live/\hologo{XeLaTeX}, sobre una plataforma Linux/Slackware.
Todo comentario, sugerencia, aviso de erratas, etc., ser\'a bienvenido a
cualquiera de las direcciones de correo electr\'onico 
\href{mailto:jvallejo@fc.uaslp.mx}{\tt jvallejo@fc.uaslp.mx},
\href{mailto:amorante@fciencias.uaslp.mx}{\tt amorante@fciencias.uaslp.mx}.
\vspace{1cm}
\begin{flushright}
San Luis Potos\'i (M\'exico) a \today
\end{flushright}

\chapter{Aritm\'etica elemental}\label{cap1}

\vspace{2.5cm}
\section{Nociones b\'asicas}

Cada una de las operaciones elementales\index{Operaciones!elementales} con n\'umeros tiene una operaci\'on inversa\index{Operaciones!inversa}. Por ejemplo,
la inversa de la suma es la resta, lo cual quiere decir que para cada n\'umero $a$ existe otro,
denotado por $-a$, tal que la suma de ambos es cero: $a+(-a)=0$. Normalmente, en lugar de escribir esta f\'ormula
se usa la notaci\'on $a-a=0$ y en ese sentido se dice que $-$ es la inversa de $+$. An\'alogamente,
para cada n\'umero no nulo $a$, existe otro, denotado por $1/a$, tal que su producto es la unidad $1$:
$$
a\cdot\frac{1}{a}=1 .
$$
Si se tiene una ecuaci\'on con una inc\'ognita (una cantidad desconocida) $x$, se puede saber su
valor `despejando'\index{Operaciones!despejar}, es decir, aplicando operaciones inversas a todos los factores que aparezcan junto
a $x$ hasta que quede aislada. Por ejemplo, en una ecuaci\'on lineal\index{Ecuaci\'on!lineal}
$$
ax+b=c,
$$
donde $a\neq 0$, primero `despejamos la $b$' sumando su inverso $-b$ 
\emph{a ambos lados de la ecuaci\'on} (pues no podemos
modificarla, as\'i que lo que hagamos en un miembro hay que hacerlo en el otro) y aplicando que $b-b=0$:
$$
ax=ax+b-b=c-b .
$$
Ahora, en la ecuaci\'on que result\'o, $ax=c-b$, `despejamos la $a$' multiplicando por su inverso:
$$
x=\frac{1}{a}ax=\frac{1}{a}(c-b).
$$
De esta manera resulta la soluci\'on final
$$
x=\frac{c-b}{a}.
$$
As\'i, si tenemos $3x-2=-1$, la soluci\'on es
$$
x=\frac{-1+2}{3}=\frac{1}{3}.
$$
N\'otese que aqu\'i $b=-2$, as\'i que su inverso respecto de la suma es $-b=-(-2)=2$.\\
Otras operaciones tambi\'en tienen su inversa, como la ra\'iz cuadrada\index{Rai@Ra\'iz!cuadrada}, que es la inversa de la
potencia al cuadrado. La diferencia aqu\'i es que al aplicar la operaci\'on inversa pueden
aparecer \emph{varias} soluciones. Consideremos por ejemplo la ecuaci\'on $x^2=4$. Aplicando
la ra\'iz a ambos miembros obtenemos
\begin{equation}\label{eq0}
\sqrt{x^2}=\sqrt{4}.
\end{equation}
Pero hay dos n\'umeros que elevados al cuadrado dan $4$: son $2$ y $-2$. Cualquiera de ellos es
un resultado v\'alido de $\sqrt{4}$, por eso decimos que \eqref{eq0} tiene dos soluciones
$$
x=\pm 2.
$$
Otra forma de expresar esto consiste en escribir
$$
\sqrt{x^2}=\pm x,
$$
aunque hay que tener mucho cuidado con expresiones de este tipo, ya que $\sqrt{x^2}$ denota un n\'umero
y no puede tener dos valores distintos a la vez. Como veremos m\'as adelante, lo correcto es escribir
$\sqrt{x^2}=|x|$, donde $|\cdot |$ denota el \emph{valor absoluto} de $x$ (v\'eanse el ejercicio \ref{ejer} y el cap\'itulo \ref{explog}).\index{Valor absoluto}

El siguiente paso es considerar una ecuaci\'on cuadr\'atica\index{Ecuaci\'on!cuadr\'atica},
como $ax^2+bx+c=0$. Si el miembro
izquierdo fuera un cuadrado perfecto\index{Cuadrado perfecto}, es decir, si se pudiera escribir como $(ax-h)^2$, la ecuaci\'on
ser\'ia $(ax-h)^2=0$, que tiene una soluci\'on muy f\'acil: tomando ra\'ices cuadradas en ambos
miembros, $ax-h=0$ (aqu\'i no salen signos $\pm$ porque $0$ s\'olo tiene una posible ra\'iz que es el
propio $0$) y despejando, ser\'ia $x=h/a$. Lo normal, sin embargo, es que nos encontremos con una
ecuaci\'on cuadr\'atica que \emph{no} es un cuadrado perfecto.\\
La forma de resolver el caso general pasa por el m\'etodo de \emph{completar cuadrados}.
\index{Me@M\'etodo!completar cuadrados}\'Este consiste en
escribir la expresi\'on $ax^2+bx+c$ en la forma $a(x-h)^2-k$, que no es un cuadrado perfecto pero que
conduce a una ecuaci\'on que se puede resolver, como veremos a continuaci\'on.\\
Veamos entonces c\'omo determinar $h$ y $k$. Si comparamos $ax^2+bx+c$ con 
$a(x-h)^2-k$, para que sean iguales
debe ser, desarrollando el cuadrado, $ax^2+bx+c = ax^2-2axh+ah^2-k$,
o sea, simplificando el t\'ermino $ax^2$,
$$
bx+c = -2axh+ah^2-k\, .
$$
Igualando los coeficientes de $x$ y el t\'ermino independiente, resulta que
\begin{align*}
b &=-2ah \\
c &=ah^2-k.
\end{align*}
De la primera de estas ecuaciones, resulta
\begin{equation}\label{eq00}
h=-\frac{b}{2a}\,,
\end{equation}
que al sustituir en la segunda nos deja (n\'otese que al elevar al cuadrado desaparece el signo, pues
$(-1)^2=1$):
$$
c=a\frac{b^2}{4a^2}-k=\frac{b^2}{4a}-k\,,
$$
de donde podemos `despejar la $k$' f\'acilmente:
\begin{equation}\label{eq000}
k=\frac{b^2}{4a} -c\,.
\end{equation}
Las ecuaciones \eqref{eq00} y \eqref{eq000} nos dicen que $ax^2+bx+c=0$ es lo mismo que
$$
a\left( x+ \frac{b}{2a} \right)^2+c-\frac{b^2}{4a}=0.
$$
Los pasos para `despejar la $x$' son los siguientes:
\begin{align*}
a\left( x+ \frac{b}{2a} \right)^2+c-\frac{b^2}{4a}=0, \\
a\left( x+ \frac{b}{2a} \right)^2=\frac{b^2}{4a}-c\,, \\
\left( x+ \frac{b}{2a} \right)^2=\frac{b^2}{4a^2}-\frac{c}{a}\,, \\
\left( x+ \frac{b}{2a} \right)^2=\frac{b^2-4ac}{4a^2}\,, \\
 x+ \frac{b}{2a} =\pm \sqrt{\frac{b^2-4ac}{4a^2}} \,.
\end{align*}
Recordando por otra parte que
\begin{align*}
\sqrt{\frac{x}{y}}&=\frac{\sqrt{x}}{\sqrt{y}}\\
\sqrt{xy}&=\sqrt{x}\sqrt{y}\,,
\end{align*}
lo anterior se escribe como
\begin{align*}
x+ \frac{b}{2a} &=\pm\frac{\sqrt{b^2-4ac}}{2a}\,,
\end{align*}
y finalmente
\begin{equation}\label{eqcuadratica}
x=\frac{-b\pm\sqrt{b^2-4ac}}{2a} \,.
\end{equation}
Esta f\'ormula es muy popular, pero no hace falta memorizarla. Lo importante es tener claro cu\'al
es el proceso de completar cuadrados y realizarlo cuando sea necesario paso a paso. La otra 
caracter\'istica importante de esta f\'ormula es que s\'olo proporciona soluciones \emph{reales} cuando
el radicando es positivo, es decir, cuando $b^2-4ac>0$. A este n\'umero se le denomina el
\emph{discriminante}\index{Discriminante} de la ecuaci\'on. Cuando es igual a $0$, s\'olo hay una soluci\'on,
$$
x=-\frac{b}{2a},
$$
mientras que cuando es negativo, las soluciones son n\'umeros complejos, como veremos
m\'as adelante.

A la hora de trabajar con varias operaciones en una misma expresi\'on, es de fundamental importancia
conocer el \emph{orden} en que se aplican, porque el resultado puede variar. Por ejemplo, una expresi\'on
como $644\div 2+10$ puede interpretarse de dos formas, con resultados bien distintos, seg\'un el orden
de las operaciones:
$$
\frac{644}{2+10}=53{.}67,
$$
o bien
$$
\frac{644}{2}+10=332.
$$
En general, si no se especifica lo contrario, cuando nos den una expresi\'on involucrando
varias expresiones se sobreentender\'a que se sigue el siguiente \emph{orden de precedencia}:
\index{Operaciones!precedencia}
\begin{enumerate}
\item Exponentes y ra\'ices.
\item Multiplicaciones y divisiones.
\item Sumas y restas.
\end{enumerate}
Por ejemplo, en la expresi\'on $135\div 15-6+5$ debe calcularse primero la divisi\'on y luego
la suma y la resta, es decir,
$$
\left( \frac{135}{15}-6\right)+5=(9-6)+5=3+5=8\,.
$$
N\'otese c\'omo se ha hecho uso de los par\'entesis ( y ) para indicar el orden de las operaciones.
En una expresi\'on que contenga par\'entesis \emph{anidados} como $(((3-5)^2 4+10)/9+3)2$, se calculan los
t\'erminos entre par\'entesis `de dentro hacia fuera':
$$
(((3-5)^2 4+10)/9+3)2=(((-2)^2 4+10)/9+3)2=((16 +10)/9+3)2=(26/9+27/9)2=(53/9)2=106/9.
$$

El conjunto de los \emph{n\'umeros naturales}\index{Num@N\'umero!natural} se denotar\'a en este libro
por $\mathbb{N}=\{1,2,3,\ldots \}$\index{N@$\mathbb{N}$}. Los \emph{n\'umeros 
enteros}\index{Num@N\'umero!entero}, consistentes en los naturales junto con sus inversos respecto
de la suma y el cero, se denotar\'an por 
$\mathbb{Z}=\{\ldots,-2,-2,0,1,2,\ldots \}$\index{N@$\mathbb{N}$}. Los
\emph{n\'umeros racionales}\index{Num@N\'umero!racional} son los cocientes de n\'umeros enteros
con denominador no nulo (teniendo en cuenta que varias fracciones pueden representar un mismo
n\'umero racional, como $2/3$ y $4/6$), denotados por
$\mathbb{Q}=\{m/n\,:m\in\mathbb{Z},n\in\mathbb{Z}-\{0\}\}$\index{Q@$\mathbb{Q}$}. Finalmente, los
\emph{n\'umeros reales}\index{Num@N\'umero!real} consisten en la uni\'on de los n\'umeros racionales
y los irracionales (como por ejemplo $\sqrt{2}$), se denotan por $\mathbb{R}$\index{R@$\mathbb{R}$}
y en este libro supondremos que son conocidos en sus propiedades m\'as elementales por el lector.

\section{Aplicaciones}

Completando el cuadrado\index{Me@M\'etodo!completar cuadrados}, la ecuaci\'on $x^2-x-1=0$ 
se puede expresar como
$$
0=\left( x-\frac{1}{2} \right)^2-\frac{1}{4}-1=\left( x-\frac{1}{2}\right)^2-\frac{5}{4}\,.
$$
De aqu\'i resulta que sus soluciones son
$$
\frac{1\pm \sqrt{5}}{2}\,.
$$
La soluci\'on con el signo positivo se conoce como el \emph{n\'umero
de oro}\index{Num@N\'umero!de oro} y se denota $\phi$, es decir,
$$
\phi =\frac{1+ \sqrt{5}}{2}\,.
$$
Como curiosidades matem\'aticas que
cumple  $\phi$ tenemos las siguientes:
\begin{enumerate}[(a)]
\item A partir de la ecuaci\'on que define al n\'umero $\phi$, $\phi^2-\phi-1=0$, 
se ve que $\phi=\sqrt{1+\phi}$, de modo que
sustituyendo repetidamente,
$$
\phi =\sqrt{1+\sqrt{1+\sqrt{1+\sqrt{1+\cdots}}}}\,.
$$
\item De nuevo a partir de la ecuaci\'on que cumple $\phi$, se tiene
$$
\phi=1+\frac{1}{\phi} \,.
$$
Entonces, sustituyendo repetidamente, se obtiene la llamada expansi\'on en \emph{fracciones continuas}
de $\phi$\index{Fracciones continuas}:
$$
\phi=\dfrac{1}{1+\dfrac{1}{1+\dfrac{1}{1+\cdots}}}
$$
\end{enumerate}
Fuera de las Matem\'aticas, el n\'umero de oro (o \emph{proporci\'on \'aurea}) ha sido ampliamente utilizado
desde la antig\"uedad para medir cuantitativamente la belleza de formas pict\'oricas o arquitect\'onicas.
\begin{center}
\includegraphics[scale=1]{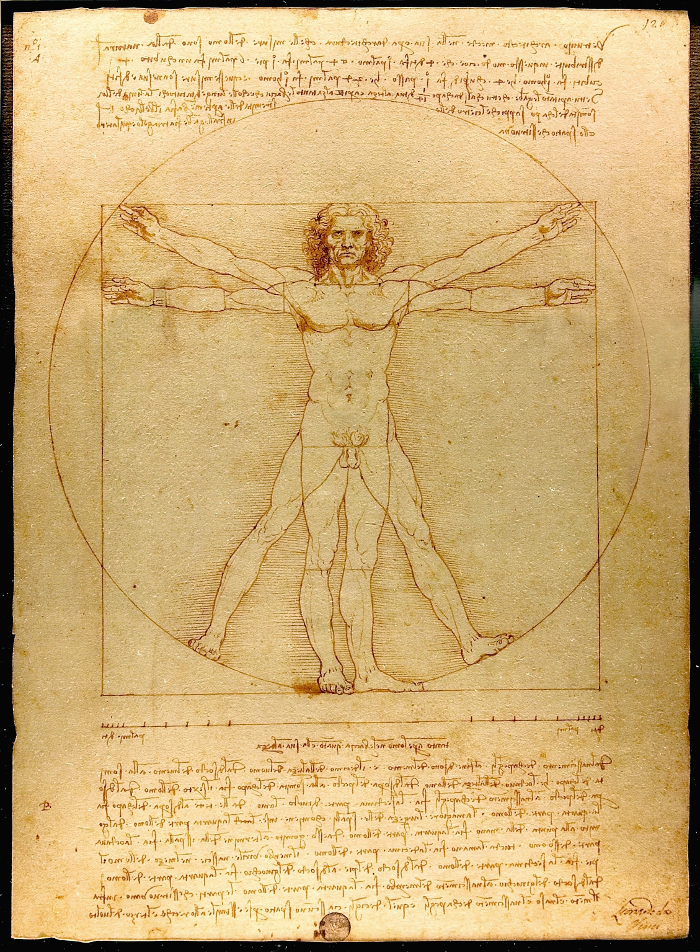} 
\end{center}
El conocido boceto de Leonardo Da Vinci\index{Da Vinci, Leonardo}, llamado \emph{El hombre de Vitrubio}, muestra las principales razones\index{Hombre de Vitrubio}
geom\'etricas que se dan entre las distintas partes del cuerpo humano.
Leonardo, como otros artistas del Renacimiento, utilizaba la proporci\'on \'aurea en sus dise\~nos. Por ejemplo,
el canon cl\'asico de belleza considera que la altura del suelo al ombligo debe ser $\phi$ veces la altura del
ombligo a la parte superior de la cabeza. En arquitectura, el Parten\'on de Atenas es una muestra de edificio
constru\'ido tomando como base $\phi$, v\'ease por ejemplo 
\url{http://www.goldennumber.net/parthenon-phi-golden-ratio/}.

Recientemente, un cirujano pl\'astico llamado Stephen Marquardt\index{Marquardt, Stephen} ha realizado unos estudios generalizando las
proporciones cl\'asicas de manera exhaustiva. As\'i, estableci\'o c\'omo (siempre seg\'un \'el) deb\'ian ser
los cocientes de la distancia entre ojos y la distancia entre orejas, la distancia de la barbilla al labio
superior y la distancia del labio a la nariz, etc. El resultado es la llamada \emph{m\'ascara de Marquardt}\index{Marquardt, Stephen!m\'ascara},
un dise\~no matem\'atico que superpuesto a una fotograf\'ia muestra qu\'e tan bien se adaptan los rasgos de
una persona al ideal de `belleza matem\'atica'. Puede ampliarse esta idea en el sitio web
Interactive Mathematics,
\url{http://www.intmath.com/numbers/math-of-beauty.php}.

Naturalmente, esto no debe tomarse como un criterio absoluto
(!`la belleza est\'a en el interior!) y, adem\'as, no est\'a exento de cr\'iticas como puede verse en
el trabajo de Erik Holland,
publicado en la revista Aesthetic Plastic Surgery y disponible en
\url{http://link.springer.com/article/10.1007/s00266-007-9080-z}. Simplemente, hay
que considerarlo como una buena explicaci\'on de lo que habitualmente se considera armonioso. 
Como un ejemplo, la siguiente figura muestra a la conocida actriz italiana
Sof\'ia Loren\index{Loren, Sof\'ia} (Sofia Villani Scicolone) con la m\'ascara de Marquardt superpuesta. 
Puede apreciarse con toda claridad la correspondencia de los rasgos con las gu\'ias 
determinadas matem\'aticamente.
\begin{center}
\includegraphics[scale=0.25]{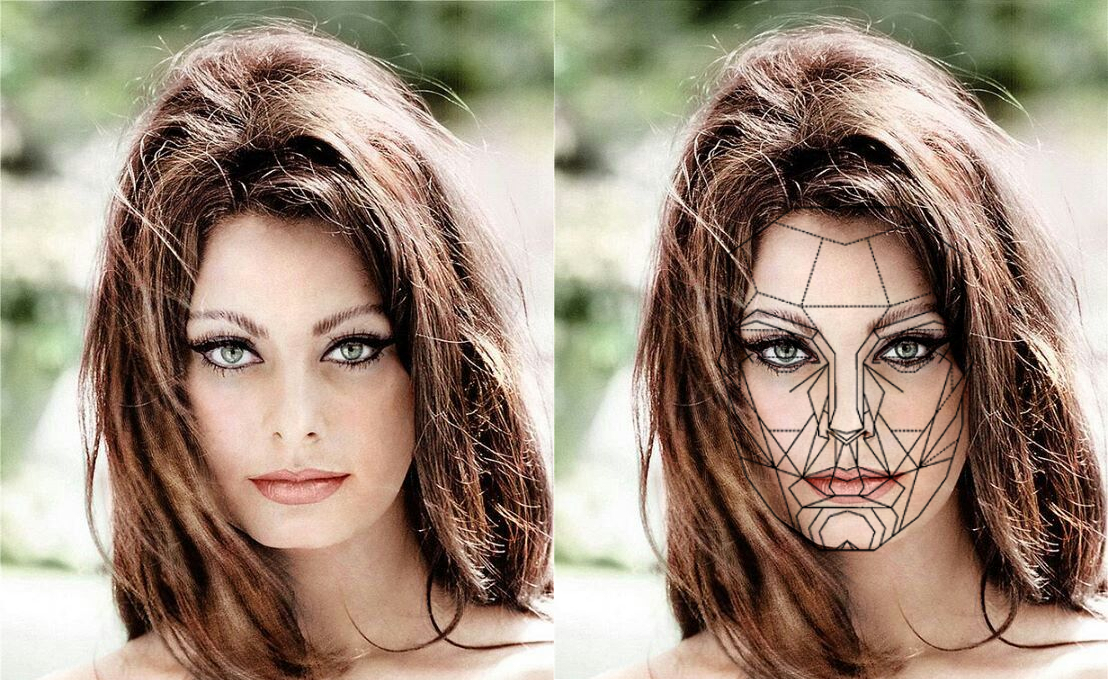} 
\end{center}
\section{Ejercicios}

\begin{enumerate}
\item\label{ejer1} Resolver las siguientes ecuaciones, en las que la inc\'ognita es $x$:
\begin{multicols}{3}
\begin{enumerate}[(a)]
\item $2x+5=-1$
\item $ax+2=0$
\item $x-\sqrt{b}=c^2$
\item $3-k=x+d$
\item $x/a-5=1$
\item $2/x+1=-3$
\item $a/x+b=1$
\item $4\sqrt{x}-1=15$
\item $\sqrt{x}+2=12$
\item $3/\sqrt{x}+b=-1$
\end{enumerate}
\end{multicols}

\item Completando el cuadrado, resolver la ecuaci\'on $x^2-5x+4=0$. Comparar el resultado con el
que proporciona la soluci\'on general \eqref{eqcuadratica}.

\item\label{parab} Una funci\'on se dice que es \emph{cuadr\'atica}\index{Funci\'on!cuadr\'atica}
si es de la forma $f(x)=ax^2+bx+c$. La gr\'afica
de una funci\'on cuadr\'atica es una curva denominada \emph{par\'abola}\index{Par\'abola}, 
cuya forma general se ve en la figura.
Reescribir la expresi\'on de la par\'abola dada por la funci\'on $f(x)=x^2-2x-3$ completando
cuadrados y a partir de la expresi\'on obtenida determinar el v\'ertice de la par\'abola y
los puntos en los que corta a los ejes coordenados.
\begin{center}
\begin{tikzpicture}[scale=0.775]
  \draw[-] (-3,0) -- (3,0) node[right] {$x$};
  \draw[-] (0,-2.5) -- (0,2.5) node[above] {$f(x)$};
  \draw[color=red] (-2,2) parabola[bend pos=0.5] bend +(0,-4) +(3,0);
\end{tikzpicture}
\end{center}

\item\label{ejer4} Repetir el ejercicio anterior para las funciones cuadr\'aticas siguientes:
\begin{multicols}{2}
\begin{enumerate}[(a)]
\item $f(x)=-x^2+2x$.
\item $f(x)=2x^2-x-1$.
\item $f(x)=3x^2-6x+7$.
\item $f(x)=\frac{1}{2}x^2-6x+1$.
\end{enumerate}
\end{multicols}

\item\label{ejer5} Realizar las siguientes operaciones.
\begin{multicols}{2}
\begin{enumerate}[(a)]
\item $(285+117)/(264-68)$
\item $\frac{2}{3}-\frac{4}{5}$
\item $(\frac{1}{2}+3)(-2-\frac{3}{2})$
\item $41-(32-21)$
\item $117+(39\div 3)$
\item $(121\div 11)-10$
\item $19+(348\div 16)$
\item $(200-135)\div 45$
\end{enumerate}
\end{multicols}

\item Colocar par\'entesis en las expresiones dadas a continuaci\'on:
\begin{multicols}{2}
\begin{enumerate}[(a)]
\item $135\div 15-6+5$ para dar $20,-2$.
\item $120\cdot 2+10-3$ para dar $247,1\, 080$.
\item $150+6\div 2$ para dar $78,153$.
\item $30-12+3\cdot 2$ para dar $24,42$.
\item $8\cdot 4+6\div 2$ para dar $35,40$.
\item $20+4-3-7$ para dar $28,14$.
\end{enumerate}
\end{multicols}

\item\label{ejer7} Realizar las siguientes operaciones.
\begin{multicols}{2}
\begin{enumerate}[(a)]
\item $10\div (10\div (10\div (10 \div (10\div 10))))$
\item $10\cdot (10\div (10\cdot \sqrt{ 10\div [10\cdot (10\div 10)]}))$
\item $3\cdot (5\div [7\cdot \{92 +5\}\cdot 4])$
\item $(9\div (7\cdot (3+7)^2\cdot 8))\cdot 5$
\item $5\cdot 4\cdot 7 \cdot (9+5+3)\cdot 5\cdot (7\div (4\cdot 5))$
\item $3\cdot ((14-9)^2\div (9-4))$
\end{enumerate}
\end{multicols}

\item Comprobar la presencia de la proporci\'on \'aurea en la vida cotidiana, midiendo:
\begin{enumerate}[(a)]
\item El cociente entre el ancho y el alto de una credencial universitaria o una tarjeta
de d\'ebito.
\item El cociente entre el ancho y el alto de un monitor de 19" (los monitores m\'as modernos,
de mayor tama\~no, utilizan otra proporci\'on debido a la necesidad de reproducir v\'ideo en
alta definici\'on. Es tambi\'en un buen ejercicio determinar esta proporci\'on y comprobar si coincide con lo anunciado por el vendedor).
\item Alguno de los art\'iculos que aparecen en la p\'agina
\href{http://www.phimatrix.com/product-design-golden-ratio/}{\tt www.phimatrix.com}.
\end{enumerate}

\end{enumerate}

\section{Uso de la tecnolog\'ia de c\'omputo}

Maxima funciona como una calculadora cient\'ifica de precisi\'on arbitraria.
Por defecto trunca los n\'umeros a 15 decimales (el valor por defecto de la
variable de entorno \texttt{fpprec} es $16$), pero este comportamiento puede
cambiarse:

\begin{maxcode}
\noindent
\begin{minipage}[t]{8ex}{\color{red}\bf
\begin{verbatim}
(%i1) 
\end{verbatim}
}
\end{minipage}
\begin{minipage}[t]{\textwidth}{\color{blue}%
\begin{verbatim}
sqrt(2),bfloat;
\end{verbatim}
}
\end{minipage}

\begin{math}\displaystyle
\parbox{8ex}{\color{labelcolor}(\%o1) }
1{.}414213562373095b0
\end{math}

\noindent
\begin{minipage}[t]{8ex}{\color{red}\bf
\begin{verbatim}
(%i2) 
\end{verbatim}
}
\end{minipage}
\begin{minipage}[t]{\textwidth}{\color{blue}
\begin{verbatim}
fpprec;
\end{verbatim}
}
\end{minipage}
\begin{math}\displaystyle
\parbox{8ex}{\color{labelcolor}(\%o2) }
16
\end{math}

\noindent
\begin{minipage}[t]{8ex}{\color{red}\bf
\begin{verbatim}
(%i3) 
\end{verbatim}
}
\end{minipage}
\begin{minipage}[t]{\textwidth}{\color{blue}
\begin{verbatim}
fpprec:50;
\end{verbatim}
}
\end{minipage}
\begin{math}\displaystyle
\parbox{8ex}{\color{labelcolor}(\%o3) }
50
\end{math}

\noindent
\begin{minipage}[t]{8ex}{\color{red}\bf
\begin{verbatim}
(%i4) 
\end{verbatim}
}
\end{minipage}
\begin{minipage}[t]{\textwidth}{\color{blue}
\begin{verbatim}
sqrt(2),bfloat;
\end{verbatim}
}
\end{minipage}
\begin{math}\displaystyle
\parbox{8ex}{\color{labelcolor}(\%o4) }
1{.}414213562373095048801688724209698078569671875377b0
\end{math}

\noindent
\begin{minipage}[t]{8ex}{\color{red}\bf
\begin{verbatim}
(%i5) 
\end{verbatim}
}
\end{minipage}
\begin{minipage}[t]{\textwidth}{\color{blue}
\begin{verbatim}
fpprec:16;
\end{verbatim}
}
\end{minipage}
\begin{math}\displaystyle
\parbox{8ex}{\color{labelcolor}(\%o5) }
16
\end{math}
\end{maxcode}

Observemos c\'omo se ha asignado el valor de una variable (\texttt{fpprec}),
mediante el uso de dos puntos `:'. Esto sirve en general para hacer asignaciones
de valores a variables\index{Maxima!variables}, las operaciones\index{Operaciones!en Maxima} 
aritm\'eticas habituales se escriben como es habitual
\begin{maxcode}
\noindent
\begin{minipage}[t]{8ex}{\color{red}\bf
\begin{verbatim}
(%i6) 
\end{verbatim}
}
\end{minipage}
\begin{minipage}[t]{\textwidth}{\color{blue}
\begin{verbatim}
x:4;
\end{verbatim}
}
\end{minipage}
\begin{math}\displaystyle
\parbox{8ex}{\color{labelcolor}(\%o6) }
4
\end{math}

\noindent
\begin{minipage}[t]{8ex}{\color{red}\bf
\begin{verbatim}
(%i7) 
\end{verbatim}
}
\end{minipage}
\begin{minipage}[t]{\textwidth}{\color{blue}
\begin{verbatim}
y:3;
\end{verbatim}
}
\end{minipage}
\begin{math}\displaystyle
\parbox{8ex}{\color{labelcolor}(\%o7) }
3
\end{math}

\noindent
\begin{minipage}[t]{8ex}{\color{red}\bf
\begin{verbatim}
(%i8) 
\end{verbatim}
}
\end{minipage}
\begin{minipage}[t]{\textwidth}{\color{blue}
\begin{verbatim}
x*y;
\end{verbatim}
}
\end{minipage}
\begin{math}\displaystyle
\parbox{8ex}{\color{labelcolor}(\%o8) }
12
\end{math}

\noindent
\begin{minipage}[t]{8ex}{\color{red}\bf
\begin{verbatim}
(%i9) 
\end{verbatim}
}
\end{minipage}
\begin{minipage}[t]{\textwidth}{\color{blue}
\begin{verbatim}
y^2*x;
\end{verbatim}
}
\end{minipage}
\begin{math}\displaystyle
\parbox{8ex}{\color{labelcolor}(\%o9) }
36
\end{math}

\noindent
\begin{minipage}[t]{8ex}{\color{red}\bf
\begin{verbatim}
(%i10) 
\end{verbatim}
}
\end{minipage}
\begin{minipage}[t]{\textwidth}{\color{blue}
\begin{verbatim}
x;
\end{verbatim}
}
\end{minipage}
\begin{math}\displaystyle
\parbox{8ex}{\color{labelcolor}(\%o10) }
4
\end{math}
\end{maxcode}

Si queremos eliminar la asignaci\'on del valor $4$ a la variable $x$, o bien
eliminar todas las asignaciones realizadas hasta el momento, podemos utilizar el comando
\texttt{remvalue} (del ingl\'es \emph{remove value}, que podr\'iamos traducir por
`eliminar el valor'), cuya sintaxis se ilustra a continuaci\'on:
\begin{maxcode}
\noindent
\begin{minipage}[t]{8ex}{\color{red}\bf
\begin{verbatim}
(%i11) 
\end{verbatim}
}
\end{minipage}
\begin{minipage}[t]{\textwidth}{\color{blue}
\begin{verbatim}
remvalue(x);
\end{verbatim}
}
\end{minipage}
\begin{math}\displaystyle
\parbox{8ex}{\color{labelcolor}(\%o11) }
[x]
\end{math}

\noindent
\begin{minipage}[t]{8ex}{\color{red}\bf
\begin{verbatim}
(%i12) 
\end{verbatim}
}
\end{minipage}
\begin{minipage}[t]{\textwidth}{\color{blue}
\begin{verbatim}
x;
\end{verbatim}
}
\end{minipage}
\begin{math}\displaystyle
\parbox{8ex}{\color{labelcolor}(\%o12) }
x
\end{math}

\noindent
\begin{minipage}[t]{8ex}{\color{red}\bf
\begin{verbatim}
(%i13) 
\end{verbatim}
}
\end{minipage}
\begin{minipage}[t]{\textwidth}{\color{blue}
\begin{verbatim}
remvalue(all);
\end{verbatim}
}
\end{minipage}
\begin{math}\displaystyle
\parbox{8ex}{\color{labelcolor}(\%o13) }
[y]
\end{math}

\noindent
\begin{minipage}[t]{8ex}{\color{red}\bf
\begin{verbatim}
(%i14) 
\end{verbatim}
}
\end{minipage}
\begin{minipage}[t]{\textwidth}{\color{blue}
\begin{verbatim}
y;
\end{verbatim}
}
\end{minipage}
\begin{math}\displaystyle
\parbox{8ex}{\color{labelcolor}(\%o14) }
y
\end{math}
\end{maxcode}

Adem\'as, Maxima tiene capacidades simb\'olicas, a diferencia de una calculadora
puede trabajar con expresiones y variables abstractas, de una forma muy parecida
a como lo hace una persona, con las notaciones matem\'aticas habituales
(incluyendo la definici\'on de funciones). Sin embargo, algo que hay que tener
presente es que, a veces, el comportamiento (aunque l\'ogico) no es intuitivo. Por
ejemplo, en una expresi\'on como

\begin{maxcode}
\noindent
\begin{minipage}[t]{8ex}{\color{red}\bf
\begin{verbatim}
(%i15) 
\end{verbatim}
}
\end{minipage}
\begin{minipage}[t]{\textwidth}{\color{blue}
\begin{verbatim}
(a+b)^4;
\end{verbatim}
}
\end{minipage}
\begin{math}\displaystyle
\parbox{8ex}{\color{labelcolor}(\%o15) }
{\left( b+a\right) }^{4}
\end{math}
\end{maxcode}

\noindent vemos que Maxima no expande la expresi\'on como esperar\'iamos\index{Maxima!expansi\'on}. 
Esto es porque
al trabajar simb\'olicamente, la expresi\'on $(a+b)^4$ tiene sentido por ella
misma y si queremos una expresi\'on equivalente constru\'ida de alguna otra
forma, debemos ped\'irselo expl\'icitamente al programa. En este caso, como lo que
queremos es una expansi\'on, podemos utilizar el comando \texttt{expand} (del verbo ingl\'es
correspondiente a `expandir'). Su uso es bastante intuitivo y no requiere de mayores explicaciones:

\begin{maxcode}
\noindent
\begin{minipage}[t]{8ex}{\color{red}\bf
\begin{verbatim}
(%i16) 
\end{verbatim}
}
\end{minipage}
\begin{minipage}[t]{\textwidth}{\color{blue}
\begin{verbatim}
expand((a+b)^4);
\end{verbatim}
}
\end{minipage}
\begin{math}\displaystyle
\parbox{8ex}{\color{labelcolor}(\%o16) }
{b}^{4}+4\,a\,{b}^{3}+6\,{a}^{2}\,{b}^{2}+4\,{a}^{3}\,b+{a}^{4}
\end{math}
\end{maxcode}

El proceso inverso ser\'ia factorizar la expresi\'on,\index{Maxima!factorizaci\'on} para lo
cual se usa el comando \texttt{factor}, tan intuitivo como el precedente:

\begin{maxcode}
\noindent
\begin{minipage}[t]{8ex}{\color{red}\bf
\begin{verbatim}
(%i17) 
\end{verbatim}
}
\end{minipage}
\begin{minipage}[t]{\textwidth}{\color{blue}
\begin{verbatim}
factor(%);
\end{verbatim}
}
\end{minipage}
\begin{math}\displaystyle
\parbox{8ex}{\color{labelcolor}(\%o17) }
{\left( b+a\right) }^{4}
\end{math}
\end{maxcode}

Aqu\'i hemos usado el s\'imbolo \texttt{ \%} para referirnos a la salida del \emph{\'ultimo}
comando. Fij\'emonos en que Maxima etiqueta autom\'aticamente cada comando tecleado,
asign\'andole un \'indice \texttt{ \%in}, y cada salida, asign\'andole un \'indice \texttt{ \%on}.
As\'i, podemos referirnos con facilidad a expresiones que ya aparecieron. Veamos
algunos ejemplos,\index{Maxima!input}\index{Maxima!output}

\begin{maxcode}
\noindent
\begin{minipage}[t]{8ex}{\color{red}\bf
\begin{verbatim}
(%i18) 
\end{verbatim}
}
\end{minipage}
\begin{minipage}[t]{\textwidth}{\color{blue}
\begin{verbatim}
valor1:sqrt(2);
\end{verbatim}
}
\end{minipage}
\begin{math}\displaystyle
\parbox{8ex}{\color{labelcolor}(\%o18) }
\sqrt{2}
\end{math}

\noindent
\begin{minipage}[t]{8ex}{\color{red}\bf
\begin{verbatim}
(%i19) 
\end{verbatim}
}
\end{minipage}
\begin{minipage}[t]{\textwidth}{\color{blue}
\begin{verbatim}
valor2:3;
\end{verbatim}
}
\end{minipage}
\begin{math}\displaystyle
\parbox{8ex}{\color{labelcolor}(\%o19) }
3
\end{math}

\noindent
\begin{minipage}[t]{8ex}{\color{red}\bf
\begin{verbatim}
(%i20) 
\end{verbatim}
}
\end{minipage}
\begin{minipage}[t]{\textwidth}{\color{blue}
\begin{verbatim}
%*(%o18)^2=valor2*(valor1)^2;
\end{verbatim}
}
\end{minipage}
\begin{math}\displaystyle
\parbox{8ex}{\color{labelcolor}(\%o20) }
6=6
\end{math}
\end{maxcode}

Las constantes matem\'aticas especiales, como $e$, $i$, $\pi$, $\phi$ etc., se denotan
en Maxima anteponi\'endolas el signo \texttt{ \%} (aunque este comportamiento se puede
modificar en las preferencias del programa) como \texttt{ \%e},  \texttt{ \%i},  \texttt{ \%pi}, 
 \texttt{ \%phi}, etc.\index{Maxima!constantes}
\begin{maxcode}
\noindent
\begin{minipage}[t]{8ex}{\color{red}\bf
\begin{verbatim}
(%i21) 
\end{verbatim}
}
\end{minipage}
\begin{minipage}[t]{\textwidth}{\color{blue}
\begin{verbatim}
%pi;
\end{verbatim}
}
\end{minipage}
\begin{math}\displaystyle
\parbox{8ex}{\color{labelcolor}(\%o21) }
\pi 
\end{math}

\noindent
\begin{minipage}[t]{8ex}{\color{red}\bf
\begin{verbatim}
(%i22) 
\end{verbatim}
}
\end{minipage}
\begin{minipage}[t]{\textwidth}{\color{blue}
\begin{verbatim}
%pi,numer;
\end{verbatim}
}
\end{minipage}
\begin{math}\displaystyle
\parbox{8ex}{\color{labelcolor}(\%o22) }
3{.}141592653589793
\end{math}

\noindent
\begin{minipage}[t]{8ex}{\color{red}\bf
\begin{verbatim}
(%i23) 
\end{verbatim}
}
\end{minipage}
\begin{minipage}[t]{\textwidth}{\color{blue}
\begin{verbatim}
%phi;
\end{verbatim}
}
\end{minipage}
\begin{math}\displaystyle
\parbox{8ex}{\color{labelcolor}(\%o23) }
\phi
\end{math}

\noindent
\begin{minipage}[t]{8ex}{\color{red}\bf
\begin{verbatim}
(%i24) 
\end{verbatim}
}
\end{minipage}
\begin{minipage}[t]{\textwidth}{\color{blue}
\begin{verbatim}
%phi,numer;
\end{verbatim}
}
\end{minipage}
\begin{math}\displaystyle
\parbox{8ex}{\color{labelcolor}(\%o24) }
1{.}618033988749895
\end{math}
\end{maxcode}

Para definir funciones, usamos el mismo s\'imbolo que en la escritura habitual,\index{Maxima!funciones}
\begin{maxcode}
\noindent
\begin{minipage}[t]{8ex}{\color{red}\bf
\begin{verbatim}
(%i25) 
\end{verbatim}
}
\end{minipage}
\begin{minipage}[t]{\textwidth}{\color{blue}
\begin{verbatim}
f(x):=x^2-x-1;
\end{verbatim}
}
\end{minipage}
\begin{math}\displaystyle
\parbox{8ex}{\color{labelcolor}(\%o25) }
\mathrm{f}\left( x\right) :={x}^{2}-x-1
\end{math}
\end{maxcode}

\noindent y podemos asignar valores haciendo
\begin{maxcode}
\noindent
\begin{minipage}[t]{8ex}{\color{red}\bf
\begin{verbatim}
(%i26) 
\end{verbatim}
}
\end{minipage}
\begin{minipage}[t]{\textwidth}{\color{blue}
\begin{verbatim}
f(1);
\end{verbatim}
}
\end{minipage}
\begin{math}\displaystyle
\parbox{8ex}{\color{labelcolor}(\%o26) }
-1
\end{math}

\noindent
\begin{minipage}[t]{8ex}{\color{red}\bf
\begin{verbatim}
(%i27) 
\end{verbatim}
}
\end{minipage}
\begin{minipage}[t]{\textwidth}{\color{blue}
\begin{verbatim}
f(%phi),numer;
\end{verbatim}
}
\end{minipage}
\begin{math}\displaystyle
\parbox{8ex}{\color{labelcolor}(\%o27) }
0{.}0
\end{math}
\end{maxcode}

Como ya hemos dicho, Maxima posee capacidades simb\'olicas. Aqu\'i vemos
como resolver la ecuaci\'on $ax-3=b$ donde la inc\'ognita es $x$ (es
importante indicar cu\'al es la inc\'ognita, pues podr\'ia ser cualquiera
de entre $x,a,b$):\index{Maxima!ecuaciones}
\begin{maxcode}
\noindent
\begin{minipage}[t]{8ex}{\color{red}\bf
\begin{verbatim}
(%i28) 
\end{verbatim}
}
\end{minipage}
\begin{minipage}[t]{\textwidth}{\color{blue}
\begin{verbatim}
solve(a*x-3=b,x);
\end{verbatim}
}
\end{minipage}
\begin{math}\displaystyle
\parbox{8ex}{\color{labelcolor}(\%o28) }
[x=\frac{b+3}{a}]
\end{math}

\noindent
\begin{minipage}[t]{8ex}{\color{red}\bf
\begin{verbatim}
(%i29) 
\end{verbatim}
}
\end{minipage}
\begin{minipage}[t]{\textwidth}{\color{blue}
\begin{verbatim}
solve(a*x-3=b,a);
\end{verbatim}
}
\end{minipage}
\begin{math}\displaystyle
\parbox{8ex}{\color{labelcolor}(\%o29) }
[a=\frac{b+3}{x}]
\end{math}
\end{maxcode}

Como en la escritura manual, es importante usar par\'entesis para indicar operaciones complicadas.
Si no se usan, Maxima utilizar\'a las reglas de precedencia de operadores,\index{Operaciones!precedencia}
\begin{maxcode}
\noindent
\begin{minipage}[t]{8ex}{\color{red}\bf
\begin{verbatim}
(%i30) 
\end{verbatim}
}
\end{minipage}
\begin{minipage}[t]{\textwidth}{\color{blue}
\begin{verbatim}
135/15-6+5;
\end{verbatim}
}
\end{minipage}
\begin{math}\displaystyle
\parbox{8ex}{\color{labelcolor}(\%o30) }
8
\end{math}

\noindent
\begin{minipage}[t]{8ex}{\color{red}\bf
\begin{verbatim}
(%i31) 
\end{verbatim}
}
\end{minipage}
\begin{minipage}[t]{\textwidth}{\color{blue}
\begin{verbatim}
135/(15-6)+5;
\end{verbatim}
}
\end{minipage}
\begin{math}\displaystyle
\parbox{8ex}{\color{labelcolor}(\%o31) }
20
\end{math}
\end{maxcode}

Para acabar esta breve introducci\'on, veamos c\'omo graficar funciones. Si ya
tenemos la funci\'on definida previamente, como en\,  \texttt{\%o25}, podemos
usar el nombre que le dimos y especificar la variable independiente y el intervalo
en que queremos el gr\'afico,\index{Maxima!gr\'aficos}
\begin{maxcode}
\noindent
\begin{minipage}[t]{8ex}{\color{red}\bf
\begin{verbatim}
(%i32) 
\end{verbatim}
}
\end{minipage}
\begin{minipage}[t]{\textwidth}{\color{blue}
\begin{verbatim}
plot2d(f(x),[x,-4,4]);
\end{verbatim}
}
\end{minipage}
\begin{math}\displaystyle
\parbox{8ex}{\color{labelcolor}(\%o32) }
\; graphics
\end{math}
\end{maxcode}

Tambi\'en podemos elegir el intervalo de ordenadas,
\begin{maxcode}
\noindent
\begin{minipage}[t]{8ex}{\color{red}\bf
\begin{verbatim}
(%i33) 
\end{verbatim}
}
\end{minipage}
\begin{minipage}[t]{\textwidth}{\color{blue}
\begin{verbatim}
plot2d(f(x),[x,-4,4],[y,-5,5]);
\end{verbatim}
}
\end{minipage}
\begin{math}\displaystyle
\parbox{8ex}{\color{labelcolor}(\%o33) }\;
graphics
\end{math}
\end{maxcode}

El comando \texttt{plot2d} muestra la gr\'afica en una ventana aparte (la del programa
Gnuplot, que es el que se usa para generar los gr\'aficos). Para que aparezca
en el documento de trabajo usamos \texttt{wxplot2d} (lo cual requiere que trabajemos con
la interfaz gr\'afica wxMaxima):\index{Maxima!wxMaxima}\index{wxMaxima}
\begin{maxcode}
\noindent
\begin{minipage}[t]{8ex}{\color{red}\bf
\begin{verbatim}
(%i34) 
\end{verbatim}
}
\end{minipage}
\begin{minipage}[t]{\textwidth}{\color{blue}
\begin{verbatim}
wxplot2d(f(x),[x,-4,4]);
\end{verbatim}
}
\end{minipage}
\begin{math}\displaystyle
\parbox{8ex}{\color{labelcolor}(\%t34) }\centering
\includegraphics[width=9cm]{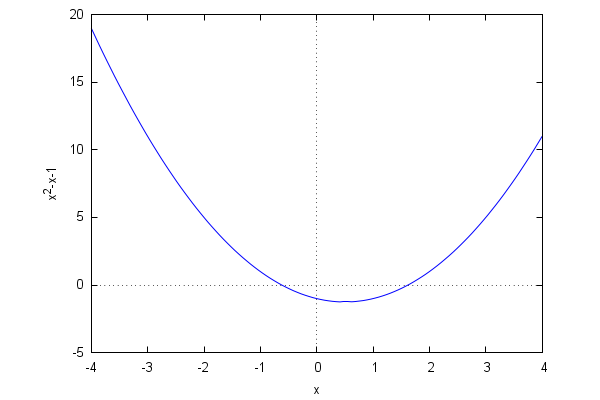}
\end{math}
\begin{math}\displaystyle
\parbox{8ex}{\color{labelcolor}(\%o34) }
\end{math}
\end{maxcode}

Se puede mejorar la est\'etica con algunos comandos opcionales. La lista de estos comandos
es muy extensa y est\'a fuera de lugar tratar de enumerarlos todos aqu\'i. El lector interesado puede
consultar la p\'agina en espa\~nol del manual en l\'inea de Maxima correspondiente a los gr\'aficos,
\url{http://maxima.sourceforge.net/docs/manual/es/maxima_12.html}. \'Unicamente mencionaremos los
elementos m\'as b\'asicos, como la \emph{leyenda} (en ingl\'es, \emph{legend}) o t\'itulo del
gr\'afico (que, como era de esperar, se especifica con el comando \texttt{legend} y el texto deseado
\emph{entre comillas}\footnote{Esto es importante: Maxima distingue como una cadena de 
texto\index{Maxima!cadena de texto} todo aquello que se escriba entre comillas.}) y el comando
\texttt{ylabel} (del ingl\'es \emph{label}, traducci\'on de `etiqueta'), que sirve para darle un
nombre al eje de las ordenadas $y$:
\begin{maxcode}
\noindent
\begin{minipage}[t]{8ex}{\color{red}\bf
\begin{verbatim}
(%i35) 
\end{verbatim}
}
\end{minipage}
\begin{minipage}[t]{\textwidth}{\color{blue}
\begin{verbatim}
wxplot2d(f(x),[x,-4,4],
[legend,"Una funcion cuadratica"],[ylabel,"f(x)"]);
\end{verbatim}
}
\end{minipage}
\begin{math}\displaystyle
\parbox{8ex}{\color{labelcolor}(\%t35) }\centering
\includegraphics[width=9cm]{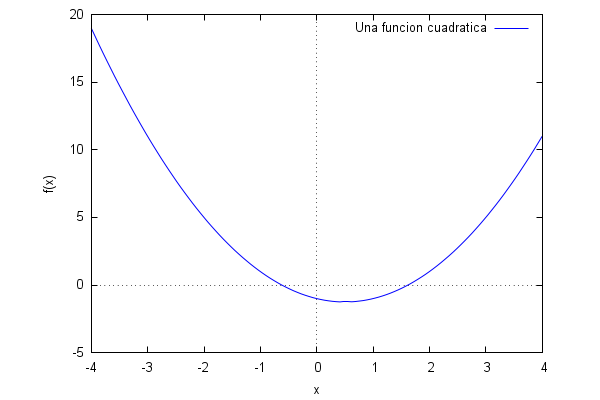}
\end{math}
\begin{math}\displaystyle
\parbox{8ex}{\color{labelcolor}(\%o35) }
\end{math}
\end{maxcode}

\section{Actividades de c\'omputo}
\begin{enumerate}[(A)]
\item Resolver el ejercicio \ref{ejer1} usando Maxima.
Comparar el resultado con el
obtenido `a mano'.
\item Graficar en Maxima las funciones del ejercicio \ref{ejer4}.
\item Realizar en Maxima las operaciones de los ejercicios \ref{ejer5} y \ref{ejer7}.
Comparar el resultado con el obtenido `a mano'.
\end{enumerate}

\chapter{Aritm\'etica modular. Divisibilidad}\label{cap2}
\vspace{2.5cm}\setcounter{footnote}{0}
\section{Nociones b\'asicas}
La aritm\'etica modular tambi\'en se conoce como `aritm\'etica del reloj':
en un reloj, cada vez que se
llega a las 12 se vuelve a contar desde 0, de modo que cada n\'umero del 1 al 12, tal como 
el 7, tiene dos interpretaciones: las 7 de la ma\~ nana y las 7 de la tarde.
En este caso, se dice que se tiene una relaci\'on de equivalencia m\'odulo 12.

En general, se dice que se tiene una \emph{relaci\'on de equivalencia m\'odulo} $m$ 
en el conjunto de los
n\'umeros enteros $\mathbb{Z}$ si `se vuelve a contar al llegar a $m$', de manera que los n\'umeros
$a\in\mathbb{Z}$ y $a+m\in\mathbb{Z}$ se consideran equivalentes. Por ejemplo, en el reloj, son
equivalentes la 1 y las 13 horas, las 2 y las 14, etc. Si la relaci\'on de equivalencia es m\'odulo
4, ser\'an equivalentes 0 y 4, 1 y 5, 2 y 6, 3 y 7. A partir de aqu\'i, los n\'umeros se repiten.
\begin{center}
\includegraphics[scale=0.25]{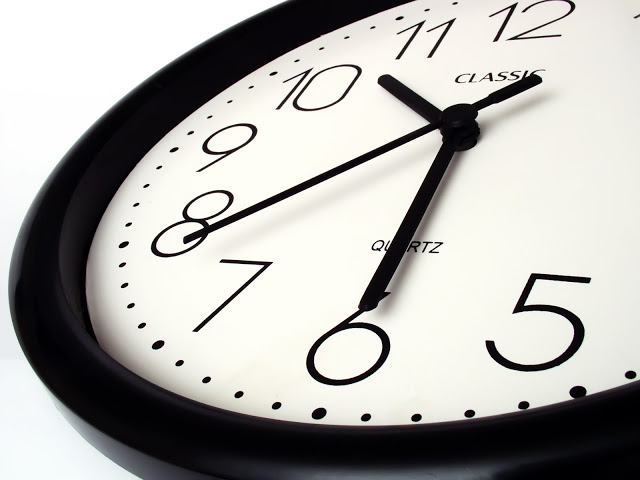} 
\end{center}
Una manera m\'as precisa de decir esto es que dos n\'umeros $a,b\in\mathbb{Z}$
son \emph{equivalentes m\'odulo} $m$
(denotado $b\equiv a\mod m$ y le\'ido `$b$ es \emph{congruente} con $a$ m\'odulo $m$')\index{Num@N\'umero!congruencia m\'odulo $m$} si $b$ 
se obtiene de $a$ sumando alg\'un m\'ultiplo de $m$, es decir, si
$$
b=a+km\,,
$$
para alg\'un $k\in\mathbb{Z}$. Y aun otra forma de expresar esto es diciendo que
la diferencia $b-a$ es un m\'ultiplo de $m$ (recordemos que el conjunto de los
m\'ultiplos de $m$ se denota $\stackrel{\circ}{m}$):
$$
b-a\in \stackrel{\circ}{m}\,,
$$
o que la diferencia $b-a$ es divisible por $m$\footnote{La notaci\'on $m|b$ se lee `$m$ divide a $b$'
y significa que $b$ es un m\'ultiplo de $m$, esto es, existe un n\'umero $p\in\mathbb{Z}-\{0\}$ tal que $b=pm$.}
$$
m|(b-a)\, .
$$
Por ejemplo, m\'odulo 13 los n\'umeros 21 y 73 son equivalentes, porque su diferencia es un m\'ultiplo
de 13:
$$
73-21=52=4\cdot 13\,.
$$

Una propiedad muy importante de la aritm\'etica modular es que
\begin{equation}\label{mods}
\mbox{si }b\equiv a\mod m,\mbox{ entonces }(b+c)\equiv (a+c)\mod m, 
\end{equation}
para cualquier $c\in\mathbb{Z}$. La demostraci\'on es muy sencilla:
si $b\equiv a\mod m$ es que existe un $k\in\mathbb{Z}$ tal que
$$
b-a=km,
$$
pero entonces, sumando y restando $c$,
$$
(b+c)-(a+c)=km,
$$
que es lo mismo que decir que $(b+c)\equiv (a+c)\mod m$.

Est\'a claro que en la pr\'actica, para saber si dos n\'umeros $a$ y $b$ son congruentes m\'odulo $m$, hay que
determinar de alguna manera si su diferencia $b-a$ es \emph{divisible} por $m$. Esto nos lleva a considerar 
algunas ideas b\'asicas sobre divisibilidad. Por ejemplo, tenemos los siguientes resultados:
\begin{list}{$\blacktriangleright$}{}
\item Si $m|a$ y $m|b$, con $m\neq 0$, entonces $m|(a+b)$ (en palabras: si un n\'umero divide a otros
dos, tambi\'en divide a su suma).
\item Si $a\geq b$, $m|a$ y $m|b$, entonces $m|(a-b)$.
\item Si $m|a$ y $k\in\mathbb{N}$, entonces $m|ka$.
\end{list}
Esto deber\'ia ser obvio. Por ejemplo, para ver que la primera afirmaci\'on es cierta, escribimos
$a=mp$ y $b=mq$ para algunos $p,q\in\mathbb{N}$, pues eso es lo que significa que $a,b$ sean
divisibles por $m$. Pero entonces $a+b=m(p+q)$ y de nuevo por definici\'on, $m|(a+b)$.

Veamos c\'omo se aplican estos resultados para construir algunos criterios elementales de divisibilidad\index{Divisibilidad!criterios}.
\begin{list}{$\blacktriangleright$}{}
\item \emph{Un n\'umero $b\in\mathbb{Z}$ es divisible por 2 si y s\'olo si su \'ultimo d\'igito es un n\'umero par}.
Esto  puede verse como sigue (consideraremos por simplicidad n\'umeros de tres d\'igitos, pero las ideas son
v\'alidas para cualquier n\'umero, s\'olo que la notaci\'on se complica sin a\~nadir nada de inter\'es).
Si tenemos un n\'umero $n$ con d\'igitos $abc$, tal como como 724 (en el que $a=7$, $b=2$, $c=4$), es que
$$
n=a\cdot 10^2+b\cdot 10^1+ c\cdot 10^0= a\cdot 10^2+b\cdot 10+c =10\cdot (a\cdot 10+b)+c .
$$
Como $2|10$, es claro que $2|10\cdot (a\cdot 10+b)$. Si $c$ es par, es que $c=2k$ para alg\'un $k\in\mathbb{N}$ y
por tanto $2|c$, luego por las propiedades que acabamos de ver $2|(10\cdot (a\cdot 10+b)+c)=n$.\\
Acabamos de ver que si $c$ es m\'ultiplo de 2, entonces 2 divide a $n=abc$. Veamos ahora el rec\'iproco,
esto es, que si 2 divide a $n=abc$, es que $c=2k$.\\
Si suponemos que $2|(10\cdot (a\cdot 10+b)+c)$, como tambi\'en $2|10\cdot (a\cdot 10+b)$ resulta que
2 divide a la resta: $2|(10\cdot (a\cdot 10+b)+c-10\cdot (a\cdot 10+b))=c$; es decir, $2|c$, que es
lo que quer\'iamos comprobar.
\item \emph{Un n\'umero $b\in\mathbb{Z}$ es divisible por 5 si y s\'olo si su \'ultimo d\'igito es un 
m\'ultiplo de 5}. La prueba de este hecho es id\'entica a la del caso del 2, cambiando 2 por 5.
\item \emph{Un n\'umero $b\in\mathbb{Z}$ es divisible por 10 si y s\'olo si su \'ultimo d\'igito es $0$}.
La prueba de este hecho es id\'entica a la del caso del 2, cambiando 2 por 10.
\item \emph{Un n\'umero $b\in\mathbb{Z}$ es divisible por 4 si y s\'olo si el n\'umero formado por sus
\'ultimos dos d\'igitos es divisible por 4}. Para verlo, de nuevo supongamos que $n=abc$, es decir,
$n=a\cdot 10^2+b\cdot 10+c$. Observando que $a\cdot 10^2=a\cdot 25\cdot 4$, tenemos que
$4|a\cdot 10^2$, as\'i que si $4|(b\cdot 10+c)$ seguro que divide a la suma $a\cdot 10^2+b\cdot 10+c$, que es $n$.\\
Rec\'iprocamente, si $4|(a\cdot 10^2+b\cdot 10+c)$, como sabemos que divide a $a\cdot 10^2$ dividir\'a a la
resta: $4|(a\cdot 10^2+b\cdot 10+c-a\cdot 10^2)=b\cdot 10+c$, es decir, 4 divide al n\'umero de d\'igitos $bc$.
\item \emph{Un n\'umero $b\in\mathbb{Z}$ es divisible por 8 si y s\'olo si el n\'umero formado por sus
\'ultimos tres d\'igitos es divisible por 8}. La prueba de este hecho es id\'entica a la del caso del 4,
cambiando 4 por 8.
\item \emph{Un n\'umero $b\in\mathbb{Z}$ es divisible por 3 si y s\'olo si la suma de sus d\'igitos es
divisible por 3}. En efecto, si escribimos el n\'umero $n=abc$ como $n=a\cdot 10^2+b\cdot 10+c$, podemos
hacer el siguiente c\'alculo:
\begin{align*}
n =& a\cdot (99+1)+b\cdot (9+1)+c \\
  =& a\cdot 99+a\cdot 1+b\cdot 9+b\cdot 1+c \\
  =& a\cdot 99+b\cdot 9 +(a+b+c) \\
  =& (a\cdot 11+b)\cdot 9 +(a+b+c).
\end{align*}
Pero como claramente $3|9$, es $3|(a\cdot 11+b)\cdot 9$ y si es $3|(a+b+c)$ entonces dividir\'a a la
suma $n=(a\cdot 11+b)\cdot 9 +(a+b+c)$.\\
Rec\'iprocamente, si $3|n$, es que $3|((a\cdot 11+b)\cdot 9 +(a+b+c))$ y como $3|(a\cdot 11+b)\cdot 9$,
tambi\'en dividir\'a a la resta: $3|((a\cdot 11+b)\cdot 9 +(a+b+c)-(a\cdot 11+b)\cdot 9)=(a+b+c)$.
\item \emph{Un n\'umero $b\in\mathbb{Z}$ es divisible por 9 si y s\'olo si la suma de sus d\'igitos es
divisible por 9}. La prueba de este hecho es id\'entica a la del caso del 3,
cambiando 3 por 9.
\item \emph{Un n\'umero $b\in\mathbb{Z}$ es divisible por 11 si y s\'olo si la suma de sus d\'igitos que ocupan
la posici\'on par menos la suma de los d\'igitos que ocupan la posici\'on impar, es 0 \'o un m\'ultiplo de 11}
(no importa si se toman los d\'igitos comenzando por la derecha o por la izquierda). Por ejemplo, el n\'umero
32 527 es divisible por 11, ya que la suma de los d\'igitos que ocupan posici\'on par (comenzando por la izquierda)
es $2+2=4$, la de los que ocupan posici\'on impar es $3+5+7=15$ y $4-15=-11$ es m\'ultiplo de 11. De hecho,
$32\, 527=11\cdot 2957$. Sin embargo el n\'umero, 1 431 no es divisible por 11, ya que $ 4+1-(1+3)=1$, que no es 0 ni
m\'ultiplo de 11. El criterio se demuestra observando que para las potencias de 10, se tiene
\begin{align*}
10^0 &=1=0\cdot 11+1 \\
10^1 &=10=1\cdot 11 -1 \\
10^2 &=100 =9\cdot 11 +1 \\
10^3 &=1000 =91\cdot 11 -1 \\
10^4 &=10000 =909\cdot 11 +1 \\
10^5 &=100000 =9091\cdot 11 -1 \\
	 &\cdots
\end{align*}
y por tanto, dado un n\'umero como $abcde=a10^4+b10^3+c10^2+d10^1+e10^0$, ocurre que
\begin{align*}
abcde &=a(909\cdot 11 +1)+b(91\cdot 11-1)+c(9\cdot 11+1)+d(1\cdot 11-1)+e(0\cdot 11+1) \\
	  &=11(909a +91b+9c+d)+a-b+c-d+e .
\end{align*}
Como, obviamente, $11(909a +91b+9c+d)$ es divisible por 11, el n\'umero $abcde$ lo ser\'a tambi\'en
si y s\'olo si lo es $a-b+c-d+e$ (\'o bien, que esta suma d\'e $0$). 
\end{list}
A partir de estos resultados se pueden obtener otros. Por ejemplo, un n\'umero es divisible por 6 si y s\'olo
si es divisible por 2 y por 3, simult\'aneamente\footnote{!`Cuidado! Esto s\'olo funciona porque 2 y 3 no tienen
divisores comunes. No es cierto, por ejemplo, que un n\'umero sea divisible por 24 si es divisible por 4 y por 6,
pues 12 cumple que $4|12$ y tambi\'en $6|12$, pero obviamente $24\nmid 12$.}. Es curioso el caso de la divisibilidad
por 7, el \'unico n\'umero entre 2 y 10 que nos falta. Existen, desde luego, criterios de divisibilidad de un
n\'umero por 7, pero son extremadamente complicados y por esa raz\'on no los veremos aqu\'i. Los interesados
pueden consultar, por ejemplo, la p\'agina del proyecto Descartes
\url{http://recursostic.educacion.es/descartes/web/materiales\_didacticos/divisibilidad/criterios\_div.htm}.

\section{Aplicaciones}

Cada libro que se publica en el mundo tiene asociado un c\'odigo identificador, su \emph{ISBN} (International
Standard Book Number, en espa\~nol: \emph{n\'umero est\'andar internacional de libro})\index{Co@C\'odigo!ISBN}. Este n\'umero se asignaba de
una forma antes del a\~no 2007 y de otra distinta a partir de esa fecha. Antes de 2007, el ISBN de un libro ten\'ia
\emph{en total} 10 d\'igitos, que se asignaban de la siguiente manera: los primeros uno a cinco d\'igitos indicaban
el idioma y pa\'is de la obra (el 0 corresponde al ingl\'es, por ejemplo). El siguiente grupo, que tambi\'en puede
tener de uno a cinco d\'igitos, identifica la editorial. Un tercer grupo, que de nuevo puede contener de uno a cinco
d\'igitos, identifica el libro concreto. El total de d\'igitos de estos tres grupos debe ser 9. Por ejemplo, un
libro publicado en ingl\'es en USA, s\'olo necesita un d\'igito en el primer grupo, el 0. Si la editorial es la
American Mathematical Society, el segundo grupo tendr\'a 4 d\'igitos, $8218$. 
En el caso concreto del libro
\emph{Principles of Functional Analysis (second edition)}, de M. Schecter, el identificador del tercer grupo tiene otras cuatro cifras, $2895$. As\'i, de los 10 d\'igitos del ISBN, los primeros 9 de este libro en particular son $0-8218-2895$. 
Otros libros de esta editorial compartir\'an los primeros 5 d\'igitos; por ejemplo, el libro 
\emph{Dirac Operators in Riemannian Geometry}, de T. Friedrich, tiene asignado un ISBN cuyos 
primeros 9 d\'igitos son $0-8218-2055$.
\begin{center}
\includegraphics[scale=1]{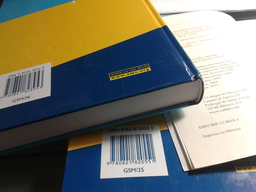} 
\end{center}
Como otro ejemplo, el libro \emph{Historia m\'inima de M\'exico}, de D. Cos\'io y otros, lo publica el Colegio de M\'exico en \'este pa\'is. Los primeros 9 d\'igitos de su ISBN son $968-12-0618$. Los tres primeros, $968$, indican que el pa\'is es
M\'exico y el idioma el espa\~nol\footnote{Puede consultarse, por ejemplo, Wikipedia: 
\url{http://es.wikipedia.org/wiki/Anexo:ISBN\_por\_pa\%C3\%ADs}.}.
El segundo grupo, 12, identifica a la editorial del Colegio de M\'exico, mientras
que el tercer grupo, $0618$, corresponde a este t\'itulo en particular, de entre todos los que edita el Colmex.

Al principio mencionamos que el ISBN (anterior a 2007) tiene 10 d\'igitos y hasta ahora s\'olo hemos hablado
de los nueve primeros. ?`Qu\'e ocurre con el \'ultimo? Se trata de un \emph{d\'igito de control}. Es decir, no
es un d\'igito que se asigne arbitrariamente, sino que se determina a partir de los nueve primeros de manera
que se pueda detectar
si el ISBN completo se escribi\'o bien o no. Es decir, el \'ultimo d\'igito est\'a dado de tal manera que si
hici\'eramos un pedido del libro \emph{Historia m\'inima de M\'exico} por internet y dij\'eramos que su ISBN es
$968-12-0618-4$ (en lugar del correcto, que es $968-12-0618-5$), la librer\'ia sabr\'ia que hay un error en
los datos y podr\'ia pedirnos que los revisemos antes de tramitar el pedido.

Este d\'igito de control se construye as\'i: dado un ISBN $a_1a_2a_3a_4a_5a_6a_7a_8a_9a_{10}$, el d\'igito
$a_{10}$ debe ser tal que
$$
a_1+2a_2+3a_3+4a_4+5a_5+6a_6+7a_7+8a_8+9a_9+10a_{10} \equiv 0\mod 11\, .
$$
En nuestro ejemplo, el ISBN ser\'a $968-12-0618-a_{10}$, donde $a_{10}$ debe cumplir que
$$
9+2\cdot 6+3\cdot 8+4\cdot 1+5\cdot 2+6\cdot 0+7\cdot 6+8\cdot 1+9\cdot 8+10a_{10}=181+10a_{10}\equiv 0\mod 11\,.
$$
Es decir: $a_{10}$ debe elegirse de manera que multiplicado por 10 y sumado a 181 d\'e un m\'ultiplo de 11.
?`C\'omo resolver este problema? En \eqref{mods} se vi\'o que se puede operar con las ecuaciones en
congruencias de la misma forma como se hace con los enteros. Entonces, tendremos que de
$181+10a_{10}\equiv 0\mod 11$ resulta $10a_{10}\equiv -181\mod 11$. Ahora bien, m\'odulo $11$, $-181$ es
equivalente a $-5$, pues $-181=11\cdot (-16)-5$ (en otras palabras, $-5$ es el resto de dividir $-181$ entre
$11$). As\'i, $ 10a_{10}\equiv -5\mod 11$. Por tanto, s\'olo hay que encontrar entre los n\'umeros congruentes
con $-5$ m\'odulo $11$ uno que sea m\'ultiplo de $10$. O sea, hay que hallar un $k\in\mathbb{Z}$ tal que
$-5+11\cdot k$ acabe en $0$. Obviamente, podemos tomar $k=5$ pues $-5+55=50$, que es m\'ultiplo de $10$.
En definitiva, el d\'igito de control para nuestro libro ha de ser tal que $10a_{10}=-5+11\cdot 5$, es
decir, $a_{10}=5$ y el ISBN completo es
$968-12-0618-5$.\\
Una observaci\'on: cuando el d\'igito de control es $10$, se escribe como X. As\'i ocurre, por ejemplo, en
el ISBN $0-486-45844-X$.

El c\'odigo ISBN que se usa a partir de 2007 tiene 13 d\'igitos, pero el principio en que se basa es el mismo
que acabamos de describir. Es interesante observar que otros c\'odigos, como los de barras que se usan en los
supermercados o los de las tarjetas de cr\'edito, tambi\'en se construyen con aritm\'etica modular.

\section{Ejercicios}

\begin{enumerate}
\setcounter{enumi}{8}

\item Determinar si los siguientes enunciados son verdaderos o falsos.
\begin{multicols}{3}
\begin{enumerate}[(a)]
\item $1\equiv 13\mod 7$
\item $1\equiv 22\mod 7$
\item $12\equiv 4\mod 3$
\item $-3 \equiv 5\mod 2$
\item $145 \equiv 122 \mod 11$
\item $145 \equiv 123 \mod 11$
\end{enumerate}
\end{multicols}

\item\label{znring} Al conjunto de los n\'umeros enteros congruentes m\'odulo $m$ \emph{no equivalentes entre si},
se le denota $\mathbb{Z}_m$\index{Zm@$\mathbb{Z}_m$}. Sus
elementos se escriben en la forma $\overline{a}\in\mathbb{Z}_m$ y en realidad son a su vez conjuntos
de n\'umeros\footnote{En terminolog\'ia matem\'atica se dice que son \emph{clases de 
equivalencia}.}:
$\overline{a}$ agrupa a todos aqu\'ellos n\'umeros que son congruentes con $a$ m\'odulo $m$.
Por ejemplo, para calcular $\mathbb{Z}_2$, empezamos con $a=0\in\mathbb{Z}$. Los n\'umeros congruentes con
$0$ m\'odulo 2 son los de la forma $0+2\cdot k$, con $k\in\mathbb{Z}$, es decir, los n\'umeros pares. Por
tanto, $\overline{0}=\{ \ldots,-4,-2,0,2,4,\ldots \}$. Tomemos ahora el primer n\'umero posterior a
$0\in\mathbb{Z}$ que no est\'e en $\overline{0}$, el 1. Los n\'umeros congruentes con
1 m\'odulo 2 son los de la forma $1+2\cdot k$, con $k\in\mathbb{Z}$, es decir, los n\'umeros impares:
$\overline{1}=\{ \ldots,-5,-3,-1,1,3,5,\ldots \}$. Entre $\overline{0}$ y $\overline{1}$ reunimos todos
los n\'umeros enteros, luego $\mathbb{Z}_2=\{ \overline{0},\overline{1}\}$.\\
Calcular $\mathbb{Z}_3$, $\mathbb{Z}_4$ y $\mathbb{Z}_5$. ?`Cu\'antos elementos tiene en general
$\mathbb{Z}_m$?

\item\label{zetam} Comprobar que si se toman dos elementos cualesquiera $\overline{m},\overline{n}\in\mathbb{Z}_6$,
y se define $\overline{m}+\overline{n}$ como la clase $\overline{m+n}\in\mathbb{Z}_6$, no importa qu\'e
elementos se elijan en $\overline{m}$ y $\overline{n}$ el resultado siempre es el mismo. Es decir, para
calcular $\overline{2}+\overline{3}$, podr\'iamos elegir cualquier elemento de
$\overline{2}=\{\ldots, -10,-4,2,8,14,\ldots\}$ (como $-4$) y cualquiera de 
$\overline{3}=\{\ldots, -9,-3,3,9,15,\ldots\}$ (como $9$) para hacer
$$
\overline{2}+\overline{3}=\overline{-4+9}=\overline{5}=\{\ldots,-7,-1,5,11,17,\ldots\}\,,
$$
y el resultado no cambiar\'a.
Hacer lo mismo si se define un producto en $\mathbb{Z}_6$ como
$$
\overline{m}\cdot \overline{n} = \overline{m\cdot n}.
$$
As\'i pues, los elementos de $\mathbb{Z}_6$ se pueden sumar y multiplicar y sigue
resultando elementos de $\mathbb{Z}_6$. Esta propiedad se resume diciendo que $\mathbb{Z}_6$ 
es un \emph{anillo algebraico}\index{Anillo algebraico}. ?`Existe un inverso para cada 
elemento respecto del producto?

\item Considerar $\mathbb{Z}_7$ y comprobar que para toda clase $\overline{m}\in\mathbb{Z}_7$ existe
un inverso respecto del producto, es decir, una clase $\overline{n}\in\mathbb{Z}_7$ tal que
$\overline{m}\cdot\overline{n}=\overline{1}$. ?`A qu\'e se debe la diferencia con $\mathbb{Z}_6$?

\item En un texto chino de matem\'aticas del siglo V dC, se encuentra el siguiente problema:
`Han Xing\index{Xing, Han} tiene un ej\'ercito. Si se forman los hombres en columnas de a 3, sobran 2 soldados.
Si se forman en columnas de a 5, sobran 3. Y si se forman en columnas de a 7, sobran 2.
?`Cu\'antos soldados tiene Han Xing? Resolverlo\footnote{Existe otro texto del a\~no 1247 dC
que da un m\'etodo general para resolver estos problemas, conocido hoy como
el \emph{teorema chino del resto}\index{Teorema!chino del resto}.}.  
\begin{center}
\includegraphics[scale=0.6]{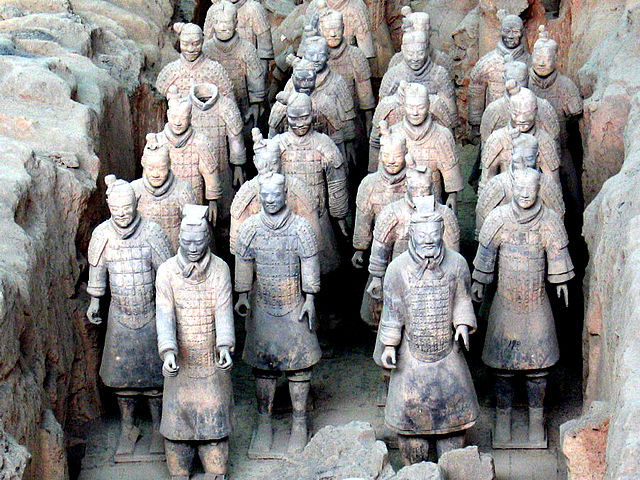} 
\end{center}

\item !`Nos vamos de cumplea\~nos!. El padrino del cumplea\~nero, un ingeniero electr\'onico
resentido con sus profesores de c\'alculo, al saber que somos estudiantes de algo que lleva el
nombre de matem\'aticas siente que ha llegado la hora de la venganza y nos propone un problema
sencillo para dejarnos en rid\'iculo delante de todos los chamacos asistentes a la fiesta. Nos
dice: `tengo una bolsa con un n\'umero de caramelos comprendido entre 20 y 50. Si los reparto
entre 5 de estos chamacos, me sobran 2. Si los reparto entre 6 me sobra 1. ?`Cu\'antos me
sobrar\'an si los reparto entre 7?'. Dejen bien alto el honor de su carrera resolviendo el
problema.
%

\item El s\'imbolo $!$ aplicado a un n\'umero se denomina \emph{factorial}. Por ejemplo,
$5!$ se lee `5 factorial' y se define por $5!=5\cdot 4\cdot 3\cdot 2\cdot 1$ (en general,
escribir\'iamos $n!=n\cdot (n-1)\cdot (n-2)\cdots 3\cdot\cdot 2\cdot 1$). Decidir cu\'al de
las siguientes afirmaciones es verdadera o falsa.
\begin{multicols}{3}
\begin{enumerate}[(a)]
\item $6|6!$
\item $5|6!$
\item $8|6!$
\item $9|6!$
\item $30|30!$
\item $40|30!$
\item $30|(30!+1)$
\end{enumerate}
\end{multicols}

\item Estudiar la divisibilidad de los siguientes n\'umeros.
\begin{multicols}{4}
\begin{enumerate}[(a)]
\item 111
\item 121
\item 221
\item 414
\end{enumerate}
\end{multicols}

\item Decidir si las siguientes afirmaciones son verdaderas o falsas
utilizando criterios de divisibilidad.
\begin{multicols}{2}
\begin{enumerate}[(a)]
\item $6|80$
\item $4|15\, 000$
\item $24|325\, 608$
\item $40|1\, 732\, 800$
\item $15| 10\, 000$
\item $12| 32\, 304$
\item $45| 13\, 075$
\item $36| 677\, 916$
\end{enumerate}
\end{multicols}

\item Temprano por la ma\~nana en una tiendita de abarrotes cualquiera (las tienditas de abarrotes, a diferencia de las grandes superficies, tienen sus precios en m\'ultiplos de 0.5\$).
Do\~na Chonita\index{Chonita, Do\~na} pone una bolsa con papas y chayotes
sobre el mostrador y dice: `llevo dos kilos de papas a 9\$ y 2 chayotes a 5\$ cada uno.
Tambi\'en llevo 3 chocolates y 6 caramelos pa' los nietecitos, pero no s\'e a cu\'anto est\'an'.
Don S\'eforo\index{S\'eforo, Don} le contesta: `de los caramelos y los chocolates son 17.50\$'. Do\~na Chonita se pone
seria: `no manche, don S\'eforo, eche bien la cuenta'. El S\'eforo vuelve a calcular y le pide
disculpas a la viejita. ?`Por qu\'e?

\item Probar que los n\'umeros de 6 d\'igitos de la forma $abcabc$ son siempre
divisibles por 13. ?`Hay alg\'un otro n\'umero que tambi\'en los divida siempre?

\item Observemos que $5!=5\cdot 4\cdot 3\cdot 2\cdot 1$ es divisible por $2,3,4$ y $5$.
\begin{enumerate}[(a)]
\item Probar que $5!+2$, $5!+3$, $5!+4$ y $5!+5$ son n\'umeros compuestos
(es decir, producto de dos o m\'as factores).
\item Encontrar $1\,000$ n\'umeros consecutivos que sean compuestos.
\end{enumerate}


\end{enumerate}

\section{Uso de la tecnolog\'ia de c\'omputo}

El sistema de \'algebra computacional Maxima puede trabajar con aritm\'etica modular directamente.
La funci\'on \texttt{mod(n,m)} devuelve el resto de dividir $n$ entre $m$. Por
ejemplo, como $7=2\cdot 3+1$, se tiene
\begin{maxcode}
\noindent
\begin{minipage}[t]{8ex}{\color{red}\bf
\begin{verbatim}
(%i1) 
\end{verbatim}
}
\end{minipage}
\begin{minipage}[t]{\textwidth}{\color{blue}
\begin{verbatim}
mod(7,3);
\end{verbatim}
}
\end{minipage}
\begin{math}\displaystyle
\parbox{8ex}{\color{labelcolor}(\%o1) }
1
\end{math}
\end{maxcode}

\noindent Otra manera de decir esto es que $7\equiv 1\mod 3$. Si queremos saber a qui\'en es
equivalente $138$ m\'odulo $6$, basta con hacer
\begin{maxcode}
\noindent
\begin{minipage}[t]{8ex}{\color{red}\bf
\begin{verbatim}
(%i2) 
\end{verbatim}
}
\end{minipage}
\begin{minipage}[t]{\textwidth}{\color{blue}
\begin{verbatim}
mod(138,6);
\end{verbatim}
}
\end{minipage}
\begin{math}\displaystyle
\parbox{8ex}{\color{labelcolor}(\%o2) }
0
\end{math}
\end{maxcode}

\noindent Este resultado se pod\'ia haber anticipado porque $138$ es un m\'ultiplo de $6$
(es divisible por $2$ por acabar en cifra par, y es divisible por $3$ por serlo
la suma de sus cifras $1+3+8=12$), en concreto $138\div 6=23$.\\

Recordando que en $\mathbb{Z}_a$ las operaciones son 
$\overline{p}+\overline{q}=\overline{p+q}$ y 
$\overline{p}\cdot\overline{q}=\overline{p\cdot q}$, podemos calcular
cosas como $11+5\mod 4$, $3\cdot 7\mod 5$, etc.
\begin{maxcode}
\noindent
\begin{minipage}[t]{8ex}{\color{red}\bf
\begin{verbatim}
(%i3) 
\end{verbatim}
}
\end{minipage}
\begin{minipage}[t]{\textwidth}{\color{blue}
\begin{verbatim}
mod(11+5,4);
\end{verbatim}
}
\end{minipage}
\begin{math}\displaystyle
\parbox{8ex}{\color{labelcolor}(\%o3) }
0
\end{math}

\noindent
\begin{minipage}[t]{8ex}{\color{red}\bf
\begin{verbatim}
(%i4) 
\end{verbatim}
}
\end{minipage}
\begin{minipage}[t]{\textwidth}{\color{blue}
\begin{verbatim}
mod(3*7,5);
\end{verbatim}
}
\end{minipage}
\begin{math}\displaystyle
\parbox{8ex}{\color{labelcolor}(\%o4) }
1
\end{math}
\end{maxcode}

Maxima proporciona el comando \texttt{inv\_mod(n,m)} para calcular el inverso
del n\'umero $n$ m\'odulo $m$:
\begin{maxcode}
\noindent
\begin{minipage}[t]{8ex}{\color{red}\bf
\begin{verbatim}
(%i5) 
\end{verbatim}
}
\end{minipage}
\begin{minipage}[t]{\textwidth}{\color{blue}
\begin{verbatim}
inv_mod(13,3);
\end{verbatim}
}
\end{minipage}
\begin{math}\displaystyle
\parbox{8ex}{\color{labelcolor}(\%o5) }
1
\end{math}
\end{maxcode}

En $\mathbb{Z}_a$ no todo elemento tiene inverso. El elemento $\overline{0}$ no tiene inverso
nunca (por ejemplo, $6$ no tiene inverso m\'odulo $3$ porque $6\equiv 0\mod 3$). En este caso
la funci\'on \texttt{inv\_mod} devuelve el mensaje \texttt{false}: 
\begin{maxcode}
\noindent
\begin{minipage}[t]{8ex}{\color{red}\bf
\begin{verbatim}
(%i6) 
\end{verbatim}
}
\end{minipage}
\begin{minipage}[t]{\textwidth}{\color{blue}
\begin{verbatim}
inv_mod(6,3);
\end{verbatim}
}
\end{minipage}
\begin{math}\displaystyle
\parbox{8ex}{\color{labelcolor}(\%o6) }
false
\end{math}
\end{maxcode}

\noindent Otro ejemplo menos trivial es que $14$ no tiene inverso m\'odulo $6$:
\begin{maxcode}
\noindent
\begin{minipage}[t]{8ex}{\color{red}\bf
\begin{verbatim}
(%i7) 
\end{verbatim}
}
\end{minipage}
\begin{minipage}[t]{\textwidth}{\color{blue}
\begin{verbatim}
inv_mod(14,6);
\end{verbatim}
}
\end{minipage}
\begin{math}\displaystyle
\parbox{8ex}{\color{labelcolor}(\%o7) }
false
\end{math}
\end{maxcode}
El factorial $m!$ de un n\'umero $m\in\mathbb{Z}$ tambi\'en es conocido
por Maxima:
\begin{maxcode}
\noindent
\begin{minipage}[t]{8ex}{\color{red}\bf
\begin{verbatim}
(%i8) 
\end{verbatim}
}
\end{minipage}
\begin{minipage}[t]{\textwidth}{\color{blue}
\begin{verbatim}
5!;
\end{verbatim}
}
\end{minipage}
\begin{math}\displaystyle
\parbox{8ex}{\color{labelcolor}(\%o8) }
120
\end{math}

\noindent
\begin{minipage}[t]{8ex}{\color{red}\bf
\begin{verbatim}
(%i9) 
\end{verbatim}
}
\end{minipage}
\begin{minipage}[t]{\textwidth}{\color{blue}
\begin{verbatim}
10!;
\end{verbatim}
}
\end{minipage}
\begin{math}\displaystyle
\parbox{8ex}{\color{labelcolor}(\%o9) }
3628800
\end{math}
\end{maxcode}
Con la aritm\'etica modular, es f\'acil saber si un n\'umero $n$ es divisible
por otro $m$, pues en ese caso, el resto de dividir $n$ entre $m$ es cero. As\'i,
$n$ es divisible por $3$ si al escribir \texttt{mod(n,3)} Maxima devuelve $0$. Por
ejemplo,
\begin{maxcode}
\noindent
\begin{minipage}[t]{8ex}{\color{red}\bf
\begin{verbatim}
(%i10) 
\end{verbatim}
}
\end{minipage}
\begin{minipage}[t]{\textwidth}{\color{blue}
\begin{verbatim}
mod(1637294754675632839,3);
\end{verbatim}
}
\end{minipage}
\definecolor{labelcolor}{RGB}{100,0,0}
\begin{math}\displaystyle
\parbox{8ex}{\color{labelcolor}(\%o10) }
1
\end{math}
\end{maxcode}
\noindent nos dice que $1\, 637\, 294\, 754\, 675\, 632\, 839$ no es divisible
por $3$, mientras que $1\, 637\, 294\, 754\, 675\, 632\, 838$ s\'i lo es:
\begin{maxcode}
\noindent
\begin{minipage}[t]{8ex}{\color{red}\bf
\begin{verbatim}
(%i11) 
\end{verbatim}
}
\end{minipage}
\begin{minipage}[t]{\textwidth}{\color{blue}
\begin{verbatim}
mod(1637294754675632838,3);
\end{verbatim}
}
\end{minipage}
\begin{math}\displaystyle
\parbox{8ex}{\color{labelcolor}(\%o11) }
0
\end{math}
\end{maxcode}
Por \'ultimo, vamos a aprovechar la ocasi\'on para introducir el concepto de \emph{listas}.
Una lista es un objeto b\'asico para Maxima: todo son listas. Una lista\index{Maxima!listas} se define dando sus
elementos encerrados entre par\'entesis cuadrados, como en

\begin{maxcode}
\noindent
\begin{minipage}[t]{8ex}{\color{red}\bf
\begin{verbatim}
(%i12) 
\end{verbatim}
}
\end{minipage}
\begin{minipage}[t]{\textwidth}{\color{blue}
\begin{verbatim}
lista1:[a,b,1,c,2];
\end{verbatim}
}
\end{minipage}
\definecolor{labelcolor}{RGB}{100,0,0}
\begin{math}\displaystyle
\parbox{8ex}{\color{labelcolor}(\%o12) }
[a,b,1,c,2]
\end{math}
\end{maxcode}
\noindent N\'otese que la lista no tiene por qu\'e consistir en n\'umeros o letras, pueden
estar combinados. Hay much\'isimos comandos en Maxima para trabajar con listas, la mayor\'ia
de los cuales son muy intuitivos:

\begin{maxcode}
\noindent
\begin{minipage}[t]{8ex}{\color{red}\bf
\begin{verbatim}
(%i13) 
\end{verbatim}
}
\end{minipage}
\begin{minipage}[t]{\textwidth}{\color{blue}
\begin{verbatim}
second(lista1);
\end{verbatim}
}
\end{minipage}
\definecolor{labelcolor}{RGB}{100,0,0}
\begin{math}\displaystyle
\parbox{8ex}{\color{labelcolor}(\%o13) }
b
\end{math}

\noindent
\begin{minipage}[t]{8ex}{\color{red}\bf
\begin{verbatim}
(%i14) 
\end{verbatim}
}
\end{minipage}
\begin{minipage}[t]{\textwidth}{\color{blue}
\begin{verbatim}
third(lista1);
\end{verbatim}
}
\end{minipage}
\definecolor{labelcolor}{RGB}{100,0,0}
\begin{math}\displaystyle
\parbox{8ex}{\color{labelcolor}(\%o14) }
1
\end{math}

\noindent
\begin{minipage}[t]{8ex}{\color{red}\bf
\begin{verbatim}
(%i15) 
\end{verbatim}
}
\end{minipage}
\begin{minipage}[t]{\textwidth}{\color{blue}
\begin{verbatim}
lista1^2;

\end{verbatim}
}
\end{minipage}
\definecolor{labelcolor}{RGB}{100,0,0}
\begin{math}\displaystyle
\parbox{8ex}{\color{labelcolor}(\%o15) }
[{a}^{2},{b}^{2},1,{c}^{2},4]
\end{math}
\end{maxcode}
Se pueden crear listas con el comando \texttt{makelist}, cuya sintaxis tambi\'en es
muy intuitiva: se da una expresi\'on que contiene un \'indice (en este ejemplo, $k$)
y se especifica el rango de valores que toma el \'indice:

\begin{maxcode}
\noindent
\begin{minipage}[t]{8ex}{\color{red}\bf
\begin{verbatim}
(%i16) 
\end{verbatim}
}
\end{minipage}
\begin{minipage}[t]{\textwidth}{\color{blue}
\begin{verbatim}
makelist(k+2,k,3,9);
\end{verbatim}
}
\end{minipage}
\definecolor{labelcolor}{RGB}{100,0,0}
\begin{math}\displaystyle
\parbox{8ex}{\color{labelcolor}(\%o16) }
[5,6,7,8,9,10,11]
\end{math}
\end{maxcode}
Naturalmente, todos los comandos que hemos visto hasta ahora se pueden combinar.
Por ejemplo, podemos obtener una lista con varios valores de n\'umeros congruentes
con $-3$ m\'odulo $8$ y ver cu\'al de ellos es divisible por $5$:

\begin{maxcode}
\noindent
\begin{minipage}[t]{8ex}{\color{red}\bf
\begin{verbatim}
(%i17) 
\end{verbatim}
}
\end{minipage}
\begin{minipage}[t]{\textwidth}{\color{blue}
\begin{verbatim}
mod(makelist(-3+8*k,k,-8,8),5);
\end{verbatim}
}
\end{minipage}
\definecolor{labelcolor}{RGB}{100,0,0}
\begin{math}\displaystyle
\parbox{8ex}{\color{labelcolor}(\%o17) }
[3,1,4,2,0,3,1,4,2,0,3,1,4,2,0,3,1]
\end{math}
\end{maxcode}
\noindent De aqu\'i se ve que los n\'umeros de la forma $-3+8k$, con $-8\leq k\leq 8$,
que son m\'ultiplos de
$5$ son los que corresponden a los \'indices $k=-4,1,6$. De hecho, esos n\'umeros son
$-3-32=-35$, $-3+8=5$ y $-3+48=45$, todos obviamente divisibles por $5$.

Se\~nalemos que en Maxima es muy sencillo generar las tablas de
composici\'on para la suma y el producto en $\mathbb{Z}_m$. Por ejemplo, para
generar la tabla de composici\'on para el producto en $\mathbb{Z}_4$, crear\'iamos
primero una funci\'on que calcule 
$\overline{i}\cdot\overline{i}=(i\cdot j)\,\mathrm{mod}4$ para cada $i,j\in\mathbb{Z}$
(notemos que los \'indices en Maxima comienzan a contar desde $1$, por lo que para
disponer de $\overline{0}$ hay que restar $1$):
\begin{maxcode}
\noindent
\begin{minipage}[t]{8ex}{\color{red}\bf
\begin{verbatim}
(%i18) 
\end{verbatim}
}
\end{minipage}
\begin{minipage}[t]{\textwidth}{\color{blue}
\begin{verbatim}
pz4[i,j]:=mod((i-1)*(j-1),4)$
\end{verbatim}
}
\end{minipage}
\end{maxcode}
\noindent Comprobemos que el resultado es el esperado. Por ejemplo,
$\overline{2}\cdot\overline{3}=(2\cdot 3)\,\mathrm{mod}_4=6\,\mathrm{mod}_4=2$,
\begin{maxcode}
\noindent
\begin{minipage}[t]{8ex}{\color{red}\bf
\begin{verbatim}
(%i19) 
\end{verbatim}
}
\end{minipage}
\begin{minipage}[t]{\textwidth}{\color{blue}
\begin{verbatim}
pz4[2,3];
\end{verbatim}
}
\end{minipage}
\definecolor{labelcolor}{RGB}{100,0,0}
\begin{math}\displaystyle
\parbox{8ex}{\color{labelcolor}(\%o19) }
2
\end{math}
\end{maxcode}
\noindent Finalmente, creamos una tabla (o matriz) que recorra los \'indices
$i,j$ desde $1$ hasta $4$ (de manera que los productos sean entre las clases
$\{\overline{0},\overline{1},\overline{2},\overline{3}\}$):
\begin{maxcode}
\noindent
\begin{minipage}[t]{8ex}{\color{red}\bf
\begin{verbatim}
(%i20) 
\end{verbatim}
}
\end{minipage}
\begin{minipage}[t]{\textwidth}{\color{blue}
\begin{verbatim}
tabla_pz4:genmatrix(pz4,4,4);
\end{verbatim}
}
\end{minipage}
\definecolor{labelcolor}{RGB}{100,0,0}
\begin{math}\displaystyle
\parbox{8ex}{\color{labelcolor}(\%o20) }
\begin{pmatrix}0 & 0 & 0 & 0\cr 0 & 1 & 2 & 3\cr 0 & 2 & 0 & 2\cr 0 & 3 & 2 & 1\end{pmatrix}
\end{math}
\end{maxcode}

\section{Actividades de c\'omputo}
\begin{enumerate}[(A)]
\item Construir las tablas de composici\'on para el producto de $\mathbb{Z}_3$, $\mathbb{Z}_5$,
$\mathbb{Z}_6$, $\mathbb{Z}_7$, $\mathbb{Z}_8$ y $\mathbb{Z}_9$. 
?`Qu\'e elementos son invertibles? ?`Podr\'ias aventurar alguna relaci\'on entre la existencia
de elementos invertibles respecto del producto en $\mathbb{Z}_m$ y el propio $m$?
\item Construir las tablas de composici\'on para la suma en $\mathbb{Z}_3$, $\mathbb{Z}_4$,
$\mathbb{Z}_5$ y $\mathbb{Z}_6$.
\item Utilizando Maxima, calcular el d\'igito de control correspondiente a los siguientes ISBN:
\begin{multicols}{2}
\begin{enumerate}[(a)]
\item $84-95594-00-?$
\item $81-297-0594-?$
\item $84-9732-471-?$
\item $970-643-644-?$
\end{enumerate}
\end{multicols}

\item Utilizando Maxima, determinar cu\'al de los siguientes ISBN es correcto y cu\'al no
\begin{multicols}{2}
\begin{enumerate}[(a)]
\item $968-511-570-9$
\item $0-140-1-2505-1$
\item $0-387-95375-2$
\item $981-01-3734-0$
\end{enumerate}
\end{multicols}

\end{enumerate}

\chapter{N\'umeros primos. Mcm y mcd}\label{cap3}
\vspace{2.5cm} \setcounter{footnote}{0}
\section{Nociones b\'asicas}
Un n\'umero primo es un n\'umero natural mayor que $1$ y que tiene exactamente dos divisores naturales\index{Num@N\'umero!primo}:
\'el mismo y el 1 (en particular, el propio 1 no es primo).
Para determinar n\'umeros primos peque\~nos se puede usar un algoritmo
conocido como la \emph{criba de Erat\'ostenes}\index{Erat\'ostenes!criba}: se escriben los n\'umeros desde el 2 hasta el mayor que se
quiere considerar (por ejemplo, el 100 si se quieren los primos menores que 100). Se van tachando los
m\'ultiplos de 2 (excepto el propio 2), pues esos no pueden ser primos. Luego se tachan los m\'ultiplos
de 3 (excepto el propio 3), y as\'i sucesivamente hasta que se llegue al tope 
elegido\footnote{En la pr\'actica no hace falta tanto,
el proceso se agota mucho antes. Concretamente, si el tope es el n\'umero $n$, de poder escribirse
$n=p\cdot q$ es claro que uno de los n\'umeros $p$ \'o $q$ debe ser menor que $\sqrt{n}$, as\'i que el proceso
se detiene como mucho al llegar a $\sqrt{n}$.} (en el ejemplo hasta llegar a los m\'ultiplos de 99). Los n\'umeros que
quedan sin tachar son precisamente los primos. En la figura se muestra el segundo paso de este proceso
para los primos menores que 100
(los m\'ultiplos de 2 est\'an marcados en naranja y los de 3 en amarillo):
\begin{center}
\begin{tikzpicture}[scale=0.7,every node/.style={minimum size=0.96cm,circle,scale=0.7}]
 \foreach \nb/\col in {2/orange,3/yellow} {%
  \pgfmathtruncatemacro{\nbi}{\nb}  
  \loop 
  \pgfmathtruncatemacro{\nbi}{\nbi+\nb}
  \testcolorednode{\nbi}{%
  \node[ball color=\col!50](\nbi) at ({mod((\nbi-1),10)},{-floor((\nbi-1)/10)})    
        {\nbi};}{} 
   \pgfmathtruncatemacro{\nbt}{\nbi+\nb}%
   \ifnum\nbt<101  \repeat}
\foreach \nb in {2,...,100} 
    {\testcolorednode{\nb}{%
       \node[ball color=gray!20](\nb) at ({mod((\nb-1),10)},{-floor((\nb-1)/10)}) {\nb};}{}}%
\end{tikzpicture}
\end{center} 

Los n\'umeros primos pueden (y deben) verse como los `ladrillos' con los
que se construyen todos los n\'umeros enteros, en el sentido de que \emph{todo n\'umero
natural se descompone de manera \'unica (salvo el orden) en factores primos}\footnote{El
nombre `primo' proviene del lat\'in `primus' en su acepci\'on de `principal' o
`primordial'. No tiene nada que ver con relaciones de parentesco. Sin embargo, las malas
traducciones frecuentemente tienen m\'as \'exito que aquellas sesudas y fundamentadas, y
dan lugar a secuelas. Por ejemplo, en Aritm\'etica encontramos n\'umeros gemelos, 
n\'umeros amigos,
etc. Hay muchos ejemplos de traducciones incorrectas que arraigaron en la terminolog\'ia
matem\'atica, entre ellas `seno' para nombrar a una de las funciones trigonom\'etricas.}\index{Factores primos}. Esto quiere decir que dado
un n\'umero natural cualquiera $n$, existen unos n\'umeros primos $p_1,\ldots ,p_k$ tales que
$$
n=p_1^{a_1}\cdots p_k^{a_k}
$$
para ciertos exponentes $a_1,\ldots,a_k$ tambi\'en naturales. Por ejemplo, 2 475 puede descomponerse as\'i:
$$
2\,475=3^2\cdot 5^2 \cdot 11 .
$$
Para hallar esta descomposici\'on, vamos dividiendo el n\'umero dado por cada unos de los primos menores
que \'el: por 2, por 3, por 5, etc. Los resultados de la divis\'on se van organizando como en la
siguiente tabla:
\begin{center}
\begin{tabular}{r|c}\label{primdes}
 2 475 & 3 \\ 
 825 & 3 \\ 
 275 & 5 \\ 
 55 & 5 \\ 
 11 & 11 \\ 
\end{tabular}
\end{center} 
de donde, contando las veces que se repite cada factor a la derecha, se ve inmediatamente
que $2\,475=3^2\cdot 5^2 \cdot 11$.

Una vez que se tiene la descomposici\'on en factores primos de dos n\'umeros, es muy f\'acil determinar
su m\'inimo com\'un m\'ultiplo (mcm)\index{Min@M\'inimo Com\'un M\'ultiplo} y su m\'aximo com\'un divisor (mcd)\index{Maximo@M\'aximo Com\'un Divisor}. El mcm de dos n\'umeros
es \emph{el producto de todos sus factores, tanto comunes como no comunes, elevado cada uno al m\'aximo
exponente} con el que aparece en la descomposici\'on de los n\'umeros. El mcd es \emph{el producto
\'unicamente de los factores comunes elevado al menor exponente} con el que aparece en la descomposici\'on
de los n\'umeros\footnote{Cuando uno de los n\'umeros es $0$ se define $\mathrm{mcd}(a,0)=a$, mientras que
$\mathrm{mcm}(a,0)=0$.}.\\
Por ejemplo, para los n\'umeros $1\,620=2^23^45$ y $1\,575=3^25^27$, tenemos que
$$
\mathrm{mcm}(1\,620,1\,575)=2^23^45^27=56\,700 ,
$$
mientras que
$$
\mathrm{mcd}(1\,620,1\,575)=3^25=45\, ,
$$
Una propiedad importante de la descomposici\'on en factores primos de dos n\'umeros $a$ y $b$ es que
$$
\mathrm{mcm}(a,b)=\frac{a\cdot b}{ \mathrm{mcd}(a,b)}\,.
$$
?`Puede el lector ver por qu\'e esto es cierto?

Una forma alternativa de calcular el m\'aximo com\'un divisor la da el algoritmo de Euclides\index{Euclides!algoritmo}.
Si se tienen dos n\'umeros $a$ y $b$, tales que $a<b$, se pueden dividir para obtener un cociente
$c$ y un resto $r$, $b=c\cdot a+r$.
Entonces, ocurre que el mcd de $(a,b)$ es el mismo que el de $(b,r)$ (?`por qu\'e?). El algoritmo
consiste en repetir el proceso, ahora con $b=c_1r+r_1$, y continuar hasta que se llega a un resto igual a
$0$. Por ejemplo, para calcular de nuevo $\mathrm{mcd}(1\,620,1\,575)$ har\'iamos:
\begin{center}
\begin{tabular}{r|l}
 $\dfrac{1620}{1575}$ cociente $0$ y resto $45$ & $\mathrm{mcd}(1620,1575)=\mathrm{mcd}(1575,45)$ \\
  & \\ 
 $\dfrac{1575}{45}$ cociente $35$ y resto $0$ & $\mathrm{mcd}(1575,45)=\mathrm{mcd}(45,0)$ 
\end{tabular}
\end{center}
luego 
$$
\mathrm{mcd}(1620,1575)=\mathrm{mcd}(1575,45)=\mathrm{mcd}(45,0)=45.
$$

\section{Aplicaciones}\label{sec32}

El c\'alculo de los factores primos en la descomposici\'on de un n\'umero entero que hemos presentado en
la p\'agina \pageref{primdes} es muy \'util desde el punto de vista te\'orico, pero no lo es tanto desde
el punto de vista pr\'actico. Supongamos que tenemos un n\'umero grande, del que incluso sabemos que es 
producto s\'olo de dos primos:
$$
7812353094515203333466690911342264649370470629995425393793599374448576999
$$
?`C\'omo podemos hallar su descomposici\'on? Tendr\'iamos que ir viendo si el n\'umero es divisible por
$2$, por $3$, por $5$, etc. Esto, incluso en la computadora m\'as r\'apida, lleva bastante tiempo. Por
ejemplo, el n\'umero precedente se ha obtenido multiplicando
$$
1837645372892083274632483893993837853
$$
y
$$
4251284393473662899058378589010198483
$$
cada uno de los cuales tiene 37 d\'igitos. El tiempo necesario para encontrar la factorizaci\'on en
la computadora m\'as r\'apida que podamos encontrar, seguramente ser\'a mayor que nuestra paciencia.\\
Es posible demostrar que el n\'umero de operaciones necesarias para factorizar un n\'umero con $n$ 
d\'igitos es del orden de $0{.}3\cdot 2^{n/2}/((n/2)\ln 2)$ para el algoritmo `de ensayo y error' que hemos descrito antes
(con $\ln$ el logaritmo natural\index{Logaritmo!natural}, v\'ease m\'as adelante el cap\'itulo \ref{explog}), y del
orden de $\exp (c \sqrt{\log n \log (\log  n)})$, donde $c$ es una constante que depende del algoritmo
(y $\log$ el logaritmo en base $10$, v\'ease de nuevo el cap\'itulo \ref{explog} m\'as adelante), para los algoritmos m\'as eficientes\footnote{Uno de los mejores que se conocen actualmente es el de 
factorizaci\'on por curvas el\'ipticas\index{Curva!el\'iptica} de Hendrik W. Lenstra\index{Lenstra, Hendrik W.}, que tiene $c=1$. Sin embargo, el de Lenstra no es el algoritmo m\'as r\'apido, t\'itulo que le corresponde a la llamada criba general de cuerpos num\'ericos.}. Teniendo en cuenta que una computadora actual es capaz de ejecutar t\'ipicamente $10^{11}$ operaciones por segundo, vemos que
la factorizaci\'on de un n\'umero de $100$ d\'igitos requiere unos $97{.}5$ segundos con el algoritmo
de ensayo y error, pero para uno de
$200$ d\'igitos, el tiempo se dispara a $5{.}5\cdot 10^{16}$ segundos, !`unos $1\,800$ millones de a\~nos!\label{time}.
Esta dificultad es la causa de que el uso de grandes n\'umeros primos sea muy \'util en criptograf\'ia.
Cuando dos personas A y B quieren tener una comunicaci\'on privada, a prueba de esp\'ias, deben
encontrar una manera de trasmitir la informaci\'on encriptada, de modo que nadie m\'as que ellos
pueda conocer el mensaje.
\begin{center}
\includegraphics[scale=0.4]{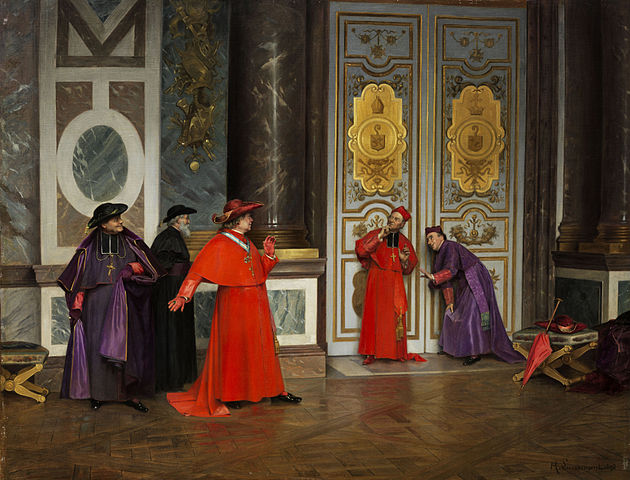} 
\end{center}
El m\'etodo de encriptaci\'on RSA\footnote{Por las iniciales de sus creadores: Rivest, Shamir, Adler.}\index{Criptograf\'ia!m\'etodo RSA}\index{Rivest, Ron}\index{Shamir, Adi}\index{Adleman, Leonard}
tiene los siguientes pasos para que Gaby\index{Gaby} pueda enviar un mensaje cifrado a 
Ale\index{Ale}
\begin{enumerate}
\item Ale comienza eligiendo un par de n\'umeros primos muy grandes, con m\'as de 200 d\'igitos, $p$ y $q$. Para
ilustrar el m\'etodo, nosotros tomaremos n\'umeros peque\~nos, como $p=11$ y $q=13$.
\item Ale calcula $N=pq$ y se lo hace llegar a Gaby. De hecho, cualquier persona tiene acceso a $N$, incluso
los posibles esp\'ias. Ale puede poner el n\'umero $N$ en su p\'agina web, por ejemplo. En nuestro caso es
$N=143$.
\item Ale, que conoce $p$ y $q$, elige un n\'umero primo relativo con $(p-1)(q-1)$, es decir, elige un
n\'umero $e$ tal que no tiene divisores comunes con $(p-1)(q-1)$. En nuestro caso, $(p-1)(q-1)=10\cdot 12=120$,
que tiene por divisores $2,3,4,5,6,8,10,12,15,20,24,30,40,60$.  Podr\'iamos elegir, por ejemplo, $e=17$.
Este n\'umero tambi\'en puede hacerse p\'ublico. De hecho, al par $(N,e)$ se le llama la \emph{clave p\'ublica}
de Ale.
\item Gaby toma su mensaje, por ejemplo \texttt{hola}, y lo traduce a un n\'umero, de acuerdo con la tabla de equivalencia del c\'odigo ASCII\index{Co@C\'odigo!ASCII} (que se puede consultar en 
\url{http://es.wikipedia.org/wiki/ASCII#Car.C3.A1cteres_imprimibles_ASCII}). En nuestro caso, resulta
`072 111 108 097'.
\item Para cada bloque de tres n\'umeros $abc$ que representa una letra, Gaby calcula el n\'umero
$abc^e\mod N$. Por ejemplo, tendr\'iamos que calcular $72^{17}\mod 143$, $111^{17}\mod 143$,
$108^{17}\mod 143$ y $97^{17}\mod 143$. Los resultados respectivos son $63$, $89$, $114$ y $15$, que
corresponden a los s\'imbolos $?$, $Y$, $r$ y $^S{_I}$ .
\item Gaby env\'ia el mensaje encriptado a Ale: \(\texttt{?Yr}^{\texttt{S}}\sb\texttt{I}\).
\item Para decodificar el mensaje, Ale calcula un n\'umero $d$ tal que
$ed\equiv 1\mod (p-1)(q-1)$, es decir, un $d$ tal que $17d\equiv 1\mod 120$. Por ejemplo,
$d=113$ (n\'otese que $17\cdot 113=1921=16\cdot 120 +1$). El par $(N,d)$ se llama la
\emph{clave privada} de Ale.
\item Ale pasa el mensaje que le envi\'o Gaby a n\'umeros, `063 089 114 015', y calcula
para cada bloque $abc$ el n\'umero
$abc^d\mod N$. Es decir, en nuestro ejemplo calcular\'iamos
$63^{113}\mod 143$, $89^{113}\mod 143$, $114^{113}\mod 143$ y $15^{113}\mod 143$.
Los resultados son, respectivamente, $72$, $ 111$, $108$ y $97$. Es decir,
recuperamos \texttt{hola}.
\end{enumerate}
N\'otese que Gaby env\'ia \emph{p\'ublicamente} su mensaje encriptado a Ale, por ejemplo,
dej\'andole una publicaci\'on en su muro de alguna red social. !`Cualquiera puede ver el 
mensaje encriptado! Sin embargo, para desencriptarlo, el esp\'ia tendr\'ia que 
conocer el n\'umero
$d$, y para eso tiene que averiguar primero $p$ y $q$, lo cual
ya hemos visto que es muy dif\'icil si $p$ y $q$ son suficientemente grandes.

Naturalmente, lo que hemos presentado aqu\'i es una versi\'on muy simplificada del
procedimiento real. Por ejemplo, hay que implementar alg\'un mecanismo para que los n\'umeros que
resulten no caigan fuera del rango de valores ASCII imprimibles (es decir, deben salir n\'umeros
entre 32 y 126). Adem\'as, aqu\'i hemos separado el mensaje en bloques de 3 en 3, una restricci\'on
que tampoco tiene el m\'etodo RSA general. 
Para una implementaci\'on m\'as realista con Maxima,
puede consultarse el art\'iculo \emph{Implementing the RSA cryptosystem with Maxima CAS},
que se puede descargar gratuitamente de la p\'agina 
\url{http://galia.fc.uaslp.mx/~jvallejo/Research.html}.

\section{Ejercicios}

\begin{enumerate}
\setcounter{enumi}{20}
\item Resolver lo siguiente. 
\begin{enumerate}[(a)]
\item Mediante el algoritmo de Erat\'ostenes\index{Erat\'ostenes} 
(la `criba'), hallar los n\'umeros primos menores que 150.
\item Se\~nalar, en la lista anterior, todas las parejas de primos gemelos (son primos gemelos aquellos que distan dos unidades).
\item Determinar si los siguientes n\'umeros son primos: 203, 1 003, 5 123, 42 534.
\item Comprobar la conjetura de Goldbach\index{Goldbach, Christian} (la conjetura
afirma que \emph{todo entero par mayor que $2$, se puede escribir como suma de dos primos})
en los n\'umeros siguientes: 50, 98, 102, 144, 296, 1 525.
\end{enumerate}

\item Calcula el m\'aximo com\'un divisor, y con \'el, el m\'{\i}nimo com\'un 
m\'ultiplo, de las siguientes parejas de n\'umeros utilizando el algoritmo de la
descomposici\'on en factores primos y el algoritmo de Euclides\index{Euclides!algoritmo}. Compara el trabajo
que requiere cada m\'etodo:
\begin{multicols}{2}
\begin{enumerate}[(a)]
\item (45, 75)
\item (102, 222)
\item (666, 1 414)
\item (20 785, 44 350)
\item (332 720 058, 212 788 065)
\end{enumerate} 
\end{multicols}

\item Factorizar en n\'umeros primos: 144, 8 162, 111 111, 9 699 690.

\item Un sat\'elite, A, en \'orbita ecuatorial, tiene un per\'{\i}odo de 180 minutos.  Otro sat\'elite, B, en \'orbita circumpolar, de 300 minutos.  Un centro de seguimiento se encuentra justo en la vertical del punto donde sus \'orbitas se cruzan (a diferentes alturas, por supuesto). 
\begin{enumerate}[(a)]
\item ?`Cada cu\'anto tiempo pasan los dos sat\'elites a la vez por encima del centro de seguimiento? Haz un diagrama de
la situaci\'on. 
\item La antena del centro env\'{\i}a autom\'aticamente se\~nales de localizaci\'on para ambos sat\'elites cada $T$ minutos.  Calcular $T$ de forma que cada sat\'elite reciba, al menos, una se\~nal en cada \'orbita, y que cada vez que coincidan los dos sobre el centro de  seguimiento reciban una se\~nal.
\end{enumerate}
\begin{center}
\includegraphics[scale=0.9]{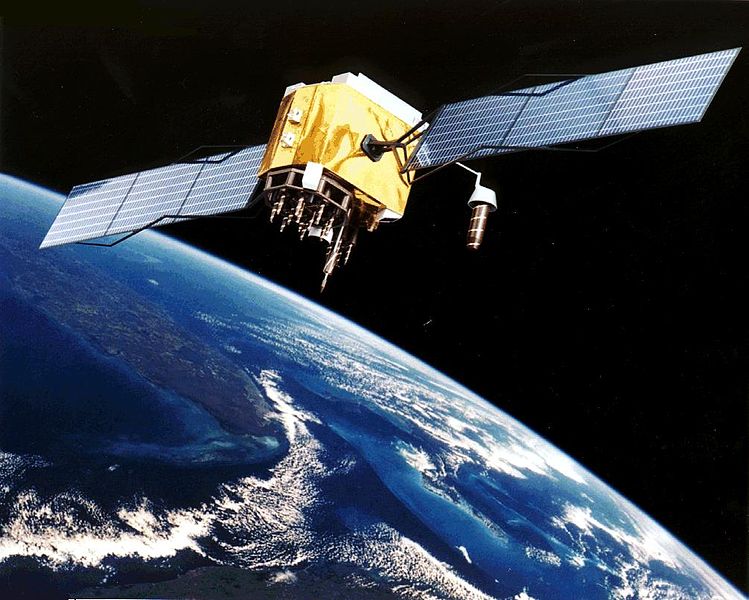} 
\end{center}
\item En un desfile conmemorativo del 16 de Septiembre va a participar una secundaria que tiene 1 000 alumnos.
El director del centro recibe las siguientes instrucciones: 
\begin{enumerate}[(a)] 
\item En una parte del recorrido han de ir formados en fila de 7.
\item En otra parte, con calles m\'as estrechas, hay que cambiar a formaci\'on de 4.
\item En la avenida principal, deben formar en filas de 11.
\item En ning\'un momento debe haber filas incompletas.
\end{enumerate} 
?`Cu\'antos alumnos de la secundaria pueden participar en el desfile?

\item La Ruta 29 y la Ruta 2 de camiones urbanos coinciden durante su recorrido en la parada de la garita de Qu\'imicas.  
Un cami\'on de la Ruta 29 tarda 84 minutos en completar su recorrido.  Otro cami\'on de la Ruta 2 emplea 54 
minutos en el suyo. Si estos dos camiones coinciden a las 12:00, ?`a qu\'e hora volver\'an a coincidir? 

\item Una f\'abrica vende tornillos grandes y peque\~nos.  Los grandes en lotes de 750 tornillos,
los peque\~nos en lotes de 1 185 tornillos. Una ferreter\'{\i}a compra 11 lotes de tornillos grandes y 5 lotes de
tornillos peque\~nos, y quiere venderlos en cajas mixtas que contengan de los dos tipos.  
Hallar todas las posibles formas de agrupar los tornillos en cajas, todas iguales.  
?`Cu\'al es el m\'aximo n\'umero de cajas mixtas que pueden hacer? ?`Y el m\'{\i}nimo? ?`Cu\'antos
tornillos de cada tipo habr\'a en cada uno de estos casos?

\item Elo\'isa\index{Elo\'isa} y Abelardo\index{Abelardo} se aman. Se conocieron en una escuela de verano de matem\'aticas en San Luis
Potos\'i, pero ella estudia en Quer\'etaro y \'el en Monterrey. Los fines de semana, para verse el
mayor tiempo posible, Abelardo toma un cami\'on que sale de Monterrey y va hacia el sur por la 57 a una
velocidad promedio de 80 km/h (teniendo en cuenta el tiempo que se pierde en las paradas). Por su parte,
Elo\'isa toma otro cami\'on desde Quer\'etaro que recorre la 57 en direcci\'on norte, a 60 km/h de
promedio (tambi\'en teniendo en cuenta las paradas). La distancia entre las dos ciudades es de
aproximadamente 700 km y ellos le piden al conductor que les deje bajar cuando sus camiones se cruzan.
\begin{enumerate}[(a)]
\item ?`Cu\'anto tiempo dura el viaje hasta que se encuentran?
\item ?`En qu\'e municipio se encontrar\'an? (puedes consultar, por ejemplo, Google maps).
\end{enumerate}

\item Recordemos que el conjunto de las clases de equivalencia m\'odulo $m$ en los enteros,
$\mathbb{Z}_m$\index{Zm@$\mathbb{Z}_m$} ten\'ia una estructura de anillo, lo que significa que las clases se pueden sumar
y multiplicar entre si. En el ejercicio \ref{zetam} se vi\'o que en $\mathbb{Z}_6$ existen elementos
que no tienen inverso respecto del producto. Probar que si $m=p$ es primo, todo elemento tiene inverso
(se dice entonces que $\mathbb{Z}_p$ es un \emph{cuerpo} \'o \emph{campo num\'erico}).

\item Descomponer el n\'umero 225 de dos formas:
\begin{enumerate}[(a)]
\item Directamente proporcional a 4 y a 5.
\item Inversamente proporcional a 4 y a 5.
\end{enumerate}

\item En un ejido se reparten las tierras de cultivo entre las familias de la comunidad de manera
directamente proporcional al n\'umero de miembros de la familia. Un lote de terreno tiene 990$m^2$ 
y hay que repartirlo entre dos familias de 5 y 7 integrantes. ?`Cu\'antos metros le corresponde a
cada familia?

\item El siguiente problema es un cl\'asico y tiene que estar en la lista,
as\'i que ah\'i va. A un peque\~no poblado cerca de Medina (en Arabia Saud\'i),
lleg\'o un peregrino con fama de sabio
en ruta hacia el sur, a La Meca. Fue recibido
por una familia de caravaneros con la proverbial hospitalidad \'arabe y durante la cena tuvo ocasi\'on de
presenciar una disputa acerca de una herencia. El anciano padre de la familia acababa de fallecer y
hab\'ia dejado instrucciones para repartir entre sus hijos un lote de camellos que \'el administraba
personalmente dentro del negocio familiar. Los camellos eran 19 y los hijos 3. El reparto que hab\'ia
dispuesto el padre, era como sigue: al hijo mayor, la mitad de los camellos. Al mediano, la cuarta parte,
y al menor la quinta parte.
\begin{center}
\includegraphics[scale=1]{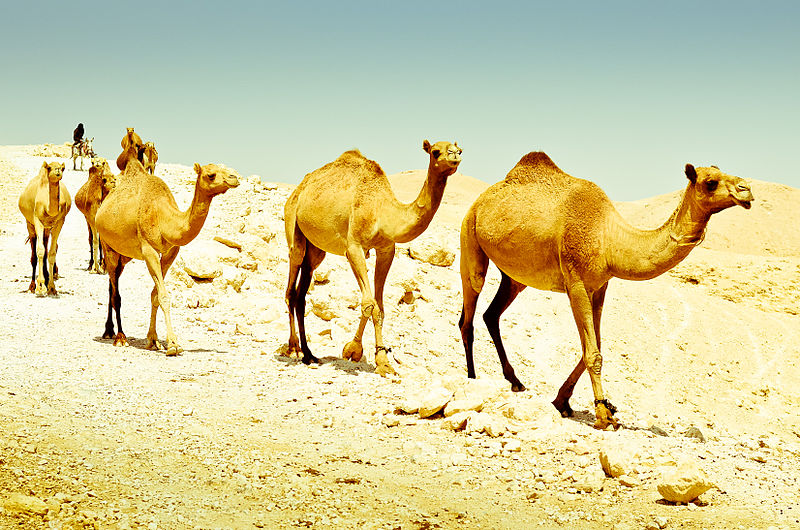} 
\end{center}
\begin{enumerate}[(a)]
\item ?`El padre plante\'o bien el reparto? Justifica la respuesta.
\item Al no ponerse de acuerdo sobre el n\'umero de camellos que le toca a cada hermano, \'estos decidieron
recurrir al forastero con fama de sabio, quien procedi\'o a realizar el reparto de la siguiente forma: tom\'o
su propio camello y lo a\~nadi\'o a los 19 que dej\'o el padre. Del nuevo total, 20, le di\'o 10 al hermano mayor,
5 al mediano y 4 al menor. Recuper\'o su camello y pidi\'o que le sirvieran un poco m\'as de t\'e con menta.
?`Puedes explicar por qu\'e funcion\'o esta manera de hacer el reparto?
\end{enumerate}

\item En un determinado juego, gana quien \emph{menos} puntos acumula al final de la partida. Los
jugadores pueden cometer faltas durante el juego, pero acumulan puntos cada vez que lo hacen. El
n\'umero de puntos es fijo e igual a 210 y se reparten de forma inversamente proporcional al n\'umero
de faltas cometidas. En una partida en que participan tres jugadores $A$, $B$ y $C$, se tiene que
el n\'umero de faltas se distribuye as\'i:
\begin{align*}
A :&\quad 4\mbox{ faltas}. \\
B: &\quad 6\mbox{ faltas}. \\
C: &\quad 12\mbox{ faltas}.
\end{align*}
?`Cu\'antos puntos le corresponden a cada jugador?

\item Las calles de la ciudad est\'an hechas un asco y el presidente municipal va a hacer una visita de 
campo para comprobar que los recursos destinados al programa emergente de bacheo se est\'an ejecutando
correctamente. El director de obras dispone de dos cuadrillas, una formada por trabajadores cualificados,
que parchean una calle en 1 hora, y otra formada por `recomendados', que parchean la misma calle en
3 horas. Si pone las dos cuadrillas a trabajar juntas, ?`en cua\'nto tiempo estar\'a parcheada la calle?

\item Una peque\~na represa para regad\'io se alimenta de dos arroyos. Uno de ellos tarda 2 d\'ias
en llenarla, mientras que otro tarda 3 d\'ias. La represa tiene una compuerta de descarga que la vac\'ia
en 6 d\'ias si se bloquea o desv\'ia el caudal de los dos arroyos. ?`Cu\'anto tiempo tarda la represa en
llenarse  si vierten los dos arroyos y la compuerta de descarga est\'a abierta?

\item Una tiendita de abarrotes vende frutos secos como botana. Compra 30kg de cacahuates a 12\$ el kilo
y los mezcla con 5kg de semillas de girasol, que tienen un precio de 4\$/kg. Vende la botana a
10\$ cada 100g. ?`Qu\'e te parece el precio?
\begin{center}
\includegraphics[scale=0.3]{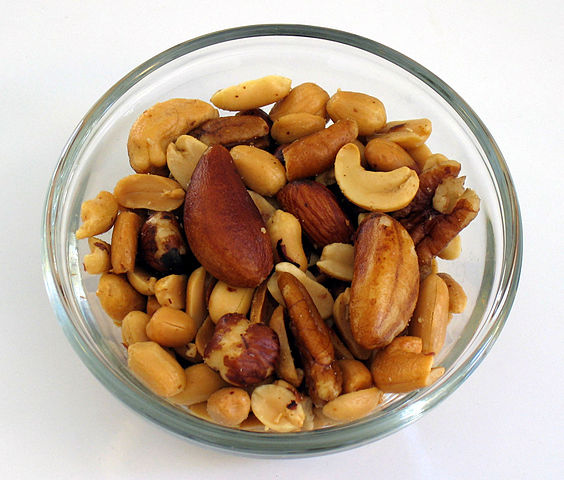}
\end{center}

\item El profe Loreto\index{Loreto, profe} no la hace con el sueldo que gana en la universidad y complementa sus ingresos
vendiendo chocolate, hecho por \'el mismo con cacao que trae de Oaxaca cuando visita a sus suegros.
El kilo de cacao puro tostado cuesta
75\$ y el de az\'ucar 15\$. ?`A qu\'e precio hay que vender el kilo de chocolate (igual cantidad de cacao que de az\'ucar)
para no perder ni ganar? El chocolate adem\'as lleva una peque\~na cantidad de canela para mejorar su sabor, que
se puede suponer que no altera el peso total, pero tiene un costo de unos 10\$ por 30g (1oz), que es lo que se suele
a\~nadir. El profe Loreto vende su chocolate con aroma de canela a 90\$/Kg. Teniendo en cuenta los gastos de 
transporte de la mercanc\'ia, costo de
elaboraci\'on (hay que moler el chocolate para mezclarlo), etc. 
?`Crees que se har\'a rico como chocolatero?  

\item El kilo de papa se vende en el mercado de abastos a unos 9\$/kg
(se puede ver informaci\'on oficial en la p\'agina \href{http://www.campomexicano.gob.mx/portal_siap/Integracion/EstadisticaDerivada/InformaciondeMercados/Mercados/snim/slpa1101.htm}{\tt www.campomexicano.gob.mx}, los precios que se dan aqu\'i eran v\'alidos a finales de
2013).
Una tonelada de aceite de palma industrial cuesta unos 750 d\'olares americanos
(se puede ver en \href{http://www.indexmundi.com/es/precios-de-mercado/?mercancia=aceite-de-palma}{\tt www.indexmundi.com}).
\begin{enumerate}[(a)]
\item Calcular el valor de producci\'on de 1kg de papas fritas (para obtenerlo, a\~nadir al precio
de la mezcla de papa y aceite un 20\% por costos de pago de sueldos, impuestos, mantenimiento, etc).
\item Una bolsa de papas fritas cuesta en un supermercado unos 9\$ y contiene 40g. ?`Crees que es un precio razonable?
\item Se llama \emph{markup}\index{Markup} de un producto al porcentaje a\~nadido a su valor de producci\'on al momento de la venta.
Suele considerarse normal un markup de un 20\% (v\'ease por ejemplo el sitio \href{http://smallbusiness.chron.com/average-profit-margin-wholesale-12941.html}{\tt smallbusiness.chron.com}). ?`Cu\'al es el markup de las papas fritas?
\item El precio de las papas fritas suele permanecer invariable por a\~nos mientras que el de los ingredientes tiene
fuertes oscilaciones. ?`Puedes explicar por qu\'e es as\'i bas\'andote en el resultado del apartado precedente?
\end{enumerate}

\item Una aleaci\'on es la mezcla de dos o m\'as metales. En joyer\'ia se suelen hacer aleaciones de un metal \emph{fino}
(como oro, plata, titanio) con otro \emph{ordinario} (cobre, n\'iquel) y se llama \emph{ley de la aleaci\'on}\index{Ley!de una aleaci\'on} a la
proporci\'on de metal fino que contiene. Por ejemplo, una aleaci\'on de plata de ley 0.12 quiere decir que tiene 12g de plata
por cada 100g de aleaci\'on. Normalmente, como unidad de ley se suele usar el \emph{quilate}: un quilate\index{Quilate} es 1/24 del peso
total, de manera que un anillo de oro de 15 quilates tiene 15/24 del peso total formado de oro. Si el anillo pesa
5g, el contenido en oro es $5\cdot 15/24 =3{.}125$g. El oro puro, naturalmente, es de 24 quilates.
\begin{enumerate}[(a)]
\item Se tienen dos aleaciones de plata, de leyes respectivas 0.75 y 0.4. ?`Qu\'e cantidad de cada una es necesaria para obtener
1kg de aleaci\'on de ley 0.5?
\item Una cadena de oro de 14 quilates contiene 12g de oro puro. ?`Cu\'anto pesa la cadena?
\end{enumerate}

\item Una b\'usqueda por internet nos convencer\'a enseguida de que hay gente dispuesta
a vender (y a comprar) oro de 32 quilates. Existe incluso una banda de rap que afirma valer m\'as que oro de 32 quilates. Pero no hay que pensar que el `oro de 32K' s\'olo existe
en la mente de habitantes del tercer mundo. En el primero, incluso hay de 40K (la captura muestra la p\'agina web de la ciudad de Colonia\index{Colonia, Alemania}, en Alemania),
\begin{center}
 \includegraphics[scale=0.5]{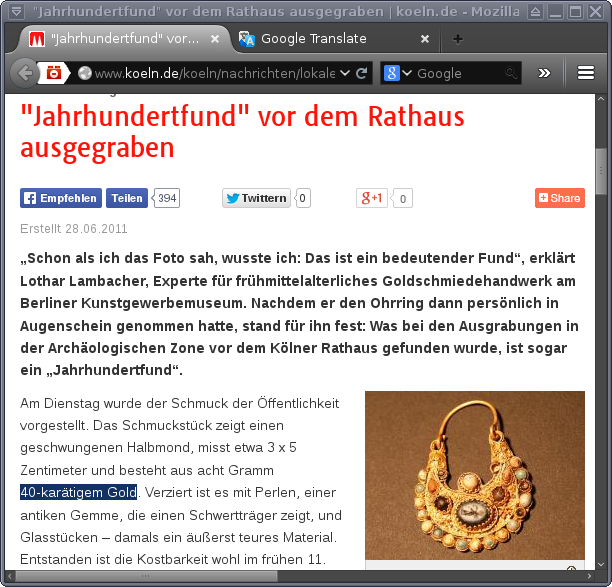} 
 \end{center}
La traducci\'on de Google muestra de qu\'e se trata: un supuesto hallazgo arqueol\'ogico.
\begin{center}
\includegraphics[scale=0.5]{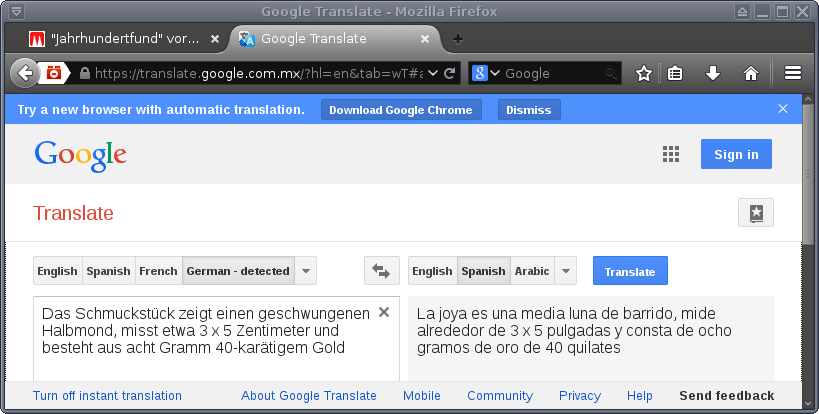} 
\end{center}
\end{enumerate}

\section{Uso de la tecnolog\'ia de c\'omputo}

Como buen programa de c\'alculo simb\'olico, Maxima puede determinar los factores 
primos\index{Factores primos} de n\'umeros relativamente peque\~nos sin problema alguno:
\begin{maxcode}
\noindent
\begin{minipage}[t]{8ex}{\color{red}\bf
\begin{verbatim}
(%i1) 
\end{verbatim}
}
\end{minipage}
\begin{minipage}[t]{\textwidth}{\color{blue}
\begin{verbatim}
factor(2475);
\end{verbatim}
}
\end{minipage}
\definecolor{labelcolor}{RGB}{100,0,0}
\begin{math}\displaystyle
\parbox{8ex}{\color{labelcolor}(\%o1) }
{3}^{2}\,{5}^{2}\,11
\end{math}

\noindent
\begin{minipage}[t]{8ex}{\color{red}\bf
\begin{verbatim}
(%i2) 
\end{verbatim}
}
\end{minipage}
\begin{minipage}[t]{\textwidth}{\color{blue}
\begin{verbatim}
factor(639287400183625434237847625432180);
\end{verbatim}
}
\end{minipage}
\definecolor{labelcolor}{RGB}{100,0,0}
\begin{math}\displaystyle
\parbox{8ex}{\color{labelcolor}(\%o2) }
{2}^{2}\,5\,19\,1998643\,841738751563495613665777
\end{math} 
\end{maxcode}

\noindent Si s\'olo nos interesa determinar si un n\'umero dado es primo o no
(sin importar la determinaci\'on de los factores), se puede hacer m\'as r\'apidamente,
sobre todo si el n\'umero es grande, mediante el uso del comando \texttt{primep}:
\begin{maxcode}
\noindent
\begin{minipage}[t]{8ex}{\color{red}\bf
\begin{verbatim}
(%i3) 
\end{verbatim}
}
\end{minipage}
\begin{minipage}[t]{\textwidth}{\color{blue}
\begin{verbatim}
primep(113);
\end{verbatim}
}
\end{minipage}
\definecolor{labelcolor}{RGB}{100,0,0}
\begin{math}\displaystyle
\parbox{8ex}{\color{labelcolor}(\%o3) }
true
\end{math}

\noindent
\begin{minipage}[t]{8ex}{\color{red}\bf
\begin{verbatim}
(%i4) 
\end{verbatim}
}
\end{minipage}
\begin{minipage}[t]{\textwidth}{\color{blue}
\begin{verbatim}
primep(2475);
\end{verbatim}
}
\end{minipage}
\definecolor{labelcolor}{RGB}{100,0,0}
\begin{math}\displaystyle
\parbox{8ex}{\color{labelcolor}(\%o4) }
false
\end{math}
\end{maxcode}

En ingl\'es, `m\'inimo com\'un m\'ultiplo' se dice `least common multiple', y
`m\'aximo com\'un divisor' es `greatest common divisor'. Las iniciales dan los
comandos correspondientes:
\begin{maxcode}
\noindent
\begin{minipage}[t]{8ex}{\color{red}\bf
\begin{verbatim}
(%i5) 
\end{verbatim}
}
\end{minipage}
\begin{minipage}[t]{\textwidth}{\color{blue}
\begin{verbatim}
lcm(1620,1575);
\end{verbatim}
}
\end{minipage}
\definecolor{labelcolor}{RGB}{100,0,0}
\begin{math}\displaystyle
\parbox{8ex}{\color{labelcolor}(\%o5) }
56700
\end{math}

\noindent
\begin{minipage}[t]{8ex}{\color{red}\bf
\begin{verbatim}
(%i6) 
\end{verbatim}
}
\end{minipage}
\begin{minipage}[t]{\textwidth}{\color{blue}
\begin{verbatim}
gcd(1620,1575);
\end{verbatim}
}
\end{minipage}
\definecolor{labelcolor}{RGB}{100,0,0}
\begin{math}\displaystyle
\parbox{8ex}{\color{labelcolor}(\%o6) }
45
\end{math}
\end{maxcode}

Finalmente, se\~nalemos que hay muchas otras funciones para trabajar con n\'umeros,
que se pueden consultar en la ayuda del programa. Una que es particularmente \'util
es \texttt{next\_prime} que, como su nombre indica,
proporciona el primer n\'umero primo a continuaci\'on de uno dado:
\begin{maxcode}
\noindent
\begin{minipage}[t]{8ex}{\color{red}\bf
\begin{verbatim}
(%i7) 
\end{verbatim}
}
\end{minipage}
\begin{minipage}[t]{\textwidth}{\color{blue}
\begin{verbatim}
next_prime(18376453728920832746324838939938376346565);
\end{verbatim}
}
\end{minipage}
\definecolor{labelcolor}{RGB}{100,0,0}
\begin{math}\displaystyle
\parbox{8ex}{\color{labelcolor}(\%o7) }
18376453728920832746324838939938376346577
\end{math}
\end{maxcode}

\section{Actividades de c\'omputo}
\begin{enumerate}[(A)]
\item Utiliza Maxima para comparar los tiempos requeridos para factorizar el n\'umero
$$
7812353094515203333466690911342264649370470629995425393793599374448576999
$$
mediante el algoritmo de ensayo y error y el de Lenstra. ?`Cu\'anto tiempo necesita el
algoritmo de Lenstra para factorizar un n\'umero de $200$ d\'igitos? (compara la respuesta con el
tiempo que toma el algoritmo de ensayo y error, visto en la p\'agina \pageref{time}).
Para finalizar, prepara un c\'odigo en Maxima que pueda hacer estas comparaciones para un n\'umero
de d\'igitos $n$ arbitrario.
\item Utilizando la clave p\'ublica de Ale\index{Ale} que se vi\'o en la secci\'on \ref{sec32},
encriptar el siguiente mensaje: \texttt{hoy?}.
\item Suponiendo que conocemos la clave privada de Ale (dada en la secci\'on Aplicaciones),
desencriptar el siguiente mensaje de Gaby\index{Gaby}: \texttt{NY[\_}.
\item ?`C\'omo habr\'ia que modificar el proceso con Maxima para poder escribir cualquier 
mensaje y asegurarnos de que el texto cifrado s\'olo contiene caracteres 
ASCII\index{Co@C\'odigo!ASCII} imprimibles?
\end{enumerate}

\chapter{Inducci\'on. Coeficientes binomiales}\label{cap4}
\vspace{2.5cm}
\section{Nociones b\'asicas}
El \emph{m\'etodo de inducci\'on}\index{Me@M\'etodo!inducci\'on} es el equivalente matem\'atico de la frase `y as\'i sucesivamente'. Proporciona
una manera de determinar propiedades de conjuntos \emph{infinitos}, que no se pueden comprobar elemento a elemento.
La idea es la siguiente: supongamos que tenemos un conjunto infinito numerable $S$, cuyos elementos escribimos como
$S=\{ s_1,s_2,s_3,\ldots \}$ (una notaci\'on m\'as adecuada matem\'aticamente es $(s_n)_{n\in\mathbb{N}}$), y una
propiedad $\mathcal{P}$ tal que
\begin{enumerate}[(a)]
\item Se cumple para un determinado elemento $s_{n_0}$ del conjunto infinito numerable $S$.
\item Cada vez que un elemento $s_n$ del conjunto $S$ cumple $\mathcal{P}$, el siguiente elemento $s_{n+1}$ 
tambi\'en la cumple.
\end{enumerate}
Entonces, podemos asegurar que todos los elementos a partir de $s_{n_0}$ tienen la propiedad $\mathcal{P}$. Si el
elemento $s_{n_0}\in S$ es el primero, es decir $n_0=1$ as\'i que $s_{n_0}=s_1$, entonces podemos asegurar que \emph{todos}
los elementos de $S$ tienen la propiedad. Intuitivamente, el proceso es como la ca\'ida de fichas de domin\'o.
\begin{center}
\includegraphics[scale=0.1]{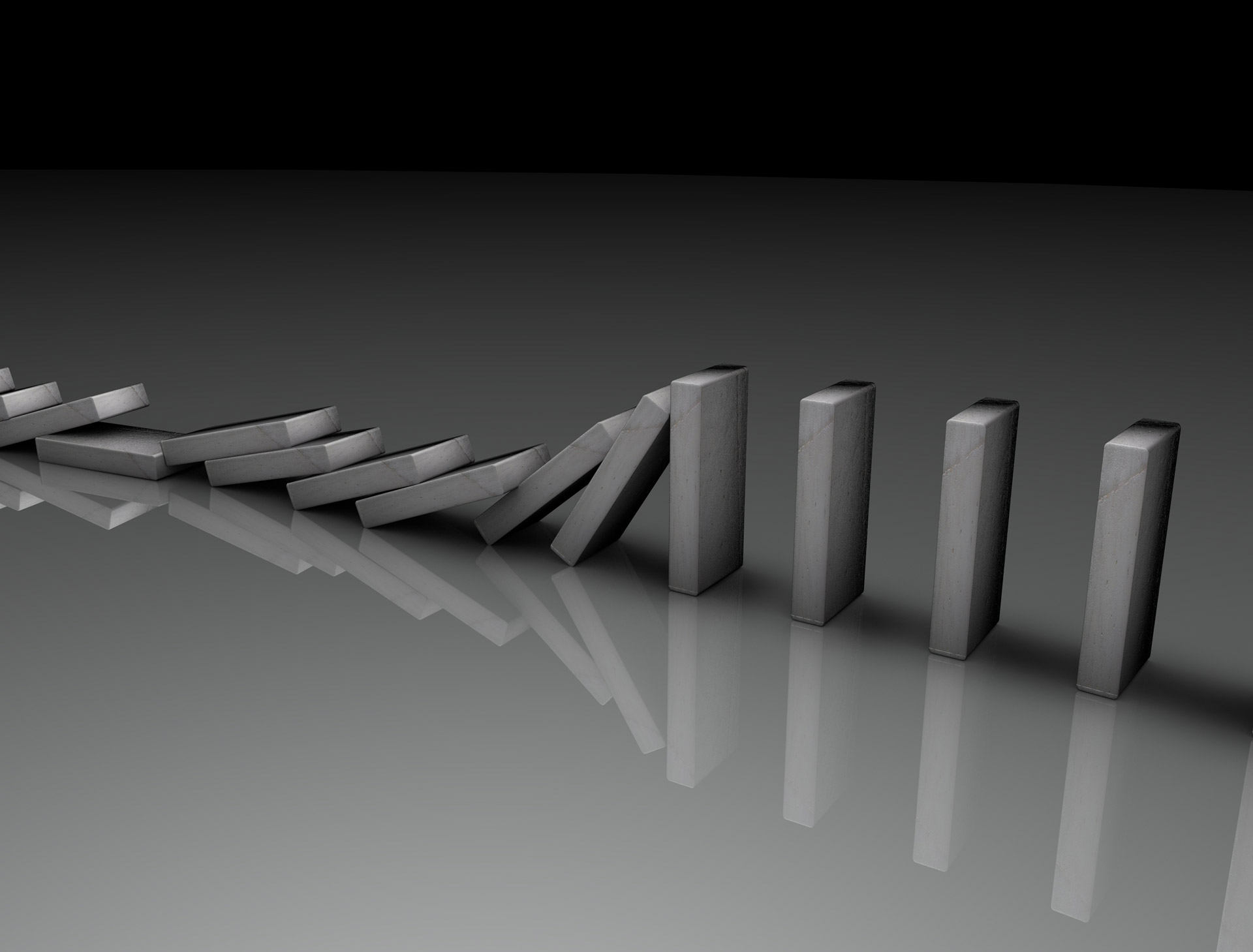} 
\end{center}
Por ejemplo, para probar que sea cual sea $n\in\mathbb{N}$ se cumple
$$
1 + 2 + 3 + \cdots + n = \frac{n(n+1)}{2},
$$
proceder\'iamos como sigue. Comprobamos que la f\'ormula propuesta es cierta para $n=1$, que en este caso es trivial:
$1=1(1+1)/2=2/2$. A continuaci\'on, \emph{suponemos} que la f\'ormula es cierta para $n=k$, es decir, suponemos que
\begin{equation}\label{eq1}
1 + 2 + 3 + \cdots + k =\frac{k(k+1)}{2} ,
\end{equation}
y vemos si tambi\'en es cierta para el siguiente $n=k+1$, lo que equivale a comprobar si es verdad que
\begin{equation}\label{eq2}
1 + 2 + 3 + \cdots + k +(k+1) =\frac{(k+1)(k+2)}{2}.
\end{equation}
Pero sustituyendo \eqref{eq1} en $1 + 2 + 3 + \cdots + k +(k+1)$ resulta:
\begin{align*}
 \color{blue}1 + 2 + 3 + \cdots + k \color{black}+k+1 &= \color{blue} \frac{k(k+1)}{2}\color{black}+k+1 \\
							&= \frac{k(k+1)}{2}+\frac{2(k+1)}{2} \\
							&= \frac{k(k+1)+2(k+1)}{2} \\
							&= \frac{(k+1)(k+2)}{2},
\end{align*}
que es precisamente \eqref{eq2}, lo que quer\'iamos probar.

En los ejercicios, usaremos el m\'etodo de inducci\'on para probar algunas propiedades de los \emph{coeficientes binomiales} y la important\'isima f\'ormula del binomio. Si
$p,q\in\mathbb{N}$ son dos n\'umeros naturales, se define $1!=1=0!$' el coeficiente binomial\index{Coeficiente binomial} de $p$ sobre $q$ es
$$
{p\choose q}=\frac{p!}{q!(p-q)!}.
$$
Por ejemplo:
$$
{7\choose 3}=\frac{7!}{4!(7-4)!}=\frac{7\cdot 6\cdot 5\cdot 4!}{4!\cdot 3!}
=\frac{7\cdot 6\cdot 5}{3!}=\frac{7\cdot 6\cdot 5}{3\cdot 2\cdot 1}=7\cdot 5=35\,.
$$

\section{Aplicaciones}
Una de las ramas de las matem\'aticas con mayores aplicaciones a nuestra vida cotidiana
es, sin lugar a dudas, la \emph{teor\'ia de la probabilidad}. Constantemente surgen
preguntas como ?`qu\'e tantas oportunidades de ganar tiene mi equipo favorito?, 
?`si voy al centro comercial $X$, me encontrar\'e con $Y$, a quien debo dinero? O incluso
otras de mayor trascendencia, del estilo ?`cu\'antas oportunidades de sobrevivir tiene
una persona con c\'ancer que se someta al tratamiento $Z$?\\
La teor\'ia de la probabilidad\index{Probabilidad} trata de contestar a estar preguntas \emph{cuantitativamente},
asignando un n\'umero entre $0$ y $1$ a cada caso que pueda darse en las situaciones
anteriores, de manera que cuanto m\'as pr\'oximo a $0$ est\'e el n\'umero asignado a un caso
concreto, menores oportunidades tendr\'a de ocurrir (o cuanto m\'as pr\'oximo est\'e a $1$,
mayores oportunidades tendr\'a de ocurrir). Para ello, la probabilidad cl\'asica se basa
en el \emph{c\'alculo de frecuencias}: la probabilidad de que un determinado suceso $S$ ocurra
se calcula como
$$
P(S)=\frac{\# \mbox{ casos favorables}}{\# \mbox{ total de casos}}\,,
$$
donde en el denominador se cuentan \emph{todos} los casos, favorables o no. 
Por ejemplo, si lanzamos un dado legal\footnote{Es decir, que no est\'a trucado: 
todas las caras se han constru\'ido con el mismo material, con las mismas dimensiones
y con la misma densidad.}
al aire, la probabilidad de que salga un $6$ es, obviamente,
$$
P(\mbox{sale }6)=\frac{1}{6}\,.
$$
Otra caracter\'istica fundamental de la probabilidad, es que la suma de todas las
probabilidades
que asignemos a \emph{todos} los posibles resultados, debe ser $1$ (es decir, la probabilidad
de que `ocurra algo' siempre es $1$). Entonces, la probabilidad de que \emph{no ocurra}
el suceso $S$, denotada $P(\overline{S})$, es $P(\overline{S})=1-P(S)$. En este ejemplo,
la probabilidad de que salga cualquier n\'umero excepto el $6$ es $5/6$.\\
Por otra parte, si lanzamos \emph{dos} dados, y nos preguntamos por la probabilidad de que la 
suma
sea $5$, la respuesta ya no es tan obvia. Podemos hacer una tabla con los posibles resultados
de los lanzamientos y contar cu\'antos de ellos suman $5$:
\begin{center}
\begin{tabular}{|c|c|c|c|c|c|}
\hline 
$(1,1)$ & $(1,2)$ & $(1,3)$ & \cellcolor{Lavender!75} $(1,4)$ & $(1,5)$ & $(1,6)$ \\ 
\hline 
$(2,1)$ & $(2,2)$ & \cellcolor{Lavender!75} $(2,3)$ & $(2,4)$ & $(2,5)$ & $(2,6)$ \\ 
\hline 
$(3,1)$ & \cellcolor{Lavender!75}$(3,2)$ & $(3,3)$ & $(3,4)$ & $(3,5)$ & $(3,6)$ \\ 
\hline 
\cellcolor{Lavender!75}$(4,1)$ & $(4,2)$ & $(4,3)$ & $(4,4)$ & $(4,5)$ & $(4,6)$ \\ 
\hline 
$(5,1)$ & $(5,2)$ & $(5,3)$ & $(5,4)$ & $(5,5)$ & $(5,6)$ \\ 
\hline 
$(6,1)$ & $(6,2)$ & $(6,3)$ & $(6,4)$ & $(6,5)$ & $(6,6)$ \\ 
\hline 
\end{tabular} 
\end{center}
Vemos que la probabilidad es ahora menor,
$$
P(\mbox{sale }5\mbox{ en }2\mbox{ dados})=\frac{4}{36}=\frac{1}{9}\,.
$$
Otro ejemplo lo tenemos en el lanzamiento al aire de una moneda (volados). Si lanzamos
$30$ veces una moneda, ?`cu\'al es la probabilidad de que salgan exactamente $13$ \'aguilas
(o cruces, o \emph{tails}, en ingl\'es). Una manera equivalente de preguntar esto es decir:
de entre todas las cadenas formadas por $0$ y $1$ con longitud $30$, ?`cu\'antas hay que
contengan exactamente $13$ unos? Como vemos, el c\'alculo de probabilidades nos lleva de
manera natural al terreno de la \emph{combinatoria}, el arte de contar. A continuaci\'on, 
veremos las f\'ormulas cl\'asicas referentes al n\'umero de maneras de combinar objetos.

Para empezar, un poco de terminolog\'ia: una \emph{combinaci\'on}\index{Combinaci\'on} de los objetos
A, B, C,..., es simplemente una lista conteniendo estos objetos, sin importar su orden.
Por ejemplo, en una barra de ensaladas nuestra preferencia podr\'ia ser una combinaci\'on
de lechuga, tomate, ma\'iz, cebolla y at\'un. No importa el orden. Sin embargo, al dar
el n\'umero de una casa lo que se necesita es una \emph{permutaci\'on}, que es una
combinaci\'on \emph{ordenada}\index{Permutaci\'on}. En este caso, una permutaci\'on de los d\'igitos
$0,1,2,\ldots ,9$, donde la permutaci\'on $524$ no equivale a la $245$.\\
Hay dos tipos de permutaciones: \emph{con repetici\'on} y \emph{sin repetici\'on}.
Por ejemplo, al dar el n\'umero de una casa admitimos el $244$, mientras que al ordenar
las clases de las materias `c\'alculo', `\'algebra' y `computaci\'on', no se permiten
repeticiones (no se puede tener la misma clase dos o m\'as veces).
\begin{list}{$\blacktriangleright$}{}
\item C\'alculo de permutaciones con repetici\'on. Si tenemos $n$ objetos para elegir
y podemos tomar $k$ de ellos con repetici\'on, ?`de cu\'antas formas lo podemos hacer?
El primer objeto lo podemos elegir de $n$ formas. Como hay repetici\'on, el segundo
objeto \emph{tambi\'en} lo podemos elegir de $n$ formas, independientemente de c\'omo
hayamos elegido el segundo. Repitiendo el argumento $k$ veces, tenemos que
$$
PR^n_k=n^k\,.
$$
Por ejemplo, ?`cu\'antas claves se pueden formar con un candado de $3$ ruedas? Cada rueda 
permite elegir un n\'umero de entre $0,1,2,\ldots ,9$ y hay $3$ de ellas, luego ser\'an
$PR^{10}_3=10^3=1\,000$ claves.
\item C\'alculo de permutaciones sin repetici\'on. En este caso se reducen las posibilidades,
una vez que una de ellas ocurre, no puede volver a darse. Entonces, si tenemos $n$ objetos de
entre los que hay que elegir $k$ sin repetirse, el primero lo podremos tomar de $n$ formas 
distintas, pero el segundo ya s\'olo de $n-1$ formas independientemente del primero, el tercero 
de $n-2$ formas independientemente de los otros dos, y as\'i
hasta el \'ultimo elemento, el $k-$\'esimo, que s\'olo se podr\'a elegir de $n-(k-1)$ maneras.
El recuento de posibilidades proporciona, por tanto,
$$
P^n_k=n\cdot (n-1)\cdot (n-(k-1))\,.
$$
Esto se puede expresar c\'omodamente mediante el uso de factoriales, pues no es 
m\'as que
$$
P^k_n=\frac{n!}{(n-k)!}\,.
$$
Por ejemplo, las placas de autom\'oviles contienen $3$ letras, que no pueden repetirse. 
?`Cu\'antas posibilidades hay para la terna? Dado que tenemos $26$ letras y se eligen $3$ sin 
repetici\'on, ser\'an
$$
P^{26}_3=\frac{26!}{23!}=26\cdot 25\cdot 24 = 15\,600\,.
$$
\end{list}
Tambi\'en con las combinaciones tenemos dos posibilidades. Es en estos casos en los que
aparecen los coeficientes binomiales\index{Coeficiente binomial}.
\begin{list}{$\blacktriangleright$}{}
\item C\'alculo de combinaciones sin repetici\'on. Se pueden contar a partir del 
correspondiente n\'umero de permutaciones, donde \emph{s\'i} importa el orden, haciendo
luego que el orden no importe. Entonces, si tenemos $n$ objetos de entre los que hay que
elegir $k$, sin que importe el orden y sin repetici\'on, podemos calcular previamente las
posibilidades en las que importa el orden, que ya conocemos
$$
\frac{n!}{(n-k)!}\,.
$$
Ahora, en cada una de esas posibilidades tenemos $k$ objetos \emph{ordenados}, pero en realidad
no nos importa el orden, de modo que estamos contando varias veces posibilidades que son la
misma. Por ejemplo, formando combinaciones de $3$ letras, nos da lo mismo cualquiera de las 
$3!$ permutaciones $UVX$, $UXV$, $XUV$, $VUX$, $VXU$, $XVU$. Si son $k$ objetos en cada 
permutaci\'on, hay $k!$ maneras de ordenarlos, todas ellas
equivalentes para formar las combinaciones. Por tanto, \'ese es el n\'umero por el que hay que
dividir para contar las combinaciones \emph{diferentes}. Llegamos as\'i a
$$
C^n_k=\frac{n!}{k!\cdot (n-k)!}={n\choose k}\, .
$$
Por ejemplo, en un popular juego de loter\'ia hay que elegir $6$ n\'umeros de entre $49$
(los que van de $1$ a $49$, de hecho)\footnote{En M\'exico\index{Me@M\'exico!loter\'ia}, la versi\'on de esta loter\'ia
funciona eligiendo $6$ n\'umeros de entre $56$.}. Si cada apuesta cuesta $\$10$, ?`cu\'anto
habr\'ia que gastar para asegurarse de sacar una combinaci\'on ganadora?
\begin{center}
\includegraphics[scale=0.275]{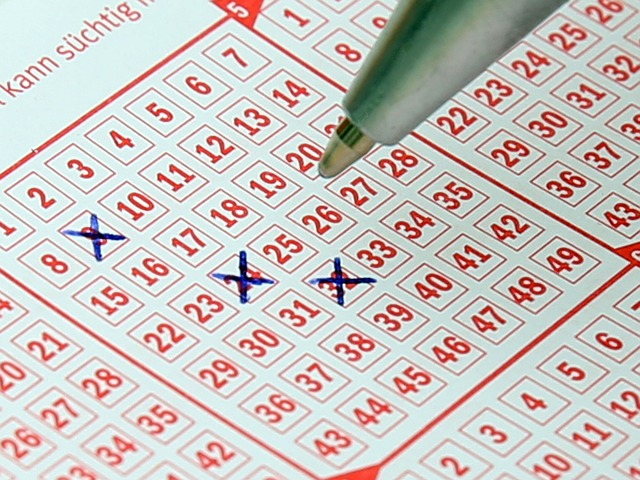}
\end{center}
Calculamos las
posibles combinaciones:
$$
C^{49}_6={49\choose 6}=\frac{49!}{6!43!}=
\frac{49\cdot 48\cdot 47\cdot 46\cdot 45\cdot 44}{6\cdot 5\cdot 4\cdot 3\cdot 2}=
 13\,983\,816\,.
$$
Para asegurarnos de que ganamos, deber\'iamos cubrir todas las posibilidades, lo que costar\'ia
$$
13\,983\,816\cdot \$10 =\$139\,838\,160\,.
$$
La probabilidad de ganar haciendo una sola apuesta es, como vimos, el n\'umero de casos
favorables (s\'olo uno) entre el n\'umero de casos posibles, es decir,
$$
P(\mbox{ganar})=\frac{1}{13\,983\,816}=0{.}00000007\,.
$$
\end{list}

\begin{list}{$\blacktriangleright$}{} 
\item C\'alculo de combinaciones con repetici\'on. \'Este es el caso m\'as complicado de todos.
Sin embargo, tambi\'en existe una t\'ecnica de c\'alculo. Para ilustrarla, recurriremos al
famoso \emph{problema de los helados de sabores}. Supongamos que estamos en una helader\'ia
donde disponen de helados de cinco sabores: nata, vainilla, chocolate, fresa y lim\'on. Los
denotamos por las letras N, V, C, F y L, respectivamente. Podemos pedir un cono con tres bolas 
de helado, que podr\'ian estar repetidas (tres bolas de fresa, por ejemplo). ?`Cu\'antas
combinaciones de sabores son posibles? Tenemos $n=5$ sabores para elegir $k=3$ de ellos, sin
que importe el orden y con posibilidad de repetici\'on. Si el mostrador de la helader\'ia
tiene los sabores en recipientes ordenados de izquierda a derecha, por ejemplo, en el orden
de arriba (NVCFL), una elecci\'on concreta como `$1$ bola de fresa y dos de vainilla' podr\'ia
describirse como una sucesi\'on de $1$s y $0$s, donde se necesitan $n-1=4$ unos para hacer de
`separadores' y $k=3$ ceros para caracterizar los sabores elegidos
\begin{center}
\begin{tabular}{c|c|c|c|c}
\hline 
N & V & C & F & L \\ 
\hline 
• & $\circ\circ$ & • & $\circ$ & • \\ 
\end{tabular} 
\end{center}
Es decir, la elecci\'on `$1$ bola de fresa y dos de vainilla' se representa por la
combinaci\'on \emph{ordenada}
$$
1\,0\,0\,1\,1\,0\,1\,.
$$
Cualquier otra elecci\'on se expresar\'a como una combinaci\'on ordenada de los $n-1=4$
separadores ($1$s) y las $k=3$ bolas ($0$s). Es decir, tenemos cadenas de $0$s y $1$s de
longitud $n+k-1$, y queremos todos los casos en que $k$ s\'imbolos son $0$.
De acuerdo con lo que hemos visto en el caso precedente, el n\'umero de combinaciones en
general ser\'a
$$
CR^n_k={n+k-1\choose k}\, .
$$
En el ejemplo, tendr\'iamos
$$
CR^5_3={7\choose 3}=\frac{7\cdot 6\cdot 5\cdot 4!}{3!\cdot 4!}=\frac{7\cdot 6\cdot 5}{6}=35\,.
$$
\end{list}

Ahora que sabemos contar, vamos c\'omo aplicar el conocimiento adquirido al c\'alculo
de probabilidades\index{Probabilidad}. En particular, vamos a fijarnos en los llamados \emph{procesos de
Bernoulli}\index{Bernoulli, Jacob}. Son procesos\footnote{O, en la terminolog\'ia probabilista, \emph{experimentos}\index{Experimento}.}
en los que el resultado\footnote{En la terminolog\'ia probabilista, 
\emph{sucesos}\index{Suceso!probabilidad}}
s\'olo puede tomar dos valores: cara o cruz, \'aguila o sello, 
$0$ o $1$, ni\~na o ni\~no, etc. Nosotros utilizaremos los s\'imbolos $+$ y $-$.
Adem\'as, en cada ensayo (sea \'este el lanzamiento de
una moneda, un nacimiento, un d\'igito impreso en una pantalla, etc.) la probabilidad
de que ocurra el suceso $+$ es $p$ y la de que ocurra el suceso $-$ es $q$, con independencia
de lo que suceda en los restantes ensayos. En este tipo de procesos, al no haber otras 
posibilidades, debe cumplirse
$$
p+q=1\,.
$$
Por ejemplo, si el experimento consiste en lanzar una moneda los sucesos
posibles tienen
$$
P(+)=\frac{1}{2}=P(-)\,.
$$
Si el experimento consiste en lanzar un dado y el suceso $+$ corresponde a `sale
un $6$', entonces $-$ es el suceso `sale cualquier n\'umero de entre $1,2,3,4,5$'
y
$$
P(+)=\frac{1}{6},\,P(-)=\frac{5}{6}\,.
$$
En el caso general, como hemos mencionado,
$$
P(+)=p,\,P(-)=q\,.
$$
Podemos definir una funci\'on que nos d\'e la probabilidad de que ocurra el suceso
$+$ o el suceso $-$. Sea $X$ el n\'umero de veces que ocurre el evento $+$. Si el
experimento s\'olo se hace una vez, claramente $X$ s\'olo puede tomar los valores
$X=x\in\{0,1\}$, con $P(X=0)=q$ (si $X=0$ es que no sali\'o $+$) y $P(X=1)=p$. La
funci\'on que describe la probabilidad de que ocurra cada uno de los valores que
puede tomar $X$, claramente se puede escribir como
$$
f(x)=P(X=x)=p^xq^{1-x}\,
$$
pues, efectivamente, $f(0)=q$ y $f(1)=p$. A una funci\'on de los resultados del
experimento, como $X$, se la denomina \emph{variable aleatoria}\index{Variable aleatoria} y a la funci\'on
$f(x)=P(X=x)$ se la denomina su \emph{funci\'on de probabilidad}\index{Probabilidad!funci\'on}.\\
Si el experimento se realiza $n$ veces, la variable aleatoria $X$ podr\'a tomar
los valores $0,1,\ldots ,n$. Si, por ejemplo, es $n=10$ y $X=4$, quiere decir que el suceso
$+$ ocurri\'o cuatro veces de un total de $10$. El proceso completo pudo haber dado como
resultados
$$
+-++---+--
$$
o quiz\'as
$$
++++------\,.
$$
De hecho, como s\'olo importa el n\'umero de veces que apareci\'o $+$ y no el orden,
el n\'umero de resultados con cuatro apariciones de $+$ es igual a
$$
{10\choose 4}\,.
$$ 
Cada uno de esos resultados conteniendo cuatro $+$, tiene una probabilidad de ocurrir dada por
$p^4q^{10-4}$ (pues cada experimento es independiente de los dem\'as y en \'el la
probabilidad de $+$ es $p$ y la de $-$ es $q$. La probabilidad de cuatro eventos $+$
es $p^4$ y la de $10-4$ eventos $-$ es $q^{10-4}$). Por tanto, la probabilidad \emph{total}
de que salgan cuatro $+$ es igual al producto del n\'umero de resultados en los que hay cuatro 
$+$ por la probabilidad de uno de esos resultados, es decir,
$$
f(4)=P(X=4)={10\choose 4}\,p^4q^{10-4}\,.
$$
El mismo razonamiento para un $X=x\in\{0,1,\ldots,n\}$ arbitrario, muestra que esta
probabilidad es
$$
f(x)=P(X=x)={n\choose x}\,p^xq^{n-x}\,.
$$

Podemos entonces contestar a una de las preguntas del inicio de secci\'on:
si lanzamos $30$ veces una moneda, ?`cu\'al es la probabilidad de que salgan
exactamente $13$ \'aguilas? Llamando $+$ al evento `sale \'aguila', tenemos
que (si la moneda no est\'a trucada),  $P(+)=1/2=P(-)$ . Si $X$ es la variable aleatoria
`n\'umero de \'aguilas en $30$ lanzamientos', lo que queremos calcular es
$$
f(13)=P(X=13)={30\choose 13}\left(\frac{1}{2}\right)^{13}\left(\frac{1}{2}\right)^{30-13}
={30\choose 13}\left(\frac{1}{2}\right)^{30}=0{.}1115\,.
$$

\section{Ejercicios}

\begin{enumerate}  
\setcounter{enumi}{40}

\item Demostrar por inducci\'on la validez de las siguientes f\'ormulas sumatorias para
$n \in \mathbb{N}$: 
\begin{enumerate}[(a)]
\item $1 + 3 + 5 + \cdots + (2n-1) = n^2$
\item  $1 + 2 + 2^2 + 2^3 + \cdots + 2^{n-1} = 2^n-1$
\item  $\displaystyle 1 + 5 + 5^2 + \cdots + 5^{n-1} = \frac{5^n-1}{4}$
\item  $\displaystyle 1^2 + 2^2 + 3^2 + \cdots + n^2 = \frac{n(n+1)(2n+1)}{6}$
\item  $\displaystyle \frac{1}{1\cdot 2} + \frac{1}{2\cdot 3} + \frac{1}{3\cdot 4} + \cdots + \frac{1}{n(n+1)} = \frac{n}{n+1}$
\item  $\displaystyle \frac{1}{1\cdot 2\cdot 3} + \frac{1}{2\cdot 3\cdot 4} + \cdots + \frac{1}{n(n+1)(n+2)} = \frac{n(n+3)}{4(n+1)(n+2)}$
\item  $\displaystyle \frac{1}{1\cdot 3} + \frac{1}{3\cdot 5} + \cdots + \frac{1}{(2n-1)(2n+1)} = \frac{n}{2n+1}$
\item  $1\cdot 4 + 4\cdot 7 + \cdots + (3n-2)(3n+1) = n(3n^2+3n-2)$
\item  $\displaystyle 1 + x + x^2 + \cdots + x^{n-1} = \frac{x^n-1}{x-1}$, $x \in \mathbb{R}-\{ 1 \}$
\item  $\displaystyle 1^4 + 2^4 + 3^4 + \cdots + n^4 = \frac{6n^5+15n^4+10n^3-n}{30}$ 
\end{enumerate}

\item\label{derivadas} Demostrar por inducci\'on la validez de las siguientes f\'ormulas para 
$n \in \mathbb{N}$: 
\begin{multicols}{3}
\begin{enumerate}[(a)]
\item  $n! \leq n^n$
\item  $n < 2^n$
\item  $\dfrac{\mathrm{d}x^n}{\mathrm{d}x} = n x^{n-1}$
\end{enumerate}
\end{multicols}

\item Calcular los siguientes coeficientes binomiales\index{Coeficiente binomial}:
\begin{multicols}{4}
\begin{enumerate}[(a)]
\item \itemequ{{7 \choose 2}}
\item \itemequ{{14 \choose 0}}
\item \itemequ{{14 \choose 13}} 
\item \itemequ{{126 \choose 125}} 
\end{enumerate}
\end{multicols}

\item Demostrar las siguientes propiedades de los coeficientes binomiales:
\begin{enumerate}[(a)]
\item  $\displaystyle {n \choose 0} + {n \choose 1} + \cdots + {n \choose n} = 2^n$ 
\item  $\displaystyle{n \choose 0} + {n \choose 2} + {n \choose 4} + \cdots  = 2^{n-1}$
\item  $\displaystyle{n \choose 1} + {n \choose 3} + {n \choose 5} + \cdots  = 2^{n-1}$
\end{enumerate}

\item Probar las siguientes propiedades:
\begin{multicols}{2}
\begin{enumerate}[(a)]
\item \itemequ{{p \choose q} = {p \choose p-q}}
\item\label{formulab} \itemequ{{p \choose q-1}+{p \choose q} = {p+1 \choose q}}
\end{enumerate}
\end{multicols}
La segunda se conoce como la \emph{regla de Pascal}\index{Pascal, Blaise}, utilizada para construir el famoso
\emph{tri\'angulo de Pascal}\index{Pascal, Blaise!tri\'angulo}. \'Este se forma de la siguiente manera: en la cima de una
pir\'amide (nivel $0$) se coloca el n\'umero $1$. En el siguiente nivel (el $1$), se
colocan dos copias del n\'umero $1$. A partir de aqu\'i, se van a\~nadiendo filas cada
una con un lugar m\'as que la anterior, y de manera que cada cifra es la suma de las dos
que tiene inmediatamente sobre ella, excepto en los extremos, donde se repite la cifra
inmediatamente superior. Es decir, empezar\'iamos con
\begin{table}[H]
\begin{center}
\begin{tabular}{c|ccc}
 $0$ &  &$1$&  \\
 $1$ & $1$& &$1$
\end{tabular}
\end{center}
\end{table}
y a\~nadir\'iamos una fila de tres elementos:
\begin{table}[H]
\begin{center} \begin{tabular}{c|ccccc}
 $0$ &   &   &$1$&   & \\
 $1$ &   &$1$&   &$1$&\\
 $2$ &$1$&   &$a$&   &$1$
\end{tabular} \end{center}
\end{table}
El coeficiente $a$ debe ser la suma de los dos que est\'an sobre \'el (a izquierda y derecha),
es decir, $1+1=2$. As\'i:
\begin{table}[H]
\begin{center} \begin{tabular}{c|ccccc}
 $0$ &   &   &$1$&   & \\
 $1$ &   &$1$&   &$1$&\\
 $2$ &$1$&   &$2$&   &$1$
\end{tabular} \end{center}
\end{table}
Continuando el procedimiento, se obtiene
\begin{table}[H]
\begin{center} 
\begin{tabular}{c|cccccccccccccccc}
 $0$ & & & & & & & &$1$& & & & & & & & \\
 $1$ & & & & & & &$1$& &$1$& & & & & & & \\
 $2$ & & & & & &$1$& &$2$& &$1$& & & & & \\
 $3$ & & & & &$1$& &$3$& &$3$& &$1$& & & \\
 $4$ & & & &$1$& &$4$& &$6$& &$4$& &$1$& \\
 $5$ & & &$1$& &$5$& &$10$& &$10$& &$5$& &$1$ \\
 $6$ & &$1$& &$6$& &$15$& &$20$& &$15$& &$6$& &$1$ \\
 $7$ & $1$& &$7$& &$21$& &$35$& &$35$& &$21$& &$7$& &$1$ 
\end{tabular} \end{center}
\end{table}
De acuerdo con la f\'ormula \eqref{formulab} del ejercicio, el tri\'angulo de Pascal se
puede escribir como
\begin{table}[H]
\begin{center} 
\begin{tabular}{c|cccccccccccccccc}
$0$ & & & & & & & &$\displaystyle {0\choose 0}$& & & & & & & & \\
$1$ & & & & & & &$\displaystyle {1\choose 0}$& &$\displaystyle {1\choose 1}$& & & & & & & \\
$2$ & & & & & &$\displaystyle {2\choose 0}$& &$\displaystyle {2\choose 1}$& &$\displaystyle {2\choose 2}$& & & & & \\
$3$ & & & & &$\displaystyle {3\choose 0}$& &$\displaystyle {3\choose 1}$& &$\displaystyle {3\choose 2}$& &$\displaystyle {3\choose 3}$& & & \\
$4$ & & & &$\displaystyle {4\choose 0}$& &$\displaystyle {4\choose 1}$& &$\displaystyle {4\choose 2}$& &$\displaystyle {4\choose 3}$& 
&$\displaystyle {4\choose 4}$&\\
$5$ & & &$\displaystyle {5\choose 0}$& &$\displaystyle {5\choose 1}$& &$\displaystyle {5\choose 2}$& &$\displaystyle {5\choose 3}$& 
&$\displaystyle {5\choose 4}$& &$\displaystyle {5\choose 5}$ \\
$6$ & &$\displaystyle {6\choose 0}$& &$\displaystyle {6\choose 1}$& &$\displaystyle {6\choose 2}$& &$\displaystyle {6\choose 3}$& &$\displaystyle {6\choose 4}$& 
&$\displaystyle {6\choose 5}$& &$\displaystyle {6\choose 6}$ 
\end{tabular} \end{center}
\end{table}
El tri\'angulo de Pascal presenta much\'isimas peculiaridades. Por ejemplo, la suma de
cada fila es una potencia de dos. La primera diagonal (desde abajo a la izquierda)
est\'a compuesta de unos, la segunda comprende los naturales, la tercera los n\'umeros 
triangulares, etc, etc.

\item Probar, por inducci\'on, la \emph{f\'ormula del binomio}\index{Fo@F\'ormula del binomio} debida a Newton\index{Newton, Isaac}:
$$
(a+b)^n= \sum^n_{k=0}{n\choose k}a^{n-k}b^k .
$$

\item Probar la siguiente f\'ormula, para $n\geq 2$ (!`no requiere inducci\'on!):
$$
a^n-b^n=(a-b)(a^{n-1}+a^{n-2}b+a^{n-3}b^2+\cdots +a^2b^{n-3}+ab^{n-2}+b^{n-1}).
$$


\item Desarrollar y simplificar las siguientes expresiones:
\begin{multicols}{3}
\begin{enumerate}[(a)]
\item  $$\left( 4x + \frac{1}{4}y\right)^5$$
\item  $$\left( \frac{x}{4y^3} - \frac{2y}{x^2}\right)^5$$
\item $$(a+b-ab)^3\,.$$
\end{enumerate}\noindent
\end{multicols}

\item ?`Cu\'antas contrase\~nas de $8$ caracteres se pueden formar con las
letras del alfabeto, may\'usculas y min\'usculas, con la condici\'on de que no haya
repeticiones?

\item En un servidor de correo electr\'onico, las contrase\~nas se asignan con
dos letras del alfabeto (min\'usculas) seguidas de $3$ d\'igitos (del $0$ al $9$). ?`Cu\'antas
contrase\~nas se pueden generar?

\item Si llamamos `palabra' a cualquier disposici\'on de letras, tenga significado o no,
?`cu\'antas palabras (de cualquier longitud) pueden formarse empleando las letras de la palabra 
`rectas' sin repeticiones?

\item En la fase final del campeonato mundial FIFA de f\'utbol participan $16$ equipos,
que se dividen en cuatro grupos A, B, C y D. Cada equipo juega un partido contra los
dem\'as de su grupo. Tras estos encuentros, el ganador del grupo A juega contra el segundo
clasificado del grupo B y el segundo del grupo A contra el primero del grupo B. An\'alogamente,
juega el ganador del grupo C contra el segundo clasificado del grupo D, y el segundo del grupo
C contra el primero del grupo D. Los ganadores de estos cuatro partidos juegan una semifinal
y los ganadores de la semifinal, la final. ?`Cu\'antos partidos se juegan en total?

\item Un restaurante ofrece $5$ entradas diferentes, $10$ platos principales y $4$
postres diferentes. Si el comensal puede comer un solo tiempo, dos o los tres seg\'un
prefiera, ?`cu\'antas comidas diferentes tiene a su disposici\'on?
\begin{center}
\includegraphics[scale=0.33]{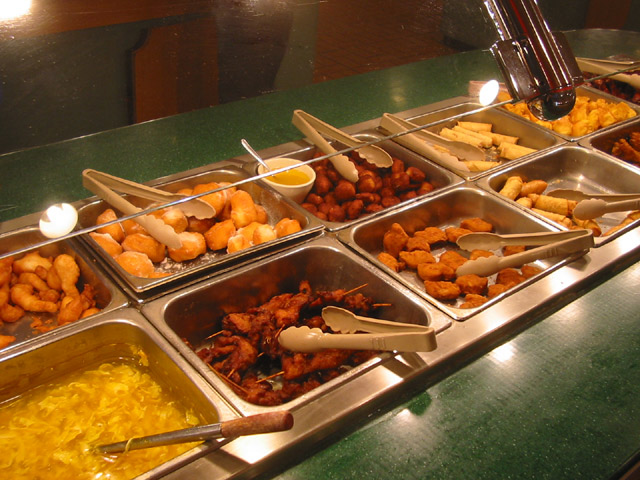} 
\end{center}

\item El restaurante del ejercicio precedente decide abrir una sucursal en la que
ofrece un buffet con cinco platos diferentes cada d\'ia. Cada comensal puede servirse
el que quiera, pero s\'olo cuatro veces (seguramente para evitar arruinarse con gente
`sin llenadera'\index{Llenadera}).
?`Cu\'antas opciones tiene un comensal?


\item Si en el experimento `nacimiento de un hijo' los resultados `ni\~no' y `ni\~na' son
igualmente probables y cada nacimiento se supone independiente de los dem\'as\footnote{Esta
es una simplificaci\'on muy dr\'astica, dado que no tiene en cuenta los detalles del proceso de 
fecundaci\'on. El \'ovulo femenino lleva \'unicamente el cromosoma X, mientras que los 
espermatozoides son portadores tanto del X como del Y. La combinaci\'on XX origina una ni\~na,
mientras que la XY da lugar a un ni\~no. La abundancia de espermatozoides con un gen u otro
depende de numerosos factores, diferentes para cada individuo.}, ?`cu\'al es
la probabilidad de que una familia con cuatro hijos est\'e formada por:
\begin{multicols}{3}
\begin{enumerate}[(a)]
\item dos ni\~nos y dos ni\~nas?,
\item tres ni\~nos y una ni\~na?,
\item cuatro ni\~nos?
\end{enumerate}
\end{multicols}


\item Se arrojan al aire cinco monedas iguales, no cargadas. Calcular las probabilidades de
obtener:
\begin{multicols}{2}
\begin{enumerate}[(a)]
\item Ninguna cara (o sello).
\item Exactamente una cara.
\item Al menos una cara.
\item No m\'as de cuatro caras.
\end{enumerate}
\end{multicols}

\item Supongamos que tenemos una caja con $n$ tornillos, de los cuales $d$ son
defectuosos.
\begin{enumerate}[(a)]
\item Determinar la probabilidad $p$ de que al extraer un tornillo, \'este sea defectuoso
\item Si se extraen $k$ tornillos con reemplazo (es decir, cada vez que se saca uno y se
comprueba si es defectuoso o no, se devuelve a la caja), los resultados de cada extracci\'on
son independientes entre s\'i, pues la situaci\'on vuelve a ser la inicial. Calcular la
probabilidad de obtener $\ell$ tornillos defectuosos de entre los $k$ extra\'idos.
\end{enumerate}

\item Un examen de respuesta m\'ultiple consta de $10$ preguntas, cada una de las cuales
tiene $4$ posibles respuestas. Si se eligen las respuestas al azar, ?`cu\'al es la 
probabilidad de acertar $6$ de ellas?
\end{enumerate}

\section{Uso de la tecnolog\'ia de c\'omputo}
Maxima implementa el c\'alculo de los coeficientes binomiales en el comando
\texttt{binomial}:
\begin{maxcode}
\noindent
\begin{minipage}[t]{8ex}{\color{red}\bf
\begin{verbatim}
(%i1) 
\end{verbatim}
}
\end{minipage}
\begin{minipage}[t]{\textwidth}{\color{blue}
\begin{verbatim}
binomial(9,4);
\end{verbatim}
}
\end{minipage}
\definecolor{labelcolor}{RGB}{100,0,0}
\begin{math}\displaystyle
\parbox{8ex}{\color{labelcolor}(\%o1) }
126
\end{math}
\end{maxcode}

\noindent Esto nos permite calcular inmediatamente el n\'umero de combinaciones
sin repetici\'on de $n$ elementos tomados de $k$ en $k$, puesto que, como sabemos,
$$
C^n_k={n\choose k}\,.
$$
Tambi\'en se puede calcular directamente el n\'umero de combinaciones con repetici\'on,
$$
CR^n_k={n+k-1\choose k}\,.
$$
Sin embargo, no existen funciones directas para calcular el n\'umero de permutaciones,
con y sin repetici\'on. Pero es muy sencillo
crear nuestras propias funciones que hagan este c\'alculo. Como ejemplo, veremos c\'omo
crear una funci\'on en Maxima para calcular el n\'umero de combinaciones con repetici\'on;
aunque acabamos de ver c\'omo se calcula, lo haremos como un sencillo ejercicio de
construcci\'on de funciones en Maxima.
La funci\'on se llamar\'a \texttt{CR} (recu\'erdese que Maxima distingue entre may\'usculas
y min\'usculas), y tomar\'a dos argumentos, $n$ y $k$, en este orden. Por tanto, la
funci\'on se invocar\'a escribiendo \texttt{CR(n,k)}.\\
El usuario puede construir sus propias funciones\index{Maxima!funciones} en Maxima mediante la funci\'on predefinida 
\texttt{block}, cuya sintaxis es
\begin{center}
\texttt{block([var1,...,vari],expr1,...,exprj)}\end{center}
Aqu\'i, \texttt{var1,...,vari} son las \emph{variables internas}\index{Maxima!variables} que se usar\'an para
definir la funci\'on (distintas a los \emph{argumentos}, que en nuestro caso ser\'an
$n$ y $k$). Podemos pensar en ellas como cantidades auxiliares que solo usaremos en
c\'alculos intermedios, para ahorrar espacio de escritura. Por su parte,
\texttt{expr1,...,exprj} son las distintas \emph{expresiones} que indican los c\'alculos
a realizar dentro de la funci\'on.\\
Todo se entiende muy f\'acilmente si estudiamos el siguiente c\'odigo, que describe
una funci\'on \texttt{foo(x,y)}, que toma dos argumentos \texttt{x} e \texttt{y}, y
calcula $\cos x^y$ utilizando la variable intermedia \texttt{var} definida como
$var=x^y$. As\'i, la \emph{expresi\'on} a calcular es simplemente 
$\cos (var)$:

\begin{maxcode}

\noindent
\begin{minipage}[t]{8ex}{\color{red}\bf
\begin{verbatim}
(%i2) 
\end{verbatim}
}
\end{minipage}
\begin{minipage}[t]{\textwidth}{\color{blue}
\begin{verbatim}
foo(x,y):=block([var],
var:x^y,
cos(var)
)$
\end{verbatim}
}
\end{minipage}

\noindent
\begin{minipage}[t]{8ex}{\color{red}\bf
\begin{verbatim}
(%i3) 
\end{verbatim}
}
\end{minipage}
\begin{minipage}[t]{\textwidth}{\color{blue}
\begin{verbatim}
foo(%pi,1);
\end{verbatim}
}
\end{minipage}
\definecolor{labelcolor}{RGB}{100,0,0}
\begin{math}\displaystyle
\parbox{8ex}{\color{labelcolor}(\%o3) }
-1
\end{math}
\end{maxcode}

\noindent Naturalmente,  en la pr\'actica no es necesario hacer esto, es mucho m\'as sencillo
\emph{no} usar la variable intermedia \texttt{var} y directamente escribir la expresi\'on
$\cos x^y$. Esto es s\'olo un ejemplo para aclarar conceptos.

\noindent Notemos c\'omo las variables auxiliares se declaran en una lista, 
y como al final de cada expresi\'on aparece una coma, excepto en la \'ultima,
que es la que devuelve la funci\'on como resultado.\\
La funci\'on que queremos es, simplemente,
\begin{maxcode}
\noindent
\begin{minipage}[t]{8ex}{\color{red}\bf
\begin{verbatim}
(%i4) 
\end{verbatim}
}
\end{minipage}
\begin{minipage}[t]{\textwidth}{\color{blue}
\begin{verbatim}
CR(n,k):=block(
binomial(n+k-1,k)
)$
\end{verbatim}
}
\end{minipage}
\noindent
\begin{minipage}[t]{8ex}{\color{red}\bf
\begin{verbatim}
(%i5) 
\end{verbatim}
}
\end{minipage}
\begin{minipage}[t]{\textwidth}{\color{blue}
\begin{verbatim}
CR(5,4);
\end{verbatim}
}
\end{minipage}
\definecolor{labelcolor}{RGB}{100,0,0}
\begin{math}\displaystyle
\parbox{8ex}{\color{labelcolor}(\%o5) }
70
\end{math}
\end{maxcode}

\section{Actividades de c\'omputo}
\begin{enumerate}[(A)]
\item Crear una funci\'on \texttt{P(n,k)} en Maxima que calcule el n\'umero de
permutaciones sin repetici\'on de $n$ elementos tomados de $k$ en $k$.
\item Crear una funci\'on \texttt{PR(n,k)} en Maxima que calcule el n\'umero de
permutaciones con repetici\'on de $n$ elementos tomados de $k$ en $k$.
\item Supongamos un proceso de Bernoulli donde, en un solo ensayo, la probabilidad
del evento $+$ es $p$ y la de $-$ es $q=1-p$. Crear una funci\'on \texttt{fprob(x,n,p)}
que proporcione la probabilidad de que en $n$ ensayos, se obtenga el evento $+$ un
n\'umero $x$ de veces. Utilizar esta funci\'on para comprobar los resultados obtenidos
en los ejercicios de este cap\'itulo.
\item Maxima dispone de comandos ya implementados para el estudio de las distribuciones de probabilidad m\'as habituales. Por ejemplo, podemos cargar el paquete \texttt{distrib},
mediante
\begin{maxcode}
\noindent
\begin{minipage}[t]{8ex}{\color{red}\bf
\begin{verbatim}
(%i6) 
\end{verbatim}
}
\end{minipage}
\begin{minipage}[t]{\textwidth}{\color{blue}
\begin{verbatim}
load("distrib")$
\end{verbatim}
}
\end{minipage}
\end{maxcode}
Este paquete dispone del comando \texttt{pdf\_binomial(x,n,p)}, que hace lo mismo que
el comando \texttt{fprob} propuesto en el \'item precedente. Comparar los resultados
de \texttt{fprob} y \texttt{pdf\_binomial}.
\end{enumerate}

\chapter{Aplicaciones y funciones}
\vspace{2.5cm}
\section{Nociones b\'asicas}
Intuitivamente, pensamos en una funci\'on como en una regla que a un n\'umero le asocia otro. Pero
el concepto es mucho m\'as amplio: cuando no tengamos conjuntos de n\'umeros, sino de cualquier tipo,
hablaremos de \emph{aplicaciones}\index{Aplicaci\'on}, y reservaremos la palabra \emph{funci\'on} para las aplicaciones\index{Funci\'on}
entre n\'umeros.

Es claro que $f(x)=x^2$ es una funci\'on, que a cada n\'umero real $x\in\mathbb{R}$ le asocia otro
n\'umero $x^2\in\mathbb{R}$, su cuadrado. O sea, que es claro lo que significa una funci\'on
$f:\mathbb{R}\to\mathbb{R}$ . Pero si tenemos el conjunto
$$
A:= \{\mbox{Alumnos de la facultad de ciencias} \}
$$
(en matem\'aticas, la notaci\'on $:=$ quiere decir `por definici\'on'), ?`qu\'e significa que
$f:A\to \mathbb{R}$ sea una aplicaci\'on? Podr\'iamos plantear algo m\'as complicado, definiendo
$$
C:=\{\mbox{\'arbol, llanta, verde, continuar, ninguno}\}
$$
y preguntando qu\'e tipo de objeto ser\'ia en este caso una aplicaci\'on $f:A\to C$ (?`una `regla' que a cada alumno
de la facultad de ciencias le asiga una llanta, o continuar, o ninguno? ?`Qu\'e regla es esa?). Para trabajar
en general, conviene recordar que dados dos conjuntos $A$ y $B$, su \emph{producto cartesiano} es el
conjunto\index{Cartesiano!producto}
$$
A\times B:=\{ (a,b):a\in A,b\in B\},
$$
es decir, el conjunto de pares ordenados donde el primer elemento del par pertenece a $A$ y el segundo a $B$.
Entonces, una aplicaci\'on $f:A\to B$ es simplemente un subconjunto de $A\times B$ con una condici\'on adicional
que discutiremos enseguida, y en caso de que
$(a,b)\in A\times B$ est\'e en el subconjunto definido por la aplicaci\'on, escribimos $b=f(a)$ y decimos que
$b$ \emph{es la imagen de} $a$\index{Imagen}.
Por ejemplo, si $A=\mathbb{R}=B$, la funci\'on $f(x)=x^2$ es una manera c\'omoda de escribir el subconjunto
de $\mathbb{R}\times\mathbb{R}$
$$
\{ (x,x^2):x\in\mathbb{R}\},
$$
algunos de cuyos elementos son $(-1,1)$, $(2,4)$, $(\pi,\pi^2)$. Al conjunto $A$ se le llama \emph{dominio}\index{Aplicaci\'on!dominio}
de la aplicaci\'on, y a $B$ el \emph{codominio}\index{Aplicaci\'on!codominio}. El \emph{rango} o \emph{recorrido}\index{Aplicaci\'on!rango} de $f$,
denotado por $\mathrm{Im}f$, es el conjunto de todas las im\'agenes de los elementos de $A$
por $f$. Claramente, es un subconjunto del codominio:
$$
\mathrm{Im}f:=\{f(a):a\in A\}\subset B\,.
$$
La condici\'on que mencionamos arriba es la siguiente: para que un subconjunto de $A\times B$ defina una 
aplicaci\'on, es necesario que para cada $a\in A$, \emph{contenga exactamente un par de la forma $(a,b)$}.
Es decir: cada elemento $a\in A$ no puede tener m\'as de una imagen.

Gr\'aficamente representamos una aplicaci\'on $f:A\to B$ escribiendo los elementos de $A$ en un eje horizontal,
los de $B$ en un eje vertical (los llamados \emph{ejes cartesianos}\index{Ejes coordenados} o \emph{ejes coordenados})\index{Cartesiano!ejes} y marcando los
pares $(a,b)$ que forman la aplicaci\'on. En el caso de tener una funci\'on, esto da lugar a la
\emph{gr\'afica de la funci\'on}\index{Funci\'on!gr\'afica} que ya hemos encontrado, por ejemplo, en el ejercicio \ref{parab}.
\begin{center}
    \begin{tikzpicture}[scale=0.8] 
        \draw [thin,gray!40] (0, 0) grid (6, 6);
        \draw [very thick] (0, 0) -- (6.5, 0);
        \draw [very thick] (0, 0) -- (0, 6.5);
        \node at (6.5, -0.5) {$A$};
        \node at (0, 7) {$B$};
        \node at (1,-0.5) {$a_1$};
        \node at (2,-0.5) {$a_2$};        
        \node at (3,-0.5) {$a_3$};        
        \node at (4,-0.5) {$a_4$};        
        \node at (5,-0.5) {$a_5$};        
        \node at (-0.5,1) {$b_1$};
        \node at (-0.5,2) {$b_2$};        
        \node at (-0.5,3) {$b_3$};        
        \node at (-0.5,4) {$b_4$};        
        \node at (-0.5,5) {$b_5$}; 
        \draw [color=blue, fill=blue] (1, 3) circle (0.1);
        \draw [color=blue, fill=blue] (2, 1) circle (0.1);
        \draw [color=blue, fill=blue] (3, 5) circle (0.1);
        \draw [color=blue, fill=blue] (4, 1) circle (0.1);
        \draw [color=blue, fill=blue] (5, 2) circle (0.1);                      
    \end{tikzpicture}
\end{center}
N\'otese que de acuerdo a la condici\'on de aplicaci\'on, si se traza una l\'inea vertical por un elemento
cualquiera $a\in A$, esa l\'inea debe contener exactamente un par $(a,b)$. As\'i, en la siguiente figura
aparecen ejemplos que no son ni aplicaci\'on (izquierda) ni funci\'on (derecha).
\begin{center}
\begin{minipage}{.5\textwidth}\begin{center}
    \begin{tikzpicture}[auto,scale=0.66] 
        \draw [thin,gray!40] (0, 0) grid (6, 6);
        \draw [very thick] (0, 0) -- (6.5, 0);
        \draw [very thick] (0, 0) -- (0, 6.5);
        \node at (6.5, -0.5) {$A$};
        \node at (0, 7) {$B$};
        \node at (1,-0.5) {$a_1$};
        \node at (2,-0.5) {$a_2$};        
        \node at (3,-0.5) {$a_3$};        
        \node at (4,-0.5) {$a_4$};        
        \node at (5,-0.5) {$a_5$};        
        \node at (-0.5,1) {$b_1$};
        \node at (-0.5,2) {$b_2$};        
        \node at (-0.5,3) {$b_3$};        
        \node at (-0.5,4) {$b_4$};        
        \node at (-0.5,5) {$b_5$}; 
        \draw [color=blue, fill=blue] (1, 4) circle (0.1);
        \draw [color=blue, fill=blue] (2, 1) circle (0.1);
        \draw [color=blue, fill=blue] (2, 2) circle (0.1);        
        \draw [color=blue, fill=blue] (3, 5) circle (0.1);
        \draw [color=blue, fill=blue] (4, 3) circle (0.1);
        \draw [color=blue, fill=blue] (5, 4) circle (0.1);                      
    \end{tikzpicture}\end{center}
\end{minipage}%
\begin{minipage}{.5\textwidth}\begin{center}
    \begin{tikzpicture}[auto,scale=0.66,x={(0,1cm)},y={(1 cm,0)}] 
        \draw [thick] (-3.25, 0) -- (3.25, 0);
        \draw [thick] (0, -3.25) -- (0, 3.25);
  		\draw [red,domain=-2.25:2.25,samples=100, thick]  plot (\x,\x*\x-2) ;  
    \end{tikzpicture}\end{center}
\end{minipage}%
\end{center}    
Existen tres tipos important\'isimos de aplicaciones $f:A\to B$.
\begin{list}{$\blacktriangleright$}{}
\item \emph{Inyectivas}\index{Aplicaci\'on!inyectiva}. Son aquellas que no repiten valores, es decir, si dos elementos
$a,b\in A$ cumplen que $f(a)=f(b)$, entonces $a=b$. Otra manera de decirlo es que si
$a\neq b$, entonces $f(a)\neq f(b)$.\\
Las gr\'aficas de las aplicaciones inyectivas no cortan m\'as de una vez a cualquier recta
\emph{horizontal} en $A\times B$ (y esto incluye la posibilidad de que \emph{no} corten a una
determinada recta).
\begin{center}
\begin{minipage}{.5\textwidth}\begin{center}
    \begin{tikzpicture}[auto,scale=0.66] 
        \draw [thin,gray!40] (0, 0) grid (6, 6);
        \draw [very thick] (0, 0) -- (6.5, 0);
        \draw [very thick] (0, 0) -- (0, 6.5);
        \node at (6.5, -0.5) {$A$};
        \node at (0, 7) {$B$};
        \node at (1,-0.5) {$a_1$};
        \node at (2,-0.5) {$a_2$};        
        \node at (3,-0.5) {$a_3$};        
        \node at (4,-0.5) {$a_4$};        
        \node at (5,-0.5) {$a_5$};        
        \node at (-0.5,1) {$b_1$};
        \node at (-0.5,2) {$b_2$};        
        \node at (-0.5,3) {$b_3$};        
        \node at (-0.5,4) {$b_4$};        
        \node at (-0.5,5) {$b_5$}; 
        \draw [color=blue, fill=blue] (1, 5) circle (0.1);
        \draw [color=blue, fill=blue] (2, 1) circle (0.1);
        \draw [color=blue, fill=blue] (3, 2) circle (0.1);
        \draw [color=blue, fill=blue] (4, 4) circle (0.1);
        \draw [color=blue, fill=blue] (5, 3) circle (0.1);                      
    \end{tikzpicture}\end{center}
\end{minipage}%
\begin{minipage}{.5\textwidth}\begin{center}
    \begin{tikzpicture}[auto,scale=0.66] 
        \draw [thick] (-3.25, 0) -- (3.25, 0);
        \draw [thick] (0, -3.25) -- (0, 3.25);
  		\draw [red,domain=-2.5:2.5,samples=100, thick]  plot (\x,{0.4*exp(\x)-2}) ;  
    \end{tikzpicture}\end{center}
\end{minipage}%
\end{center}
Estudiemos como ejemplo la funci\'on $f:\mathbb{R}\to\mathbb{R}$ dada por $f(x)=x^2+3$, cuyo dominio es $\mathbb{R}$ y recorrido\footnote{La notaci\'on $[a,b]$ significa `todos los n\'umeros $x$
tales que $a\leq x\leq b$'. Si se escribe, por ejemplo, $]a,b]$, se quiere decir
`todos los n\'umeros $x$ tales que $a<x\leq b$', y una expresi\'on como $[a,+\infty [$
significa `todos los $x$ tales que $a\leq x<+\infty$. Los conjuntos de la forma $[a,b]$ se
conocen como \emph{intervalos}, abiertos o cerrados seg\'un se incluya o no sus extremos.} $[3,+\infty [\subset \mathbb{R}$. La
inyectividad implica que, dado un n\'umero en el dominio $a\in\mathbb{R}$, la ecuaci\'on $f(a)=f(b)$ tiene una
soluci\'on o ninguna. En nuestro caso, esta ecuaci\'on es
$$
a^2+3=b^2+3,
$$
de donde resulta que $a^2=b^2$, una ecuaci\'on que, al tomar ra\'ices, tiene dos soluciones $b=\pm a$. Por
tanto, conclu\'imos que $f(x)=x^2+3$ no es inyectiva. Esto tambi\'en puede verse gr\'aficamente, notando
que $f(x)$ es una funci\'on cuadr\'atica\index{Funci\'on!cuadr\'atica} y recordando c\'omo son las gr\'aficas de estas funciones (ejercicio
\ref{parab}).\\
Por el contrario, una funci\'on \emph{lineal}\index{Funci\'on!lineal} como $f(x)=mx+n$ (con $m\neq 0$) \emph{siempre} es inyectiva: si fuese
$f(a)=f(b)$, tendr\'iamos
$$
ma+n=mb+n,
$$
de donde, sumando $-n$, $ma=mb$ y, al dividir por $m\neq 0$ resulta $a=b$. Este resultado tambi\'en se deduce de
la forma de la gr\'afica de la funci\'on.

\item \emph{Sobreyectivas}\index{Aplicaci\'on!sobreyectiva}. Son aquellas que recorren todos los valores de su codominio. En otras palabras, si $b\in B$, existe un $a\in A$ tal que $b=f(a)$. La gr\'afica de una aplicaci\'on sobreyectiva corta al menos una vez a 
\emph{cualquier} recta horizontal en $A\times B$.
\begin{center}
\begin{minipage}{.5\textwidth}\begin{center}
    \begin{tikzpicture}[auto,scale=0.66] 
        \draw [thin,gray!40] (0, 0) grid (6, 6);
        \draw [very thick] (0, 0) -- (6.5, 0);
        \draw [very thick] (0, 0) -- (0, 6.5);
        \node at (6.5, -0.5) {$A$};
        \node at (0, 7) {$B$};
        \node at (1,-0.5) {$a_1$};
        \node at (2,-0.5) {$a_2$};        
        \node at (3,-0.5) {$a_3$};        
        \node at (4,-0.5) {$a_4$};        
        \node at (5,-0.5) {$a_5$};        
        \node at (-0.5,1) {$b_1$};
        \node at (-0.5,2) {$b_2$};        
        \node at (-0.5,3) {$b_3$};        
        \node at (-0.5,4) {$b_4$};        
        \node at (-0.5,5) {$b_5$}; 
        \draw [color=blue, fill=blue] (1, 3) circle (0.1);
        \draw [color=blue, fill=blue] (2, 1) circle (0.1);
        \draw [color=blue, fill=blue] (3, 2) circle (0.1);
        \draw [color=blue, fill=blue] (4, 5) circle (0.1);
        \draw [color=blue, fill=blue] (5, 4) circle (0.1);                      
    \end{tikzpicture}\end{center}
\end{minipage}%
\begin{minipage}{.5\textwidth}\begin{center}
    \begin{tikzpicture}[auto,scale=0.66,yscale=0.25] 
        \draw [thick] (-3.25, 0) -- (3.25, 0);
        \draw [thick] (0, -12.25) -- (0, 12.25);
  		\draw [red,domain=-2.75:2.75,samples=100, thick]  plot (\x,{\x*\x*\x-3*\x}) ;  
    \end{tikzpicture}\end{center}
\end{minipage}%
\end{center}
Consideremos por ejemplo la funci\'on $f(x)=x^3$. Su dominio y codominio coinciden con
$\mathbb{R}$ (n\'otese que al ser la potencia de $x$ impar, puede tomar valores positivos
y negativos). Para que sea sobreyectiva, cualquier n\'umero real del codominio, $b\in\mathbb{R}$,
debe ser la imagen de un cierto elemento del dominio, $a\in\mathbb{R}$, de modo que $b=f(a)$.
Si conocemos $b$ podemos determinar $a$ resolviendo la ecuaci\'on
$$
b=a^3,
$$
de donde $a=\sqrt[3]{b}$. En general, un n\'umero real tiene $n$ ra\'ices $n-$\'esimas\index{Rai@Ra\'iz!$n$-\'esima}, de las cuales
algunas pueden ser complejas, como veremos m\'as adelante. En este caso, de las 3 ra\'ices de $b$,
seguro existe una real\footnote{Esto es as\'i porque las ra\'ices complejas aparecen por pares (si hay una ra\'iz
compleja, su conjugada tambi\'en es ra\'iz)
de modo que de las 3 ra\'ices, o las 3 son reales o dos son complejas y la otra real. No hay m\'as posibilidades.},
que es la soluci\'on buscada.

\item \emph{Biyectivas}\index{Aplicaci\'on!biyectiva}, las que son simult\'aneamente inyectivas y sobreyectivas. Cada elemento del dominio
tiene una sola imagen y cada elemento del codominio tiene una sola preimagen. La gr\'afica de la aplicaci\'on
corta \emph{exactamente} una vez a cada recta horizontal y a cada recta vertical en $A\times B$.\\
Las funciones de $\mathbb{R}$ en $\mathbb{R}$ lineales, $f(x)=mx+n$, son biyectivas, mientras que las cuadr\'aticas $f(x)=ax^2+bx+c$ nunca lo son (?`puede el lector ver por qu\'e?). Las
aplicaciones biyectivas de la recta $\mathbb{R}$ o del 
\emph{plano cartesiano}\index{Plano!cartesiano}\index{R2@$\mathbb{R}^2$} $\mathbb{R}^2=\{(x,y)\in\mathbb{R}\}$
se llaman \emph{transformaciones}\index{Transformaciones} (de la recta o del plano, respectivamente). As\'i, por ejemplo, $f(x)=mx+n$ es una transformaci\'on 
lineal\index{Transformaci\'on!lineal} de la recta.
\begin{center}
\begin{minipage}{.5\textwidth}\begin{center}
    \begin{tikzpicture}[auto,scale=0.66] 
        \draw [thin,gray!40] (0, 0) grid (6, 6);
        \draw [very thick] (0, 0) -- (6.5, 0);
        \draw [very thick] (0, 0) -- (0, 6.5);
        \node at (6.5, -0.5) {$A$};
        \node at (0, 7) {$B$};
        \node at (1,-0.5) {$a_1$};
        \node at (2,-0.5) {$a_2$};        
        \node at (3,-0.5) {$a_3$};        
        \node at (4,-0.5) {$a_4$};        
        \node at (5,-0.5) {$a_5$};        
        \node at (-0.5,1) {$b_1$};
        \node at (-0.5,2) {$b_2$};        
        \node at (-0.5,3) {$b_3$};        
        \node at (-0.5,4) {$b_4$};        
        \node at (-0.5,5) {$b_5$}; 
        \draw [color=blue, fill=blue] (1, 2) circle (0.1);
        \draw [color=blue, fill=blue] (2, 1) circle (0.1);
        \draw [color=blue, fill=blue] (3, 3) circle (0.1);
        \draw [color=blue, fill=blue] (4, 5) circle (0.1);
        \draw [color=blue, fill=blue] (5, 4) circle (0.1);                      
    \end{tikzpicture}\end{center}
\end{minipage}%
\begin{minipage}{.5\textwidth}\begin{center}
    \begin{tikzpicture}[auto,scale=0.66,yscale=0.25] 
        \draw [thick] (-3.25, 0) -- (3.25, 0);
        \draw [thick] (0, -12.25) -- (0, 12.25);
  		\draw [red,domain=-2.75:2.75,samples=100, thick]  plot (\x,{0.55*\x*\x*\x}) ;  
    \end{tikzpicture}\end{center}
\end{minipage}%
\end{center}
\end{list}

Un concepto esencial cuando se trabaja con aplicaciones es el de \emph{composici\'on}\index{Aplicaci\'on!composici\'on}. La
composici\'on de dos aplicaciones no es m\'as que el resultado de hacerlas actuar sucesivamente
sobre un elemento. Si tenemos aplicaciones $f:A\to B$ y $g:B\to C$, donde \emph{el dominio de la
segunda contiene al rango de la primera}, podemos hacerlas actuar as\'i:
$$
A\stackrel{f}{\longrightarrow}B\stackrel{g}{\longrightarrow}C 
$$
El resultado es la composici\'on de $g$ con $f$ (en este orden), denotada $g\circ f$:
$$
A\xrightarrow{\; g\circ f\; }C
$$
Por ejemplo, si $f(x)=x^2+1$ y $g(x)=\frac{x+1}{x-1}$, la composici\'on $g\circ f$ ser\'ia
la funci\'on dada para todo real distinto del $0$ por
$$
(g\circ f)(x)=g(\color{blue}f(x)\color{black})
=\frac{\color{blue}f(x)\color{black}+1}{\color{blue}f(x)\color{black}-1}
=\frac{\color{blue}x^2+1\color{black}+1}{\color{blue}x^2+1\color{black}-1}=\frac{x^2+2}{x^2}.
$$
Es imprescindible observar que la composici\'on de funciones \emph{no es conmutativa}, aun cuando los
dominios y codominios sean los adecuados (?`qu\'e tiene que cumplirse para que tenga sentido plantear
la igualdad $g\circ f=f\circ g$?). En el
ejemplo que estamos considerando se ve claramente que $g\circ f\neq f\circ g$, pues $f\circ g$ es la
funci\'on definida para todo real distinto de $1$ por
$$
(f\circ g)(x)=f(g(x))=\left( \frac{x+1}{x-1} \right)^2+1
=\frac{(x+1)^2}{(x-1)^2}+\frac{(x-1)^2}{(x-1)^2}
=\frac{x^2+2x+1+x^2-2x+1}{(x-1)^2}
=\frac{2x^2+2}{(x-1)^2}
$$

Una propiedad fundamental de las aplicaciones biyectivas es que admiten una \emph{aplicaci\'on
inversa}. Dada una aplicaci\'on $f:A\to B$, su inversa\index{Aplicaci\'on!inversa} es una aplicaci\'on $f^{-1}:B\to A$ (con
el dominio y el codominio intercambiados) tal que la composici\'on de $f$ con $f^{-1}$ \emph{en
cualquier orden} da la aplicaci\'on identidad, es decir,
\begin{align}\begin{split}\label{inversas}
(f^{-1}\circ f)(a)=a \\
(f\circ f^{-1})(b)=b,
\end{split}\end{align}
para cualesquiera $a\in A$, $b\in B$.

Las funciones lineales\index{Funci\'on!lineal}, por ser biyectivas, tienen inversa.
Si consideramos $f(x)=2x-1$, podemos hallar su inversa escribiendo $y=f(x)=2x-1$ y calculando $x$
como funci\'on de $y$\footnote{El razonamiento con todo detalle es el siguiente: como $f^{-1}\circ f=\mathrm{Id}$,
aplicando $f^{-1}$ a ambos miembros de $y=f(x)$, resulta $f^{-1}(y)=x$.}
$$
y+1=2x \Rightarrow x=\frac{y+1}{2}\, ,
$$
(el s\'imbolo $\Rightarrow$ se lee `implica que' y es de uso com\'un en matem\'aticas), de donde, en
la notaci\'on habitual,
$$
f^{-1}(x)=\frac{x+1}{2}\, .
$$
No est\'a de m\'as comprobar que estas funciones cumplen la condici\'on de ser inversas una de la otra,
\eqref{inversas}. Por una parte, tenemos que
$$
(f^{-1}\circ f)(x)=f^{-1}(f(x))=f^{-1}(2x-1)=\frac{2x-1+1}{2}=\frac{2x}{2}=x,
$$
y por otra:
$$
(f\circ f^{-1})(x)=f(f^{-1}(x))=f\left( \frac{x+1}{2}\right) =2\frac{x+1}{2}-1
=x+1-1=x.
$$

Gr\'aficamente, se puede obtener la representaci\'on de $f^{-1}$ sin m\'as que hacer una
reflexi\'on de la gr\'afica de $f$ respecto de la recta $y=x$ (bisectriz del primer cuadrante),
pues esta acci\'on se corresponde con el intercambio de los ejes $x$ e $y$.

Por \'ultimo, introduciremos un tipo especial de funci\'on que aparecer\'a con frecuencia.
Una \emph{funci\'on polin\'omica}\index{Funci\'on!polin\'omica} o \emph{polinomio}\index{Polinomio} de
\emph{grado} $k$\index{Polinomio!grado} es una de la forma
$$
f(x)=a_kx^k +a_{k-1}x^{k-1}+\cdots +a_1x+a_0
$$
donde los n\'umeros reales $a_n$ (con $n\in\{0,\ldots,k\}$) se denominan los \emph{coeficientes}
\index{Polinomio!coeficientes} del polinomio y $a_0$ en particular, el \emph{t\'ermino independiente}.
Por ejemplo, $f(x)=-3x^2+1$ es un polinomio de grado $2$ cuyo t\'ermino independiente es $1$.

\section{Aplicaciones}
La principal utilidad de las funciones es que, si se usan para modelar fen\'omenos
naturales, permiten hacer \emph{predicciones}, que es el objetivo de cualquier
disciplina cient\'ifica. Consideremos por ejemplo el lanzamiento de un objeto, tal
como una piedra o un misil, o el salto de un atleta. Mediante principios f\'isicos
elementales se llega a que la funci\'on de describe la posici\'on $x$ sobre el suelo
en funci\'on del tiempo es
$$
x(t)=v_0 t\, \cos\alpha\,,
$$
donde $v_0$ es la velocidad inicial y $\alpha$ es el \'angulo de lanzamiento.\\
Por su parte, la altura sobre el suelo como funci\'on del tiempo est\'a dada por
la funci\'on
$$
y(t)=v_0 t\, \sin\alpha -\frac{1}{2}gt^2\,,
$$
donde $g=9{.}80$m/s$^2$ es la aceleraci\'on debida a la gravedad en la superficie
terrestre (el peso de un cuerpo de masa $m$ en la superficie terrestre es pues $P=mg$).
A partir de la funci\'on $y(t)$, podemos calcular el tiempo que el proyectil o el
saltador estar\'an en el aire, pues ser\'a el valor $t_0$ para el cual $y(t_0)=0$
(es decir, toca de nuevo el suelo). Planteando la ecuaci\'on, resulta, 
$$
v_0 t_0\, \sin\alpha -\frac{1}{2}gt^2_0=0\,,
$$
que se puede poner como
$$
t_0\left(v_0\,\sin\alpha -\frac{1}{2}gt_0\right) =0\,.
$$
De aqu\'i es obvio que las soluciones son $t_0=0$, que corresponde al momento del
lanzamiento o del salto, y
$$
t_0=\frac{2v_0\,\sin\alpha}{g}\,,
$$
que corresponde al tiempo de vuelo.
\begin{center}
\includegraphics[scale=0.4]{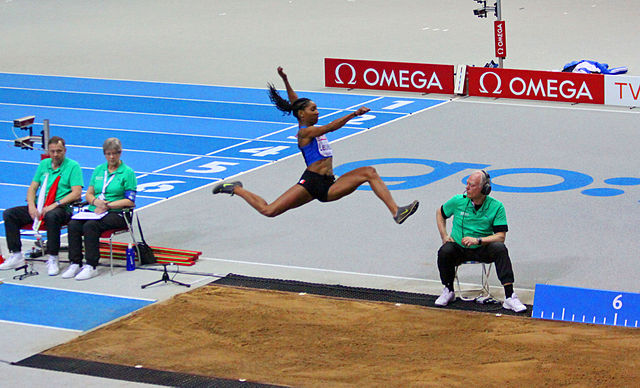} 
\end{center}
Con estos datos, podemos \emph{predecir} el lugar donde caer\'a el proyectil o el saltador. 
Ser\'a el valor de la posici\'on que corresponde a $t=t_0$, es decir, $x(t_0)$. Sustituyendo:
$$
x(t_0)=v_0\frac{2}{g}v_0\,\sin\alpha\cos\alpha = \frac{v^2_0}{g}\,\sin (2\alpha)\,,
$$
donde hemos hecho uso de una identidad trigonom\'etrica ($2\sin\alpha\cos\alpha =\sin 2\alpha$)
que veremos m\'as adelante (en el cap\'itulo \ref{chaptrig}).\\
A la vista de esta expresi\'on, vemos que hay dos posibilidades para mejorar la distancia
alcanzada: hacer la velocidad inicial $v_0$ tan grande como se pueda (eso es lo que hace el
saltador de longitud con la carrera previa al salto) y tratar de que la inclinaci\'on inicial
$\alpha$ est\'e tan cerca de $45\degree$ como se pueda (pues, como veremos en el cap\'itulo
\ref{chaptrig}, \'este es el valor para el cual $\sin 2\alpha$ alcanza su valor m\'aximo,
que es $1$). Sin embargo, aunque el segundo m\'etodo es factible para los lanzamientos de
proyectiles, no lo es en el caso de saltos de longitud. Los saltadores profesionales no
se acercan a $\alpha =45\degree$, sino a $\alpha =22\degree$. El motivo son las limitaciones
f\'isicas del cuerpo humano: para un atleta, no es igual de f\'acil invertir la energ\'ia
de la carrera en lograr un \'angulo inicial de $22\degree$ que uno de $45\degree$, debido
a que en el segundo caso la inercia de la carrera le impide al saltador colocar los m\'usculos
de la pierna en una posici\'on \'optima para darse impulso hacia arriba. Un an\'alisis detallado
puede verse en el trabajo de Linthorne, Guzm\'an y Bridgett en la revista Journal of Sports 
Sciences, disponible en
\url{http://www.tandfonline.com/doi/abs/10.1080/02640410400022011#.U_ANPitqnIE}. 

\section{Ejercicios}

\begin{enumerate}  
\setcounter{enumi}{58}

\item Si $f(t)=t+(1/t)$, probar que $f(a)=f(1/a)$ para todo $a\neq 0$.

\item Si $g(x)=x^2-x+1$, probar que $g(a)=g(1-a)$ para cualquier $a\in\mathbb{R}$.

\item A veces, al dar una funci\'on, no se especifica su dominio. En estos casos, se supone que se trata
del \emph{dominio natural}, es decir, el mayor subconjunto de $\mathbb{R}$ en el que la expresi\'on de la
funci\'on tiene sentido. Calcular los dominios naturales de las siguientes funciones.
\begin{multicols}{2}
\begin{enumerate}[(a)]
\item $f(x)=\sqrt{x-1}$
\item $f(x)=\sqrt{1-x}$
\item $g(x)=\sqrt{1-x^2}$
\item $g(x)=1/\sqrt{x-1}$
\item $f(t)=t/(t-1)$
\item $h(t)=t/(t^2-1)$
\item $g(x)=x^3-3x^2$
\item $h(u)=u^3-3/u^2$
\item $f(x)=x/(x^2+1)$
\end{enumerate}
\end{multicols}

\item En los siguientes ejercicios, se dan unas funciones con sus respectivos dominios y se pide
determinar el rango $\mathrm{Im}f$.
\begin{enumerate}[(a)]
\item $f(x)=x^4$ con dominio $\mathbb{R}$.
\item $f(x)=1+x^2$, con dominio $\mathbb{R}$.
\item $f(x)=\sqrt{x-1}$, con dominio el conjunto $\{x\in\mathbb{R}:x\geq 1\}$.
\item $f(x)=1/x$, con dominio $\mathbb{R}-\{0\}$.
\item $f(x)=1/\sqrt{x}$, con dominio $\mathbb{R}^+$.
\item $f(x)=2/(3+x^2)$, con dominio $\mathbb{R}$.
\end{enumerate}

\item Determinar si las curvas que se presentan a continuaci\'on son o no las gr\'aficas de alguna funci\'on
$y=f(x)$. En caso afirmativo, decidir si se trata de una funci\'on inyectiva, sobreyectiva o biyectiva.
\begin{center}
\begin{minipage}{.5\textwidth}(a)\begin{center}
    \begin{tikzpicture}[auto,scale=0.6]
        \draw [thick] (-3.25, 0) -- (3.25, 0);
        \draw [thick] (0, -3.25) -- (0, 3.25);
        \begin{scope}[rotate=-90]
  		\draw [red,domain=-2.75:2.75,samples=100, thick]  plot (\x,{0.3*\x*\x*\x-3*\x}) ;
  		\end{scope}                       
    \end{tikzpicture}\end{center}
\end{minipage}%
\begin{minipage}{.5\textwidth}(b)\begin{center}
    \begin{tikzpicture}[auto,scale=0.6,yscale=0.25] 
        \draw [thick] (-3.25, 0) -- (3.25, 0);
        \draw [thick] (0, -13.5) -- (0, 13.5);
  		\draw [red,domain=-2.75:2.75,samples=100, thick]  plot (\x,{\x*\x*\x-3*\x}) ;  
    \end{tikzpicture}\end{center}
\end{minipage}%
\end{center}
\begin{center}
\begin{minipage}{.5\textwidth}(c)\begin{center}
    \begin{tikzpicture}[auto,scale=0.6]
        \draw [thick] (-3.25, 0) -- (3.25, 0);
        \draw [thick] (0, -3.25) -- (0, 3.25);
  		\draw [red, thick]  (-3.2,-3)--(-1,0.5) ;
  		\draw [red, thick]  (-1,0.5)--(1,0.5) ;
  		\draw [red, thick]  (1,0.5)--(3.25,3.2) ;  		  		
    \end{tikzpicture}\end{center}
\end{minipage}%
\begin{minipage}{.5\textwidth}(d)\begin{center}
    \begin{tikzpicture}[auto,scale=0.6] 
        \draw [thick] (-3.25, 0) -- (3.25, 0);
        \draw [thick] (0, -3.25) -- (0, 3.25);
  		\draw [red, thick]  (-3.25,-3)--(0.5,0.5) ;
  		\draw [red,thick]   (0.5,0.5) -- (0.5,2);
  		\draw [red,thick]   (0.5,2)--(3.25,3);
    \end{tikzpicture}\end{center}
\end{minipage}%
\end{center}
\begin{center}
\begin{minipage}{.5\textwidth}(e)\begin{center}
    \begin{tikzpicture}[auto,scale=0.6]
        \draw [thick] (-3.25, 0) -- (3.25, 0);
        \draw [thick] (0, -3.25) -- (0, 3.25);
         \begin{scope}[rotate=180]       
  \draw[red,domain=-3.25:2,samples=100, thick] plot(\x,{1.5*(1-exp(\x-1))*(1-exp(\x-1))-1.5}) ;
  		\end{scope}  	  		
    \end{tikzpicture}\end{center}
\end{minipage}%
\begin{minipage}{.5\textwidth}(f)\begin{center}
    \begin{tikzpicture}[auto,scale=0.6] 
        \draw [thick] (-3.25, 0) -- (3.25, 0);
        \draw [thick] (0, -3.25) -- (0, 3.25);
        \begin{scope}[rotate=-90]
  		\draw [red,domain=-2.75:2.75,samples=100, thick]  plot (\x,{0.7*\x*\x-2)}) ;
  		\end{scope}  
    \end{tikzpicture}\end{center}
\end{minipage}%
\end{center}

\item Determinar si las curvas que se presentan a continuaci\'on son o no la gr\'afica de 
alguna funci\'on $y=f(x)$ \emph{biyectiva}. En caso afirmativo, construir la gr\'afica de la
funci\'on inversa.
\begin{center}
\begin{minipage}{.5\textwidth}(a)\begin{center}
    \begin{tikzpicture}[auto,scale=0.6]
        \draw [thick] (-3.25, 0) -- (3.25, 0);
        \draw [thick] (0, -3.25) -- (0, 3.25);
  		\draw [red,domain=-2.75:2.75,samples=100, thick]  plot (\x,{exp(1/(0.85*\x+3))-2}) ;
    \end{tikzpicture}\end{center}
\end{minipage}%
\begin{minipage}{.5\textwidth}(b)\begin{center}
    \begin{tikzpicture}[auto,scale=0.6,yscale=0.25] 
        \draw [thick] (-3.25, 0) -- (3.25, 0);
        \draw [thick] (0, -13.5) -- (0, 13.5);
  		\draw [red,domain=-2.75:3.25,samples=100, thick]  plot (\x,{1.2*(\x*\x-2*\x-3)}) ;  
    \end{tikzpicture}\end{center}
\end{minipage}%
\end{center}
\begin{center}
\begin{minipage}{.5\textwidth}(c)\begin{center}
    \begin{tikzpicture}[auto,scale=0.6]
        \draw [thick] (-3.25, 0) -- (3.25, 0);
        \draw [thick] (0, -3.25) -- (0, 3.25);
	    \draw[red,domain=0:3.25,samples=100, thick]  plot (\x,{0.5*sqrt(\x*\x*\x)}) ;
	    \draw[red,domain=0:3.25,samples=100, thick]  plot (\x,{-0.5*sqrt(\x*\x*\x)}) ;	     
 		  		
    \end{tikzpicture}\end{center}
\end{minipage}%
\begin{minipage}{.5\textwidth}(d)\begin{center}
    \begin{tikzpicture}[auto,scale=0.6] 
        \draw [thick] (-3.25, 0) -- (3.25, 0);
        \draw [thick] (0, -3.25) -- (0, 3.25);        
		\begin{scope}[rotate=90]
	    \draw[red,domain=0:3.25,samples=100, thick]  plot (\x,{0.5*sqrt(\x*\x*\x)}) ;
	    \draw[red,domain=0:3.25,samples=100, thick]  plot (\x,{-0.5*sqrt(\x*\x*\x)}) ;		
		\end{scope}
    \end{tikzpicture}\end{center}
\end{minipage}%
\end{center}
\begin{center}
\begin{minipage}{.5\textwidth}(e)\begin{center}
    \begin{tikzpicture}[auto,scale=0.6,xscale=5]
        \draw [thick] (-0.65, 0) -- (0.65, 0);
        \draw [thick] (0, -3.25) -- (0, 3.25);
  \draw[red,domain=-0.65:0.65,samples=100, thick] plot(\x,{10*\x*\x*\x}) ;
    \end{tikzpicture}\end{center}
\end{minipage}%
\begin{minipage}{.5\textwidth}(f)\begin{center}
    \begin{tikzpicture}[auto,scale=0.6] 
        \draw [thick] (-3.25, 0) -- (3.25, 0);
        \draw [thick] (0, -3.25) -- (0, 3.25);
  		\draw [red,domain=0.25:3.2,samples=100, thick]  plot (\x,{0.8/\x}) ;
  		\draw [red,domain=-3.2:-0.25,samples=100, thick]  plot (\x,{0.8/\x}) ;  		
    \end{tikzpicture}\end{center}
\end{minipage}%
\end{center}

\item\label{ejer65} En los siguientes ejercicios, se pide determinar si la funci\'on dada posee inversa.
En caso afirmativo, determinar expl\'icitamente $f^{-1}$, dando su dominio y rango.
\begin{enumerate}[(a)]
\item $f(x)=1/x^2$, para $x\neq 0$.
\item $f(x)=1/x^2$, para $x>0$.
\item $f(x)=1/x^2$, para $x<0$.
\item $f(x)=1/x^3$, para $x\neq 0$.
\item $f(x)=1/(x-2)$, para $x\neq 2$.
\item $f(x)=\sqrt{x-2}$, para $x\geq 2$.
\end{enumerate}

\item Consideremos la funci\'on
$$
f(x)=\frac{1}{(x-1)^2},
$$
con dominio los reales $x>1$. Calcular su inversa $f^{-1}$ y especificar su dominio y rango.

\item\label{ejer67} Calcular las composiciones $f\circ g$ y $g\circ f$ para los siguientes pares de funciones.
\begin{enumerate}[(a)]
\item $f(x)=x^5$ y $g(x)=x^2+1$.
\item $f(x)=1/x$ y $g(x)=x/(x+1)$.
\item $f(x)=ax+b$ y $g(x)=x^2-1$.
\end{enumerate}

\item\label{ejer68} Una funci\'on $f(x)$ se dice que es \emph{creciente}\index{Funci\'on!creciente} si para cada par de elementos en su dominio
$a,b$, tales que $a\leq b$, se tiene que $f(a)\leq f(b)$. Cuando las desigualdades son estrictas, se
dice que la funci\'on es \emph{estrictamente creciente}. De la misma forma, se dice que $f(x)$ es
\emph{decreciente}\index{Funci\'on!decreciente} si para cada par $a\leq b$ es $f(a)\geq f(b)$ (y se dice \emph{estrictamente
decreciente} si las desigualdades son estrictas). El dominio de una funci\'on puede dividirse en
subconjuntos donde la funci\'on sea creciente o decreciente.\\
Para las siguientes funciones, determinar d\'onde son crecientes y d\'onde decrecientes.
\begin{multicols}{2}
\begin{enumerate}[(a)]
\item $f(x)=-2x+7$
\item $f(x)=\frac{x}{3}-1$
\item $f(x)=(x+1)^2$
\item $f(x)=x^3+1$
\item $f(x)=1/x$, $x\neq 0$
\item $f(x)=\sqrt{x-1}$, $x\geq 1$
\item $f(x)=1/(x-1)^2$, $x\neq 1$
\item $f(x)=1/(1+x^2)$
\end{enumerate}
\end{multicols}

\item\label{ejer69} Una funci\'on se dice que es \emph{par}\index{Funci\'on!par} si para cualquier $x$ en su dominio, ocurre que
$-x$ tambi\'en lo est\'a y, adem\'as,
$$
f(-x)=f(x).
$$
La funci\'on se dice que es \emph{impar}\index{Funci\'on!impar} si se cumple
$$
f(-x)=-f(x).
$$
En los ejercicios que siguen, determinar si la funci\'on es par, impar o ninguna de las dos.
\begin{multicols}{2}
\begin{enumerate}
\item $f(x)=2x$
\item $f(x)=2x-1$
\item $f(x)=x^4-3x^2+1$
\item $f(x)=x^3-3x+1$
\item $f(x)=x^3+5x$
\item $f(x)=x/(x^2+1)$
\item $f(x)=x^2/(x^2+1)$
\item $f(x)$ definida como
$$
f(x)=
\begin{cases}
1-x,&\mbox{ si }x\geq 0 \\
1+x,&\mbox{ si }x\leq 0.
\end{cases}
$$
\end{enumerate}
\end{multicols}

\end{enumerate}

\section{Uso de la tecnolog\'ia de c\'omputo}
En Maxima se puede calcular la composici\'on de dos funciones. La manera de hacerlo
que vamos a describir no es la mejor, pero es la m\'as directa y para las
tareas sencillas en que la usaremos, no dar\'a problemas.
Comenzamos definiendo las funciones con las que vamos a trabajar,

\begin{maxcode}
\noindent
\begin{minipage}[t]{8ex}{\color{red}\bf
\begin{verbatim}
(%i1) 
\end{verbatim}
}
\end{minipage}
\begin{minipage}[t]{\textwidth}{\color{blue}
\begin{verbatim}
f(x):=x^2+1$
\end{verbatim}
}
\end{minipage}

\noindent
\begin{minipage}[t]{8ex}{\color{red}\bf
\begin{verbatim}
(%i2) 
\end{verbatim}
}
\end{minipage}
\begin{minipage}[t]{\textwidth}{\color{blue}
\begin{verbatim}
g(x):=cos(2*x)$
\end{verbatim}
}
\end{minipage}
\end{maxcode}

La composici\'on se puede construir con la notaci\'on habitual
mediante los siguientes comandos:
\begin{maxcode}
\noindent
\begin{minipage}[t]{8ex}{\color{red}\bf
\begin{verbatim}
(%i3) 
\end{verbatim}
}
\end{minipage}
\begin{minipage}[t]{\textwidth}{\color{blue}
\begin{verbatim}
gof(x):=g(f(x))$
\end{verbatim}
}
\end{minipage}

\noindent
\begin{minipage}[t]{8ex}{\color{red}\bf
\begin{verbatim}
(%i4) 
\end{verbatim}
}
\end{minipage}
\begin{minipage}[t]{\textwidth}{\color{blue}
\begin{verbatim}
gof(u);
\end{verbatim}
}
\end{minipage}
\definecolor{labelcolor}{RGB}{100,0,0}
\begin{math}\displaystyle
\parbox{8ex}{\color{labelcolor}(\%o4) }
\mathrm{cos}\left( 2\,\left( {u}^{2}+1\right) \right) 
\end{math}

\noindent
\begin{minipage}[t]{8ex}{\color{red}\bf
\begin{verbatim}
(%i5) 
\end{verbatim}
}
\end{minipage}
\begin{minipage}[t]{\textwidth}{\color{blue}
\begin{verbatim}
fog(x):=f(g(x))$
\end{verbatim}
}
\end{minipage}

\noindent
\begin{minipage}[t]{8ex}{\color{red}\bf
\begin{verbatim}
(%i6) 
\end{verbatim}
}
\end{minipage}
\begin{minipage}[t]{\textwidth}{\color{blue}
\begin{verbatim}
fog(u);
\end{verbatim}
}
\end{minipage}
\definecolor{labelcolor}{RGB}{100,0,0}
\begin{math}\displaystyle
\parbox{8ex}{\color{labelcolor}(\%o6) }
{\mathrm{cos}\left( 2\,u\right) }^{2}+1
\end{math}
\end{maxcode}

Anal\'iticamente es claro que $g\circ f(u)\neq f\circ g(u)$, pero
esto se v\'e a\'un m\'as dr\'asticamente comparando los gr\'aficos.
N\'otese como es posible superponer varios gr\'aficos en Maxima, basta
con escribir las funciones a representar en una lista\index{Maxima!lista}:
\begin{maxcode}
\noindent
\begin{minipage}[t]{8ex}{\color{red}\bf
\begin{verbatim}
(%i7) 
\end{verbatim}
}
\end{minipage}
\begin{minipage}[t]{\textwidth}{\color{blue}
\begin{verbatim}
wxplot2d([gof(y),fog(y)],[y,-3,3]);
\end{verbatim}
}
\end{minipage}
\definecolor{labelcolor}{RGB}{100,0,0}
\begin{math}\displaystyle
\parbox{8ex}{\color{labelcolor}(\%t7) }\centering
\includegraphics[width=9cm]{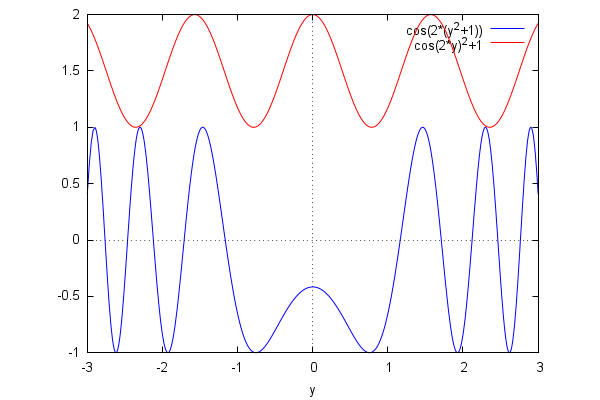}
\end{math}
\begin{math}\displaystyle
\parbox{8ex}{\color{labelcolor}(\%o7) }
\end{math}
\end{maxcode}

Para estudiar la paridad de una funci\'on, recurrimos a lo siguiente.
La funci\'on $h(x)$ es par si $h(x)=h(-x)$, que es equivalente a
$h(x)-h(-x)=0$. El comando \texttt{zeroequiv(}\emph{expr}\texttt{,x)} comprueba si
la expresi\'on \emph{expr}, dependiente de $x$, es equivalente a cero
tras aplicar simplificaciones sint\'acticas y algebraicas, por lo que,
aplicado a $h(x)-h(-x)$, devolver\'a \emph{true} si la funci\'on es par.
An\'alogamente, para estudiar la imparidad comprobaremos si
$h(x)+h(-x)$ es equivalente a cero o no. Los siguientes comandos ilustran
esta t\'ecnica:
\begin{maxcode}
\noindent
\begin{minipage}[t]{8ex}{\color{red}\bf
\begin{verbatim}
(%i8) 
\end{verbatim}
}
\end{minipage}
\begin{minipage}[t]{\textwidth}{\color{blue}
\begin{verbatim}
h(x):=x/(x^2+1)$
\end{verbatim}
}
\end{minipage}

\noindent
\begin{minipage}[t]{8ex}{\color{red}\bf
\begin{verbatim}
(%i9) 
\end{verbatim}
}
\end{minipage}
\begin{minipage}[t]{\textwidth}{\color{blue}
\begin{verbatim}
zeroequiv(h(x)-h(-x),x);
\end{verbatim}
}
\end{minipage}
\definecolor{labelcolor}{RGB}{100,0,0}
\begin{math}\displaystyle
\parbox{8ex}{\color{labelcolor}(\%o9) }
false
\end{math}
\end{maxcode}

El resultado obtenido (\emph{false}) nos indica que la diferencia $h(x)-h(-x)$ no es cero.
Lo comprobamos pidi\'endole a Maxima que nos diga expl\'icitamente qu\'e vale esa diferencia,

\begin{maxcode}
\noindent
\begin{minipage}[t]{8ex}{\color{red}\bf
\begin{verbatim}
(%i10) 
\end{verbatim}
}
\end{minipage}
\begin{minipage}[t]{\textwidth}{\color{blue}
\begin{verbatim}
h(x)-h(-x);
\end{verbatim}
}
\end{minipage}
\definecolor{labelcolor}{RGB}{100,0,0}
\begin{math}\displaystyle
\parbox{8ex}{\color{labelcolor}(\%o10) }
\frac{2\,x}{{x}^{2}+1}
\end{math}
\end{maxcode}
Para la imparidad procederemos de manera an\'aloga, esta vez estudiando la suma
$h(x)+h(-x)$:
\begin{maxcode}
\noindent
\begin{minipage}[t]{8ex}{\color{red}\bf
\begin{verbatim}
(%i11) 
\end{verbatim}
}
\end{minipage}
\begin{minipage}[t]{\textwidth}{\color{blue}
\begin{verbatim}
zeroequiv(h(x)+h(-x),x);
\end{verbatim}
}
\end{minipage}
\definecolor{labelcolor}{RGB}{100,0,0}
\begin{math}\displaystyle
\parbox{8ex}{\color{labelcolor}(\%o11) }
true
\end{math}
\end{maxcode}

Por \'ultimo, veamos el aspecto de $h(x)$. El gr\'afico muestra claramente
su car\'acter de funci\'on impar:
\begin{maxcode}
\noindent
\begin{minipage}[t]{8ex}{\color{red}\bf
\begin{verbatim}
(%i12) 
\end{verbatim}
}
\end{minipage}
\begin{minipage}[t]{\textwidth}{\color{blue}
\begin{verbatim}
wxplot2d(h(x),[x,-5,5]);
\end{verbatim}
}
\end{minipage}
\definecolor{labelcolor}{RGB}{100,0,0}
\begin{math}\displaystyle
\parbox{8ex}{\color{labelcolor}(\%t12) }\centering
\includegraphics[width=9cm]{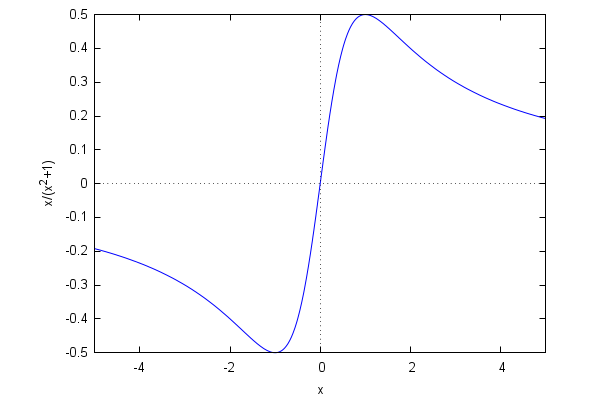} 
\end{math}
\begin{math}\displaystyle
\parbox{8ex}{\color{labelcolor}(\%o12) }
\end{math}
\end{maxcode}

En GeoGebra tambi\'en es muy sencillo comprobar si una funci\'on dada es par o impar.
Tomando de nuevo la funci\'on $h(x)=x/(1+x^2)$ como ejemplo, la representamos en GeoGebra
introduciendo \texttt{h(x):=x/(1+x\^\ 2)} en la entrada de comandos (la ventana inferior). Al pulsar
\texttt{Enter} en el teclado se
visualizar\'a la gr\'afica de la funci\'on. Para ver si es par, seleccionamos la herramienta
\texttt{Simetr\'ia Axial}\index{GeoGebra!simetr\'ia axial} del men\'u, que requiere que seleccionemos primero el objeto a reflejar
(en este caso, picamos sobre la gr\'afica de la funci\'on) y luego el eje de simetr\'ia (en este
caso, el eje de ordenadas). Si la funci\'on fuera par, el resultado ser\'ia la misma gr\'afica
que ya ten\'iamos, lo que no se observa en nuestro ejemplo.
\begin{center}
\includegraphics[scale=0.4]{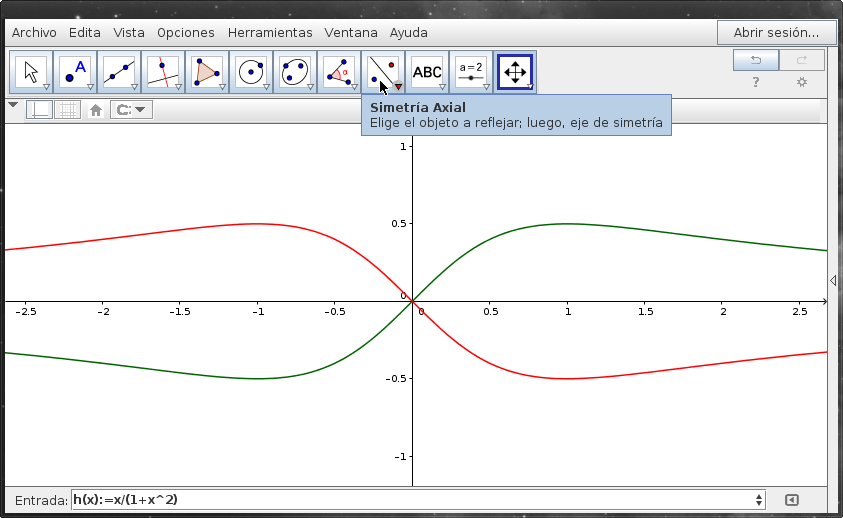} 
\end{center}
Un proceso similar se sigue en el caso de querer estudiar si la funci\'on es impar, el \'unico cambio
es que tras representar $h(x)$ se grafica tambi\'en la recta $y=x$ (escribiendo por ejemplo
\texttt{k(x):=x} en la entrada de comandos) y en la herramienta \texttt{Simetr\'ia Axial} se
selecciona esta recta como el eje de simetr\'ia.
\begin{center}
\includegraphics[scale=0.4]{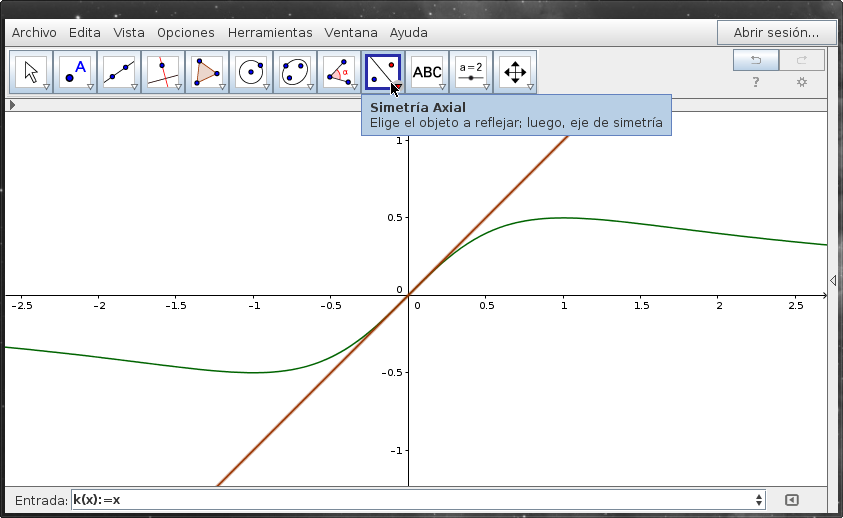} 
\end{center}

\section{Actividades de c\'omputo}
\begin{enumerate}[(A)]
\item Representar gr\'aficamente las funciones del ejercicio \ref{ejer65} y discutir
la existencia de la funci\'on inversa en cada caso.
\item Calcular con Maxima las composiciones de funciones del ejercicio \ref{ejer67}.
\item Representar gr\'aficamente las funciones del ejercicio \ref{ejer68} y
discutir su crecimiento, comparando los resultados con los obtenidos te\'oricamente.
\item Estudiar con Maxima la paridad de las funciones del ejercicio \ref{ejer69}.
Representarlas gr\'aficamente y comparar el gr\'afico con el resultado obtenido.
\item Recordemos que la distancia alcanzada por un proyectil lanzado con velocidad
inicial $v_0$ y \'angulo $\alpha$, viene dada por
$$
 \frac{v^2_0}{g}\,\sin (2\alpha)\,.
$$
Para una misma velocidad inicial de $v_0=20$m/s (y aproximando $g\simeq 10$m/s$^2$),
representar en un mismo gr\'afico
la distancia alcanzada en funci\'on del \'angulo $\alpha$ y estudiar gr\'aficamente
cu\'ando se alcanza la m\'axima distancia\footnote{Para usar las funciones trigonom\'etricas
de Maxima con \'angulos medidos en grados, hay que hacer una
peque\~na modificaci\'on y escribir el seno de $x$ como \texttt{sin(\%pi*x/180)}.}.
Considerar \'angulos entre $0\degree$ y $90\degree$. Repetir el ejercicio pero ahora
con distintas velocidades, entre $v_0=5$m/s y $v_0=50$m/s (variando de $5$m/s en $5$m/s).
%
%
%
%
%
\end{enumerate}

\addtocontents{toc}{\protect\newpage}
\chapter{Exp, log y valor absoluto}\label{explog}
\vspace{2.5cm}\setcounter{footnote}{0}
\section{Nociones b\'asicas}
Las funciones exponenciales probablemente sean las m\'as integradas en la vida cotidiana.
Todo el mundo ha escuchado alguna vez la expresi\'on `crecimiento exponencial' para referirse a
algo que crece muy deprisa. Aunque las funciones exponenciales no son las que presentan
el crecimiento m\'as r\'apido en el mundo de las matem\'aticas\footnote{El lector interesado puede
buscar informaci\'on sobre las funciones superfactoriales, la funci\'on de Ackermann o el
`castor ocupado' (busy beaver).}, ciertamente modelan algunos de los fen\'omenos de crecimiento m\'as 
espectaculares en la naturaleza.

Una \emph{funci\'on exponencial}\index{Funci\'on!exponencial} es una de la forma $f(x)=a^x$. Al n\'umero $a$ se le denomina
\emph{base} y $x$ es el \emph{exponente}. En particular, cuando se habla de \emph{la} funci\'on
exponencial, sin especificar la base, se entiende la funci\'on $f(x)=e^x$, donde la base est\'a dada
por la \emph{constante (o n\'umero) de Euler}\index{Euler, Leonhard!n\'umero}\index{Num@N\'umero!e@$e$} $e=2{.}7182\ldots$, uno de los n\'umeros m\'as importantes de las 
matem\'aticas\footnote{Los `n\'umeros sagrados de las matem\'aticas' son $0,1,\pi,i,e$. 
Existe entre ellos una preciosa relaci\'on,
descubierta precisamente por Euler: $e^{i\pi}+1=0$\index{Euler, Leonhard!f\'ormula}. Si es absolutamente imprescindible que se tat\'uen
alguna parte del cuerpo, esta f\'ormula es un buen candidato para el dise\~no:
\url{https://thedailyomnivore.net/2010/08/16/scarification/}.}

Para apreciar la rapidez de crecimiento de las funciones exponenciales, nada mejor que recordar la leyenda
del ajedrez\index{Ajedrez}, recogida por un poeta persa en los comienzos del siglo X: un matem\'atico hind\'u, Sissa\index{Sissa},
cre\'o un juego (el que hoy conocemos como ajedrez) que consider\'o lo suficientemente interesante como
para mostr\'arselo a su se\~nor. \'Este reconoci\'o que, en efecto, el juego era de lo m\'as interesante;
tras explicarle las reglas, el matem\'atico jug\'o la primera partida contra el se\~nor y le derrot\'o en
4 movimientos (una partida conocida como el \emph{mate pastor}\index{Mate pastor}, 
porque en versiones posteriores del poema,
de corte m\'as rom\'antico, era un pastorcillo el que jugaba contra un rey). El se\~nor, agradecido por
aquella fuente inagotable de entretenimiento, le prometi\'o a Sissa que le dar\'ia aquello que le pidiera,
a lo que aqu\'el replic\'o que se conformaba con que le diese la cantidad de trigo que se obtuviese de la
siguiente manera. El tablero de ajedrez consta de $8\times 8=64$ casillas blancas y negras alternadas. Por
la primera casilla, un grano de trigo. Por la segunda, $2^1=2$ granos. Por la tercera $2^2=4$ granos y
as\'i hasta agotar las casillas.
\begin{center}
\includegraphics[scale=0.4]{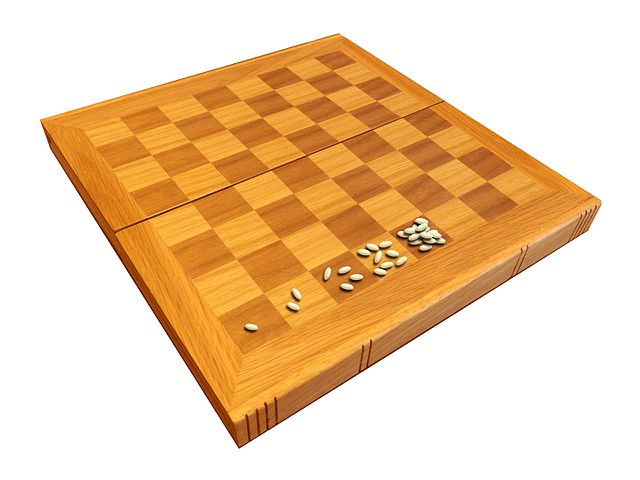} 
\end{center}
El se\~nor, desencantado por aquella petici\'on tan humilde que parec\'ia desde\~nar
su riqueza, le dijo a sus administradores que le dieran a Sissa el pu\~nado de trigo que 
ped\'ia y se retir\'o. Al d\'ia siguiente tuvo necesidad de recurrir a sus administradores 
para un asunto econ\'omico y los mand\'o llamar, pero ellos a\'un estaban ocupados calculando 
el n\'umero de granos de trigo y pensando en d\'onde los conseguir\'ian.

La funci\'on que determina el n\'umero de granos es una exponencial de base 2 y exponente $x-1$,
$f(x)=2^{x-1}$. Los valores que resultan al ir colocando granos en las casillas del tablero son
\begin{align*}
f(1)=&2^0=1 \\
f(2)=&2^1=2 \\
f(3)=&2^2=4 \\
&\cdots 
\end{align*}
Por ejemplo, al final de la primera fila (que tiene 8 casillas) tendremos
$f(8)=2^7=128$, una cantidad rid\'icula, desde luego. Al final de la segunda fila ser\'an
$f(16)=2^{15}=32\, 768$, una cantidad mayor pero todav\'ia bastante reducida.
El peso de un grano de trigo es aproximadamente
0.05 gramos, as\'i que esto equivale a poco m\'as de un kilo y medio de trigo (1 638.4 gramos). 
Al final de la tercera fila, las cosas se ponen m\'as interesantes:
$f(24)=2^{23}=8\, 388\, 608$, unos 420 kilos de trigo. Para la cuarta fila,
los administradores ya empezaban a ponerse nerviosos:
$f(32)=2^{31}=2\, 147\, 483\, 648$, que son unas 107 toneladas.
Para no hacerlo largo, aqu\'i tenemos lo que corresponde s\'olo a la 
\'ultima casilla (a lo que habr\'ia que a\~nadir lo de las dem\'as):
$$
f(64)=2^{63}\simeq 9{.}2\cdot 10^{18},
$$
o sea, unos $4{.}6\cdot 10^{14}$ kilos de trigo, m\'as de 400 mil millones de toneladas (la producci\'on
de trigo en la India en 2012 fue de unos 96 millones de toneladas\footnote{Datos tomados
de Wikipedia: \url{http://en.wikipedia.org/wiki/International_wheat_production_statistics}}.
Toda la India cultivando trigo
durante 4 000 a\~nos dar\'ia para el trigo de esta \'ultima casilla).

Hasta aqu\'i hemos utilizado valores naturales en el exponente, con lo cual est\'a claro el significado
de la exponencial; por ejemplo,
$$
2^6=2\cdot 2\cdot 2\cdot 2\cdot 2\cdot 2 = 64.
$$
Para valores racionales positivos del exponente, $p/q >0$, la funci\'on exponencial se define como
\begin{equation}\label{exp1}
a^{p/q}=\sqrt[q]{a^p}=\left( \sqrt[q]{a} \right)^p,
\end{equation}
 para valores racionales negativos,
\begin{equation}\label{exp2}
a^{-p/q}=1/a^{p/q} .
\end{equation}
Si el exponente $x$ es irracional, como por ejemplo $x=\sqrt{2}$, hay que utilizar un proceso de l\'imite
que se estudia en cursos m\'as avanzados de c\'alculo. La idea es que todo n\'umero real, incluyendo los
irracionales, se puede expresar como l\'imite de una sucesi\'on de n\'umeros racionales que lo
aproximan; en el caso de $\sqrt{2}$, tendr\'iamos la sucesi\'on de racionales
$\{1,1{.}4,1{.}41,1{.}414,1{.}4142,\ldots\}$. As\'i, por ejemplo, para calcular $3^{\sqrt{2}}$ formar\'iamos la
sucesi\'on $\{ 3^1,3^{1{.}4},3^{1{.}41},\ldots\}$ que da
$$
\{ 3, 4{.}655, 4{.}707,\ldots\},
$$
lo cual nos va aproximando al valor de $3^{\sqrt{2}}\simeq 4{.}7288$.

Con esto, podemos graficar las funciones exponenciales $f(x)=a^x$. Su aspecto para valores
$a>1$, es\footnote{Aqu\'i aparece de nuevo el valor $e\simeq 2{.}7182...$ del
\emph{n\'umero de Euler}.}
\begin{center}
\begin{tikzpicture}
\draw[thick,->] (-2,0) -- (2.5,0) node[right]{$x$};
\draw[thick,->] (0,0) -- (0,5) node[above]{$y$};
\draw[red,thick,domain=-2:2.15] plot (\x,{pow(2,\x)}) node[above]{$2^x$};
\draw[orange,thick,domain=-2:0.65] plot (\x,{pow(10,\x)}) node[above]{$10^x$};
\draw[blue,thick,domain=-2:1.5] plot (\x,{exp(\x)}) node[above]{$e^x$};
\end{tikzpicture}
\end{center}
A la vista de este gr\'afico, es claro que cualitativamente el comportamiento de las funciones exponenciales
es el mismo; la base \'unicamente determina la velocidad de crecimiento, pero este presenta unas 
caracter\'isticas comunes:
\begin{list}{$\blacktriangleright$}{}
\item La gr\'afica pasa por el punto $(0,1)$.
\item El dominio es todo $\mathbb{R}$ y el rango $\mathbb{R}^+$ (los n\'umeros reales \emph{positivos}\index{R+@$\mathbb{R}^+$}, es decir $\mathbb{R}^+=\{x\in\mathbb{R}:x>0\}$.
\item Son funciones crecientes.
\item Hacia la derecha se extienden hacia $\infty$ y hacia la izquierda tienden a $0$.
\end{list}
Para las funciones con base $0<a<1$, el comportamiento es distinto. Gr\'aficamente:
\begin{center}
\begin{tikzpicture}
\draw[thick,->] (-2.5,0) -- (2,0) node[right]{$x$};
\draw[thick,->] (0,0) -- (0,5) node[above]{$y$};
\draw[red,thick,domain=-2.15:1.9] plot (\x,{pow(0.5,\x)}) 
node[below]at (-2.25,5.5) {$\left( \frac{1}{2} \right)^x$};
\draw[orange,thick,domain=-1.24:1.9] plot (\x,{pow(0.3,\x)}) 
node[below]at (-1.4,5.5) {$\left( \frac{1}{3} \right)^x$};
\draw[blue,thick,domain=-0.93:1.9] plot (\x,{pow(0.2,\x)}) 
node[below]at (-0.7,5.5) {$\left( \frac{1}{5} \right)^x$};
\end{tikzpicture}
\end{center}
Las caracter\'isticas que se observan, son:
\begin{list}{$\blacktriangleright$}{}
\item La gr\'afica pasa por el punto $(0,1)$.
\item El dominio es todo $\mathbb{R}$ y el rango $\mathbb{R}^+$.
\item Son funciones decrecientes.
\item Hacia la derecha tienden a $0$ y hacia la izquierda se extienden hacia $\infty$.
\end{list}
Las reglas de c\'alculo con funciones exponenciales son muy sencillas y semejantes a las
de las funciones potenciales. Para cualquier base positiva $a>0$ y cualesquiera exponentes
$x,y\in\mathbb{R}$, se cumple:
\begin{list}{$\blacktriangleright$}{}
\item $a^xa^y=a^{x+y}$,
\item 
$\dfrac{a^x}{a^y}=a^{x-y}$,
\item $(a^x)^y=a^{xy}$,
\item $a^0=1$,
\item $a^{-x}=\dfrac{1}{a^x}$.
\end{list}
La deducci\'on de estas reglas es inmediata para el caso de valores naturales de los exponentes. A
partir de aqu\'i, se deducen para exponentes racionales usando \eqref{exp1} y \eqref{exp2}. Por
ejemplo,
$$
a^{p/q}a^{r/s}=a^{ps/qs}a^{rq/qs}=\sqrt[qs]{a^{ps}}\sqrt[qs]{a^{rq}}=\sqrt[qs]{a^{ps+rq}}
=a^{(ps+rq)/qs}=a^{p/q+r/s}.
$$
De aqu\'i es inmediato (usando \eqref{exp2}) que
$$
\frac{a^{p/q}}{a^{r/s}}=a^{p/q}a^{-r/s}=a^{p/q-r/s}.
$$
Las restantes propiedades se siguen entonces sin dificultad. Tomando $x=y$ en la f\'ormula que 
acabamos de deducir:
$$
a^0=a^{x-x}=\frac{a^x}{a^x}=1 ,
$$
y tomando $x=0$:
$$
a^{-y}=a^{0-y}=\frac{a^0}{a^y}=\frac{1}{a^y} .
$$

Las funciones logar\'itmicas\index{Funci\'on!logar\'itmica} se definen como las inversas de las exponenciales. M\'as concretamente,
dado un n\'umero real $b>0$ diferente de la unidad, $b\neq 1$, el \emph{logaritmo en base $b$ del
n\'umero}\index{Logaritmo} $x>0$ es el n\'umero $y=\log_b x$ tal que
\begin{equation}\label{logar}
b^y=x .
\end{equation}
Un caso particular importante es el logaritmo inverso de la funci\'on exponencial (con base $e$). Es
el n\'umero $y=\log_e x$ tal que $e^y=x$, y se denota $\ln x$ en lugar de $\log_e x$ y se llama
logaritmo \emph{natural} o \emph{neperiano}\index{Logaritmo!natural}\index{Logaritmo!neperiano},
por su descubridor John Napier\index{Napier, John}.\\
Fij\'emonos en que s\'olo podemos calcular logaritmos de n\'umeros positivos $x>0$. Esto es as\'i
porque cualquier logaritmo es funci\'on inversa\index{Funci\'on!inversa} de una exponencial, y \'estas tienen rango
$\mathbb{R}^+$. Por contra, el rango de cualquier logaritmo es todo $\mathbb{R}$.

As\'i como las exponenciales crecen muy deprisa, los logaritmos (que son sus inversas) lo hacen muy
lentamente\index{Funci\'on!creciente}. Gr\'aficamente, el comportamiento de las funciones logar\'itmicas es:
\begin{center}
    \begin{tikzpicture}[auto,scale=0.85] 
     \draw [thick,->] (-2, 0) -- (4.25, 0) node[right]{$x$};
     \draw [thick,->] (0, -3.25) -- (0, 3) node[right]{$y$};
     \draw [orange,domain=0.05:4.2,samples=100,thick] plot(\x,{ln(\x)}) node[right]{$\ln (x)$};
  	 \draw [red,domain=0.05:4.2,samples=100,thick] plot(\x,{log10(\x)}) node[right]{$\mathrm{log}_{10}(x)$};  		 
  	 \draw [blue,domain=0.1:4.2,samples=100,thick] plot(\x,{log2(\x)}) node[right]{$\mathrm{log}_{2}(x)$};
    \end{tikzpicture}
\end{center}
Las caracter\'isticas m\'as importantes deducidas de las gr\'aficas de los logaritmos con base $b>1$ son
las siguientes:
\begin{list}{$\blacktriangleright$}{}
\item La gr\'afica pasa por el punto $(1,0)$ (?`por qu\'e?).
\item El dominio es $\mathbb{R}^+$ y el rango $\mathbb{R}$.
\item Son funciones crecientes.
\item Al acercarnos a $0$, las gr\'aficas tienden a $-\infty$, y hacia la derecha tienden a $\infty$.
\end{list}

Los logaritmos tienen unas reglas de c\'alculo que se deducen de las de las funciones exponenciales.
Son muy sencillas; para todos $x,y\in\mathbb{R}^+$, se cumple:
\begin{list}{$\blacktriangleright$}{}
\item $log_b (xy)=log_b(x)+\log_b(y)$.
\item $log_b\left(\dfrac{x}{y}\right)=log_b(x) -\log_b(y)$.
\item $log_b(x^y)=y\log_b(x)$, para cualquier $y\in\mathbb{R}$.
\item $\log_b 1=0$.
\item $\log_b\left(\dfrac{1}{x}\right)=-\log_b (x)$.
\item (Cambio de base) $log_b(x)=\log_a(x)\cdot \log_b(a)$, para $a>1$.
\end{list}
Veamos alguna de ellas. Dados dos n\'umeros $x,y\in\mathbb{R}^+$, pongamos $r=\log_b x$, $s=\log_b y$,
de modo que seg\'un la definici\'on \eqref{logar} es $b^r=x$ y $b^s=y$. Entonces, por las propiedades 
de las exponenciales:
$$
xy=b^rb^s=b^{r+s},
$$
y de nuevo por la definici\'on \eqref{logar} esto implica que
$$
\log_b(xy)=r+s=\log_b x+\log_b y.
$$
Por otra parte, si hacemos
$$
\frac{x}{y}=\frac{b^r}{b^s}=b^{r-s},
$$
resulta que, por \eqref{logar},
$$
log_b(x/y)=r-s=\log_bx-\log_by .
$$
?`Puede el lector probar las restantes propiedades?

El \emph{valor absoluto}\index{Valor absoluto} es una funci\'on con dominio $\mathbb{R}$ y codominio $\mathbb{R}^+$, definida
como
$$
|x|=
\begin{dcases}
x &\mbox{ si }x\geq 0 \\
-x &\mbox{ si }x\leq 0 .
\end{dcases}
$$
Su gr\'afica es la siguiente:
\begin{center}
    \begin{tikzpicture}[auto,scale=1] 
     \draw [->,thick] (-3, 0) -- (3, 0) node[right]{$x$};
     \draw [->,thick] (0, -0.5) -- (0, 3) node[right]{$y$};
  	 \draw [red,domain=-2.9:2.9,samples=100,thick] plot(\x,{abs(\x)});  		
    \end{tikzpicture}
\end{center}
La importancia del valor absoluto radica en que puede verse como una \emph{distancia}\index{Distancia} en la
recta real. Dado un n\'umero cualquiera $x\in\mathbb{R}$, $|x|$ da la distancia a la que se
encuentra del origen $0\in\mathbb{R}$.

Ampliando esta idea, podemos decir que el valor absoluto $|x-y|$ mide la separaci\'on
que hay entre $x$ e $y$ independientemente de la posici\'on que ocupen en la recta real.
De esta forma, tendremos una funci\'on que \emph{a cada par de puntos} en la recta, le
asigna otro n\'umero, la \emph{distancia} entre ellos. Esta funci\'on distancia ser\'a
del tipo $d:\mathbb{R}\times\mathbb{R}\to\mathbb{R}$, y su actuaci\'on se escribir\'a
$$
d(x,y)=\vert x-y\vert\,.
$$
Veamos las propiedades de esta funci\'on distancia.
\begin{list}{$\blacktriangleright$}{}
\item La distancia entre dos n\'umeros siempre es mayor o igual que cero.
\end{list}
Esto es obvio por la propia definici\'on de valor absoluto, que siempre es un n\'umero
mayor o igual que cero. De hecho, el \'unico n\'umero cuyo valor absoluto es cero es el
propio $0$, por lo que tenemos otra propiedad.
\begin{list}{$\blacktriangleright$}{}
\item La distancia entre dos n\'umeros $x,y\in\mathbb{R}$ es cero si y s\'olo si $x=y$.
\end{list}
Geom\'etricamente esperar\'iamos que la distancia entre 3 y 5 fuera la misma que la
que hay entre -5 y -3, y as\'i ocurre:
$$
|5-3|=|2|=2=-(-2)=|-2|=|-5-(-3)| .
$$
Esta propiedad es cierta para cualquier par de n\'umeros, pues $|-1|=1$ y
$|x-y|=|(-1)(y-x)|=|-1||y-x|$ (el elctor probar\'a sin dificultad que 
$|x\cdot y|=|x||y|$).
\begin{list}{$\blacktriangleright$}{}
\item La distancia entre dos n\'umeros $x,y\in\mathbb{R}$ cumple la \emph{propiedad
de simetr\'ia}
$$
d(x,y)=\vert x-y\vert =\vert y-x\vert =d(y,x)\,.
$$
\end{list}
Finalmente, tenemos la llamada \emph{desigualdad triangular}\index{Desigualdad triangular}, una de las propiedades
m\'as importantes para el estudio de los n\'umeros reales.
\begin{list}{$\blacktriangleright$}{}
\item La distancia entre dos n\'umeros $x,y\in\mathbb{R}$ cumple que
$$
d(x,y)=\vert x-y\vert \leq\vert x-z\vert +\vert z-y\vert =d(x,z)+d(z,y)\,,
$$
para cualquier otro $z\in\mathbb{R}$. 
\end{list}
Esta es una propiedad que se deduce
inmediatamente de la desigualdad triangular para el valor absoluto (v\'ease el
ejercicio \ref{triangdes}).

\section{Aplicaciones}
En esta secci\'on vamos a ver la relaci\'on existente entre los logaritmos y
una manera de medir la informaci\'on contenida en un sistema f\'isico o l\'ogico,
a trav\'es de una magnitud conocida como \emph{entrop\'ia}\footnote{La motivaci\'on a la 
definici\'on de entrop\'ia que damos a continuaci\'on, est\'a basada en unas notas no 
publicadas de Jos\'e Vicente Arnau C\'ordoba\index{Arnau, Jos\'e Vicente}, de la Universitat de Val\`encia, que constituyen en nuestra opini\'on la mejor introducci\'on elemental a este fascinante concepto.}\index{Entrop\'ia}.

Supongamos una fuente $F$ de informaci\'on de cualquier tipo (puede ser un peri\'odico,
una emisi\'on de radio, un programa de c\'omputo,
una cadena de ADN, etc.). Esta fuente construir\'a su informaci\'on a partir de $n$
caracteres distintos ($26$ letras en el caso del peri\'odico o la emisi\'on de radio,
los dos caracteres $0$ y $1$ en el caso del programa de c\'omputo, los caracteres
$A$, $C$, $G$ y $T$ en el caso del ADN). Llamaremos \emph{mensaje de longitud} $k$ a
una sucesi\'on aleatoria de $k$ caracteres de entre los $n$ de que dispone la fuente,
posiblemente repetidos.\\
Notemos que en esta definici\'on prescindimos del sentido sem\'antico del mensaje. Es decir,
un mensaje como `akenrbfyde' es tan v\'alido como `ornitorrinco'\footnote{De hecho, se puede
asignar un significado a mensajes aparentemente absurdos, por ejemplo usando acr\'onimos
como `ONU', o notaciones como la del ajedrez, ` f3 e5 g4?? Qh4 '.}.\\
Una caracter\'istica esencial de los mensajes emitidos por $F$ es que \emph{aportan tanta m\'as
informaci\'on cuanto menos previsible son}\index{Informaci\'on}, independientemente de su longitud. Por ejemplo,
un paciente obtiene mucha m\'as informaci\'on del mensaje `las pruebas dieron negativo' que
del mensaje `los resultados no son concluyentes', a pesar de que el segundo es m\'as largo que el
primero.\\
Al tratar con fuentes de informaci\'on, hay que distinguir dos casos. En el primero s\'olo
conocemos el n\'umero $N$ de mensajes que puede emitir la fuente, sin mayores datos sobre
cada mensaje individual. En el segundo caso, se puede asignar a cada mensaje individual
un `peso relativo' frente a los restantes (por ejemplo, su probabilidad de aparici\'on).
\begin{list}{$\blacktriangleright$}{}
\item Mensajes equiprobables \\
En este caso, es sencillo definir una medida de la informaci\'on global que aporta la fuente
$F$. Si puede emitir $N$ mensajes de $k$ caracteres, pondr\'iamos $\tilde{I}(F)=N$. Las
ventajas de esta medida es que es calculable, es positiva y aumenta cuando lo hace $N$, como 
cabr\'ia esperar. Sin embargo, no cumple una propiedad esencial que esperamos de cualquier
medida: lo que se conoce como \emph{aditividad}\index{Aditividad}. Esto quiere decir que esperamos que la
informaci\'on total de dos fuentes, $I(F_1+F_2)$, sea la suma de las informaciones de cada
una por separado $I(F_1)+I(F_2)$. Notemos que si $F_1$ puede emitir $N_1$ mensajes, de longitud 
$k_1$, y $F_2$ puede emitir $F_2$ mensajes, de longitud $k_2$, al unirlas para que emitan 
mensajes de longitud $k_1+k_2$ (por ejemplo, yuxtaponiendo los mensajes, uno a continuaci\'on
del otro) se obtienen $N_1\cdot N_2$ posibilidades, no $N_1+N_2$.\\
?`C\'omo solucionar este problema? La soluci\'on la dan los logaritmos. Si en lugar de usar
$\tilde{I}(F)$ definimos la \emph{medida de informaci\'on global}\index{Informaci\'on!global} de la fuente $F$ como
$$
I(F)=\log (N)\, ,
$$
obtenemos una medida que tiene \emph{casi} todas las propiedades deseables:
como debe ser $N\geq 1$
(en otro caso, no podr\'iamos hablar de `fuente de informaci\'on, pues no habr\'ia mensajes),
se tiene $I(F)\geq 0$. Adem\'as, por las propiedades del logaritmo, si los mensajes de
$F_1$ y $F_2$ son todos \emph{distintos}, ser\'a $I(F_1+F_2)=I(F_1)+I(F_2)$. En general,
sin embargo, habr\'a mensajes que puedan ser emitidos por las dos fuentes, con lo cual
$$
I(F_1+F_2)\leq I(F_1)+I(F_2)\,.
$$  
Esta propiedad se conoce como \emph{subaditividad}\index{Subaditividad}.\\
Hemos mencionado que esta medida tiene \emph{casi} todas las propiedades deseables. ?`Qu\'e
es lo que falta? La \emph{normalizaci\'on}.
Nos gustar\'ia que la informaci\'on m\'inima no trivial
producida por la fuente tuviera una medida $1$, pero esto no ocurre con la definici\'on que
hemos dado. La causa es que la informaci\'on m\'inima \emph{no trivial} que 
puede emitir una fuente no se
da cuando $N=1$ (pues en ese caso como solo hay un mensaje posible, ya sabemos lo que la fuente
va a emitir, luego la informaci\'on es cero). Esta informaci\'on m\'inima se da con 
$N=2$, as\'i que supongamos una fuente que emite dos mensajes, ambos equiprobables. En este caso
ser\'ia, con nuestra definici\'on, $I(F)=\log (2)$. Para que esto valga $1$, debemos tomar
el logaritmo \emph{en base} $2$, y en esto consiste la normalizaci\'on, en fijar la base del
logaritmo\index{Logaritmo} con que se trabaja (v\'ease el ejercicio \ref{baselog}).\\
Con todo esto, podemos definir la \emph{medida de informaci\'on normalizada}\index{Informaci\'on!normalizada}
de la fuente $F$ como
$$
I(F)=\log_2 (N)\, .
$$
En este caso, la unidad de informaci\'on es un \emph{bit}\index{Bit} (un d\'igito binario)\footnote{En
ingl\'es, `\textbf{b}inary dig\textbf{it}'.}: $1\,\mbox{bit}=\log_2(2)$. Por ejemplo,
consideremos como fuente de informaci\'on el lanzamiento de una moneda. Los mensajes
posibles son $+$ y $-$ (cara o cruz, sello o \'aguila, head o tail), cada uno con la
misma probabilidad igual a $1/2$. De acuerdo con lo expuesto, en el lanzamiento de una
moneda la fuente puede emitir $N=2$ mensajes y por tanto, la medida de informaci\'on
normalizada de este proceso es
$$
I(F)=\log_2(2)=1\,,
$$
es decir, un bit de informaci\'on.\\
Si el proceso consta del lanzamiento de la moneda $k$ veces, se forman mensajes de
logitud $k$ y la fuente puede emitir $2^k$ de ellos. Por tanto, la informaci\'on
obtenida con cada resultado es ahora
$$
I(F)=\log_2(2^k)=k\,,
$$
como cabr\'ia esperar.

\item Mensajes con peso \\
Supongamos que tenemos una fuente de informaci\'on capaz de emitir $N$ mensajes distintos
$m_i$ (con $i\in\{1,2\ldots ,N\}$), a los que asignaremos probabilidades de aparici\'on 
diferentes, $p_i$. Estas probabilidades\index{Probabilidad}, como sabemos, deben cumplir que
$$
\sum^N_{i=1}p_i=1\,.
$$
Por otra parte, notemos que si $0<p_i<1$, es $\log_2(p_i)<0$, de manera que cuando las
$p_i$ representan probabilidades\footnote{Si un mensaje tuviera $p_i=0$ no
podr\'ia emitirse, por lo que podemos tomar $0<p_i$. Y si alguno tuviera $p_i=1$, ser\'ia
el \'unico que se emitir\'ia y toda la discusi\'on carecer\'ia de sentido, por lo que
tambi\'en podemos tomar $p_i<1$.}, $-\log_2(p_i)>0$.\\
Una idea central en teor\'ia de la informaci\'on, es que \emph{cada mensaje puede considerarse
como una microfuente}, cuya medida de informaci\'on est\'a dada por su probabilidad de 
aparici\'on,
$$
I(m_i)=-\log_2(p_i),\quad i\in\{0,1,\ldots,N\}\,.
$$
La justificaci\'on de introducir esta aparente complicaci\'on es que cuando $p_i$ se
aproxima a cero (en cualquier logaritmo), su informaci\'on asociada $I(p_i)$ crece
indefinidamente. Es decir, le damos sentido matem\'atico a la caracter\'istica que mencionamos 
de los mensajes, que aportan tanta m\'as informaci\'on cuanto menos predecibles son.\\
Ahora, podemos definir la informaci\'on global\index{Informaci\'on!global} de la fuente como la media ponderada de
la informaci\'on de cada microfuente:
$$
I(F)=\frac{\sum^N_{i=1}p_i\cdot I(m_i)}{\sum^N_{i=1}p_i}\,,
$$
o bien, por ser $\sum^N_{i=1}p_i=1$,
$$
I(F)=-\sum^N_{i=1}p_i\log_2(p_i)\,.
$$
La medida de informaci\'on expresada en esta forma (para el caso en que cada mensaje
tiene una probabilidad diferente de ser emitido) se denomina \emph{entrop\'ia}\index{Entrop\'ia}, del
griego \emph{en-} y $\tau\rho$\emph{\'o}$\pi$o$\varsigma$ (\emph{tropos}, en rotaci\'on).
Intuitivamente, se habla de la entrop\'ia como de una \emph{medida del desorden}: en una
habitaci\'on desordenada, por ejemplo, conocer d\'onde se encuentra determinado objeto
aporta mucha m\'as informaci\'on que si se tiene todo `en su lugar'.\\
Veamos un ejemplo. La siguiente secuencia es muy popular por ser el comienzo de la secuencia 
del ADN\index{ADN} (por sus siglas \textbf{\'a}cido \textbf{d}esoxirribo\textbf{n}ucleico, en ingl\'es DNA) 
de la bacteria E. Coli, concretamente sus 32 primeros pares de bases:
$$
AGCTTTTCATTCTGACTGCAACGGGCAATATG
$$
Se puede obtener de la p\'agina \url{http://www.ecogene.org/old/SequenceDownload.php}.
En el di\'alogo que aparece, seleccionamos la versi\'on U00096.2 y el rango de bases 1 a 32, 
dejando la opci\'on referente al complemento en `No'. Aparecer\'a un fichero de texto para 
descargar conteniendo la secuencia que mostramos arriba.
Es sabido que a partir del ADN se generan las prote\'inas\index{Prote\'inas}; en concreto, cada bloque de tres bases de entre A, C, G, T (denominado un \emph{cod\'on}\index{Cod\'on}) codifica un amino\'acido\index{Amino\'acidos}, que son los constituyentes primarios de las prote\'inas; (por ejemplo, TTT codifica el amino\'acido fenilalanina\index{Fenilalanina}, uno de los $9$ amino\'acidos esenciales para el ser humano\footnote{La ausencia de fenilalanina es fatal para el ser humano, pero su exceso tambi\'en lo es. Cuando la fenilalanina no se metaboliza adecuadamente, debido a una enfermedad conocida como fenilcetonuria
\index{Fenilcetonuria}, se acumula en el cuerpo causando graves da\~nos al sistema 
nervioso central. La fenilalanina se usa en conjunci\'on con el \'acido asp\'artico
para obtener el popular edulcorante aspartame\index{Aspartame}, habitual en las bebidas refrescantes `light'.}).
Podemos estimar la probabilidad con que aparece cada
cod\'on a partir de nuestro fragmento de ADN, contando:
\begin{center}
\begin{tabular}{|cccccccccc|}
\hline 
AAT & 1 & AAC & 1 & ACG & 1 & ACT & 1 & AGC & 1 \\ 
ATA & 1 & ATG & 1 & ATT & 1 & CAA & 2 & CAT & 1 \\ 
CGG & 1 & CTG & 2 & CTT & 1 & GAC & 1 & GCA & 2 \\ 
GCT & 1 & GGC & 1 & GGG & 1 & TTC & 2 & TAT & 1 \\ 
TCA & 1 & TCT & 1 & TGA & 1 & TGC & 1 & TTT & 2 \\ 
\hline 
\end{tabular} 
\end{center}
As\'i, por ejemplo, la probabilidad del cod\'on AAT es $1/30$, mientras que la del
TTC es $2/30$ (n\'otese que, en el caso del cod\'on TTT, se cuentan dos apariciones
porque hay un segmento TTTT que lo contiene dos veces, \textbf{TTT}T y
T\textbf{TTT}). La entrop\'ia de este fragmento particular de ADN es
$$
-\sum^{30}_{i=1}p_i\log_2(p_i)=-\left( 5\cdot \frac{2}{30}\log_2\left(\frac{2}{30}\right)
+20\cdot \frac{1}{30}\log_2\left(\frac{1}{30}\right)\right) =4{.}5736\mbox{ bits}\,.
$$
\begin{center}
\includegraphics[scale=1]{E_coli.jpg} 
\end{center}
Aunque ir m\'as all\'a en este tema queda fuera del objetivo de estas notas,
debemos resaltar que este tipo de an\'alisis es muy
\'util en gen\'omica\index{Gen\'omica}. Por ejemplo, es posible probar que la m\'axima entrop\'ia
para un fragmento se obtiene cuando todos los codones son equiprobables\footnote{La prueba de 
este hecho es muy sencilla cuando se dispone de la maquinaria del \emph{c\'alculo
diferencial}\index{C\'alculo diferencial} (lo que deber\'ia ser un est\'imulo para su estudio). En general, sea $\mathbf{p}$
el vector de probabilidades de todos los mensajes, es decir, 
$\mathbf{p}=(p_1,\ldots ,p_N)$ con $p_i$
la probabilidad del mensaje $m_i$. Para probar que \emph{de todas las fuentes capaces de
emitir $N$ mensajes distintos, la que tiene todos sus mensajes equiprobables ($p_i=m_i$)
es la de informaci\'on m\'axima}, basta con ver que $I(\mathbf{p})=-\sum_i p_i\log_2(p_i)$
alcanza un m\'aximo cuando $p_i=1/N$ para todo $i\in\{1,\ldots,N\}$. Cambiando de base $2$
a base $N$ en el logaritmo y exigiendo que la derivada parcial de $I(\mathbf{p})$ respecto
a cada $p_j$ se anule, resulta:
$$
\frac{\partial I}{\partial p_j}=\frac{-1}{\log_N(2)}(\log_N(p_j)+1)=0\,,
$$
y esto ocurre si y s\'olo si $-\log_N(p_j)=1$, esto es, $p_j=1/N$.}  y su valor
es, entonces, $\log_2(64)=6\mbox{ bits}$. As\'i, podr\'iamos decir que el fragmento 
del ADN de la E. Coli que hemos elegido presenta poca informaci\'on, pues su entrop\'ia est\'a 
relativamente lejana al valor m\'aximo. Esto es algo
que est\'a en consonancia con el hecho de que una buena parte del ADN carece de 
funci\'on codificadora, lo que se conoce como \emph{ADN basura}, aunque esto no 
significa que una porci\'on de ADN con poca entrop\'ia no sea codificante ni tampoco que
no sea interesante el estudio del ADN basura\index{ADN!basura} (un excelente repaso de este tema se da en el
art\'iculo de Scientific American `Hidden treasures in junk DNA', disponible en 
\url{https://www.scientificamerican.com/article/hidden-treasures-in-junk-dna/}).
\end{list}

Las bases de la teor\'ia de la informaci\'on, tal y como las hemos expuesto aqu\'i, fueron
enunciadas por C. Shannon\index{Shannon, Claude}, un matem\'atico estadounidense del Instituto de Estudios Avanzados
de Princeton que cooperaba con los laboratorios Bell, en 1948. La teor\'ia tiene un alcance
extraordinario, en ramas tan dispares como la sociolog\'ia o la mec\'anica estad\'istica, y es
uno de los pilares del desarrollo de las ciencias de la computaci\'on actuales. Todo ello
gracias a la modesta contribuci\'on de los logaritmos, que fueron inventados por J. Napier\index{Napier, John}
(un astr\'onomo) en 1614, para poder hacer c\'alculos de productos de n\'umeros naturales
con rapidez.

\section{Ejercicios}

\begin{enumerate}  
\setcounter{enumi}{69}

\item Para cada apartado, representar \emph{cualitativamente} las funciones que se indican
en un mismo gr\'afico:
\begin{multicols}{2}
\begin{enumerate}[(a)]
\item $f(x)=2^x$ y $g(x)=2^{-x}$
\item $f(x)=3^x$ y $g(x)=e^x$
\item $f(x)=\ln x,\; g(x)=\log_2 x$ y $h(x)=\log_{10}x$
\item $f(x)=e^x$ y $g(x)=\ln x$
\end{enumerate}
\end{multicols}

\item\label{baselog} Probar la propiedad de cambio de base de los logaritmos: 
para cualesquiera $x>0$, $a,b>1$,  
$$
log_b(x)=\log_a(x)\cdot \log_b(a) .
$$

\item\label{ejer72} Resolver las ecuaciones y sistemas siguientes:
\begin{multicols}{3}
\begin{enumerate}[(a)]
\item $e^{3x}=2$
\item $e^x + 12e^{-x}=7$
\item $5^x =7$
\item $2^{x-3}=16$
\item $7^{5x+1}=49^{3x+2}$
\item $3^{2-x^2}=3$
\item $4^{2x+1}=(1/2)^{3x+5}$
\item $9^x - 4\cdot 3^x + 3=0$
\item
$
\begin{cases}
2^x - 3^{y-1}=5 \\
2^{x+1}+8\cdot 3^y =712
\end{cases}
$
\item
$
\begin{cases}
2^x +2^y =24 \\
2^{x+y}=128
\end{cases}
$
\end{enumerate}
\end{multicols}

\item\label{ejer73} Resolver las ecuaciones siguientes:
\begin{multicols}{3}
\begin{enumerate}[(a)]
\item $log_x 8=3$
\item $log_2 x=7$
\item $log_3 x=-1$
\item $log_{27} x=-{\displaystyle \frac{ 2}{ 3}}$
\item $\ln x^4 =-1$ 
\item $\ln(\ln x)=2$ 
\item $\log_2 x^2 +\log_2 x=4$ 
\item $\log_2 2^{4x}=20$
\item $\log_3 x+\log_3 (x-2)=1$
\item $\log_8 3+\frac{1}{2}\log_8 25=\log_8 x$ 
\end{enumerate}
\end{multicols}

\item Cuando sea posible, escribir las expresiones siguientes como una \'unica exponencial:
\begin{multicols}{4}
\begin{enumerate}[(a)]
\item $e^2 \cdot e^{x+1}\cdot e^{y-3}$
\item $2^{3x}\cdot 3^x$
\item $\sqrt{4^{3x-1}}$
\item $\frac{{\displaystyle 5}^2}{\displaystyle 25^{3x-2}}$
\end{enumerate}
\end{multicols}

\item Cuando sea posible, escribir las siguientes expresiones como un \'unico logaritmo:
\begin{multicols}{2}
\begin{enumerate}[(a)]
\item $ \ln x+2\ln y$
\item $ \log (x+1)-\frac12 \log x-\log (x^2 -1)$
\item $ \log_3 (x-2)-\log_3 (x+2)$
\item $ 3\ln x+2\ln y-4\ln z$
\end{enumerate}
\end{multicols}

\item La intensidad de los terremotos se puede medir de varias formas. Una de las m\'as populares
(aunque no la mejor t\'ecnicamente) es la escala de Richter-Gutenberg, que se basa en las mediciones
hechas por un tipo de sism\'ografo denominado de Wood-Anderson.\index{Richter, Charles Francis}
\begin{center}
\includegraphics[scale=0.2]{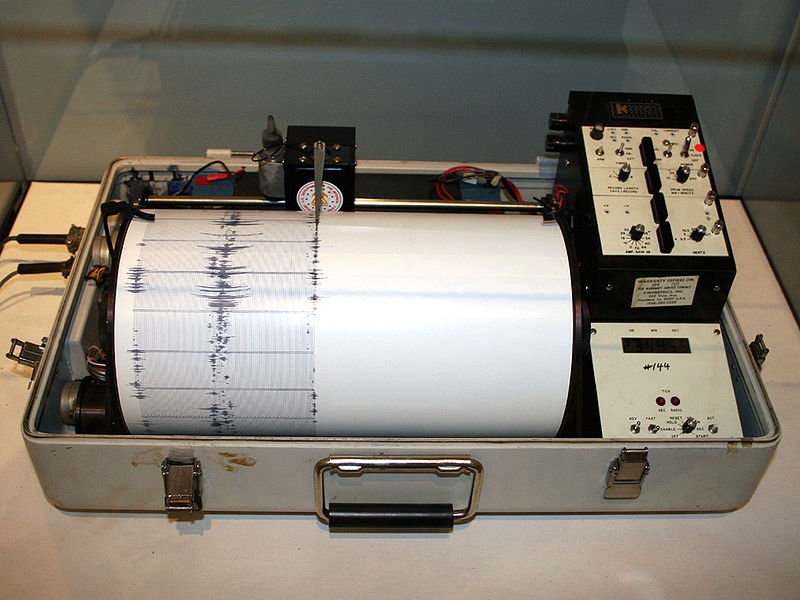} 
\end{center}
La escala asigna un n\'umero $M_L$ a cada terremoto bas\'andose en la lectura de la amplitud m\'axima de las
oscilaciones registradas por el sism\'ografo, $A$, y una funci\'on (determinada emp\'iricamente) que
depende de la distancia $\delta$ de la estaci\'on s\'ismica al centro del epicentro, $A_0(\delta)$\footnote{Una
de las limitaciones fundamentales de la escala de Richter es que no es adecuada para estaciones que est\'en
alejadas m\'as de 600km del epicentro, por la dificultad de determinar $A_0(\delta)$ en esos casos.}:
$$
M_L=\log_{10}\left( \frac{A}{A_0(\delta)} \right) .
$$
Terremotos con un valor inferior a 2 en la escala Richter se registran continuamente (se denominan
microsismos). Los que tienen un valor de 4.5 se detectan en estaciones s\'ismicas por todo el globo
terrestre.
Lo destacable de la escala de Richter es que es \emph{logar\'itmica} y no lineal. Eso hace que un
peque\~no aumento en el valor de la escala corresponda a terremotos con capacidades destructivas completamente
distintas.
\begin{enumerate}[(a)]
\item El 28 de Julio de 1957 se registr\'o un terremoto en M\'exico, con su epicentro cerca de
Acapulco (Gro.), que alcanz\'o un valor de 7.9 en la escala Richter (datos del United States Geological
Survey, en \url{http://earthquake.usgs.gov/earthquakes/?source=sitenav}). Oficialmente 
caus\'o 68 v\'ictimas. El 19 de Septiembre de 1985 se registr\'o
otro terremoto, con epicentro en la costa de Michoac\'an, que afect\'o gravemente a la ciudad de M\'exico,
donde caus\'o m\'as de 6 000 muertos seg\'un cifras oficiales. La magnitud de este 
terremoto fue de 8.1. Calcular qu\'e tanto m\'as intenso fue este terremoto que el de Acapulco.
\item El terremoto m\'as intenso del que se tiene registro ocurri\'o en Chile, el 22 de Mayo de 1960.
Alcanz\'o una magnitud de 9.5. Calcular cu\'antas veces mayor fue la intensidad de este terremoto que el
de M\'exico en 1985 (puede suponerse que los datos corresponden a una estaci\'on que equidista de M\'exico
y Chile).
\end{enumerate}

\item Si un terremoto registra una magnitud de 6.2 en la escala de Richter, ?`qu\'e magnitud registrar\'a
otro tal que la amplitud m\'axima de sus ondas s\'ismicas sea 4 veces mayor?

\item En 1859, un granjero originario de Inglaterra residente en Australia solt\'o 12 parejas de conejos
dom\'esticos con el fin de que se reprodujeran para poder cazarlos luego. Los conejos se aplicaron
concienzudamente a la tarea que se les hab\'ia encomendado y para 1920, la poblaci\'on de 
conejos en Australia se estimaba en unos 10 000 millones (datos de la RFA: 
\url{http://www.rabbitfreeaustralia.com.au/rabbits/the-rabbit-problem/}).
Suponiendo que el crecimiento de los conejos fue exponencial (del tipo $f(t)=ae^{bt}$),
estimar los par\'ametros del modelo.
\begin{center}
\includegraphics[scale=0.62]{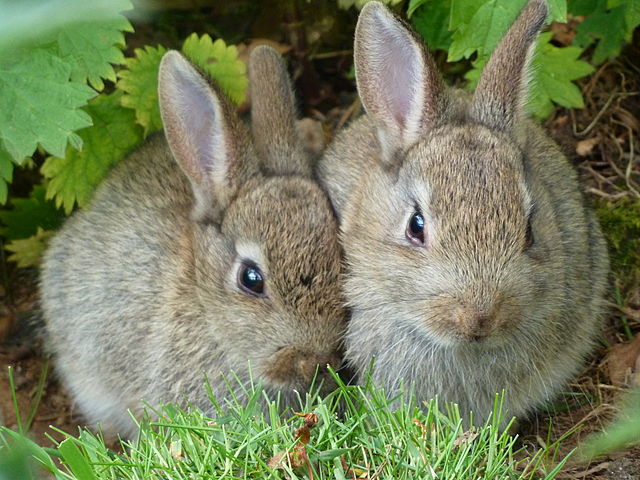} 
\end{center}%

\item El Salmo 139 del Antiguo Testamento (vers\'iculos 17-18) dice que los pensamientos de
Jehov\'a\index{Jehov\'a} sobre cada uno de los hombres son numerosos como los granos de arena (v\'ease
\url{http://www.biblegateway.com/passage/?search=Salmos+139&version=RVR1995}). Esta alegor\'ia
suele utilizarse con un desenfado que va de la mano de la ignorancia num\'erica, como
sin\'onimo de infinitud. Quiz\'as
sea interesante calcular el n\'umero de granos de arena (se supone que de playa) que hay en la 
tierra. Incluso tirando a lo alto\footnote{Arqu\'imedes\index{Arqu\'imedes}, uno de los matem\'aticos m\'as grandes de la historia, trat\'o de estimar el n\'umero de granos de arena que llenar\'ian el universo en su
obra `El arenario'. Para
ello, desarroll\'o una notaci\'on y nomenclatura para los grandes n\'umeros. Su conclusi\'on es
que el n\'umero de granos que llenan el universo es $10^{63}$, partiendo del c\'alculo err\'oneo
de que el radio del universo es de unos $2$ a\~nos luz. Sin embargo, resulta curioso observar que
el n\'umero de nucleones contenidos en esos $10^{63}$ granos de arena es aproximadamente
$10^{80}$, el n\'umero de Eddington\index{Eddington, Arthur}, que es el n\'umero estimado de nucleones en el universo
visible actual.}.
\begin{enumerate}
\item La tierra tiene un radio de 6 370km. Supuesta esf\'erica, calcula su \'area en metros cuadrados.
\item Como es sabido, de esa superficie s\'olo el 30\% es tierra firme, el resto es agua. Ahora vamos a calcular
qu\'e porcentaje de esa superficie corresponde a playas. Para ello, extrapolaremos datos de M\'exico (esto se
justifica por la gran extensi\'on de M\'exico). Seg\'un los datos oficiales (v\'ease
\url{http://en.wikipedia.org/wiki/List_of_countries_by_length_of_coastline}), M\'exico tiene
unos 2 000 000km$^2$ de superficie y aproximadamente 10 000km de costa\index{Me@M\'exico!playas}. Suponiendo que una playa t\'ipica tiene una
anchura de 20m y que toda la costa es playa, calcular el porcentaje de playa sobre la superficie total.
Aplicar el resultado para estimar el \'area de las playas en el planeta.
\item Un grano medio de arena puede suponerse tambi\'en esf\'erico, con un radio de 0.5$mm$
(datos: \url{http://en.wikipedia.org/wiki/Orders_of_magnitude_(volume)}). Suponiendo que
el espesor en arena de una playa t\'ipica es de 1m, calcular el n\'umero de granos que contienen
las playas de la tierra.
\item El cerebro de una persona adulta tiene $10^{14}$ sinapsis\index{Sinapsis} neuronales \emph{en promedio}\footnote{Las
sinapsis se est\'an modificando continuamente, seg\'un estudios recientes
 sobre 
plasticidad neuronal, v\'ease \url{http://www.sciencedaily.com/releases/2013/10/131010205325.htm}, por eso hablamos de promedio.}
(datos de la Universidad de Washington \url{http://faculty.washington.edu/chudler/facts.html}). Suponiendo
(a falta de una definici\'on mejor de pensamiento) que cada sinapsis da lugar a un `pensamiento' y que cada una 
se activa unas 1 000 veces al d\'ia, determinar en cu\'anto tiempo una persona adulta promedio tiene 
un n\'umero de `pensamientos' igual al del \'item precedente.
\end{enumerate}

\item Para cada apartado, representar las funciones que se indican en un mismo gr\'afico:
\begin{multicols}{2}
\begin{enumerate}[(a)]
\item $f(x)=x$, $g(x)=-x $, $h(x)=|x|$
\item $ f(x)=x^2 -2$, $g(x)=|x^2 -2|$
\item $ f(x)=x+1$, $g(x)=|x+1|$
\item $ f(x)=\ln x$, $g(x)=\ln |x|$
\end{enumerate}
\end{multicols}

\item\label{ejer} Razonar si las siguientes igualdades son ciertas:
\begin{multicols}{3}
\begin{enumerate}[(a)]
\item $|-3|=3 $
\item $ |27|=27$
\item $ |a^2|=a$
\item $ |-a|=a$
\item $ |a|=a$
\item $ \sqrt{a^2}=a$
\item $ \sqrt{a^2}=|a|$
\item $ \sqrt[3]{a^3}=a$
\item $|\sqrt{2}+\sqrt{3}-3|=\sqrt{2}+\sqrt{3}-3 $
\end{enumerate}
\end{multicols}

\item Resolver las ecuaciones siguientes:
\begin{multicols}{4}
\begin{enumerate}[(a)]
\item $|x-3|=x-3 $
\item $|x-3|=4$
\item $|x^2 -6x+8|=3$
\item $|x^2-6x+5|=x^2-6x+5 $
\end{enumerate}
\end{multicols}

\item\label{triangdes} Probar la desigualdad triangular del valor absoluto:
para todos $x,y\in\mathbb{R}$,
$$
|x+y|\leq |x|+|y|\,.
$$

\item\label{ejer84} Resolver las inecuaciones siguientes:
\begin{multicols}{2}
\begin{enumerate}[(a)]
\item $|x|<1 $
\item $|x+1|\leq 3 $
\item $|x|\geq 1 $
\item $|x+1|>3 $
\item $|1-x|<4x+1 $
\item $ \frac{{\displaystyle |x-1|}}{{\displaystyle x}}\leq 0$
\item $(1+x)^2 \geq |1-x^2| $
\item $ |34+21x-x^2|\leq -1$
\end{enumerate}
\end{multicols}

\item Comprobar que si se tiene una fuente $F$, capaz de emitir $N$ mensajes
\emph{equiprobables}, la medida de informaci\'on normalizada se puede obtener
como una caso particular de la entrop\'ia.

\item Un conocido refr\'an afirma que `una imagen vale m\'as que mil palabras'. Para
decidir si esto es cierto o no, consid\'erese la imagen ofrecida por la pantalla de
una computadora actual, con una resoluci\'on de $1920\times 1080$ pixels y 32 bits 
de color.\index{Informaci\'on!imagen y palabras}
\begin{enumerate}[(a)]
\item Asimilando un s\'imbolo a cada estado de un p\'ixel, determinar el n\'umero
de mensajes que puede emitir esta fuente. Calcular su entrop\'ia.
\item Comparar el resultado con el obtenido por un locutor profesional de radio,
con un vocabulario de unas $1024$ palabras (sup\'ongase que todas equiprobables,
aunque esta suposici\'on es obviamente sobresimplificadora) y capaz de pronunciar
$50$ de ellas por minuto.
\end{enumerate}

\item Supongamos que la fuente $F_p$ puede emitir los mensajes $m_1$, $m_2$ con probabilidades
respectivas $p_1=1/2=p_2$. Sea $F_q$ otra fuente que puede emitir $8$ mensajes
$n_1,\ldots,n_8$ con probabilidades $q_j=1/8$, $j\in\{1,2\ldots ,8\}$. Se pide calcular
la informaci\'on global de las siguientes fuentes:
\begin{enumerate}[(a)]
\item $F_{\mathrm{eq}}$, que puede emitir un mensaje de $F_p$ o de $F_q$ con la misma
probabilidad.
\item $F_{2,8}$, que puede emitir un mensaje de $F_p$ con probabilidad $2/10$ o uno de
$F_q$ con probabilidad $8/10$.
\item $F_{8,2}$, que puede emitir un mensaje de $F_p$ con probabilidad $8/10$ o uno de
$F_q$ con probabilidad $2/10$.
\end{enumerate}
\end{enumerate}

\section{Uso de la tecnolog\'ia de c\'omputo}
Maxima tiene implementadas todas las herramientas que necesitamos
para trabajar con las funciones que hemos tratado en este cap\'itulo.
Por ejemplo, las potencias se indican mediante el circunflejo $\wedge$:
\begin{maxcode}
\noindent
\begin{minipage}[t]{8ex}{\color{red}\bf
\begin{verbatim}
(%i1) 
\end{verbatim}
}
\end{minipage}
\begin{minipage}[t]{\textwidth}{\color{blue}
\begin{verbatim}
3^2;
\end{verbatim}
}
\end{minipage}
\definecolor{labelcolor}{RGB}{100,0,0}
\begin{math}\displaystyle
\parbox{8ex}{\color{labelcolor}(\%o1) }
9
\end{math}
\end{maxcode}
\noindent Notemos que las reglas de composici\'on de los exponentes
est\'an implementadas por defecto,
\begin{maxcode}
\noindent
\begin{minipage}[t]{8ex}{\color{red}\bf
\begin{verbatim}
(%i2) 
\end{verbatim}
}
\end{minipage}
\begin{minipage}[t]{\textwidth}{\color{blue}
\begin{verbatim}
a^2*a^3;
\end{verbatim}
}
\end{minipage}
\definecolor{labelcolor}{RGB}{100,0,0}
\begin{math}\displaystyle
\parbox{8ex}{\color{labelcolor}(\%o2) }
{a}^{5}
\end{math}
\end{maxcode}
Por supuesto, podemos definir funciones m\'as complejas
combinando otras ya conocidas,
\begin{maxcode}
\noindent
\begin{minipage}[t]{8ex}{\color{red}\bf
\begin{verbatim}
(%i3) 
\end{verbatim}
}
\end{minipage}
\begin{minipage}[t]{\textwidth}{\color{blue}
\begin{verbatim}
f(x):=3^(2*x-1);
\end{verbatim}
}
\end{minipage}
\definecolor{labelcolor}{RGB}{100,0,0}
\begin{math}\displaystyle
\parbox{8ex}{\color{labelcolor}(\%o3) }
\mathrm{f}(x):=3^{2x-1}
\end{math}

\noindent
\begin{minipage}[t]{8ex}{\color{red}\bf
\begin{verbatim}
(%i4) 
\end{verbatim}
}
\end{minipage}
\begin{minipage}[t]{\textwidth}{\color{blue}
\begin{verbatim}
f(2);
\end{verbatim}
}
\end{minipage}
\definecolor{labelcolor}{RGB}{100,0,0}
\begin{math}\displaystyle
\parbox{8ex}{\color{labelcolor}(\%o4) }
27
\end{math}
\end{maxcode}
El logaritmo que Maxima tiene implementado es el \emph{natural},
donde la base es el n\'umero $e\simeq 2{.}7172$,
\begin{maxcode}
\noindent
\begin{minipage}[t]{8ex}{\color{red}\bf
\begin{verbatim}
(%i5) 
\end{verbatim}
}
\end{minipage}
\begin{minipage}[t]{\textwidth}{\color{blue}
\begin{verbatim}
log(%e^4);
\end{verbatim}
}
\end{minipage}
\definecolor{labelcolor}{RGB}{100,0,0}
\begin{math}\displaystyle
\parbox{8ex}{\color{labelcolor}(\%o5) }
4
\end{math}
\end{maxcode}
\noindent Podemos definir el logaritmo en una base cualquiera muy f\'acilmente,
mediante la f\'ormula de cambio de base (notemos que la nueva funci\'on no puede
llamarse \texttt{log(x,n)}, pues ese nombre ya est\'a reservado para la funci\'on
\texttt{log} de Maxima):
\begin{maxcode}
\noindent
\begin{minipage}[t]{8ex}{\color{red}\bf
\begin{verbatim}
(%i6) 
\end{verbatim}
}
\end{minipage}
\begin{minipage}[t]{\textwidth}{\color{blue}
\begin{verbatim}
Log(x,n):=log(x)/log(n);
\end{verbatim}
}
\end{minipage}
\definecolor{labelcolor}{RGB}{100,0,0}
\begin{math}\displaystyle
\parbox{8ex}{\color{labelcolor}(\%o6) }
\mathrm{Log}\left( x,n\right) :=\frac{\mathrm{log}\left( x\right) }{\mathrm{log}\left( n\right) }
\end{math}

\noindent
\begin{minipage}[t]{8ex}{\color{red}\bf
\begin{verbatim}
(%i7) 
\end{verbatim}
}
\end{minipage}
\begin{minipage}[t]{\textwidth}{\color{blue}
\begin{verbatim}
Log(8,2);
\end{verbatim}
}
\end{minipage}
\definecolor{labelcolor}{RGB}{100,0,0}
\begin{math}\displaystyle
\parbox{8ex}{\color{labelcolor}(\%o7) }
\frac{\mathrm{log}\left( 8\right) }{\mathrm{log}\left( 2\right) }
\end{math}

\noindent
\begin{minipage}[t]{8ex}{\color{red}\bf
\begin{verbatim}
(%i8) 
\end{verbatim}
}
\end{minipage}
\begin{minipage}[t]{\textwidth}{\color{blue}
\begin{verbatim}
Log(8,2),numer;
\end{verbatim}
}
\end{minipage}
\definecolor{labelcolor}{RGB}{100,0,0}
\begin{math}\displaystyle
\parbox{8ex}{\color{labelcolor}(\%o8) }
3{.}0
\end{math}
\end{maxcode}
Por defecto, Maxima no aplica las reglas de simplificaci\'on para los logaritmos,
pero este comportamiento puede modificarse ajustando
el valor de la variable de entorno \texttt{logexpand}.
Los valores posibles son \texttt{true} (que hace que $\log (a^b)=b\log a$),
\texttt{all} (que, adem\'as, transforma $\log (ab)$ en $\log a+\log b$),
\texttt{super} (que a\~nade la regla $\log (a/b)=\log a-\log b$) y
\texttt{false} (que no efect\'ua simplificaci\'on alguna):

\begin{maxcode}
\noindent
\begin{minipage}[t]{8ex}{\color{red}\bf
\begin{verbatim}
(%i9) 
\end{verbatim}
}
\end{minipage}
\begin{minipage}[t]{\textwidth}{\color{blue}
\begin{verbatim}
log(x*y);
\end{verbatim}
}
\end{minipage}
\definecolor{labelcolor}{RGB}{100,0,0}
\begin{math}\displaystyle
\parbox{8ex}{\color{labelcolor}(\%o9) }
\mathrm{log}\left( x\,y\right) 
\end{math}

\noindent
\begin{minipage}[t]{8ex}{\color{red}\bf
\begin{verbatim}
(%i10) 
\end{verbatim}
}
\end{minipage}
\begin{minipage}[t]{\textwidth}{\color{blue}
\begin{verbatim}
logexpand:super;
\end{verbatim}
}
\end{minipage}
\definecolor{labelcolor}{RGB}{100,0,0}
\begin{math}\displaystyle
\parbox{8ex}{\color{labelcolor}(\%o10) }
super
\end{math}

\noindent
\begin{minipage}[t]{8ex}{\color{red}\bf
\begin{verbatim}
(%i11) 
\end{verbatim}
}
\end{minipage}
\begin{minipage}[t]{\textwidth}{\color{blue}
\begin{verbatim}
log(x*y);
\end{verbatim}
}
\end{minipage}
\definecolor{labelcolor}{RGB}{100,0,0}
\begin{math}\displaystyle
\parbox{8ex}{\color{labelcolor}(\%o11) }
\mathrm{log}\left( y\right) +\mathrm{log}\left( x\right) 
\end{math}

\noindent
\begin{minipage}[t]{8ex}{\color{red}\bf
\begin{verbatim}
(%i12) 
\end{verbatim}
}
\end{minipage}
\begin{minipage}[t]{\textwidth}{\color{blue}
\begin{verbatim}
log(m/n);
\end{verbatim}
}
\end{minipage}
\definecolor{labelcolor}{RGB}{100,0,0}
\begin{math}\displaystyle
\parbox{8ex}{\color{labelcolor}(\%o12) }
\mathrm{log}\left( m\right) -\mathrm{log}\left( n\right) 
\end{math}
\end{maxcode}
\noindent Otro comando \'util para trabajar con logaritmos
es \texttt{logcontract}, que realiza las simplificaciones
en sentido contrario a \texttt{logexpand},
\begin{maxcode}
\noindent
\begin{minipage}[t]{8ex}{\color{red}\bf
\begin{verbatim}
(%i13) 
\end{verbatim}
}
\end{minipage}
\begin{minipage}[t]{\textwidth}{\color{blue}
\begin{verbatim}
logcontract(2*log(a));
\end{verbatim}
}
\end{minipage}
\definecolor{labelcolor}{RGB}{100,0,0}
\begin{math}\displaystyle
\parbox{8ex}{\color{labelcolor}(\%o13) }
\mathrm{log}\left( {a}^{2}\right) 
\end{math}
\end{maxcode}
La funci\'on exponencial con base $e$, tiene su propia notaci\'on,
\texttt{exp}:
\begin{maxcode}
\noindent
\begin{minipage}[t]{8ex}{\color{red}\bf
\begin{verbatim}
(%i14) 
\end{verbatim}
}
\end{minipage}
\begin{minipage}[t]{\textwidth}{\color{blue}
\begin{verbatim}
exp(a);
\end{verbatim}
}
\end{minipage}
\definecolor{labelcolor}{RGB}{100,0,0}
\begin{math}\displaystyle
\parbox{8ex}{\color{labelcolor}(\%o14) }
{e}^{a}
\end{math}

\noindent
\begin{minipage}[t]{8ex}{\color{red}\bf
\begin{verbatim}
(%i15) 
\end{verbatim}
}
\end{minipage}
\begin{minipage}[t]{\textwidth}{\color{blue}
\begin{verbatim}
exp(log(a));
\end{verbatim}
}
\end{minipage}
\definecolor{labelcolor}{RGB}{100,0,0}
\begin{math}\displaystyle
\parbox{8ex}{\color{labelcolor}(\%o15) }
a
\end{math}

\noindent
\begin{minipage}[t]{8ex}{\color{red}\bf
\begin{verbatim}
(%i16) 
\end{verbatim}
}
\end{minipage}
\begin{minipage}[t]{\textwidth}{\color{blue}
\begin{verbatim}
log(exp(a));
\end{verbatim}
}
\end{minipage}
\definecolor{labelcolor}{RGB}{100,0,0}
\begin{math}\displaystyle
\parbox{8ex}{\color{labelcolor}(\%o16) }
a
\end{math}

\noindent
\begin{minipage}[t]{8ex}{\color{red}\bf
\begin{verbatim}
(%i17) 
\end{verbatim}
}
\end{minipage}
\begin{minipage}[t]{\textwidth}{\color{blue}
\begin{verbatim}
exp(a+b);
\end{verbatim}
}
\end{minipage}
\definecolor{labelcolor}{RGB}{100,0,0}
\begin{math}\displaystyle
\parbox{8ex}{\color{labelcolor}(\%o17) }
{e}^{b+a}
\end{math}
\end{maxcode}
Por \'ultimo tenemos el valor absoluto, que en Maxima se
escribe como \texttt{abs},
\begin{maxcode}
\noindent
\begin{minipage}[t]{8ex}{\color{red}\bf
\begin{verbatim}
(%i18) 
\end{verbatim}
}
\end{minipage}
\begin{minipage}[t]{\textwidth}{\color{blue}
\begin{verbatim}
abs(3);
\end{verbatim}
}
\end{minipage}
\definecolor{labelcolor}{RGB}{100,0,0}
\begin{math}\displaystyle
\parbox{8ex}{\color{labelcolor}(\%o18) }
3
\end{math}

\noindent
\begin{minipage}[t]{8ex}{\color{red}\bf
\begin{verbatim}
(%i19) 
\end{verbatim}
}
\end{minipage}
\begin{minipage}[t]{\textwidth}{\color{blue}
\begin{verbatim}
abs(-3);
\end{verbatim}
}
\end{minipage}
\definecolor{labelcolor}{RGB}{100,0,0}
\begin{math}\displaystyle
\parbox{8ex}{\color{labelcolor}(\%o19) }
3
\end{math}
\end{maxcode}
\noindent A la hora de trabajar con el valor absoluto, es muy conveniente
utilizar el comando \texttt{assume} que, como su nombre indica, le dice a
Maxima que una determinada variable est\'a restringida a ciertos valores,
\begin{maxcode}
\noindent
\begin{minipage}[t]{8ex}{\color{red}\bf
\begin{verbatim}
(%i20) 
\end{verbatim}
}
\end{minipage}
\begin{minipage}[t]{\textwidth}{\color{blue}
\begin{verbatim}
assume(u>1);
\end{verbatim}
}
\end{minipage}
\definecolor{labelcolor}{RGB}{100,0,0}
\begin{math}\displaystyle
\parbox{8ex}{\color{labelcolor}(\%o20) }
[u>1]
\end{math}

\noindent
\begin{minipage}[t]{8ex}{\color{red}\bf
\begin{verbatim}
(%i21) 
\end{verbatim}
}
\end{minipage}
\begin{minipage}[t]{\textwidth}{\color{blue}
\begin{verbatim}
abs(1-u);
\end{verbatim}
}
\end{minipage}
\definecolor{labelcolor}{RGB}{100,0,0}
\begin{math}\displaystyle
\parbox{8ex}{\color{labelcolor}(\%o21) }
u-1
\end{math}
\end{maxcode}

\noindent Esto cobra especial relevancia al tratar con n\'umeros
entre $0$ y $1$, para los cuales su cuadrado es menor que el propio
n\'umero, como se puede apreciar gr\'aficamente\footnote{Anal\'iticamente, si
$0<x<1$, al multiplicar ambos miembros por $x$ resulta $x^2<x$.},
\begin{maxcode}
\noindent
\begin{minipage}[t]{8ex}{\color{red}\bf
\begin{verbatim}
(%i22) 
\end{verbatim}
}
\end{minipage}
\begin{minipage}[t]{\textwidth}{\color{blue}
\begin{verbatim}
wxplot2d([x,x^2],[x,-0.5,1.5]);
\end{verbatim}
}
\end{minipage}
\definecolor{labelcolor}{RGB}{100,0,0}
\begin{math}\displaystyle
\parbox{8ex}{\color{labelcolor}(\%t22) }\centering
\includegraphics[width=9cm]{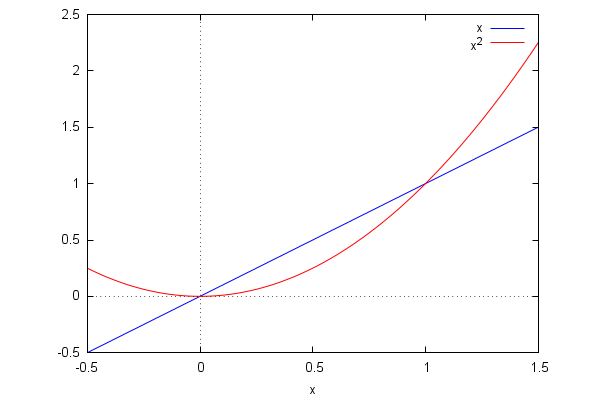} 
\end{math}
\begin{math}\displaystyle
\parbox{8ex}{\color{labelcolor}(\%o22) }
\end{math}
\end{maxcode}
\noindent La manera de trabajar en estos casos consiste en usar \texttt{assume}:
\begin{maxcode}
\noindent
\begin{minipage}[t]{8ex}{\color{red}\bf
\begin{verbatim}
(%i23) 
\end{verbatim}
}
\end{minipage}
\begin{minipage}[t]{\textwidth}{\color{blue}
\begin{verbatim}
assume(0<v,v<1);
\end{verbatim}
}
\end{minipage}
\definecolor{labelcolor}{RGB}{100,0,0}
\begin{math}\displaystyle
\parbox{8ex}{\color{labelcolor}(\%o23) }
[v>0,v<1]
\end{math}

\noindent
\begin{minipage}[t]{8ex}{\color{red}\bf
\begin{verbatim}
(%i24) 
\end{verbatim}
}
\end{minipage}
\begin{minipage}[t]{\textwidth}{\color{blue}
\begin{verbatim}
abs(v^2-v);
\end{verbatim}
}
\end{minipage}
\definecolor{labelcolor}{RGB}{100,0,0}
\begin{math}\displaystyle
\parbox{8ex}{\color{labelcolor}(\%o24) }
v-v^2
\end{math}
\end{maxcode}

\section{Actividades de c\'omputo}
\begin{enumerate}[(A)]
\item Resolver el ejercicio \ref{ejer72} utilizando Maxima. Comparar el resultado
obtenido con el te\'orico.
\item Resolver el ejercicio \ref{ejer73} utilizando Maxima. Comparar el resultado
obtenido con el te\'orico.
\item Resolver el ejercicio \ref{ejer84} gr\'aficamente con Maxima, representando
en una misma figura las funciones involucradas y estudiando sus puntos de corte.
Comparar el resultado obtenido con el te\'orico.
\end{enumerate}

\chapter{Trigonometr\'ia}\label{chaptrig}
\vspace{2.5cm}
\section{Nociones b\'asicas}
La trigonometr\'ia es una de las ramas venerables de la geometr\'ia. Su aparici\'on se remonta
a Babilonia y el antiguo Egipto, 4\,000 adC, \'intimamente relacionada con la agrimensura y la
astronom\'ia.
B\'asicamente, podemos caracterizarla como el estudio de la geometr\'ia de tri\'angulos y circunferencias: la palabra \emph{trigonometr\'ia} proviene de \emph{-tri} (tres) y el griego
$-\gamma$o$\nu$o$\varsigma$ (\emph{-gonos}, \'angulo), y significa precisamente
`medida de tres \'angulos'.
Para nuestros prop\'ositos, un \emph{tri\'angulo}\index{Tri\'angulo} puede definirse como la regi\'on \emph{acotada}
del plano resultante de la intersecci\'on de los semiplanos determinados por tres rectas no paralelas.
\begin{center}
    \begin{tikzpicture}[auto,scale=0.66]
	\coordinate (A) at (-6,1);
	\coordinate (B) at (4,1);
	\coordinate (C) at (-1,-3);
	\coordinate (D) at (2,3);
	\coordinate (E) at (2,-2);
	\coordinate (F) at (-5,1.5);
	\coordinate (P1) at (0,-1);
	\coordinate (P2) at (-4,1);
	\coordinate (P3) at (1,1);	
	\draw [ultra thin,draw=white,fill=yellow,opacity=0.2](P1)--(P2)--(P3)--cycle;
    \draw[thick,red] (A) -- (B);
	\draw[thick,blue] (C) -- (D);
	\draw[thick,orange] (E) -- (F);
	     
    \end{tikzpicture}
\end{center}
Una \emph{circunferencia}\index{Circunferencia} de centro
$P\equiv (a,b)$ y radio $r$, denotada $C_r(P)$, es el conjunto de puntos del plano que equidistan de $P$,
donde la distancia entre dos puntos $(a,b),(x,y)\in\mathbb{R}^2$ es la \emph{euclidea}\index{Euclides!distancia}\index{Distancia!euclidea}:
$$
\mathrm{d}((x,y);(a,b))=\sqrt{(x-a)^2+(y-b)^2}.
$$
De aqu\'i resulta inmediato que la ecuaci\'on que satisfacen los puntos $(x,y)$ de $C_r(P)$ es
$$
(x-a)^2+(y-b)^2=r^2 .
$$
Sin embargo, en esta secci\'on nos ocuparemos b\'asicamente de las propiedades de las circunferencias que son
independientes de su representaci\'on en coordenadas. Por ejemplo, utilizaremos que la
longitud de una circunferencia de radio $r$ es $2\pi r$ (un resultado que \emph{no} pertenece al \'ambito
de la geometr\'ia euclidea, sino de la geometr\'ia m\'etrica, e involucra el uso del c\'alculo diferencial). De hecho, uno e los mayores descubrimientos de la geometr\'ia griega fu\'e el
hecho de que el cociente entre la longitud de \emph{cualquier} circunferencia y su di\'ametro
es una constante, el n\'umero $\pi$, con un valor aproximado (determinado por Arqu\'imedes\index{Arqu\'imedes}) de $\pi\simeq 3{.}141592$\index{Num@N\'umero!pi@$\pi$}.

Un \emph{\'angulo}\index{An@\'Angulo} es cada una de las porciones de plano determinadas por dos semirrectas
con origen com\'un. N\'otese que cada par de semirrectas define \emph{dos} \'angulos, el
\emph{exterior} y el \emph{interior}. El \'angulo interior es el formado por la regi\'on
del plano que se recorre al
realizar el menor giro que lleva una recta en la otra. La \'unica ambig\"uedad que se presenta es en el
caso del \'angulo llano\index{An@\'Angulo!llano} (formado por dos semirrectas colineales opuestas), pero en este caso el \'angulo
externo y el interno coinciden, por lo que no se presenta problema alguno.
\begin{center}
\begin{tikzpicture}
    \draw[red,thick] (0,0) -- (30:3);
	\draw[blue,thick] (0,0) -- (3,0) ;
    \draw (1,0) arc (0:30:1);
    \node[] at (15:1.25)  {$\alpha$};
\end{tikzpicture}
\end{center}
La magnitud de un \'angulo mide la separaci\'on
entre las rectas. Convencionalmente (como herencia de la matem\'atica babilonia cuyo sistema de 
numeraci\'on era sexagesimal, con base 60), el \'angulo que le corresponde a una circunferencia,
es decir, el \'angulo externo formado por dos rectas coincidentes, tiene una magnitud de 360$\degree$,
de modo que cualquiera de los \'angulos \emph{llanos}\index{An@\'Angulo!llano} determinados por un \emph{di\'ametro} (es decir,
el segmento de cualquier recta que pasa por el origen de la circunferencia comprendido entre sus
intersecciones con ella) vale 180$\degree$. Dos rectas perpendiculares determinan cuatro \'angulos
iguales que, por tanto, tienen una magnitud de 90$\degree$ cada uno. Se denominan \'angulos \emph{rectos}\index{An@\'Angulo!recto}.  

Obviamente, cada tri\'angulo determina tres \'angulos, uno en cada v\'ertice (los puntos de intersecci\'on
de las rectas que definen el tri\'angulo). Las propiedades b\'asicas de los tri\'angulos se pueden deducir 
f\'acilmente de las generales de la geometr\'ia euclidea. Por ejemplo, sabiendo que los \'angulos alternos internos\index{An@\'Angulo!alterno interno}\footnote{Son los \'angulos que cumplen: (i) ambos son internos, es decir, comprendidos entre
las paralelas, (ii) se encuentran en semiplanos distintos con respecto a la recta transversal, (iii)
no son adyacentes. Por su parte, los \'angulos adyacentes\index{An@\'Angulo!adyacente} son los que comparten una semirrecta (y, por tanto,
el origen).} determinados por la intersecci\'on de dos rectas paralelas con una tercera transversal
(o sea, no paralela) a ellas son iguales,
\begin{center}
\begin{tikzpicture}
    \draw[thick] (-4,1) -- (3,1);
    \draw[thick] (-4,-1) -- (3,-1);
    \draw[thick,gray] (-0.5,1) arc (180:210:0.5);
    \node[] at (-0.75,0.75)  {$\gamma$};
    \draw[thick,gray] (-3,-1) arc (0:30:0.5);
    \node[] at (-2.75,-0.8)  {$\gamma$};    
    \draw[orange,very thick] (-4,-1.28) -- (1.5,1.85);
\end{tikzpicture}
\end{center}
se puede probar el important\'isimo resultado v\'alido para todo
tri\'angulo: la suma de las magnitudes de sus \'angulos es igual a un \'angulo llano, es decir, 
180$\degree$. La demostraci\'on, usando dos rectas paralelas, es puramente gr\'afica:\\
\begin{center}
\begin{tikzpicture}
    \draw[thick] (-4,1) -- (3,1);
    \draw[thick] (-4,-1) -- (-3.5,-1);
    \draw[thick] (2,-1) -- (3,-1);
    \draw[thick,gray] (-0.5,1) arc (180:360:0.5);
    \node[] at (0,0.25)  {$\beta$};
    \node[] at (-0.75,0.75)  {$\gamma$};
    \node[] at (0.75,0.75)  {$\alpha$};           
    \draw[thick,gray] (-3,-1) arc (0:30:0.5);
    \node[] at (-2.75,-0.8)  {$\gamma$};    
    \draw[thick,gray] (1.5,-1) arc (180:132:0.5);
     \node[] at (1.25,-0.8)  {$\alpha$};                 
    \draw[orange,very thick] (-3.5,-1) -- (2,-1);
    \draw[orange,very thick] (-3.5,-1) -- (0,1);
    \draw[orange,very thick] (0,1) -- (2,-1);
\end{tikzpicture}
\end{center}

Un tipo particularmente interesante de tri\'angulos es el formado por aquellos que tienen un \'angulo
recto. Por este motivo se denominan tri\'angulos \emph{rect\'angulos}\index{Tri\'angulo!rect\'angulo}. Como consecuencia del resultado
precedente, la suma de los otros dos \'angulos debe ser $\alpha +\beta =90\degree$, en particular
$0<\alpha ,\beta <90\degree$.

Para un tri\'angulo rect\'angulo como el de la figura,
\begin{center}
    \begin{tikzpicture}[scale=0.75]
	   \draw[-,very thick] (4,0) -- (4,4) node[midway, right] {$a$};
	   \draw[-,very thick] (0,0) -- (4,0) node[midway, below] {$b$};
	   \draw[-,very thick] (0,0) --node[above]{$c$} (4,4) ;
	   \draw[thick] (0.75,0) arc (0:45:0.75);
	   \node[] at (20:1.15) {$\alpha$};
	   \draw[thick] (3.5,3.5) arc (225:270:0.75);
	   \node[] at (37.5:4.5) {$\beta$};	   
	\end{tikzpicture}
\end{center}
\noindent se definen las \emph{razones trigonom\'etricas}\index{Razones trigonom\'etricas}
de sus \'angulos: seno\index{Seno}, coseno\index{Coseno} y tangente\index{Tangente}, respectivamente,
\begin{align}
\sin \alpha &= \frac{a}{c} \nonumber \\
\cos\alpha &= \frac{b}{c} \label{trigo}\\
\tan \alpha &= \frac{a}{b} \nonumber.
\end{align}
N\'otese que, por definici\'on,
$$
\tan \alpha = \frac{\sin\alpha}{\cos\alpha}\,.
$$
An\'alogamente, para el \'angulo $\beta$ se tiene
\begin{align*}
\sin \beta &= \frac{b}{c} \\
\cos\beta &= \frac{a}{c} \\
\tan \beta &= \frac{b}{a}.
\end{align*}
Com\'unmente se dice que `el (co)seno de un \'angulo es la longitud del cateto opuesto (contiguo)
dividida por la longitud de la hipotenusa', mientras que `la tangente de un \'angulo es el seno dividido
por el coseno'. Para el \'angulo recto, estas razones se definen directamente como
\begin{align*}
\sin 90\degree &= 1 \\
\cos 90\degree &= 0 \\
\tan 90\degree &= \infty .
\end{align*}
Los tri\'angulos rect\'angulos cumplen el conocid\'isimo \emph{teorema de Pit\'agoras}\index{Pit\'agoras!teorema}:
\begin{equation}\label{pitagoras}
a^2+b^2=c^2 ,
\end{equation}
cuya demostraci\'on, acorde con el origen etimol\'ogico del t\'ermino\footnote{La
palabra \emph{teorema}\index{Teorema} proviene directamente del griego: el lexema `thea' equivale al del
verbo \emph{mirar}, y el morfema `ema' equivale al espa\~nol \emph{-ci\'on}, es decir,
\emph{efecto de}, \emph{resultado de}. As\'i, un \emph{teorema} es el resultado de mirar
algo (de mirar una \emph{demostraci\'on}\index{Demostraci\'on}, en este caso).}, es puramente visual:
\begin{center}
\begin{tikzpicture}[scale=0.67]
\tkzDefPoint(0,0){A} \tkzDefPoint(7,0){B}
\tkzDefSquare(A,B) \tkzGetPoints{C}{D}
\tkzDrawPolygon(B,C,D,A)
\tkzDefPoint(3,0){I} \tkzDefPoint(7,3){J}
\tkzDefSquare(I,J)
\tkzGetPoints{K}{L}
\tkzDrawPolygon(I,J,K,L) \tkzFillPolygon[color=red!50](A,I,L)
\tkzFillPolygon[color=blue!30](I,B,J)\tkzFillPolygon[color=yellow!20](K,D,L)
\tkzFillPolygon[color=olive!50](K,C,J)
\end{tikzpicture}\mbox{ }
\begin{tikzpicture}[scale=0.67]
\tkzDefPoint(0,0){A} \tkzDefPoint(7,0){B}
\tkzDefSquare(A,B) \tkzGetPoints{C}{D}
\tkzDefPoint(3,0){I} \tkzDefPoint(7,3){J}
\tkzDefSquare(A,I) \tkzGetPoints{M}{L}
\tkzDefSquare(M,J) \tkzGetPoints{C}{P}
\tkzDrawPolygon(A,I,M,L) \tkzDrawPolygon(M,J,C,P)
\tkzDrawSegments(I,J D,M I,B B,J P,D D,L)
\tkzFillPolygon[color=red!50](M,J,I) \tkzFillPolygon[color=yellow!20](D,M,L)
\tkzFillPolygon[color=blue!30](I,B,J) \tkzFillPolygon[color=olive!50](D,P,M)
\end{tikzpicture} 
\end{center}
La identidad \eqref{pitagoras} en t\'erminos de las razones trigonom\'etricas \eqref{trigo} se
traduce en
$$
\sin^2\alpha +\cos^2\alpha =1,
$$
relaci\'on que es v\'alida para cualquier \'angulo en un tri\'angulo rect\'angulo.
N\'otese que esto implica las propiedades
\begin{equation}\label{cota1}
|\cos\alpha| \leq 1\qquad \mbox{ y }\qquad |\sin\alpha|\leq 1\,.
\end{equation}

Por otra parte, con el debido respeto hacia los babilonios\index{Babilonios}, el dividir la circunferencia en
360 partes (esto es, asignarle
unos \'angulos de 0$\degree$ y 360$\degree$ a dos semirrectas coincidentes), es algo muy arbitrario\footnote{Aunque no carente de sentido. El hecho de que $360$ sea un
n\'umero con una gran cantidad de divisores de peque\~na magnitud (su factorizaci\'on
prima es $360=2^3\, 3^2\, 5$) hace que sea f\'acil calcular con una gran variedad de
valores de \'angulos predeterminados, como $12\degree$, $15\degree$, $30\degree$, etc.}. 
Desde un punto de vista
formal, es mucho m\'as conveniente asignarle a la circunferencia un \'angulo relacionado con su longitud:
en lugar de 360$\degree$ le asignamos $2\pi$ radianes. Un \emph{radi\'an}\index{Radi\'an} es el \'angulo que
determina una longitud de arco de circunferencia igual al radio de la misma, de manera que la circunferencia
tiene un \'angulo de $2\pi$ radianes independientemente de su radio.
\begin{center}
\begin{tikzpicture}
\draw[thick] (0,0) circle(2);
\draw[gray] (0,0) -- (15:2);
\draw[gray] (0,0) -- node[black,midway,above]{$r$}  (60:2);
\draw[very thick,red] (1.93,0.52) arc (15:60:2);
\draw[very thick,red] (0.48,0.11) arc (15:60:0.5);
\node[] at (37:0.7) {$\theta$};
\node[] at (37:2.4) {$r\theta$};
\end{tikzpicture}
\end{center}
Como una circunferencia tiene un \'angulo de $2\pi\simeq 6{.}28$ radianes, un radi\'an expresado en grados
es $1\,\mathrm{rad}\simeq 360\degree/6{.}28\simeq 57{.}3\degree$. En general, para pasar de grados a radianes y 
viceversa, se usa una relaci\'on de proporcionalidad (regla de tres). Por ejemplo, si deseamos pasar
$60\degree$ a radianes, planteamos
$$
60\cdot 1\degree =60\cdot\frac{2\pi}{360}\mathrm{rad}=\frac{2\pi}{6}\mathrm{rad}=\frac{\pi}{3}\mathrm{rad}.
$$

Mediante el uso de una circunferencia de radio unidad, se puede extender la definici\'on de las razones
trigonom\'etricas ($\sin$, $\cos$, $\tan$) para \'angulos $\alpha$ de cualquier magnitud real, no s\'olo
los estrictamente comprendidos entre $0\degree$ y $90\degree$. En otras palabras, se pueden considerar
$\sin$, $\cos$, $\tan$ como funciones reales.\\
Para ello, comenzamos observando que, por ser el radio $r=1$, las razones trigonom\'etricas tienen la
interpretaci\'on geom\'etrica siguiente:
\begin{center}
\begin{tikzpicture}[scale=0.63]
	   \draw[thick] (0,0) circle(4);
	   \draw[-] (-4.5,0) -- (4.5,0) ;
       \draw[-] (0,-4.5) -- (0,4.5) ;
	   \draw[-] (0,0) -- (4.25,4.25) node[midway,above,left] {$r=1$} ;
	   \draw[-,very thick,red] (4,0) -- (4,4) node[midway, right] {$\tan \alpha$};
	   \draw[-,very thick,blue] (2.83,0) -- (2.83,2.83);
	   \draw[-,very thick,orange] (0,0) -- (2.83,0) node[midway, below] {$\cos \alpha$};
	   \draw[thick] (0.75,0) arc (0:45:0.75);
	   \node[] at (22.5:1.15) {$\alpha$};
	   \node[blue] at (2,0.95) {$\sin \alpha$};
\end{tikzpicture}
\end{center}     
En particular, la interpretaci\'on de $\tan\alpha$ como la longitud del segmento indicado, se debe
a un resultado elemental de geometr\'ia euclidea, el \emph{teorema de Tales}\index{Tales!teorema} referente a tri\'angulos
semejantes:
\begin{center}
    \begin{tikzpicture}[scale=0.75]
	   \draw[-,very thick,red] (4,-0.027) -- (4,4) node[midway, right] {$B$};
	   \draw[-,very thick,blue] (2.83,0) -- (2.83,2.83) node[midway, right] {$b$};
	   \draw[-,very thick,orange] (0,0.025) -- (2.83,0.025) node[midway, above] {$a$};
	   \draw[-,very thick,teal] (0,-0.025) -- (4,-0.025) node[midway, below] {$A$};
	   \draw[-] (0,0) -- (4,4) ;
     \end{tikzpicture}
\end{center}
El teorema afirma que se cumple la relaci\'on
$$
\frac{b}{a}=\frac{B}{A} .
$$
Tiene sentido entonces definir las \emph{funciones}\index{Funci\'on!trigonom\'etrica}
$\sin x,\cos x,\tan x$ para cualquier n\'umero $x$
simplemente tomando m\'odulo $2\pi$ (y por tanto olvidando el n\'umero de 'vueltas enteras'); es decir, si por ejemplo
$x=25$, para calcular $\sin 25$ podemos pensar en $25=7\pi + 3{.}01$ ( $25\equiv 3{.}01\mod 2\pi$, tres vueltas y media m\'as $3{.}01$ radianes), y definir $\sin 25=-\sin 3{.}01$. N\'otese que, por 
construcci\'on, las funciones trigonom\'etricas
$\sin x,\cos x,\tan x$  son \emph{peri\'odicas}\index{Funci\'on!peri\'odica}: sus valores se repiten a intervalos de $2\pi$ (es decir,
$\sin (x+2\pi)=\sin x$, etc). Leibniz\index{Leibniz, Gottfried Wilhelm} fue el primero en considerar las magnitudes
trigonom\'etricas como funciones.

De lo expuesto hasta ahora se deducen algunas propiedades importantes.
Por ejemplo, observemos que si 
$0<\alpha <\pi/2$, entonces
$$
0<\sin\alpha <\alpha <\tan\alpha\,,
$$
como es evidente a partir de las definiciones y su interpretaci\'on gr\'afica sobre la
circunferencia \emph{unidad} (de modo que la longitud del arco es igual al \'angulo)
\begin{center}
\begin{tikzpicture}[scale=0.63]
	   \draw[thick] (0,0) circle(4);
	   \draw[-] (-4.5,0) -- (4.5,0) ;
       \draw[-] (0,-4.5) -- (0,4.5) ;
	   \draw[-] (0,0) -- (4.25,4.25) node[midway,above,left] {$r=1$} ;
	   \draw[-,very thick,red] (4,0) -- (4,4) node[midway, right] {$\tan \alpha$};
	   \draw[-,very thick,blue] (2.83,0) -- (2.83,2.83);
	   \draw[-,very thick,orange] (4,0) arc (0:45:4);
	   \node[orange] at (22.5:3.5) {$\alpha$};
	   \draw[] (0.75,0) arc (0:45:0.75);	   
	   \node[blue] at (2,0.95) {$\sin \alpha$};
\end{tikzpicture}
\end{center}
Si son dos los \'angulos, y satisfacen $0<\beta<\alpha<\pi/2$, 
entonces resulta (dividiendo las correspondientes ecuaciones)
$$
\frac{\sin\alpha}{\sin\beta}<\frac{\alpha}{\beta}<\frac{\tan\alpha}{\tan\beta}\,.
$$

Otras propiedades son las siguientes: si $0<\alpha <90\degree$,
\begin{list}{$\blacktriangleright$}{}
\item $\sin (180\degree-\alpha )= \sin\alpha $
\item $\cos (180\degree-\alpha )= -\cos\alpha$.
\end{list}
Gr\'aficamente estas relaciones resultan obvias, como se aprecia en la figura:
\begin{center}
    \begin{tikzpicture}[scale=0.63]
	   \draw[thick] (0,0) circle(4);
	   \draw[-] (-4.5,0) -- (4.5,0) ;
       \draw[-] (0,-4.5) -- (0,4.5) ;
	   \draw[-] (0,0) -- (4.25,4.25) ;
	   \draw[-] (0,0) -- (-4.25,4.25) ;
	   \draw[-,very thick,blue] (2.83,0) -- (2.83,2.83) node[midway, right] {$\sin \alpha$};
	   \draw[-,very thick,orange] (0,0) -- (2.83,0) node[midway, below] {$\cos \alpha$};
	   \draw[thick] (0.75,0) arc (0:45:0.75);
	   \node[] at (22.5:1) {$\alpha$};
	   \draw[thick] (0.6,0) arc (0:135:0.6);
	   \node[] at (100:1) {$180-\alpha$};
	   \draw[-,very thick,blue] (-2.83,0) -- (-2.83,2.83) node[midway, left] {$\sin(180- \alpha)$};
	   \draw[-,very thick,orange] (0,0) -- (-2.83,0) node[midway, below] {$\cos (180-\alpha)$};
     \end{tikzpicture}
\end{center}
Tambi\'en es posible derivar las f\'ormulas para las razones trigonom\'etricas de la suma y diferencia de
\'angulos. En la figura se muestra la construcci\'on para el caso de la suma
(esta preciosa demostraci\'on visual se debe al matem\'atico persa Ab\={u} al-Waf\={a}\index{Abu al-Wafa},
quien vivi\'o aproximadamente entre los a\~nos 940-988):
\begin{list}{$\blacktriangleright$}{}
\item $\sin (\alpha +\beta )= \sin\alpha\cos\beta +\cos\alpha\sin\beta$
\item $\cos (\alpha +\beta )= \cos\alpha\cos\beta -\sin\alpha\sin\beta$.
\end{list}
\begin{center}
    \begin{tikzpicture}
\draw[-,thick,color=blue] (-4,0) -- (4,0) node[midway,below]{$\cos\alpha\cos\beta$} ;
\draw[-,thick,color=red] (4,0) --node[sloped,midway,below]{$\sin\alpha\cos\beta$} (4,4) ;
\draw[-] (-4,0) -- (4,4) node[midway,above,sloped]{$\cos \beta$} ;
\draw[-] (4,4) -- (2,8) node[midway,above,sloped]{$\sin\beta$} ;
\draw[-] (-4,0) -- (2,8) node[midway,above=8pt]{$1$} ;
\draw[-]  (2,8) -- (2,0) ;
\draw[-,thick,color=blue]  (2,4) -- (4,4) node[midway,above]{$\sin\alpha\sin\beta$};
\draw[-]  (-4,0) -- (-4,8) node[midway,above,sloped,red]{$\sin (\alpha +\beta )$};
\draw[-]  (-4,8) -- (2,8) node[midway,above,blue]{$\cos (\alpha +\beta)$};
\draw[-,thick,color=red] (2,4) -- (2,8) node[midway,above,sloped]{$\cos\alpha\sin\beta$} ;
\draw[thick] (2,7.25) arc (-90:-62:0.75)  ;
\node[] at (2.25,7.05)  {$\alpha$};

\draw[-,thick] (-9,0) -- (-5,0);
\draw[-,thick] (-5.25,0) -- (-5.25,0.25);
\draw[-,thick] (-5.25,0.25) -- (-5,0.25);
\draw[-,thick]  (-5,0) -- (-5,2);
\draw[-,thick]  (-9,0) -- (-5,2);
\draw[-,thick]  (-5,2) -- (-6,4);
\draw[-,thick]  (-5.25,1.875) -- (-5.375,2.125);
\draw[-,thick]  (-5.375,2.125) -- (-5.125,2.250);

\draw[-,thick]  (-9,0) -- (-6,4);
\draw[thick] (-8.25,0) arc (0:27:0.75)  ;
 \node[] at (-1.75:-8)  {$\alpha$};
\draw[thick] (-7.5,0) arc (0:46:1.75)  ;
 \node[] at (-8.75:-7.45)  {$\alpha+\beta$};
\draw[thick] (-8,0.485) arc (0:57:0.57)  ;
 \node[] at (-4:-8.25)  {$\beta$};
\end{tikzpicture}  
\end{center}

Como casos particulares de estas expresiones, tenemos las importantes f\'ormulas
\begin{align*}
\sin (2\alpha )& = 2\sin\alpha \cos\alpha\,, \\
\cos (2\alpha )& = \cos^2\alpha -\sin^2\alpha \,.
\end{align*}

A efectos pr\'acticos, conviene recordar la siguiente tabla de valores particulares: 
\begin{center}
\begin{tabular}{ccccc}
$\alpha$ & $\mathrm{rad}$ & $\sin\alpha$ & $\cos\alpha$ & $\tan\alpha$ \\ 
\hline 
\rule{0pt}{3ex} $0\degree$ & $0$ & $0$ & $1$ & $0$ \\ 
$30\degree$ & $\pi/6$ & $1/2$ & $\sqrt{3}/2$ & $\sqrt{3}/3$ \\ 
$45\degree$ & $\pi/4$ & $\sqrt{2}/2$ & $\sqrt{2}/2$ & $1$ \\ 
$60\degree$ & $\pi/3$ & $\sqrt{3}/2$ & $1/2$ & $\sqrt{3}$ \\ 
$90\degree$ & $\pi/2$ & $1$ & $0$ & $\infty$ \\ 
\hline 
\end{tabular} 
\end{center}

Para finalizar, diremos algunas palabras acerca de tri\'angulos arbitrarios,
no necesariamente rect\'angulos (como los \emph{escalenos}\index{Tri\'angulo!escaleno}, aqu\'ellos con 
todos sus lados de distinta longitud, o los \emph{is\'osceles}\index{Tri\'angulo!is\'osceles}, con dos 
lados iguales entre s\'i y distintos al tercero). Para un tri\'angulo como el de la figura,
\begin{center}
    \begin{tikzpicture}[scale=0.66]
	   \draw[-,very thick] (3,0) -- (4,4) node[midway, right] {$a$};
	   \draw[-,very thick] (0,0) -- (3,0) node[midway, below] {$b$};
	   \draw[-,very thick] (0,0) --node[above]{$c$} (4,4) ;
	   \draw[thick] (0.75,0) arc (0:45:0.75);
	   \node[] at (20:1.15) {$\alpha$};
	   \draw[thick] (3.5,3.5) arc (225:256:0.75);
	   \node[] at (40:4.5) {$\beta$};
	   \draw[thick] (2.5,0) arc (180:80:0.5);
	   \node[] at (16:2.75) {$\gamma$};	   	   	   
	\end{tikzpicture} 
\end{center}

\noindent se cumplen las llamadas \emph{leyes del coseno}\index{Ley!de los cosenos}:
\begin{align}\label{leyescoseno}
a^2& = b^2+c^2-2bc\cos\alpha\nonumber \\
b^2& = a^2+c^2-2ac\cos\beta \\
c^2& = a^2+b^2 -2ab\cos\gamma\nonumber\,.
\end{align}
Es importante recalcar que cuando uno de los \'angulos es rect\'angulo, su coseno
vale cero y entonces estas f\'ormulas se reducen a las ya conocidas.
Tambi\'en se tienen otras relaciones involucrando el seno (\emph{ley de los senos})\index{Ley!de los senos}:
\begin{equation}\label{leyeseno}
\frac{a}{\sin\alpha}=\frac{b}{\sin\beta}=\frac{c}{\sin\gamma}\,.
\end{equation}

Con las f\'ormulas vistas en el cap\'itulo, es posible resolver cualquier tri\'angulo 
conociendo tres de sus elementos, bien sean dos \'angulos y la longitud
de un lado, los tres lados, etc.

\section{Aplicaciones}

\subsection{C\'alculo de distancias astron\'omicas}
Desde los tiempos m\'as remotos, el hombre ha tenido necesidad de
conocer el movimiento de los astros. Desde luego, los dos m\'as importantes
son el Sol y la Luna, nuestro sat\'elite natural. Entre ambos determinan
dos de los eventos de mayor influencia en la vida sobre la Tierra: las estaciones (el Sol)
y las mareas (el Sol y la Luna).

El conocimiento acerca del Sol y la Luna se vi\'o obstaculizado por el hecho de que
(al menos en tiempos antiguos) no se ten\'ia acceso directo a ellos. S\'olo podemos
utilizar m\'etodos indirectos y aqu\'i es donde entra en juego la trigonometr\'ia.
Aristarco de Samos\index{Aristarco} (ca. 310-230 adC) fue un destacado matem\'atico griego que realiz\'o
algunos c\'alculos referentes al Sol y la Luna. Su observaci\'on fundamental es que
durante la media luna, el Sol, la Luna y la Tierra forman un \'angulo recto, como se puede
apreciar en la figura siguiente.
\begin{center}
    \begin{tikzpicture}[scale=4]
        \def\rS{0.2}                                

        \def\Earthangle{-14}                         
        \def\rE{0.1}                                
        \def\eE{0.25}                               
        \pgfmathsetmacro\bE{sqrt(1-\eE*\eE)}        

        \def\Moonangle{59}                         
        \pgfmathsetmacro\rM{.7*\rE}                 
        \pgfmathsetmacro\aM{2.5*\rE}                
        \def\eM{0.4}                                
        \pgfmathsetmacro\bM{\aM*sqrt(1-\eM*\eM)}    
        \def\offsetM{30}                            

        \pgfmathdeclarefunction{f}{1}{%
            \pgfmathparse{
                ((-\eE+cos(#1))<0) * ( 180 + atan( \bE*sin(#1)/(-\eE+cos(#1)) ) ) 
                +
                ((-\eE+cos(#1))>=0) * ( atan( \bE*sin(#1)/(-\eE+cos(#1)) ) ) 
            }
        }

        \pgfmathdeclarefunction{d}{1}{%
            \pgfmathparse{ sqrt((-\eE+cos(#1))*(-\eE+cos(#1))+\bE*sin(#1)*\bE*sin(#1)) }
        }

        \draw[very thin,color=gray!75] (0,0) ellipse (1 and \bE);
        
        \shade[
            top color=yellow!70,
            bottom color=red!70,
            shading angle={45},
            ] ({sqrt(1-\bE*\bE)},0) circle (\rS);

        \pgfmathsetmacro{\radiation}{100*(1-\eE)/(d(\Earthangle)*d(\Earthangle))}
        \colorlet{Earthlight}{yellow!\radiation!blue}
        \shade[%
            top color=Earthlight,%
            bottom color=blue,%
            shading angle={90+f(\Earthangle)},%
        ] ({cos(\Earthangle)},{\bE*sin(\Earthangle)}) circle (\rE);

        \draw[very thin,color=gray!75,rotate around={{\offsetM}:({cos(\Earthangle)},{\bE*sin(\Earthangle)})}]
            ({cos(\Earthangle)},{\bE*sin(\Earthangle)}) ellipse ({\aM} and {\bM});
        \shade[
            top color=black!60,
            bottom color=white,
            shading angle={0},
        ]   ({cos(\Earthangle)+\aM*cos(\Moonangle)*cos(\offsetM)-\bM*sin(\Moonangle)*sin(\offsetM)},%
            {\bE*sin(\Earthangle)+\aM*cos(\Moonangle)*sin(\offsetM)+\bM*sin(\Moonangle)*cos(\offsetM)}) circle (\rM);
        
        \draw[thin,color=black] (0.25,0) -- 
 ({cos(\Earthangle)+\aM*cos(\Moonangle)*cos(\offsetM)-\bM*sin(\Moonangle)*sin(\offsetM)},%
            {\bE*sin(\Earthangle)+\aM*cos(\Moonangle)*sin(\offsetM)+\bM*sin(\Moonangle)*cos(\offsetM)});
        \draw[thin,color=black] (0.25,0) -- ({cos(\Earthangle)},{\bE*sin(\Earthangle)});                    
        \draw[thin,color=black]
 ({cos(\Earthangle)+\aM*cos(\Moonangle)*cos(\offsetM)-\bM*sin(\Moonangle)*sin(\offsetM)},%
  {\bE*sin(\Earthangle)+\aM*cos(\Moonangle)*sin(\offsetM)+\bM*sin(\Moonangle)*cos(\offsetM)})      
   -- ({cos(\Earthangle)},{\bE*sin(\Earthangle)});
\end{tikzpicture}
\end{center}
\noindent Aristarco consider\'o los segmentos $\overline{TL}$, que une la Tierra con la Luna, y
$\overline{TS}$, que une la Tierra con el Sol, y propuso una manera de determinar su longitud
relativa. Consideremos el siguiente diagrama:
\begin{center}
\begin{tikzpicture}
\coordinate (O) at (-3,0);
\coordinate (A) at (3,0);
\coordinate (B) at (3,-2);
\draw (O)--(A)--(B)--cycle;
\tkzLabelSegment[above=2pt](O,A){$\overline{SL}$}
\tkzLabelSegment[below=2pt](O,B){$\overline{TS}$}
\tkzLabelSegment[right=2pt](A,B){$\overline{TL}$}
\tkzMarkAngle[fill= Gainsboro,size=0.85cm,opacity=.4](B,O,A)
\tkzLabelAngle[pos = 1.15](A,O,B){$\varphi$}
\tkzMarkAngle[fill= Gainsboro,size=0.65cm,opacity=.4](A,B,O)
\tkzLabelAngle[pos = 0.95](A,B,O){$\psi$} 
\end{tikzpicture}
\end{center}
Es claro que
$$
\overline{TS}=\frac{\overline{TL}}{\sin\varphi}\,.
$$
Para evaluar num\'ericamente este resultado, se necesitaba una estimaci\'on de uno de los
dos \'angulos $\varphi$ \'o $\psi$ (por tratarse de un tri\'angulo rect\'angulo, 
el segundo ser\'ia igual a
$\pi/2$ menos el primero). Aristarco consider\'o que un ciclo lunar completo dura $29{.}5$
d\'ias, y que el doble del \'angulo $\psi$ subtiende el tramo de la \'orbita comprendido
entre la media luna menguante y la media luna creciente:\\
\medskip
\begin{center}
    \begin{tikzpicture}[scale=4]
        \def\rS{0.2}                                

        \def\Earthangle{0}                         
        \def\rE{0.1}                                
        \def\eE{0.25}                               
        \pgfmathsetmacro\bE{sqrt(1-\eE*\eE)}        

        \def\Moonangle{78.5}                         
        \def\Moonanglem{-133.5}                         
        \pgfmathsetmacro\rM{.7*\rE}                 
        \pgfmathsetmacro\aM{2.5*\rE}                
        \def\eM{0.4}                                
        \pgfmathsetmacro\bM{\aM*sqrt(1-\eM*\eM)}    
        \def\offsetM{30}                            

        \pgfmathdeclarefunction{f}{1}{%
            \pgfmathparse{
                ((-\eE+cos(#1))<0) * ( 180 + atan( \bE*sin(#1)/(-\eE+cos(#1)) ) ) 
                +
                ((-\eE+cos(#1))>=0) * ( atan( \bE*sin(#1)/(-\eE+cos(#1)) ) ) 
            }
        }

        \pgfmathdeclarefunction{d}{1}{%
            \pgfmathparse{ sqrt((-\eE+cos(#1))*(-\eE+cos(#1))+\bE*sin(#1)*\bE*sin(#1)) }
        }

        
        \shade[
            top color=yellow!70,
            bottom color=red!70,
            shading angle={45},
            ] ({sqrt(1-\bE*\bE)},0) circle (\rS);

        \pgfmathsetmacro{\radiation}{100*(1-\eE)/(d(\Earthangle)*d(\Earthangle))}
        \colorlet{Earthlight}{yellow!\radiation!blue}
        \shade[%
            top color=Earthlight,%
            bottom color=blue,%
            shading angle={90+f(\Earthangle)},%
        ] ({cos(\Earthangle)},{\bE*sin(\Earthangle)}) circle (\rE);

        \shade[
            top color=black!60,
            bottom color=white,
            shading angle={0},
        ]   ({cos(\Earthangle)+\aM*cos(\Moonangle)*cos(\offsetM)-\bM*sin(\Moonangle)*sin(\offsetM)},%
            {\bE*sin(\Earthangle)+\aM*cos(\Moonangle)*sin(\offsetM)+\bM*sin(\Moonangle)*cos(\offsetM)}) circle (\rM);
            
        \shade[
            top color=black!60,
            bottom color=white,
            shading angle={0},
        ]   ({cos(\Earthangle)+\aM*cos(\Moonanglem)*cos(\offsetM)-\bM*sin(\Moonanglem)*sin(\offsetM)},%
            {\bE*sin(\Earthangle)+\aM*cos(\Moonanglem)*sin(\offsetM)+\bM*sin(\Moonanglem)*cos(\offsetM)}) circle (\rM);            
        
        \draw[thin,color=black] (0.25,0) -- 
 ({cos(\Earthangle)+\aM*cos(\Moonangle)*cos(\offsetM)-\bM*sin(\Moonangle)*sin(\offsetM)},%
            {\bE*sin(\Earthangle)+\aM*cos(\Moonangle)*sin(\offsetM)+\bM*sin(\Moonangle)*cos(\offsetM)});
        \draw[thin,color=black] (0.25,0) -- ({cos(\Earthangle)},{\bE*sin(\Earthangle)});                    
        \draw[thin,color=black]
 ({cos(\Earthangle)+\aM*cos(\Moonangle)*cos(\offsetM)-\bM*sin(\Moonangle)*sin(\offsetM)},%
  {\bE*sin(\Earthangle)+\aM*cos(\Moonangle)*sin(\offsetM)+\bM*sin(\Moonangle)*cos(\offsetM)})      
   -- ({cos(\Earthangle)},{\bE*sin(\Earthangle)});
   
       \draw[thin,color=black] (0.25,0) -- 
 ({cos(\Earthangle)+\aM*cos(\Moonanglem)*cos(\offsetM)-\bM*sin(\Moonanglem)*sin(\offsetM)},%
            {\bE*sin(\Earthangle)+\aM*cos(\Moonanglem)*sin(\offsetM)+\bM*sin(\Moonanglem)*cos(\offsetM)});
        \draw[thin,color=black] (0.25,0) -- ({cos(\Earthangle)},{\bE*sin(\Earthangle)});                    
        \draw[thin,color=black]
 ({cos(\Earthangle)+\aM*cos(\Moonanglem)*cos(\offsetM)-\bM*sin(\Moonanglem)*sin(\offsetM)},%
  {\bE*sin(\Earthangle)+\aM*cos(\Moonanglem)*sin(\offsetM)+\bM*sin(\Moonanglem)*cos(\offsetM)})      
   -- ({cos(\Earthangle)},{\bE*sin(\Earthangle)});   
\end{tikzpicture}
\end{center}
\bigskip
\noindent El tiempo entre dos medias lunas dura aproximadamente $14{.}25$ d\'ias,
as\'i que suponiendo
que la velocidad de la Luna es uniforme, podemos determinar $2\psi$ mediante una simple
proporcionalidad,
$$
2\psi =2\pi\frac{14{.}25}{29{.}5}=3{.}035\, \mbox{rad}\simeq 174\degree\,,
$$
que conduce a $\psi \simeq 87\degree\simeq 1{.}52\,\mbox{rad}$. En consecuencia, 
$\varphi\simeq 3\degree\simeq 0{.}052\,\mbox{rad}$. Sustituyendo, encontramos la famosa
estimaci\'on de Aristarco,
$$
\overline{TS}\simeq 19\overline{TL}\,.
$$
As\'i, la distancia entre el Sol y la Tierra ser\'ia unas $19$ veces mayor que la que hay de la Tierra a la Luna. El valor de esta estimaci\'on no reside en el n\'umero concreto que proporciona (que es una aproximaci\'on muy pobre: hoy sabemos que la distancia real es del orden de $\overline{TS}\simeq 390\overline{TL}$), sino en el m\'etodo empleado. 
Es un buen ejemplo de
c\'omo las matem\'aticas pueden porporcionar informaci\'on sobre partes del universo
f\'isico que no son directamente accesibles al experimento.

\subsection{Navegaci\'on a\'erea por deriva}
Vamos a ver una aplicaci\'on m\'as actual de la trigonometr\'ia, relacionada con
el c\'alculo de trayectorias en aeron\'autica. Cuando un piloto va a realizar un
trayecto entre dos puntos, como por ejemplo el aeropuerto de San Luis Potos\'i 
(c\'odigo ICAO\index{ICAO}\footnote{Siglas de \emph{International
Civil Aviation Organization}.} MMSP) y el de la Ciudad de M\'exico
(MMMX), debe tener en cuenta que el trayecto no se
realizar\'a en l\'inea recta si hay viento. El procedimiento para corregir
la trayectoria es el siguiente:
\begin{enumerate}
\item En primer lugar, el piloto marca en su carta de navegaci\'on\index{Navegaci\'on!carta} el origen y destino,
trazando un segmento entre ellos.
El \'angulo formado por este segmento y la direcci\'on
norte-sur es el \emph{curso verdadero}\index{Curso verdadero} (\emph{true course}, en ingl\'es). Se puede
determinar este \'angulo directamente sobre la carta de navegaci\'on porque \'esta
trae marcas con los meridianos locales a espacios de una milla n\'autica. En nuestro
caso, el curso verdadero es de $143\degree$.
\item Se obtiene la informaci\'on sobre el viento a lo largo de la trayectoria.
Usualmente, \'esta la provee el aeropuerto de destino por medio de un reporte
metereol\'ogico (METAR report)\index{METAR}, que contiene todos los datos relevantes. Para
el prop\'osito de esta secci\'on, podemos acudir a la p\'agina
\href{http://skyvector.com}{\tt http://skyvector.com}, que nos proporciona de manera
gratuita una carta de navegaci\'on bastante completa junto con el reporte METAR
que necesitamos.
\begin{center}
\includegraphics[scale=0.46]{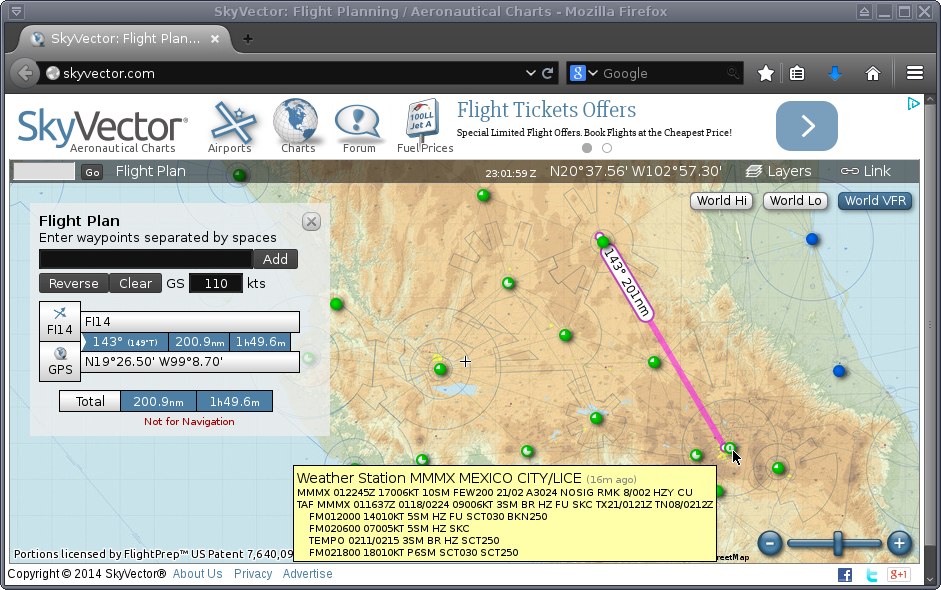} 
\end{center}
\'Este lo podemos ver picando en el aeropuerto de inter\'es,
es la primera l\'inea del desplegable que aparece. Otra posibilidad, mucho m\'as fiable,
consiste en acudir a la
p\'agina (dependiente del Servicio Nacional de 	Climatolog\'ia de los Estados Unidos) \href{http://www.aviationweather.gov/adds/metars/}{\tt www.aviationweather.gov}
e introducir el c\'odigo ICAO del aeropuerto que nos interese.
\begin{center}
\includegraphics[scale=0.3]{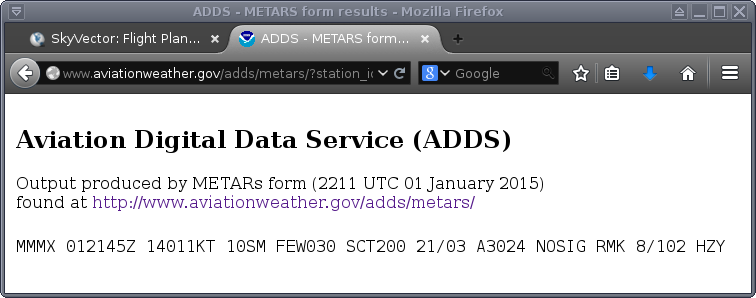} 
\end{center}
El primer bloque del reporte METAR es el c\'odigo del aeropuerto (MMMX). El segundo
contiene la fecha (d\'ia 01) y la hora ($21:45$ Zul\'u).
El tercer bloque es el que nos interesa, aqu\'i est\'a la informaci\'on sobre el 
viento. Los tres primeros n\'umeros ($140$) dan la direcci\'on del viento en grados, los dos
siguientes ($11$) junto con la abreviatura KT nos dice que su velocidad es de $11$ nudos
(en ingl\'es, \emph{knots}). Como un nudo equivale a $1{.}852$Km/h, la velocidad del
viento en estas unidades es de unos $20{.}4$Km/h. Los siguientes bloques no los
necesitaremos\footnote{Por ejemplo, el siguiente ($10$SM) proporciona la visibilidad en millas.}.
\item Ya estamos en condiciones de calcular el \emph{rumbo verdadero}\index{Rumbo verdadero} (en ingl\'es,
\emph{true heading}), que es el que corrige la acci\'on del viento. El \'angulo 
$\theta$ entre el curso verdadero y el rumbo verdadero se denomina 
\emph{\'angulo de correcci\'on del viento} o \emph{\'angulo de deriva}
(en ingl\'es, \emph{drift angle}).\index{An@\'Angulo!de deriva}
\begin{center}
\begin{tikzpicture}
\draw[->] (0,0) -- (0,3.5)node[right]{N};
\draw[] (-1.875,1.03) arc (180:141:0.5);
\node[] at (147:2.4)  {$\theta$};
\draw[->,red,thick] (-4,3) -- (0,0)node[midway,above right]{rumbo verdadero};
\draw[->,orange,thick] (-3.75,2.06) -- (-4,3)node[midway,left]{viento};
\draw[->,blue,thick] (-3.75,2.06) -- (0,0)node[midway,below left]{curso verdadero};
\draw[color=black,fill=black] (-3.75,2.06) circle(0.05)node[below left]{posici\'on actual};
\end{tikzpicture}
\end{center}
Si, una vez en aire, el aeroplano se dirige al destino con un \'angulo igual al de deriva
seguir\'a el curso verdadero, marcado en la carta de navegaci\'on.
\item El \'angulo de deriva se encuentra teniendo en cuenta las longitudes de los lados en
el tri\'angulo precedente, que vienen dadas por las velocidades del viento ($11$kts)
y la velocidad en aire del aeroplano, que seguir\'a a lo largo del rumbo verdadero 
(esta velocidad viene en su manual y nosotros la supondremos de unos 
$120$kts\footnote{Por ejemplo la Cessna Skyhawk, v\'ease
\url{http://cessna.txtav.com/en/piston/cessna-skyhawk\#_model-specs}.}).
Resulta as\'i el llamado \emph{tri\'angulo de viento}:\index{Tri\'angulo!de viento}
\begin{center}
\begin{tikzpicture}
\draw[->] (0,0) -- (0,3.5)node[right]{N};
\draw[] (0,0) -- (1.5,0);
\draw[] (-1.875,1.03) arc (180:141:0.5);
\draw[thin,gray!50] (-3.75,2.06) -- (-2.5,2.06);
\draw[] (-3.6,2.06) arc (0:108:0.15);
\node[] at (147:2.4)  {$\theta$};
\node[] at (45:1.4)  {$143\degree$};
\draw[->,red,thick] (-4,3) -- (0,0)node[midway,above right]{$120$};
\draw[->,orange,thick] (-3.75,2.06) -- (-4,3)node[midway,left]{$11$};
\draw[->,blue,thick] (-3.75,2.06) -- (0,0);
\draw[->] (1,0) arc (0:150:1);
\draw[color=black,fill=black] (-3.75,2.06) circle(0.05);
\node[] at (144:4)  {$140\degree$};
\end{tikzpicture}
\end{center}
\end{enumerate}

En este ejemplo, un sencillo c\'alculo trigonom\'etrico nos dice que se tienen las
relaciones adicionales

\begin{center}
\begin{tikzpicture}
\draw[->] (0,0) -- (0,3.5)node[right]{N};
\draw[] (-1.5,0) -- (1.5,0);
\draw[thin,gray!50] (-3.75,2.06) -- (-2.5,2.06);
\draw[] (-3.6,2.06) arc (0:108:0.15);
\node[] at (45:1.4)  {$143\degree$};
\draw[->,red,thick] (-4,3) -- (0,0)node[midway,above right]{$120$};
\draw[->,orange,thick] (-3.75,2.06) -- (-4,3)node[midway,left]{$11$};
\draw[->,blue,thick] (-3.75,2.06) --node[midway,below left]{$gs$} (0,0);
\draw[->] (1,0) arc (0:150:1);
\draw[color=black,fill=black] (-3.75,2.06) circle(0.05);
\node[] at (144:4)  {$140\degree$};
\draw[->] (-1,0) arc (180:150:1);
\node[] at (170:1.4)  {$37\degree$};
\draw[->] (-3.1,2.06) arc (0:-37:0.5);
\node[] at (147:3.25)  {$37\degree$};
\end{tikzpicture}
\end{center}

El lado restante, que nos da la velocidad del aeroplano con respecto a un observador
en el suelo $gs$ (\emph{ground speed}, en ingl\'es), se obtiene f\'acilmente de la ley de
los cosenos \eqref{leyescoseno}:\index{Ley!de los cosenos}
$$
120^2=11^2+gs^2-2\cdot 11\cdot gs\cdot \cos 177\,,
$$
ecuaci\'on cuadr\'atica cuya soluci\'on positiva (la que tiene sentido) es
$$
gs=109{.}01\mbox{kts}\,.
$$
En unidades m\'as comunes, esto es una velocidad de aproximadamente $200$Km/h.
La distancia San Luis-M\'exico es de unas $200$mn (millas n\'auticas), es
decir, unos $370{.}5$Km, por lo que el tiempo estimado de vuelo ser\'a de
$1{.}85$ horas (unos $110$ minutos). Con
estos datos, el piloto ya puede calcular el combustible que se va a consumir
en el trayecto.

Calculemos ahora el \'angulo de deriva $\theta$ con el que hay que dirigir
el aeroplano hacia el destino. En este caso, lo m\'as c\'omodo es aplicar la ley
de los senos \eqref{leyeseno}:\index{Ley!de los senos}
$$
\frac{120}{\sin 177}=\frac{11}{\sin\theta}\,,
$$
de donde
$$
\theta=\arcsin \left(\frac{11}{2293.19} \right)\simeq 0{.}27\,.
$$
N\'otese que en este caso, casi no se requiere ajuste, como cab\'ia esperar
del hecho de que el viento sopla pr\'acticamente en la misma direcci\'on del vuelo. En
este caso el efecto del viento es simplemente el de reducir la velocidad efectiva del
aeroplano de $120$ a $109$kts. El c\'alculo que hemos hecho, sin embargo, es aplicable
a cualquier situaci\'on.

Finalmente, observemos que en una situaci\'on de vuelo real a\'un hay algunos c\'alculos adicionales que realizar, como por ejemplo el ajuste del norte geogr\'afico (con el que
hemos trabajado) al magn\'etico que indicar\'ia una br\'ujula sobre el aeroplano. Estos
c\'alculos son muy sencillos y se dejan al cuidado del lector interesado como ejercicio.

\section{Ejercicios}

\begin{enumerate}  
\setcounter{enumi}{87}

\item Los siguientes \'angulos est\'an dados en grados. Pasarlos a radianes:
\begin{multicols}{3}
\begin{enumerate}[(a)]
\item $360\degree $
\item $ -30\degree$
\item $15\degree $
\item $150\degree $
\item $ 90\degree$
\item $ 45\degree$
\end{enumerate}
\end{multicols}

\item Los siguientes \'angulos est\'an dados en radianes. Pasarlos a grados:
\begin{multicols}{3}
\begin{enumerate}[(a)]
\item $\pi /6$ rad
\item $ \pi /5$ rad
\item $-5\pi /3$ rad
\item $\pi $ rad
\item $1 $ rad
\item $3\pi /7 $ rad
\end{enumerate}
\end{multicols}

\item Dado un \'angulo arbitrario, calcular $\tan (\pi+\alpha )$. 
Haciendo un diagrama apropiado, determinar tambi\'en:
\begin{enumerate}[(a)]
\item $\sin (\pi/2+\alpha )$
\item $\cos (\pi/2+\alpha )$
\item $\tan (\pi/2+\alpha )$.
\end{enumerate}

\item Representar gr\'aficamente las razones trigonom\'etricas principales 
($\sin$, $\cos$, $\tan$) de los \'angulos siguientes:
\begin{multicols}{4}
\begin{enumerate}[(a)]
\item $ \pi /3$
\item $3\pi /4 $
\item $7\pi /6 $
\item $3\pi /2 $
\end{enumerate}
\end{multicols}

\item Calcular las razones trigonom\'etricas principales 
($\sin$, $\cos$, $\tan$) de los \'angulos siguientes:
\begin{multicols}{2}
\begin{enumerate}[(a)]
\item $\pi /3,4\pi /3,-7\pi /3$
\item $ \pi /6,5\pi /6,19\pi /6$
\item $\pi /4,-\pi /4,27\pi /4 $
\item $ \pi /2,7\pi,12\pi ,51\pi /2$
\end{enumerate}
\end{multicols}

\item Representar gr\'aficamente los \'angulos correspondientes y dar su valor en radianes:
\begin{multicols}{3}
\begin{enumerate}[(a)]
\item $ \cos \alpha =\sqrt{3}/2 $
\item $ \sin \alpha =-\sqrt{2}/2 $
\item $ \tan \alpha =2$
\item $ \cos \alpha =3/5$
\item $ \sin \alpha =4/5$
\end{enumerate}
\end{multicols}

\item Simplificar las siguientes expresiones
\begin{multicols}{2}
\begin{enumerate}[(a)]
\item $5\sin\theta -2\sin(-\theta)$
\item $4 \cos\theta- \cos(-\theta)$
\item $4\cos(\frac{\pi}{2}-\theta)-3\sin\theta$
\item $4 \cos\theta + 2 \cos(-\theta)$
\item $\tan(-\theta) + \tan\theta$
\item $\cos(-\alpha) \cos\alpha - \sin(-\alpha) \sin \alpha$
\item $1+\tan\alpha /\cot\alpha$
\item $\tan^2\theta /(1+\sec\theta)$
\item $\tan\theta\sqrt{1-\sin^2\theta} /\sin\theta$
\end{enumerate}
\end{multicols}

\item Probar las siguientes relaciones
\begin{multicols}{2}
\begin{enumerate}[(a)]
\item $(\sin\alpha +\cos\alpha )^2 =1+\sin (2\alpha)$
\item $\cos^4\alpha -\sin^4\alpha =\cos (2\alpha)$
\item $\csc \theta -\sin\theta =\cot\theta\cos\theta$
\item $1/(\tan\theta +\cot\theta ) =\sin\theta\cos\theta$
\item $\cos\alpha\sec\alpha /(1+\tan^2\alpha )=\cos^2\alpha$
\item $\sin 3\alpha +\sin\alpha =2\sin (2\alpha)\cos\alpha$
\end{enumerate}
\end{multicols}

\item Para cada apartado, representar gr\'aficamente las funciones que se indican
(las funciones $\arcsin$, $\arccos$ y $\arctan$ son las inversas de $\sin$, $\cos$
y $\tan$, respectivamente):
\begin{multicols}{3}
\begin{enumerate}[(a)]
\item $ f(x)=\sin x,g(x)=\cos x$
\item $ f(x)=\sin x,g(x)=\cos 2x$
\item $ f(x)=\tan x$
\item $ f(x)=\arcsin x$
\item $ f(x)=\arccos x$
\item $ f(x)=\arctan x$
\end{enumerate}
\end{multicols}

\item Las funciones trigonom\'etricas son de uso frecuente en F\'isica para modelar
fen\'omenos oscilatorios, que ejecutan el llamado \emph{movimiento arm\'onico simple}.\index{Movimiento arm\'onico simple}
El ejemplo m\'as sencillo es el de una onda descrita por una funci\'on como
$$
\psi (t)=A\sin (wt +\phi)\,.
$$
Aqu\'i, el par\'ametro $A$ representa la \emph{amplitud}\index{Amplitud} de la onda, $w$ es su 
\emph{frecuencia}\index{Frecuencia} y $\phi$ la \emph{fase}\index{Fase}. La siguiente figura, en la que se representan
las ondas $f(x)=\sin x$, $g(x)=\sin (x/2)$ y $h(x)=2\sin (x/2 +1)$, muestra el efecto de
estos par\'ametros.
\begin{center}
\begin{tikzpicture}[line cap=round,line join=round,>=triangle 45,x=1.0cm,y=1.0cm]
\draw[color=black] (-3.5,0.0) -- (11.5,0.0);
\foreach \x in {-pi,-pi/2,pi/2,pi,3*pi/2,2*pi,5*pi/2,3*pi}
\draw[shift={(\x,0)},color=black,fill=black] circle(0.05);
\draw[fill=black,color=black] (-pi,0)node[below]{$-\pi$};
\draw[fill=black,color=black] (-pi/2,0)node[below]{$-\pi/2$};
\draw[fill=black,color=black] (pi/2,0)node[below]{$\pi/2$};
\draw[fill=black,color=black] (pi,0)node[below]{$\pi$};
\draw[fill=black,color=black] (3*pi/2,0)node[below]{$3\pi/2$};
\draw[fill=black,color=black] (2*pi,0)node[below]{$2\pi$};
\draw[fill=black,color=black] (5*pi/2,0)node[below]{$5\pi/2$};
\draw[fill=black,color=black] (3*pi,0)node[below]{$3\pi$};
\draw[color=black] (0,-2.35) -- (0,2.35);
\foreach \y in {-2,-1,1,2}
\draw[shift={(0,\y)},color=black,fill=black] circle(0.05) node[left] {$\y$};
\draw[color=black] (0pt,-10pt) node[right] {$0$};
\clip(-3.5,-2.6) rectangle (13.45,2.35);
\draw[color=red,thick,smooth,samples=100,domain=-3.5:11.45] plot(\x,{sin(((\x))*180/pi)})node[below]{$f(x)$};
\draw[color=blue,thick,smooth,samples=100,domain=-3.5:11.45] plot(\x,{sin(((\x)/2.0)*180/pi)})node[above]{$g(x)$};
\draw[color=orange,thick,smooth,samples=100,domain=-3.5:11.45] plot(\x,{2.0*sin(((\x)/2.0+1.0)*180/pi)})node[above]{$h(x)$};
\end{tikzpicture}
\end{center}
La amplitud\index{Ondas!amplitud} $A$ representa la mitad de diferencia entre las alturas m\'axima y m\'inima que
alcanza la onda. La frecuencia\index{Ondas!frecuencia} $w$ da el n\'umero de oscilaciones completas en un
intervalo de longitud $2\pi$, y la fase $\phi$ mide el desplazamiento de los picos
respecto a los que tendr\'ia la funci\'on seno. Teniendo en cuenta estas interpretaciones,
aplicables a cualquier fen\'omeno peri\'odico, proponer un modelo que se ajuste a las
siguientes observaciones de datos de temperaturas mensuales en David City (Nebraska, USA)\index{David City}\footnote{Datos reales obtenidos de \href{http://www.usclimatedata.com/climate/david-city/nebraska/united-states/usne0135}{\tt www.usclimatedata.com}. Se han redondeado al primer
d\'igito decimal por simplicidad.}
\begin{center}
\begin{tabular}{|c|c|c|c|c|c|c|c|c|c|c|c|c|}
\hline 
\rule[-1ex]{0pt}{2.5ex} • & Ene & Feb & Mar & Abr & May & Jun & Jul & Ago & Sep & Oct & Nov & Dic \\ 
\hline 
\rule[-1ex]{0pt}{2.5ex} 2004 & 18.7 & 20.9 & 41.8 & 51.7 & 61.3 & 66.3 & 71.3 & 69.1 & 69.2 & 53.5 & 40.2 & 27.7 \\ 
\hline 
\rule[-1ex]{0pt}{2.5ex} 2005 & 18 & 31.9 & 38.5 & 52.4 & 60.5 & 72.9 & 77.2 & 74.3 & 68.7 & 53.9 & 40.6 & 22.6 \\ 
\hline 
\rule[-1ex]{0pt}{2.5ex} 2006 & 34.7 & 27 & 36.5 & 53.9 & 62.5 & 72.3 & 77.3 & 73.5 & 60.7 & 48.8 & 38.2 & 30.7 \\ 
\hline 
\rule[-1ex]{0pt}{2.5ex} 2007 & 20.7 & 19.6 & 44.5 & 47.1 & 64.6 & 70.7 & 76.2 & 76.3 & 65.2 & 55.3 & 38 & 19.3 \\ 
\hline 
\end{tabular} 
\end{center}
Utilizar el modelo obtenido para predecir la temperatura en Septiembre de 2008 y comparar el 
resultado con la temperatura promedio real de ese mes, disponible en la p\'agina web
\href{http://www.usclimatedata.com/climate/david-city/nebraska/united-states/usne0135/2008/9}{\tt www.usclimatedata.com/climate/david-city/nebraska/united-states/usne0135/2008/9}.

\item Probar las leyes de los cosenos \eqref{leyescoseno}.

\item Probar la ley de los senos \eqref{leyeseno}.

\item Resolver las siguientes ecuaciones (para $0\leq\theta\leq 2\pi$):
\begin{multicols}{2}
\begin{enumerate}[(a)]
\item $2\sin^2\theta = \sin\theta$
\item $3\cos^2 \theta =2\cos\theta +1$
\item $2\sin^2\theta =3\sin\theta +5$
\item $6\cos^2\theta +5\cos\theta =-1$ 
\end{enumerate}
\end{multicols}

\item El profe Hern\'an\index{Hern\'an, profe} anda por ah\'i proponiendo este ejercicio a todo el que se encuentra por los
pasillos. Resu\'elvanlo y as\'i, si se topan con \'el, se lo pueden quitar de encima r\'apidamente
d\'andole la respuesta\footnote{Por si les urge esto \'ultimo, la respuesta es $60/37\simeq 1{.}62$.}.
Se trata de calcular el lado del cuadrado inscrito en la figura, donde el 
tri\'angulo es rect\'angulo de lados $3$, $4$ y $5$.
\begin{center}
    \begin{tikzpicture}[scale=1.5]
	   \draw[-] (-2,0) -- (2,0) ;
       \draw[-] (-2,0) -- (-2,3) ;
	   \draw[-] (-2,3) -- (2,0) ;
	   \draw[-,very thick,red] (-0.7,0) -- (-2,0.97) ;
       \draw[-,very thick,red] (-2,0.97) -- (-1.03,2.26) ;
	   \draw[-,very thick,red] (-0.7,0) -- (0.27,1.3) ;
	   \draw[-,very thick,red] (-1.03,2.26) -- (0.27,1.3) ;	   
     \end{tikzpicture}
\end{center}

\item Se tienen dos circunferencias de radios $a<b$, separadas por una distancia
$d>a+b$. Se traza una recta tangente a ambas circunferencias como se muestra en la
figura. Calcular el \'angulo $\alpha$ que forma esta recta con la que une los centros
de las circunferencias.
\begin{center}
\begin{tikzpicture}[line cap=round,line join=round,>=triangle 45]
\clip(-5.06,-3.5) rectangle (6.46,2.84);
\draw [shift={(-0.1,0.0)},color=gray,fill=gray,fill opacity=0.1] (-0.9,0) -- (0:0) arc (0:32:0.8) -- cycle;
\draw [thick,color=red] (-3.0,0.0) circle (1.0);
\draw [thick,color=red] (3.0,0.0) circle (2.0);
\draw [thick,color=blue,domain=-5.06:6.46] plot(\x,{(-2.66-2.6712*\x)/-4.639});
\draw (-3.0,0.0)-- (3.0,0.0);
\draw (-3.0,0.0)-- (-2.55,0.89);
\draw (3.0,0.0)-- (3.88,1.80);

\draw [fill=black] (-3.0,0.0) circle (0.05);
\draw [fill=black] (3.0,0.0) circle (0.05);
\draw[color=black] (0.16,-0.24) node {$d$};
\draw[color=black] (-2.44,0.48) node {$a$};
\draw[color=black] (3.78,0.92) node {$b$};
\draw[] (0.1,0.32) node {$\alpha$};

\end{tikzpicture}
\end{center}

\end{enumerate}

\section{Uso de la tecnolog\'ia de c\'omputo}
Vamos a ver c\'omo visualizar las funciones trigonom\'etricas usando GeoGebra.\index{GeoGebra!deslizador}
Comenzamos creando una circunferencia de radio unidad (fijando un punto arbitrario
y usando el comando \texttt{Circunferencia dado su Centro y Radio}). Sobre esa circunferencia
seleccionamos un punto, que se ir\'a moviendo sobre ella y al que le iremos calculando
el coseno y el seno de su radiovector\index{Radiovector} (el segmento que lo une con el centro de la
circunferencia). Por conveniencia, elegimos el punto sobre el di\'ametro horizontal
(si se toma el origen de la circunferencia en $A\equiv (0,0)$, este punto es simplemente
el $B\equiv (1,0)$). Creamos adem\'as un deslizador de tipo \'angulo, con un tama\~no
relativamente grande (por ejemplo, $200$px). El siguiente paso es la clave de esta
construcci\'on: seleccionamos el comando \texttt{\'Angulo dada su Amplitud}, tomando como
referencia el segmento $\overline{BA}$ formado por el centro de la circunferencia y
el punto elegido y estableciendo el sentido antihorario.
Como \'angulo \emph{se toma el $\theta$ del deslizador} (esto autom\'aticamente
renombrar\'a el deslizador a $\theta_1$).
\begin{center}
\includegraphics[scale=0.4]{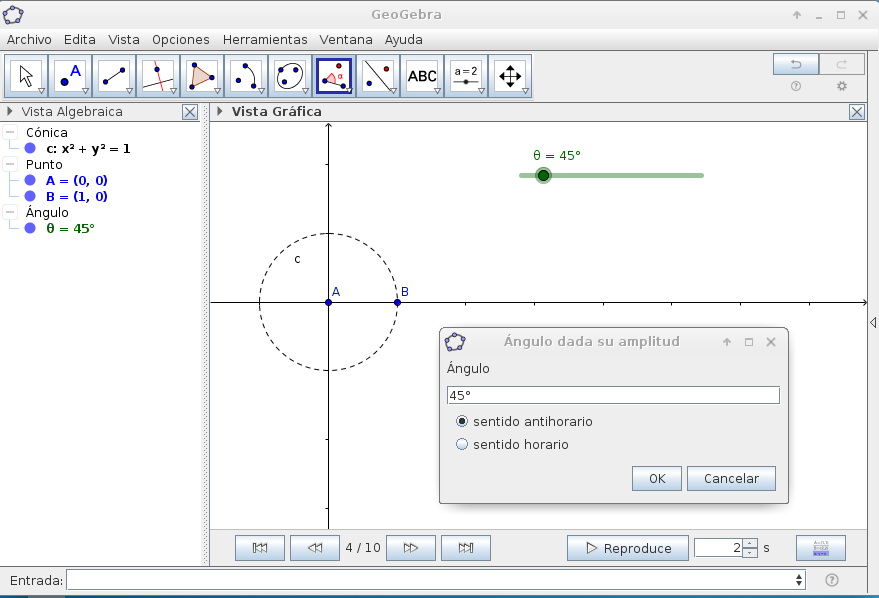}  
\end{center}
Esto generar\'a un nuevo punto $B'$ que es el `desplazado' de $B$
un \'angulo indicado por $\theta$. Para una mejor visualizaci\'on, seleccionamos
\texttt{Arco de Circunferencia con Centro entre dos Puntos} y marcamos el arco
$\widearc{BB'}$, d\'andole un color llamativo (digamos, rojo). Para
que el gr\'afico se vea menos abigarrado, se pueden eliminar algunas etiquetas y
renombrar los puntos (por ejemplo, ocultando la etiqueta de $A$ y $B$ y renombrando
el punto $B'$ a $P$); tambi\'en se puede a\~nadir color al radiovector de $P$, por
ejemplo con un gris suave.

Ahora creamos el punto que, en su movimiento, trazar\'a la gr\'afica de la funci\'on
seno. Este punto tendr\'a coordenadas $(d,y(P))$, donde $d$ es la longitud del arco
$\widearc{BP}$ y el n\'umero $y(P)$ es la segunda coordenada del $P$ (su altura).
Para ello, simplemente introducimos \texttt{Q=(d,y(P))} en la ventana de input. A
continuaci\'on, en la ventana gr\'afica seleccionamos el punto reci\'en creado $Q$
y modificamos sus propiedades para que deje traza y tenga un color apropiado (en
el ejemplo, un naranja).

Para trazar la gr\'afica del coseno, repetimos el proceso, creando un nuevo punto
$R=(d,x(P))$, de manera que su ordenada es la proyecci\'on sobre el eje de abcisas
de $P$, una proyecci\'on que por encontrarse $P$ en la circunferencia unidad sabemos que coincide con el coseno del \'angulo de su radiovector.
La figura muestra el resultado final, ya con algunos retoques est\'eticos.

\begin{center}
\includegraphics[scale=0.5]{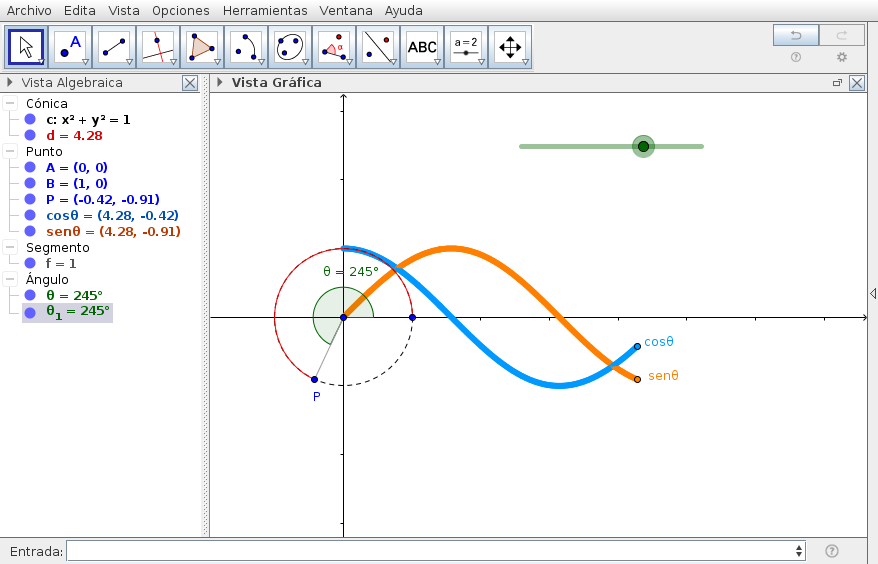} 
\end{center}

De esta construcci\'on resulta aparente una de las principales propiedades de
las funciones seno y coseno, la acotaci\'on \eqref{cota1}.

\section{Actividades de c\'omputo}
\begin{enumerate}[(A)]
\item Crear un programa en Maxima que acepte como entrada dos listas: una lista de dos elementos \texttt{l1:[l1[1],l1[2]]} conteniendo el curso verdadero \texttt{l1[1]} medido
en grados entre un origen y un destino junto con la velocidad del aeroplano \texttt{l1[2]},
y una segunda lista \texttt{l2:[l2[1],l2[2]]} cuyas entradas sean la direcci\'on del
viento en grados y su velocidad. El programa debe devolver el \'angulo de deriva
$\theta$ y la velocidad en tierra del aeroplano, $gs$.
\item Dado un tri\'angulo cualquiera, existe una circunferencia exterior a \'el
que pasa por sus v\'ertices (se la denomina \emph{circunferencia circunscrita}\index{Circunferencia!circunscrita} y
a veces con el anglicismo \emph{circunc\'irculo}).
\begin{enumerate}
\item Justificar por qu\'e el punto que se encuentra sobre la intersecci\'on de las
mediatrices de las aristas (llamado \emph{circuncentro}\index{Circuncentro}) es equidistante de los v\'ertices.
\item Utilizar este hecho para construir en GeoGebra la circunferencia circunscrita
a un tri\'angulo cualquiera.
\end{enumerate}
\item An\'alogamente, en un tri\'angulo cualquiera existe una circunferencia
interior a \'el que es tangente a sus tres aristas, denominada \emph{circunferencia
inscrita}\index{Circunferencia!inscrita} o \emph{inc\'irculo}.
\begin{enumerate}
\item Justificar por qu\'e el punto obtenido considerando la intersecci\'on
de las bisectrices de los \'angulos (llamado \emph{incentro}\index{Incentro}) equidista de las aristas.
\item Utilizar este hecho para construir en GeoGebra la circunferencia inscrita
a un tri\'angulo cualquiera.
\item Comprobar, usando los comandos de GeoGebra, que el radio de la circunferencia
inscrita satisface
$$
r=\frac{2a}{p}\,,
$$
donde $a$ es el \'area del tri\'angulo y $p$ su per\'imetro.?`Qu\'e relaci\'on se
satisfar\'a en un tri\'angulo equil\'atero?
\end{enumerate} 

\end{enumerate}

\chapter{Geometr\'ia anal\'itica}\label{cap8}
\vspace{2.5cm}
\section{Nociones b\'asicas}

La geometr\'ia anal\'itica se basa en el uso de \emph{coordenadas} para describir las figuras geom\'etricas y sus
propiedades. Su invenci\'on se remonta al S. XVII y se debe al genio de R. Descartes\index{Descartes, Ren\'e} (Cartesius en su forma latina,
de donde deriva el adjetivo \emph{cartesiano})\index{Cartesiano} y G. W. Leibniz\index{Leibniz, Gottfried Wilhelm}. Aqu\'i nos centraremos exclusivamente en el caso bidimiensional,
aunque el estudio del caso $n-$dimensional es completamente an\'alogo.

Los ejes cartesianos\index{Cartesiano!ejes} en el plano consisten en dos rectas perpendiculares (es decir, formando cuatro \'angulos rectos)
sobre las que se toma un intervalo arbitrario al que se asigna longitud $1$. Esto permite, mediante el teorema de
Tales, asignar divisiones sobre cada una de las rectas que correspondan a los n\'umeros racionales y, por aproximaci\'on,
a cualquier n\'umero real. El plano se puede describir entonces como el producto cartesiano\index{Cartesiano!producto} $\mathbb{R}\times\mathbb{R}$, \'o $\mathbb{R}^2$\index{R2@$\mathbb{R}^2$},
y sus puntos se identifican con los pares \emph{ordenados} de n\'umeros reales. El par $(x,y)$ que le corresponde al punto $P$
del plano, constituye sus \emph{coordenadas cartesianas}\index{Cartesiano!coordenadas}. Al n\'umero $x$ se le llama la \emph{abcisa}\index{Abcisa} de $P$, y a $y$ su
\emph{ordenada}\index{Ordenada}.
\begin{center}
\begin{tikzpicture}
\draw[very thin,color=gray!20,step=0.5] (-1.15,-1.15) grid (1.15,1.15);
\draw[-] (-1.5,0) -- (1.5,0) node[right] {$x$};
\draw[-] (0,-1.5) -- (0,1.5) node[above] {$y$};
\draw[very thick,color=red] (0.5,0) node[below]{$1$} -- (0.5,1) node[above=6pt,right] {$P(1,2)$};
\draw[very thick,color=red] (0,1) node[left]{$2$} -- (0.5,1);
\draw[fill=red,color=red] (0.5,1) circle (0.05);
\end{tikzpicture}
\end{center}

Una \emph{curva algebraica}\index{Curva!algebraica} en el plano es el lugar geom\'etrico de los puntos que satisfacen una ecuaci\'on entre sus
coordenadas de la forma 
$$
F(x,y)=0\,,
$$
donde $F$ es una funci\'on polin\'omica\index{Funci\'on!polin\'omica} de sus argumentos (es decir, un polinomio en $x$ y en $y$). El grado del polinomio
$F$ se denomina el \emph{grado de la curva}.
Cuando $F$ tiene grado $1$ (es decir, es \emph{lineal} en $x$ e $y$) se dice que la curva es una \emph{recta}\index{Recta}, como por ejemplo $F(x,y)=2x-2y+1=0$,
que escribir\'iamos en la forma habitual (despejando la ordenada)
$$
y=x+\frac{1}{2}\,.
$$
\begin{center}
\begin{tikzpicture}
\draw[very thin,color=gray!20,step=0.5] (-2.15,-2.15) grid (2.15,2.15);
\draw[-] (-2.15,0) -- (2.15,0) node[right] {$x$};
\draw[-] (0,-2.15) -- (0,2.15) node[above] {$y$};
\draw[very thick,color=red] (-1.5,-2) -- (2,1.5) node[left=5pt]{$y=x-\frac{1}{2}$};
\end{tikzpicture}
\end{center}

La ecuaci\'on m\'as general de una recta es $F(x,y)=ax+by+c=0$, pero a veces (?`cu\'ando?) puede ponerse en la forma
$y=mx+n$, que se denomina la \emph{ecuaci\'on expl\'icita}\index{Recta!ecuaci\'on expl\'icita}. En esta forma, $m$ tiene la interpretaci\'on de la
\emph{pendiente}\index{Recta!pendiente} de la recta, en tanto que $n$ es la \emph{ordenada en el origen}, es decir, la altura a la que la recta
corta al eje de ordenadas.

Cuando se tienen dos rectas en forma expl\'icita, digamos $y=mx+n$, $y=ax+b$, es muy sencillo determinar su \emph{posici\'on
relativa}:
\begin{list}{$\blacktriangleright$}{}
\item Las rectas son paralelas si $m=a$.
\item Las rectas son perpendiculares si $m=-1/a$.
\item Las rectas se cortan obl\'icuamente en otro caso.
\end{list}
A veces resulta conveniente determinar la recta que pasa por dos puntos dados, $(x_0,y_0)$, $(x_1,y_1)$. Un razonamiento
elemental usando el teorema de Tales\index{Tales!teorema} (que se deja al cuidado del lector como ejercicio) muestra que su ecuaci\'on es\index{Recta!por dos puntos}
$$
\frac{x-x_0}{x_1-x_0}=\frac{y-y_0}{y_1-y_0}\,.
$$

Si $F$ es un polinomio cuadr\'atico, es decir, una expresi\'on del tipo
$F(x,y)=ax^2+bxy+cy^2+dx+ey+f$ (con $a,b,c,d,e,f$ n\'umeros reales), se dice que la correspondiente curva es una \emph{c\'onica}\index{Curva!c\'onica}. A continuaci\'on repasamos las c\'onicas
m\'as importantes (por simplicidad, estudiaremos \'unicamente c\'onicas \emph{no giradas}, para
las cuales el t\'ermino cruzado $xy$ est\'a ausente, $b=0$).
\begin{enumerate}[(a)]
\item \emph{Elipses}\index{Elipse}\\
Son las c\'onicas dadas por una ecuaci\'on de la forma
$$
F(x,y)=ax^2+cy^2+dx+ey+f=0\, ,
$$
donde $a,c>0$: 
\begin{center}
\begin{tikzpicture}
\draw[very thin,color=gray!20,step=0.5] (-2.15,-2.15) grid (2.15,2.15);
\draw[-] (-2.15,0) -- (2.15,0) node[right] {$x$};
\draw[-] (0,-2.15) -- (0,2.15) node[above] {$y$};
\draw[color=red,thick]  (-1,0.5) ellipse (0.5 and 1.5);
\draw[color=red,thick]  (1,1) circle (0.5);
\end{tikzpicture}
\end{center}
que incluyen como caso particular a las circunferencias\index{Circunferencia}, cuando $a=c$. Una circunferencia se dice que est\'a en
\emph{forma can\'onica} cuando se expresa como suma de cuadrados:
$$
(x-x_0)^2+(y-y_0)^2=r^2\, ,
$$
para unos ciertos n\'umeros reales $x_0$,$y_0$, $r$. En ese caso, las coordenadas $(x_0,y_0)$ representan el centro
de la circunferencia, mientras que $r$ da su radio. De la misma forma, una elipse se encuentra en forma can\'onica si se escribe como
$$
\frac{(x-x_0)^2}{a^2}+\frac{(y-y_0)^2}{b^2}=1\,.
$$
El paso de la ecuaci\'on en forma polinomial $ax^2+cy^2+dx+ey+f=0$ a la forma
can\'onica se puede hacer f\'acilmente completando cuadrados, como se explic\'o
en el cap\'itulo \ref{cap1}. En las actividades de c\'omputo se pide implementar
este proceso en Maxima.\\
En el caso de tener la elipse en forma can\'onica,
la \emph{excentricidad}\index{Elipse!excentricidad} $e$ de la elipse mide cu\'anto se aleja
\'esta de ser una
circunferencia. Claramente, esta propiedad est\'a relacionada con los coeficientes
$a$ y $b$, pues la circunferencia aparece precisamente cuando $a=b=r$. Al mayor de los n\'umeros
$a,b$ se le denomina \emph{semieje mayor}, y al menor \emph{semieje menor}. La excentricidad se define
por
$$
e=\sqrt{1-\left(\frac{b}{a}\right)^2}\,.
$$
El punto $(x_0,y_0)$ es el centro\index{Elipse!centro} de la elipse, y los puntos
$F_1=(x_0-ea,y_0)$, $F_2=(x_0+ea,y_0)$ son los \emph{focos}\index{Elipse!focos}.
La figura muestra estos elementos.
\begin{center}
\begin{tikzpicture}
\draw[very thin,color=gray!20,step=0.5] (-3.15,-2.15) grid (3.15,2.15);
\draw[-] (-3.15,0) -- (3.15,0) node[right] {$x$};
\draw[-] (0,-2.15) -- (0,2.15) node[above] {$y$};
\draw[color=red,thick]  (-0.5,0.5) ellipse (2.5 and 1);
\draw[color=red,fill=red] (-0.5,0.5) circle (0.05)node[below]{$(x_0,y_0)$};
\draw[color=blue] (-3,0.5) --node[midway,above]{$a$} (-0.5,0.5);
\draw[color=blue] (-0.5,0.5) --node[midway,left]{$b$} (-0.5,1.5);
\draw[color=red,fill=red] (1.79,0.5) circle (0.05)node[below]{$F_2$};
\draw[color=red,fill=red] (-2.79,0.5) circle (0.05)node[below]{$F_1$};
\end{tikzpicture}
\end{center}
N\'otese que $0\leq e\leq 1$; el caso $e=0$ corresponde a una circunferencia
($a=b$) y el $e=1$ a una elipse degenerada a un segmento de recta ($b=0$).\\
Los puntos sobre la elipse tienen la propiedad de que la suma de sus distancias
a los focos es una constante, que obviamente debe ser igual al doble del semieje
mayor $2a$ (para verlo basta con tomar el punto $P$ exactamente en el extremo del semieje):
\begin{center}
\begin{tikzpicture}
\draw[color=red,thick]  (-0.5,0.5) ellipse (2.5 and 1);
\draw[color=red,fill=red] (1.79,0.5) circle (0.05)node[below]{$F_2$};
\draw[color=red,fill=red] (-2.79,0.5) circle (0.05)node[below]{$F_1$};
\draw[color=blue] (-2.79,0.5) --node[midway,above]{$F_1P$} (-1,1.48)node[above]{$P$} 
--node[midway,below]{$PF_2$} (1.79,0.5);
\draw[color=orange] (-2.79,0.5) --node[midway,below]{$F_1Q$} (1,1.3)node[above]{$Q$} 
--node[midway,right]{$QF_2$} (1.79,0.5);
\end{tikzpicture}
\end{center}
$$
F_1P+PF_2=2a=F_1Q+QF_2\,.
$$
Como consecuencia de esta f\'ormula y del teorema de Pit\'agoras\index{Pit\'agoras!teorema},
tomando el punto $P$ como el extremo del semieje \emph{menor}
resulta que, si $F$ denota la distancia de uno de los focos al centro,
\begin{equation}\label{formelip}
F^2+b^2=a^2\,.
\end{equation}
Desde un punto de vista f\'isico, si una superficie tiene una secci\'on el\'iptica
los rayos de luz que emanan desde uno de los focos, tras rebotar en la superficie
pasan por el otro foco. Lo mismo ocurre con las ondas de sonido. El National Statuary 
Hall\index{National Statuary Hall},
una sala del Capitolio de los Estados Unidos, tiene una c\'upula de secci\'on el\'iptica
que presenta este efecto\footnote{Una leyenda urbana cuenta que el congresista Quincy 
Adams\index{Adams, Quincy} fing\'ia quedarse dormido sobre una mesa colocada sobre uno de los focos para escuchar
lo que comentaban los integrantes del partido contrario en el otro extremo de la sala,
aunque esto es imposible porque la c\'upula el\'iptica se a\~nadi\'o al edificio cuando Quincy ya hab\'ia fallecido.}, lo mismo que la Whispering Gallery\index{Whispering Gallery} en la Grand Central Terminal
de Nueva York. Hay muchos ejemplos m\'as en distintos lugares del mundo, como la Sala de los Secretos\index{Sala de los Secretos} en la Alhambra de Granada o la Bas\'ilica de San 
Pedro en Ciudad del Vaticano.\\
Como una aplicaci\'on m\'as de esta propiedad de la elipse, tenemos los \emph{billares
el\'ipticos}\index{Billar}. En ellos, pueden darse diferentes situaciones:
\begin{list}{$\blacktriangleright$}{}
\item Si el primer lanzamiento es tal que la bola pasa sobre uno de los focos, rebotar\'a
e ir\'a al otro, repitiendo el proceso indefinidamente (mientras la bola tenga energ\'ia).
Se observa que la trayectoria tiende a identificarse con el semieje mayor.
\begin{center}
\includegraphics[scale=0.4]{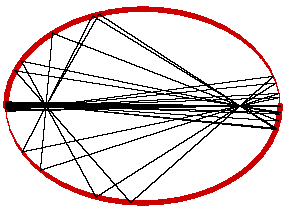} 
\end{center}
\item Si el lanzamiento no pasa entre los focos, las trayectorias determinan tangencialmente una elipse interior a la del billar. Como casos excepcionales, pueden presentarse
trayectorias peri\'odicas poligonales.
\item Cuando el primer lanzamiento pasa entre los focos, la trayectoria (tras
numerosos rebotes) define tangencialmente otra c\'onica, una hip\'erbola, cuyas propiedades veremos a continuaci\'on.
\begin{center}
\includegraphics[scale=0.4]{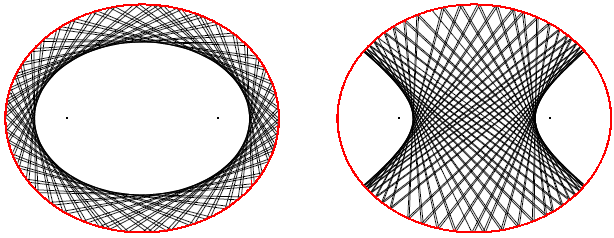} 
\end{center}
\end{list}

\item \emph{Hip\'erbolas}\index{Hip\'erbola}\\
Tienen por ecuaci\'on
$$
F(x,y)=ax^2-cy^2+dx+ey+f=0 \,.
$$
La forma normal de una hip\'erbola es semejante a la de la elipse (y tambi\'en se hace
completando cuadrados), salvo un signo,
y se presenta en dos variantes:
\begin{align}
\frac{(x-x_0)^2}{a^2}-\frac{(y-y_0)^2}{b^2}&= 1 \label{hiphor}\\
\frac{(y-y_0)^2}{b^2}-\frac{(x-x_0)^2}{a^2}&= 1\label{hipver}\,.
\end{align}
El punto $(x_0,y_0)$ es el \emph{centro}\index{Hip\'erbola!centro} de la hip\'erbola en ambos casos.
En el primer caso \eqref{hiphor}, la
hip\'erbola se dice \emph{horizontal}; la gr\'afica de
la curva es sim\'etrica respecto de las rectas $y=y_0$, que corta a dicha
gr\'afica en los puntos $x=-a$ y $x=a$, y la $x=x_0$, que no corta a la gr\'afica. 
En el segundo caso \eqref{hipver}, la hip\'erbola se dice \emph{vertical}, y ahora es la 
recta $y=y_0$ la que no corta a la gr\'afica, mientras que la recta
$x=x_0$ lo hace en los puntos $y=-b$, $y=b$.\\
La siguiente figura muestra una hip\'erbola horizontal centrada en el origen:
\begin{center}
\begin{tikzpicture}
\draw[very thin,color=gray!20,step=0.5] (-2.15,-2.15) grid (2.15,2.15);
\draw[-] (-2.15,0) -- (2.15,0) node[right] {$x$};
\draw[-] (0,-2.15) -- (0,2.15) node[above] {$y$};
\draw[color=red,thick] plot[variable=\t,samples=100,domain=-60:60] ({sec(\t)},{tan(\t)}) ;
\draw[color=red,thick] plot[variable=\t,samples=100,domain=-60:60] ({-sec(\t)},{tan(\t)}) ;
\draw[color=red,fill=red] (-1,0) circle (0.05)node[below right]{$-a$};
\draw[color=red,fill=red] (1,0) circle (0.05)node[below left]{$a$};
\end{tikzpicture}
\end{center}
La hip\'erbola tambi\'en tiene focos y excentricidad. Los focos de una hip\'erbola
horizontal\index{Hip\'erbola!focos} se encuentran en los puntos $F_1=(x_0-|c|,y_0)$, $F_2=(x_0+|c|,y_0)$, donde $c^2=a^2+b^2$ (para el caso de una hip\'erbola vertical, en $(x_0,y_0\pm |c|)$),
y la excentricidad\index{Hip\'erbola!excentricidad} se define
como $e=c/a$. Notemos que en una hip\'erbola siempre es $c>a$, luego $e>1$ (a diferencia de
lo que ocurre en la elipse). Cuanto mayor es $e$, m\'as `recta' es la hip\'erbola (en el 
l\'imite $e\to\infty$, la hip\'erbola degenera en dos rectas paralelas).\\
La propiedad geom\'etrica que tienen los puntos de la hip\'erbola es que la 
\emph{diferencia} entre sus distancias a los focos, es una constante ($2a$ para una
hip\'erbola horizontal y $2b$ para una hip\'erbola vertical). La figura muestra el caso
de una hip\'erbola vertical.
\begin{center}
\begin{tikzpicture}
\draw[-] (-3,0) -- (3,0) node[right] {$x$};
\draw[-] (0,-3) -- (0,3) node[above] {$y$};
\draw[color=red,thick] plot[variable=\t,samples=100,domain=-70:70] ({tan(\t)},{sec(\t)}) ;
\draw[color=red,thick] plot[variable=\t,samples=100,domain=-70:70] ({tan(\t)},{-sec(\t)}) ;
\draw[color=red,fill=red] (0,-1) circle (0.05)node[below right]{$-b$};
\draw[color=red,fill=red] (0,1) circle (0.05)node[above right]{$b$};
\draw[color=red,fill=red] (0,-2.41) circle (0.05)node[below right]{$F_1$};
\draw[color=red,fill=red] (0,2.41) circle (0.05)node[above right]{$F_2$};
\draw[color=blue] (0,-2.41) --node[midway,below left]{$F_1P$}
(-1,-1.41)node[above left]{$P$} --node[midway,left]{$PF_2$} (0,2.41);
\draw[color=orange] (0,-2.41) --node[midway,below]{$F_1Q$}
(2.19,-2.41)node[above right]{$Q$} --node[midway,above right]{$QF_2$} (0,2.41);
\end{tikzpicture}
\end{center} 
$$
PF_2-F_1P=2b=QF_2-F_1Q\,.
$$

\item \emph{Par\'abolas}\index{Par\'abola}\\
Otro ejemplo de curva de segundo grado que es una c\'onica lo da la par\'abola, cuya ecuaci\'on general es
$$
F(x,y)=ax^2+dx+ey=0.
$$
De acuerdo con lo visto en el cap\'itulo \ref{cap1}, una par\'abola siempre puede
escribirse (por el m\'etodo de completar cuadrados) en la forma normal 
$$
y=A(x-h)^2-k\,,
$$
(en la literatura anglosajona tambi\'en se denomina \emph{forma v\'ertice} o \emph{vertex
form}), en la que resulta inmediato localizar el v\'ertice de la par\'abola, situado en
el punto $(h,-k)$. El par\'ametro $A$ determina si la par\'abola se abre `hacia arriba'
(cuando $A>0$) o `hacia abajo' (cuando $A<0$), y su valor absoluto $|A|$ determina c\'omo
es esta apertura (mayor cuanto m\'as grande es $|A|$). 
La par\'abola es sim\'etrica respecto de la recta vertical $x=h$.\\
Tambi\'en existe un \emph{foco}\index{Par\'abola!foco} $F$ para las par\'abolas: se encuentra situado sobre el eje
de simetr\'ia $x=h$ a una distancia $1/4A$ del v\'ertice. Tiene la propiedad de que los 
rayos de luz incidentes sobre la par\'abola \emph{paralelos al eje de simetr\'ia},
rebotan y pasan por el foco. Despu\'es, rebotan de nuevo sobre la par\'abola y salen
paralelos al eje de simetr\'ia.
\begin{center}
\begin{tikzpicture}
\draw[color=red,thick] plot[variable=\t,samples=100,domain=-2:4] ({\t},{0.5*(\t -1)^2-2});
\draw[color=red,fill=red] (1,-2) circle (0.05)node[below]{$(h,-k)$};
\draw[color=red,fill=red] (1,-1) circle (0.05)node[below right]{$F$};
\draw[color=orange,thick] (1,-2) -- (1,2.5);
\draw[color=orange,thick] (-1,2.5) -- (-1,0) -- (2,-1.5) -- (2,2.5);
\draw[color=orange,thick] (3,2.5) -- (3,0) -- (0,-1.5) -- (0,2.5);
\draw[color=orange,thick] (0.5,2.5) -- (0.5,-1.875) -- (3,2.5);
\draw[color=orange,thick] (1.5,2.5) -- (1.5,-1.875) -- (-1,2.5);
\end{tikzpicture}
\end{center}
Esta propiedad es el fundamento de las \emph{antenas parab\'olicas}, cuya secci\'on transversal tiene forma de par\'abola. El receptor de se\~nal se coloca 
precisamente en el foco, de manera que las ondas recibidas por la antena se concentran
en \'el. De este tipo son tanto las antenas dom\'esticas para la recepci\'on de se\~nales
por sat\'elite como antenas sofisticadas de uso cient\'ifico y militar.
  \begin{center}
  \includegraphics[scale=0.85]{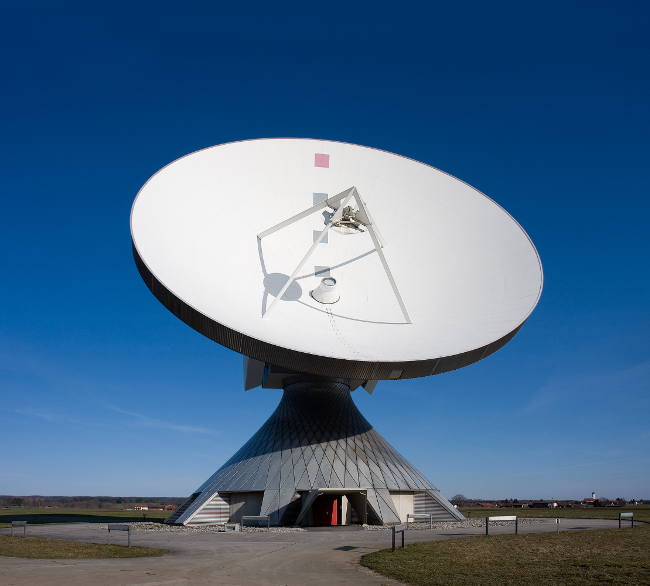} 
  \end{center}
El fen\'omeno tambi\'en se aprovecha a la inversa: los faros de los autom\'oviles
tienen una secci\'on parab\'olica y el foco de iluminaci\'on en el foco geom\'etrico, de 
modo que los rayos de luz salen paralelos y no se dispersan.\\
Otro elemento geom\'etrico fundamental en una par\'abola es su \emph{directriz}. Se
trata de una recta $d$ perpendicular al eje de simetr\'ia, 
a una distancia del v\'ertice igual
a la que tiene el foco, pero en sentido opuesto. Los puntos de la par\'abola cumplen que
su distancia al foco es igual a su distancia a la directriz.
\begin{center}
\begin{tikzpicture}
\draw[color=red,thick] plot[variable=\t,samples=100,domain=-1.25:3.25] ({\t},{0.5*(\t -1)^2-2});
\draw[color=red,fill=red] (1,-1.5) circle (0.05)node[above]{$F$};
\draw[color=blue] (-1.5,-2.5) --node[midway,below left]{$d$} (3.5,-2.5);
\draw[color=orange] (0,-1.5)node[below left]{$P$} -- (1,-1.5);
\draw[color=orange] (0,-1.5) -- (0,-2.5)node[below]{$P'$};
\draw[color=orange] (1.5,-1.875)node[below right]{$Q$} -- (1,-1.5);
\draw[color=orange] (1.5,-1.875) -- (1.5,-2.5)node[below]{$Q'$};
\draw[color=orange] (3,0)node[above left]{$R$} -- (1,-1.5);
\draw[color=orange] (3,0) -- (3,-2.5)node[below]{$R'$};
\draw[color=orange,fill=orange] (0,-1.5) circle (0.05);
\draw[color=orange,fill=orange] (1.5,-1.875) circle (0.05);
\draw[color=orange,fill=orange] (3,0) circle (0.05);
\end{tikzpicture}
\end{center}
\noindent Es decir, en la figura precedente se cumplen las relaciones:
$$
PP'= PF \qquad ,\qquad QQ'= QF\qquad ,\qquad RR'= RF\,.
$$
\end{enumerate}

Finalizamos esta secci\'on con unas breves nociones sobre la degeneraci\'on de curvas. 
Una curva algebraica se dice que es degenerada\index{Curva!degenerada} cuando el polinomio que la define \emph{no es irreducible}\index{Polinomio!irreducible}.
Esto ocurre bien sea porque la curva se reduce a un punto (o al vac\'io), como en el caso de
$F(x,y)=x^2+y^2=0$ (cuyo lugar geom\'etrico consta \'unicamente de $\{(0,0)\}$), 
o bien porque se puede factorizar, es decir,
existen unos polinomios $f$ y $g$ tales que $F(x,y)=f(x,y)g(x,y)$. En tal caso,
la curva tiene dos \emph{componentes}\index{Curva!componentes}, dadas por los lugares geom\'etricos $f(x,y)=0$ 
y $g(x,y)=0$, respectivamente.\\
Por ejemplo, una c\'onica degenerada es
$$
F(x,y)=x^2-y^2=0\,.
$$
Notemos que $x^2-y^2=(x+y)(x-y)$, de modo que la curva est\'a compuesta de 
dos l\'ineas rectas que se intersectan en un punto.
\begin{center}
\begin{tikzpicture}
\draw[very thin,color=gray!20,step=0.5] (-2.15,-2.15) grid (2.15,2.15);
\draw[-] (-2.15,0) -- (2.15,0) node[right] {$x$};
\draw[-] (0,-2.15) -- (0,2.15) node[above] {$y$};
\draw[very thick,color=red] (-2,-2) -- (2,2) ;
\draw[very thick,color=red] (-2,2) -- (2,-2) ;
\end{tikzpicture}
\end{center}

\section{Aplicaciones}
Vamos a ver una situaci\'on cotidiana, alejada de cuestiones tecnol\'ogicas o
cient\'ificas avanzadas, que involucra el uso de la geometr\'ia anal\'itica. En el
mercado pueden encontrarse chimeneas\index{Chimenea} dom\'esticas que son muy \'utiles en zonas de
clima fr\'io.
\begin{center}
\includegraphics[scale=0.25]{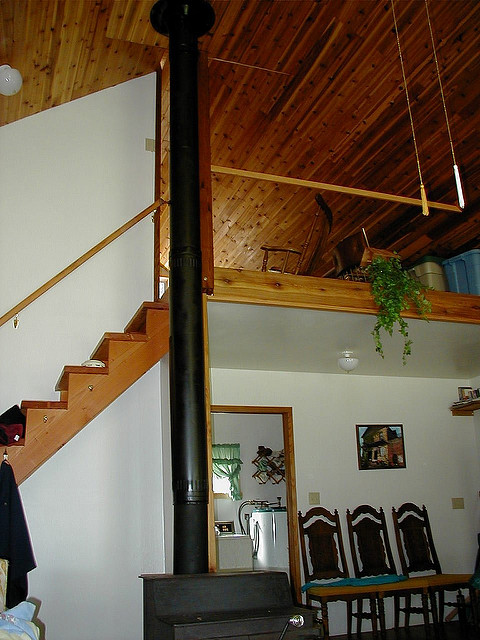} 
\end{center}
Como se aprecia en la figura, es necesario realizar una perforaci\'on en el techo para
hacer pasar el tubo que conduce el humo al exterior. Dado que la secci\'on del 
tubo de la chimenea es circular, habitualmente se procede realizando un corte con esa
forma.
\begin{center}
\includegraphics[scale=0.4]{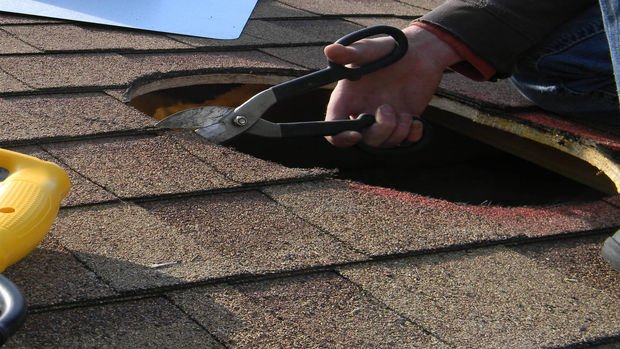} 
\end{center}
Sin embargo, la secci\'on de un cilindro por un plano inclinado \emph{no} es una
circunferencia, sino una elipse\index{Elipse}. Es f\'acil comprobar esto experimentalmente
llenando un vaso cil\'indrico con agua e inclin\'andolo. La principal consecuencia
de esto es que al pasar el tubo por el agujero, quedan espacios a los lados que es
necesario tapar para evitar p\'erdidas de calor, lo que conlleva trabajo y costos
adicionales, adem\'as de que el resultado final es poco est\'etico.\\
Se plantea entonces el problema de determinar la forma (el\'iptica) del corte que
hay que practicar en un techo que tiene una inclinaci\'on de \'angulo $\alpha$
para que pase un tubo que supondremos de secci\'on circular con radio $R$.

\begin{center}
\begin{tikzpicture}[MyPersp,font=\large,scale=0.63]
	\def\h{2.5}
	\def\a{35}
	\def\aa{35}
	\pgfmathparse{\h/tan(\a)}
  	\let\b\pgfmathresult

	\coordinate (A) at (2,\b,0);
	\coordinate (B) at (-2,\b,0);
	\coordinate (C) at (-2,-1.5,{(1.5+\b)*tan(\a)});
	\coordinate (D) at (2,-1.5,{(1.5+\b)*tan(\a)});
	\coordinate (E) at (2,-1.5,0);
	\coordinate (F) at (-2,-1.5,0);
	\draw[color=gray](-2,-0.575,2.95) -- (C);
	\draw[color=gray] (E) -- (-0.1,-1.5,0);	
	\draw[color=gray] (-0.1,-1.5,{(1.5+\b)*tan(\a)}) -- (C);
	\draw[color=gray] (-0.1,-1.5,0) -- (F);
	\draw[color=gray] (F)--(C);	
	\draw[color=gray] (B)--(-2,-0.575,0);  
	\draw[color=gray] (-2,-0.575,0)-- (F);
	\draw[fill,color=gray!25,opacity=0.45] (A)--(B)--(C)--(D)--cycle;
	\draw[]  (A) --(B);
	\draw[] (B) -- (-2,-0.575,2.95);		
	\coordinate[](centro) at (0,0,-0.55);

	\node[draw,
cylinder,rotate=90, minimum height=11cm,minimum width=4.1cm,scale=0.63, 
cylinder uses custom fill, cylinder end fill=blue!25, cylinder body fill=blue!40,
opacity=0.5,
shape aspect=8, right=0cm of centro](){};	

	\draw (A) --(D);
	\draw (D) -- (-0.1,-1.5,{(1.5+\b)*tan(\a)});
	\draw (A)--(E) -- (D); 

	\foreach \t in {20,40,...,360}
		\draw[gray,dashed] ({cos(\t)},{sin(\t)},0)
      --({cos(\t)},{sin(\t)},{2.0*\h});
	\draw[blue,very thick,opacity=0.5] (1,0,0) 
		\foreach \t in {5,10,...,360}
			{--({cos(\t)},{sin(\t)},0)}--cycle;
	\draw[blue,very thick] (1,0,{2*\h}) 
		\foreach \t in {10,20,...,360}
			{--({cos(\t)},{sin(\t)},{2*\h})}--cycle;
	\fill[blue!15,draw=blue,very thick,opacity=0.5]
     (0,1,{\h-tan(\a)})
		\foreach \t in {5,10,...,360}
			{--({sin(\t)},{cos(\t)},{-tan(\a)*cos(\t)+\h})}--cycle;
\end{tikzpicture}
\end{center}
La vista lateral de esta figura nos ofrece las relaciones geom\'etricas que necesitamos
para resolver el problema. Tomaremos un sistema de coordenas cartesianas con origen en el 
centro de la elipse y sus ejes coordenados dirigidos seg\'un los semiejes de la elipse.
\begin{center}
    \begin{tikzpicture}[scale=0.75]
	   \draw[-] (4,-0.027) -- (4,2) node[midway, right] {$2R\sin\alpha$};
	   \draw[thick,color=blue] (0,-1) -- (0,3);
	   \draw[thick,color=blue] (4,-1) -- (4,3);	   
	   \draw[-] (0,-0.025) -- (4,-0.025) node[midway, below] {$2R$};
	   \draw[-,very thick] (-2,-1) --node[above]{$L$} (6,3) ;
	   \draw[<-] (0.05,0.25) -- (1.75,1.1);
	   \draw[->] (2.25,1.35) -- (3.95,2.2);
	   \draw[thick] (0.9,-0.05) arc (0:30:0.9);
	   \node[] at (12:1.25) {$\alpha$};
	   \draw[color=gray!50] (10,-2.5) -- (10,4.5)node[black,above right]{$y$};
	   \draw[color=gray!50] (8,1) -- (12,1)node[black,below right]{$x$};	   
	   \draw[color=red,thick]  (10,1) ellipse (1.5 and 3);
	   \draw[color=red,fill=red] (8.5,1) circle (0.06)node[black,above left]{$-R$};
	   \draw[color=red,fill=red] (11.5,1) circle (0.06)node[black,above right]{$R$};
	   \draw[color=red,fill=red] (10,4) circle (0.06)node[black,above left]{$S$};
	   \draw[color=red,fill=red] (10,-2) circle (0.06)node[black,below left]{$-S$};
	   \draw[color=red,fill=red] (10,3.25) circle (0.06)node[black,below right]{$F$};
	   \draw[color=red,fill=red] (10,-1.25) circle (0.06)node[black,above right]{$-F$};	      
	      
	\end{tikzpicture}
\end{center}
El semieje menor claramente mide $R$. Llamando $\pm F$ a los focos y
$S$ al semieje mayor, se cumple $L=2S$. 
Aplicando el teorema de Pit\'agoras al techo, $4S^2 = 4R^2\sin^2\alpha +4R^2$,
de donde despejando
\begin{equation}\label{eqs}
S^2=R^2(1+\sin^2\alpha )\,.
\end{equation}
Por otra parte, la ecuaci\'on de la elipse ser\'a
$\frac{x^2}{R^2}+\frac{y^2}{S^2}=1$,
de manera que, por \eqref{formelip},
\begin{equation}\label{eqs2}
F^2+R^2=S^2\,.
\end{equation}
Pero de \eqref{eqs} y \eqref{eqs2} se obtiene
$F^2=S^2-R^2=R^2(1+\sin^2\alpha -1)=R^2\sin^2\alpha$, luego
\begin{equation}\label{focos}
F=\pm R\sin\alpha\,.
\end{equation}
Adem\'as, de \eqref{eqs} sabemos que
\begin{equation}\label{eql}
L=2S=2R\sqrt{1+\sin^2\alpha}\,.
\end{equation}
Con esto resulta muy f\'acil dibujar la elipse sobre el techo; basta con tomar un hilo
de longitud \eqref{eql}, fijar sus extremos con dos clavos en los valores de $F$ y $-F$
dados por \eqref{focos} y trazar la elipse manteniendo tirante el hilo.
\begin{center}
\includegraphics[scale=0.3]{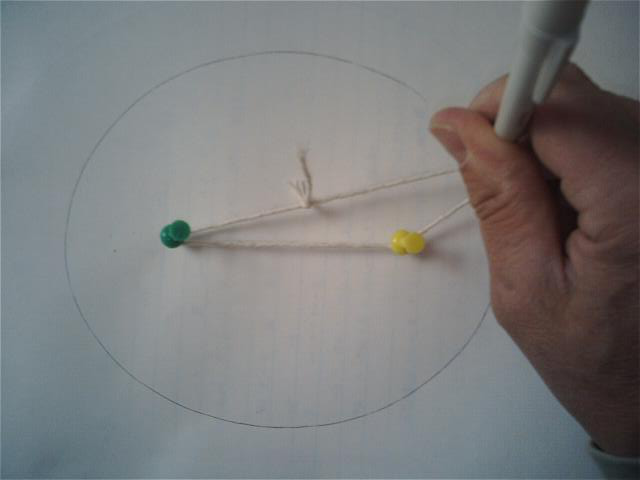} 
\end{center}
\section{Ejercicios}

\begin{enumerate}  
\setcounter{enumi}{102}

\item Dar la ecuaci\'on en la forma $y=mx+n$ para cada una de las rectas siguientes:
\begin{enumerate}[(a)]
\item Recta que pasa por los puntos $(4,-2)$ y $(1,7)$.
\item Recta con pendiente $3$ y ordenada en el origen $4$.
\item Recta que pasa por los puntos $(-1,0)$ y $(0,3)$.
\item Recta que pasa por $(2,-3)$ y es paralela al eje $Ox$.
\item Recta que pasa por $(2,3)$ y aumenta cuatro unidades por cada unidad que crece $x$.
\item Recta que pasa por $(-2,2)$ y disminuye dos unidades por cada unidad que crece $x$.
\item Recta que pasa por $(3,-4)$ y es paralela a la recta $5x-2y=4$.
\item Recta que pasa por el origen y es paralela a la recta $y=2$.
\item Recta que pasa por $(-2,5)$ y es perpendicular a la recta $4x+8y=3$.
\item Recta que pasa por el origen y es perpendicular a la recta $3x-2y=1$.
\item Recta que pasa por $(2,1)$ y es perpendicular a la recta $x=2$.
\item Recta que pasa por el origen y es bisectriz del \'angulo determinado por los semiejes positivos.
\end{enumerate}

\item Hallar las pendientes y las ordenadas en el origen de las rectas siguientes. Tambi\'en,
dar las coordenadas de un punto sobre la recta que sea diferente del $(0,b)$:
\begin{multicols}{3}
\begin{enumerate}[(a)]
\item $ y=3x-2$
\item $ 2x-5y=3$
\item $ y=4x-3$
\item $ y=-3$
\item $ {\displaystyle \frac{y}{2}}+{\displaystyle \frac{x}{3}}=1$
\end{enumerate}
\end{multicols} 

\item
\begin{enumerate}[(a)]
\item Determinar $k$ sabiendo que el punto $(3,k)$ se encuentra sobre la recta de pendiente
$m=-2$ que pasa por el punto $(2,5)$.
\item El punto $(3,-2)$, ?` Est\'a sobre la recta que pasa por los puntos $(8,0)$ y $(-7,-6)$?
\item Determinar si el tri\'angulo de v\'ertices $(7,-1)$, $(10,1)$ y $(6,7)$ es rect\'angulo.
\item Determinar si los puntos $(8,0)$, $(-1,-2)$, $(-2,3)$ y $(7,5)$ son v\'ertices de un paralelogramo.
\item Determinar $k$ de manera que los puntos $A\equiv (7,3)$, $B\equiv (-1,0)$ y $C\equiv (k,-2)$
sean los v\'ertices de un tri\'angulo rect\'angulo con \'angulo recto en $B$.
\end{enumerate}

\item Para cada par de rectas, determinar si son paralelas, perpendiculares o ninguna de ellas:
\begin{multicols}{2}
\begin{enumerate}[(a)]
\item $ y=3x+2,y=3x-2$
\item $ y=2x-4,y=3x+5$
\item $ 3x-2y=5,2x+3y=4$
\item $ 6x+3y=1,4x+2y=3$
\item $ x=3,y=-4$
\item $ 5x+4y=1,4x+5y=2$
\item $ x=-2,x=7$
\end{enumerate}
\end{multicols}

\item Hallar las ecuaciones can\'onicas de las circunferencias que satisfacen las condiciones siguientes:
\begin{enumerate}[(a)]
\item Tiene centro $(4,-1)$ y radio $1$.
\item Tiene centro $(-2,-2)$ y radio $5\sqrt{2}$.
\item Tiene centro $(-2,3)$ y pasa por $(3,-2)$.
\item Tiene centro $(6,1)$ y pasa por el origen.
\end{enumerate}

\item Identificar y representar las gr\'aficas de las ecuaciones siguientes, especificando
sus elementos geom\'etricos (focos, simetr\'ias, excentricidad, etc.):
\begin{multicols}{2}
\begin{enumerate}[(a)]
\item $ x^2+y^2+16x-12y+10=0$
\item $ x^2+y^2-4x+5y+10=0$
\item $ x^2+y^2+x-y=0$
\item $ 4x^2+4y^2+8y-3=0$
\item $ x^2+y^2-x-2y+3=0$
\item $ x^2+y^2+\sqrt{2}x-2=0$
\end{enumerate}
\end{multicols}

\item Representar, en un mismo gr\'afico, las par\'abolas $y=x^2$, $y=-x^2$, $x=y^2$,
$x=-y^2$, identificando los puntos de corte entre ellas.

\item Representar, en un mismo gr\'afico, las par\'abolas $y=x^2+1$, $y=-x^2-1$, $x=y^2+1$,
$x=-y^2-1$, identificando los puntos de corte entre ellas.

\item Representar, en un mismo gr\'afico, las par\'abolas $2y=x^2$, $2y=(x+1)^2$, $2y-2=(x-2)^2$,
identificando los puntos de corte entre ellas.

\item Identificar y representar las gr\'aficas de las ecuaciones siguientes, especificando
sus elementos geom\'etricos (focos, simetr\'ias, excentricidad, etc.):
\begin{multicols}{2}
 
\begin{enumerate}[(a)]
\item $ y^2-x^2 =1$
\item $ 25x^2 +36y^2 =900$
\item $ 2x^2-y^2 =4$
\item $ xy=4$
\item $ 4x^2 + 4y^2 =1$
\item $ 8x=y^2$
\item $ 10y=x^2$
\item $ 4x^2 +9y^2 =16$
\item $xy=-1$
\item $3y^2 -x^2 =9$
\item $4x^2 + 9y^2 - 48x + 72y + 144 = 0$
\item $4x^2 - 5y^2 + 40x - 30y - 45 = 0$
\item $(y - 3)^2 = 8(x- 5)$
\item $(x + 3)^2 = -20(y- 1) $
\item $x^2-4y^2+4=0 $
\end{enumerate}
\end{multicols}

\end{enumerate}

\section{Uso de la tecnolog\'ia de c\'omputo}
Vamos a ver la construcci\'on geom\'etrica (con regla y comp\'as) y las\index{GeoGebra!construcciones con regla y comp\'as}
propiedades de las c\'onicas con GeoGebra.\index{GeoGebra!c\'onicas}
\subsection{La elipse}
Seguiremos el protocolo que se enumera a continuaci\'on:
\begin{enumerate}
\item Seleccionamos la vista de \texttt{Geometr\'ia B\'asica} y con la herramienta \texttt{Nuevo Punto}
marcamos los que ser\'an los focos de la elipse, renombr\'andolos como $F_1$ y
$F_2$.
\item Con la herramienta \texttt{Circunferencia dados su Centro y uno de sus Puntos} constru\'imos
una circunferencia con centro en $F_2$ y que pase por un punto $A$, de tal forma que 
contenga en su \emph{interior} a ambos focos. 
\item Seleccionamos un punto sobre la circunferencia, digamos $B$, y constru\'imos el radio 
$F_2B$.
\item Con la herramienta \texttt{Mediatriz} constru\'imos la mediatriz del segmento $F_1B$,
renombr\'andola como $m$.
\item Con la herramienta \texttt{Intersecci\'on de Dos Objetos} determinamos el punto de 
intersecci\'on de $F_2B$ con $m$, renombr\'andolo como $P$.
\item Picando con el bot\'on derecho del rat\'on en el punto $P$ marcamos la opci\'on
\texttt{Activa Rastro}.
\item Desplazando el punto $B$ por la circunferencia, ?`qu\'e figura traza el punto 
$P$? La respuesta puede justificarse usando resultados concocidos de geometr\'ia elemental 
(puede ser de ayuda trazar los segmentos auxiliares $F_1B$ y $F_1P$). 
\item Resulta interesante repetir el experimento con distintas posiciones de $F_1,\ F_2$ y 
$A$, manteniendo a $F_1$ dentro de la circunferencia.
\end{enumerate}
La siguiente figura ilustra este proceso:
\begin{center}
\includegraphics[scale=0.4]{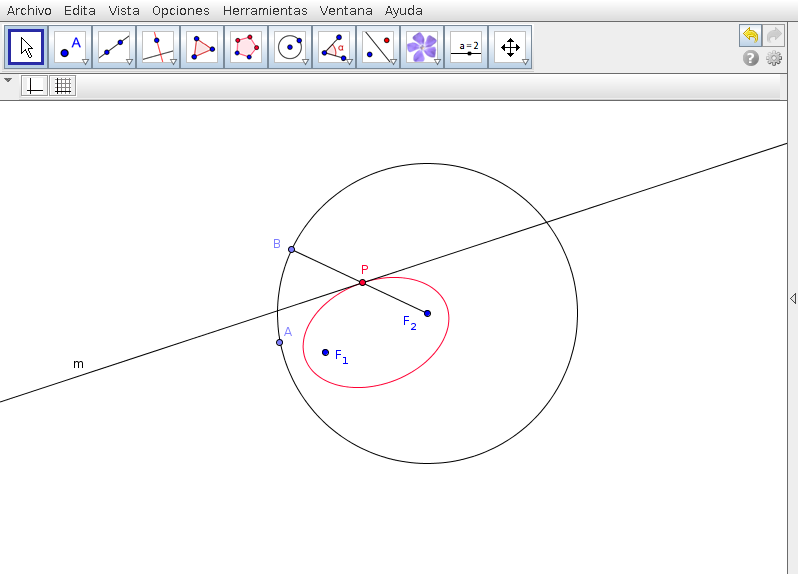}
\end{center}

\subsection{La hip\'erbola}
La hip\'erbola tambi\'en puede construirse con regla y comp\'as siguiendo un procedimiento
similar al de la elipse:
\begin{enumerate}
\item Seleccionamos la vista de \texttt{Geometr\'ia B\'asica} y marcamos los que ser\'an los focos 
de la hip\'erbola, renombr\'andolos como $F_1$ y $F_2$.
\item Constru\'imos una circunferencia con centro en $F_2$ y que pase por un punto $A$, de tal forma que el foco $F_1$ quede en el \emph{exterior} de ella. 
\item Seleccionamos un punto sobre la circunferencia, digamos $B$, y trazamos la recta 
$F_2B$, renombr\'andola como $n$.
\item Constru\'imos la mediatriz del segmento $F_1B$, renombr\'andola como $m$.
\item Determinamos el punto de intersecci\'on de $m$ con $n$, renombr\'andolo como $P$.
\item Activamos el rastro del punto $P$.
\item Desplazando el punto $B$ por la circunferencia, ?`qu\'e figura traza el punto $P$? La 
respuesta puede justificarse usando resultados concocidos de geometr\'ia elemental (puede 
ser de ayuda trazar los segmentos auxiliares $F_1B$ y $F_1P$). 
\item Resulta interesante repetir el experimento con distintas posiciones de $F_1,\ F_2$ y $A$, manteniendo a $F_1$ fuera de la circunferencia.
\end{enumerate}
Todo este proceso en GeoGebra se ve como sigue:
   \begin{center}
\includegraphics[scale=0.4]{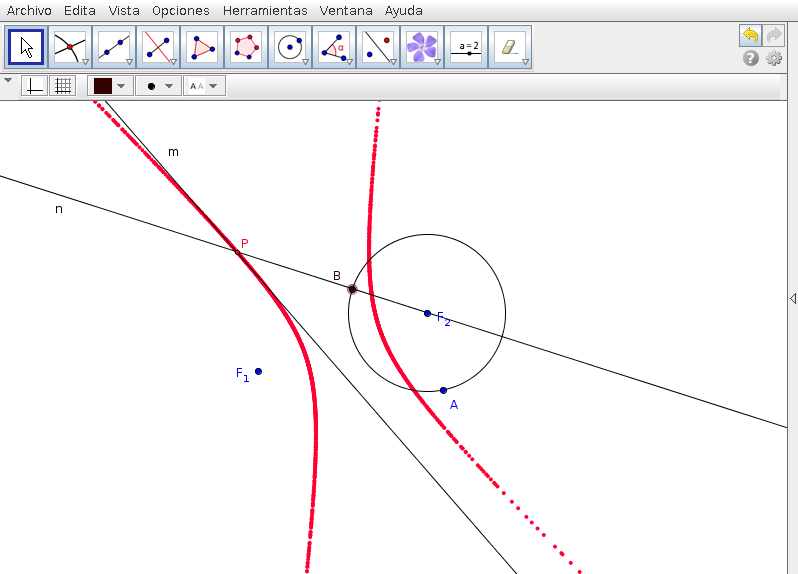}
\end{center}

\subsection{La par\'abola}
Por \'ultimo, veamos el caso de la par\'abola siguiendo las ideas presentadas para
la elipse y la hip\'erbola.
\begin{enumerate}
\item Seleccionamos la vista de \texttt{Geometr\'ia B\'asica}; con la herramienta
\texttt{Recta que pasa por Dos Puntos} construimos una recta renombr\'andola como
$d$. Marcamos un punto fuera de $d$ y lo renombramos como $F$ (estos elementos
ser\'an la directriz y  el foco de la par\'abola, respectivamente). 
\item Seleccionamos un punto sobre la recta, digamos $C$, y construimos la mediatriz del 
segmento $FC$, renombr\'andola como $m$.
\item Trazamos la perpendicular a $d$ que pasa por $C$, renombr\'andola como $n$.
\item Localizamos el punto de intersecci\'on de $m$ con $n$, renombr\'andolo como $P$.
\item Activamos el rastro del punto $P$.
\item Desplazando el punto $C$ a lo largo de la recta, ?`qu\'e figura traza el punto $P$?
De nuevo, la respuesta se puede justificar usando resultados concocidos de geometr\'ia 
elemental (puede ser de ayuda trazar los segmentos auxiliares $FC$ y $FP$). 
\item Una vez m\'as, es interesante repetir el experimento con distintas posiciones de
$F$ y $d$.
\end{enumerate}
El proceso completo en GeoGebra se ve como sigue:
\begin{center}
\includegraphics[scale=0.4]{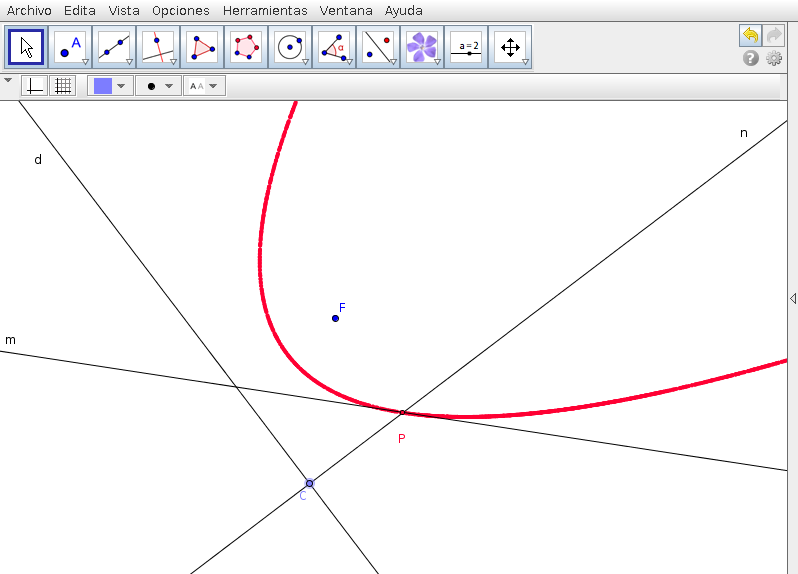}
\end{center}

\section{Actividades de c\'omputo}
\begin{enumerate}[(A)]
\item Construir en Maxima una funci\'on que reduzca una c\'onica dada por la
ecuaci\'on $Ax^2+By^2+Cx+Dy+E=0$ con $A\neq 0\neq B$ (una elipse o una hip\'erbola)
a su forma est\'andar
$$
\frac{(x-x_0)^2}{a}\pm\frac{(y-y_0)^2}{b}=1\,.
$$
La siguiente captura de pantalla muestra la salida de una tal funci\'on implementada por los
autores.
\begin{center}
\includegraphics[scale=0.4]{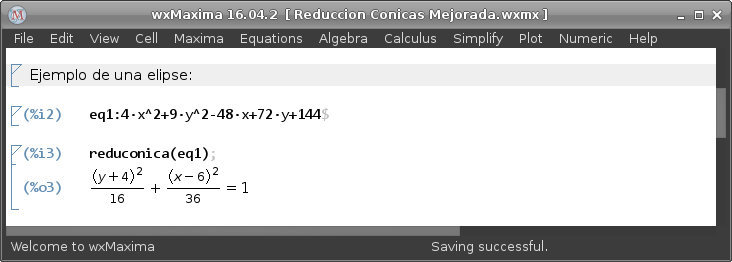} 
\end{center}

\item Construir en Maxima una funci\'on que reduzca una c\'onica dada por la
ecuaci\'on $Ax^2+Bx+C+Dy=0$ (o la $Ay^2+By+C+Dx=0$) a su forma est\'andar.

\item Las \emph{propiedades de reflexi\'on de la elipse} se pueden visualizar con la siguiente construcci\'on:
\begin{enumerate}
\item Crear dos deslizadores $a$ y $b$ que tomen valores positivos, digamos con $1$ como m\'inimo y $5$ como m\'aximo, y ajustarlos con los valores $a=3$ y $b=2$, por
ejemplo.
\item En la ventana de entrada de \'ordenes, escribir la ecuaci\'on est\'andar de la elipse
cuyos ejes son los segmentos  que van de $(-a,0)$ a $(a,0)$ y de $(0,-b)$ a $(0,b)$.
\item Con la \'orden \texttt{Foco} determinar los focos de la elipse, renombr\'andolos
como $F_1$ y $F_2$.
\item Seleccionar un punto sobre la elipse, renombr\'andolo como $P$.
\item Para simular un haz de luz o una onda sonora que emana del foco $F_1$, trazar la recta $F_1P$, renombr\'andola como $l$.
\item Con la \'orden \texttt{Tangente} trazar la recta tangente a la elipse que pasa por el
punto $P$, renombr\'andola como $t$.
\item Con la \'orden \texttt{Perpendicular} trazar la recta normal a la elipse que pasa por $P$,
renombr\'andola como $n$.
\item Con la \'orden \texttt{Refleja} trazar la imagen de $l$ reflejada en $n$.
\end{enumerate}
Al desplazar el punto $P$ a lo largo de la elipse, ?`qu\'e sucede con la reflexi\'on $l'$ de 
$l$? Repitir el experimento con distintos valores de $a$ y de $b$ y explicar el resultado.
   \begin{center}
\includegraphics[scale=0.4]{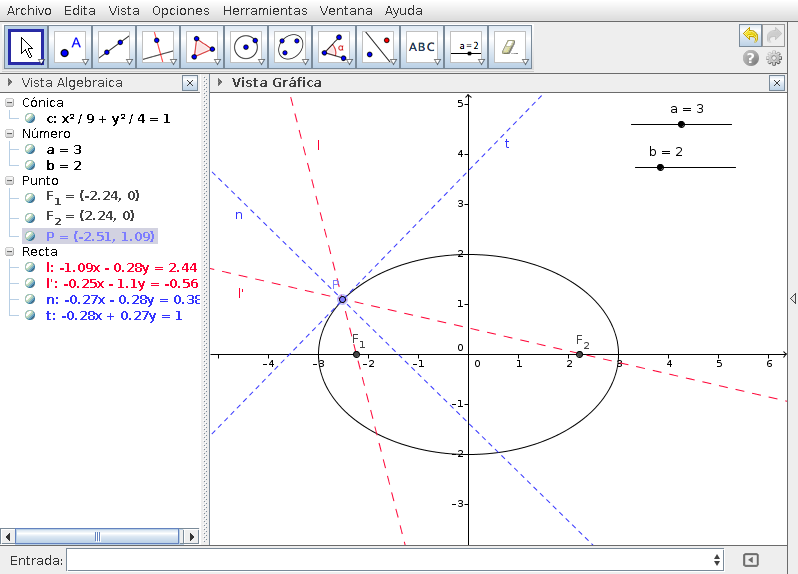}
   \end{center}
   
\item Estudiar las \emph{propiedades de reflexi\'on de la par\'abola} siguiendo
el esquema:
\begin{enumerate}
\item Crear un deslizador $p$ con $-5$ como m\'inimo y $5$ como m\'aximo, y ajustarlo
con $p=2$, por ejemplo.
\item Escribir la ecuaci\'on est\'andar de la par\'abola con v\'ertice en el origen y foco en $(0,p)$.
\item Determinar el foco de la par\'abola, renombr\'andolo como $F$.
\item Seleccionar un punto sobre la par\'abola, renombr\'andolo como $P$.
\item Para simular un haz de luz o una onda sonora que emana del foco $F$,
trazar la recta $FP$, renombr\'andola como $l$.
\item Trazar la recta tangente a la par\'abola que pasa por el punto $P$,
renombr\'andola como $t$.
\item Trazar la recta normal a la par\'abola que pasa por $P$, renombr\'andola como $n$.
\item Trazar la imagen de $l$ reflejada en $n$.
\end{enumerate}
Al desplazar el punto $P$ a lo largo de la par\'abola, ?`que sucede con la reflexi\'on
$l'$ de $l$? Repetir el experimento con distintos valores de $p$ y explicar el resultado.
   \begin{center}
\includegraphics[scale=0.4]{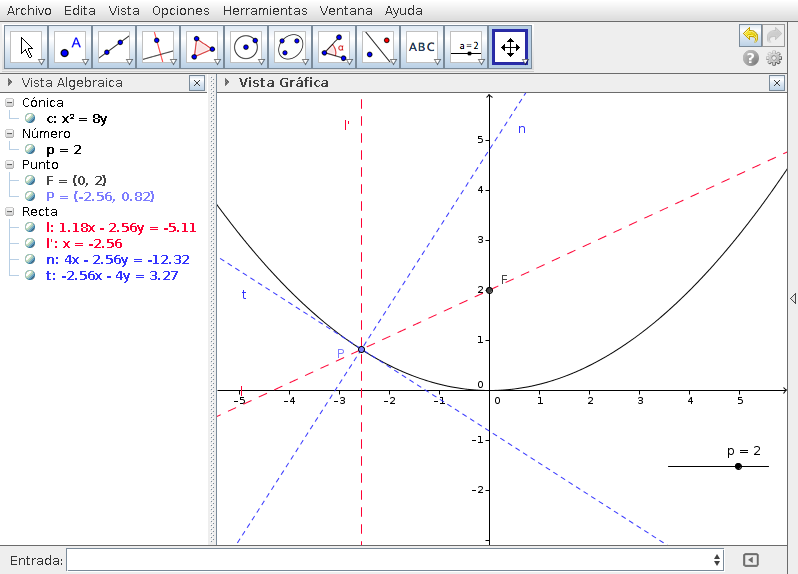}
   \end{center}
   
\item En esta actividad, se propone estudiar las \emph{propiedades de reflexi\'on de la hip\'erbola}. Para ello, procederemos de manera similar al caso de la elipse. 
\begin{enumerate}
\item Crear dos deslizadores $a$ y $b$ que tomen valores positivos, digamos con
$1$ como m\'inimo y $5$ como m\'aximo, y ajustarlos con los valores $a=3$ y $b=2$,
por ejemplo.
\item Escribir la ecuaci\'on est\'andar de la hip\'erbola cuyo rect\'angulo central tiene 
puntos medios $(-a,0)$, $(a,0)$, $(0,-b)$ y $(0,b)$.
\item Determinar los focos de la hip\'erbola, renombr\'andolos como $F_1$ y $F_2$.
\item La hip\'erbola divide al plano en 3 regiones: marcar un punto $O$ en la regi\'on que 
no contiene a ninguno de los focos.
\item Para simular un haz de luz o una onda sonora que emana del punto $O$ en la direcci\'on 
de $F_1$, trazar la recta $F_1P$, renombr\'andola como $l$.
\item Determinar el punto de intersecci\'on de $l$ con la rama de la hip\'erbola que contiene a $F_1$, renombr\'andolo como $P$.
\item Trazar la recta tangente a la hip\'erbola que pasa por el punto $P$, renombr\'andola como $t$.
\item Trazar la recta normal a la hip\'erbola que pasa por $P$, renombr\'andola como $n$.
\item Trazar la imagen de $l$ y la de $O$ reflejados en $n$.
\item Con la herramienta \texttt{Vector entre Dos Puntos} trace los vectores que van de $O$ a $P$
y de $P$ a $O'$, renombr\'andolos como $u$ y $u'$, respectivamente.
\end{enumerate}
Al desplazar el punto $O$ en la regi\'on que no contiene a ninguno de los focos,
?`que sucede con la reflexi\'on $l'$ de $l$? Repetir el experimento con distintos
valores de $a$ y de $b$ y explicar el resultado.   
   \begin{center}
\includegraphics[scale=0.33]{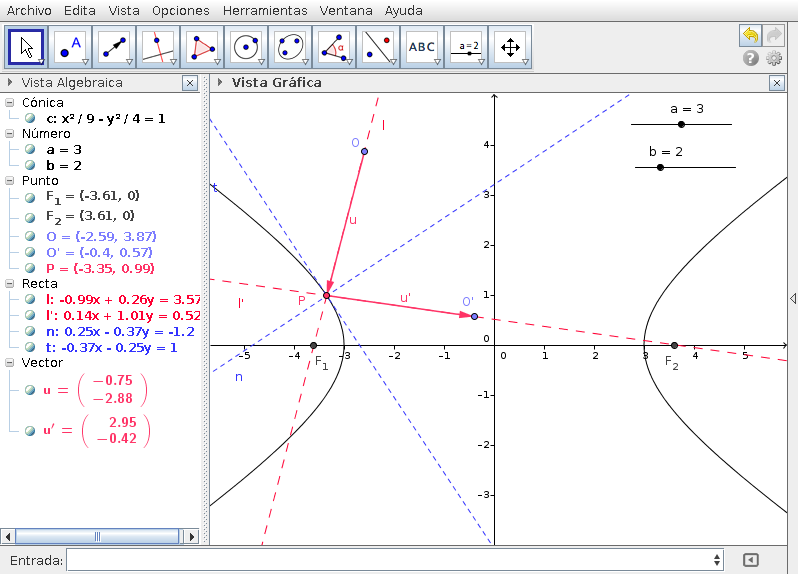}
   \end{center}
Para disminuir el tama\~no de las antenas parab\'olicas sin reducir su potencia, se 
utilizan las llamadas \emph{antenas de Cassegrain}\index{Cassegrain, Laurent}, que constan de un reflector primario 
parab\'olico, que comparte su (\'unico) foco con uno de los dos focos de un reflector 
secundario hiperb\'olico. En la figura se muestra una imagen\footnote{Disponible
libremente en los archivos fotogr\'aficos de la NASA, \url{https://solarsystem.nasa.gov/galleries/prototype-voyager} y tambi\'en en Wikimedia
\url{https://upload.wikimedia.org/wikipedia/commons/7/7d/Voyager_Spacecraft_During_Vibration_Testing_-_GPN-2003-000008.jpg}.} 
de la antena compuesta del sat\'elite Voyager, en la que mediante ensayo y error se han
superimpuesto\index{GeoGebra!im\'agenes}, respectivamente, una par\'abola y una hip\'erbola en la vista lateral de los
reflectores primario y secundario de la antena (por supuesto, con un cierto error por la perspectiva de la fotograf\'ia utilizada). Repetir la construcci\'on, a\~nadiendo los trazos necesarios para comprobar que, efectivamente, un rayo vertical que incide en el reflector primario, pasa por el otro foco de la hip\'erbola.
\begin{center}
\includegraphics[scale=0.37]{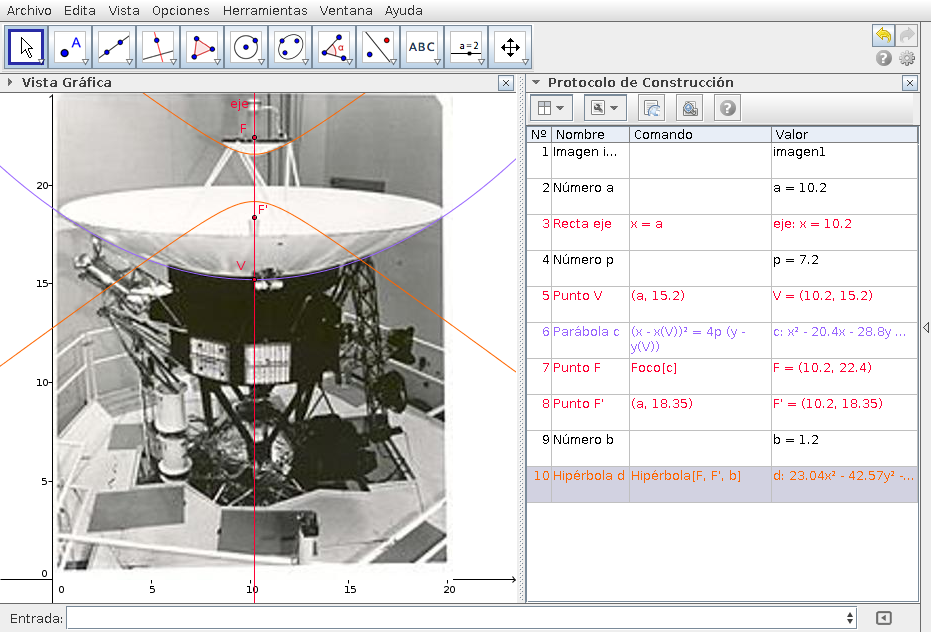}
\end{center}

\end{enumerate}

\chapter{N\'umeros complejos}
\vspace{2.5cm}\setcounter{footnote}{0}
\section{Nociones b\'asicas}

Desde un punto de vista abstracto, llamamos `n\'umeros' a cualquier clase de objetos en los que se haya definido
un par de operaciones con las propiedades de la suma y el producto de los n\'umeros reales\footnote{Tales clases
de objetos se denominan \emph{anillos algebraicos}.\index{Anillo algebraico}}. En este cap\'itulo y
en los restantes, veremos tres ejemplos de tales `n\'umeros' generalizados: los complejos,
los polinomios y las matrices.

Los n\'umeros complejos surgen de la necesidad de resolver ecuaciones como $x^2+1=0$. No hay ning\'un n\'umero real
$x$ que la cumpla, dado que equivale a $x^2=-1$ y el cuadrado de cualquier n\'umero real, positivo o negativo, es siempre
mayor o igual que cero. De forma sorprendente, se puede definir un par de operaciones \emph{en los puntos del plano
cartesiano}\index{R2@$\mathbb{R}^2$}\index{Plano!cartesiano} $\mathbb{R}^2=\mathbb{R}\times\mathbb{R}$ tales que esta ecuaci\'on tenga soluci\'on. Si $(a,b)$, $(c,d)$
son puntos de $\mathbb{R}^2$, definimos su \emph{suma}\index{Operaciones!con complejos} como
$$
(a,b)+(c,d):=(a+c,b+d)\,.
$$ 
No hay nada raro aqu\'i, la suma es componente a componente; pero el \emph{producto} es m\'as extra\~no:
\begin{equation}\label{multc}
(a,b)\cdot (c,d):=(ac-bd,ad+bc)\, .
\end{equation}
Por ejemplo, el producto de $(1,-2)$ y $(3,4)$ se calcula como
$$
(1,-2)\cdot (3,4)=(1\cdot 3-(-2)\cdot 4,1\cdot 4+(-2)\cdot 3)=(3+8,4-6)=(11,-2)\,.
$$
Aunque estas operaciones puedan resultar poco familiares, reproducen las de los n\'umeros reales. Podemos identificar
cualquier n\'umero real $x\in\mathbb{R}$ con el punto del plano $(x,0)$, pues \'este se localiza en la recta de
abcisas en la misma posici\'on que ocupa $x$ en la recta real.
\begin{center}
    \begin{tikzpicture}
	   \draw[-,very thick,orange] (-3.5,0) -- (3.5,0) ;
       \draw[-] (0,-1.5) -- (0,1.5)  ;
	   \draw[-,very thick,orange] (-12,0) -- (-5,0) ;
	   \draw[-,fill,blue] (-10,0) circle (0.075) node[above,black]{$x$} ;
	   \draw[-,fill,blue] (-1.5,0) circle (0.075) node[above,black]{$(x,0)$} ;
	   \draw[->] (-4.75,0.5) to[out=30,in=150]  (-3.75,0.5);
 \end{tikzpicture}
\end{center}
Entonces, si multiplicamos $(x,0)$ por $(y,0)$ obtenemos $(xy,0)$, que es el punto del plano que corresponde
al n\'umero real $xy$. En otras palabras, el producto de puntos en el plano generaliza el de los n\'umeros
reales, de la misma forma que el producto de n\'umeros racionales generaliza el de n\'umeros enteros.

De acuerdo al p\'arrafo precedente, al n\'umero real $-1$ le corresponde el punto $(-1,0)$. Nos preguntamos si existir\'a alg\'un punto del plano $(a,b)\in\mathbb{R}^2$ tal que
$$
(a,b)^2=(a,b)\cdot (a,b)=(-1,0)\,.
$$
La respuesta es que s\'i: consideremos $(a,b)=(0,1)$ y calculemos de acuerdo con \eqref{multc}:
$$
(0,1)^2=(0,1)\cdot (0,1)=(0\cdot 0-1\cdot 1,0\cdot 1+1\cdot 0)=(-1,0)\,.
$$
De esta manera, !`con este producto resulta que $(0,1)$ es la ra\'iz cuadrada de $-1\equiv (-1,0)$!.

Para poner \'enfasis en que veremos los puntos del plano como `n\'umeros' vamos a utilizar la notaci\'on
$\mathbb{C}=\mathbb{R}^2$\index{C@$\mathbb{C}$} y llamaremos a los puntos del plano \emph{n\'umeros complejos}\index{Num@N\'umero!complejo}.
Como ya hemos se\~nalado, podemos considerar que todo n\'umero real es tambi\'en un complejo identificando
$x\in\mathbb{R}$ con $(x,0)\in\mathbb{C}$ (de la misma forma que cada entero $a\in\mathbb{Z}$ se identifica con
el racional $a/1\in\mathbb{Q}$). Los n\'umeros complejos que est\'an sobre el eje de ordenadas tambi\'en
tienen una notaci\'on especial. Comenzamos definiendo el \emph{n\'umero} $\mathbf{i}$\index{Num@N\'umero!i@$\mathbf{i}$} como
$$
\mathbf{i}:=(0,1),
$$
y observando que un n\'umero complejo sobre el eje de ordenadas, $(0,b)$, puede ponerse como
$$
(0,b)=b\cdot (0,1)=b\cdot \mathbf{i}\,,
$$
pues
$$
(0,b)=(b\cdot 0-0\cdot 1,b\cdot 1+0\cdot 0)=(b,0)\cdot (0,1)=b\cdot (0,1)\,.
$$
A partir de ahora omitiremos el punto $\cdot$ para denotar el producto. Con lo que acabamos de ver, se tiene
que los n\'umeros complejos se pueden escribir como
$$
(a,b)=(a,0)+(0,b)=a+b(0,1)=a+b\mathbf{i}\,.
$$
A $a+b\mathbf{i}$ se le llama la forma \emph{rectangular} o forma \emph{bin\'omica} del n\'umero complejo $(a,b)$\index{Num@N\'umero!complejo!forma bin\'omica}
(la forma $(a,b)$ es la \emph{cartesiana}). Al n\'umero real $a$ se
la llama la \emph{parte real}\index{Num@N\'umero!complejo!parte real} de $(a,b)=a+b\mathbf{i}$, y al n\'umero real $b$, su \emph{parte compleja} o \emph{imaginaria}\index{Num@N\'umero!complejo!parte imaginaria}. Es
com\'un escribir $z=a+b\mathbf{i}$ y $\mathrm{Re}z=a$, $\mathrm{Im}z=b$.

Los n\'umeros complejos est\'an \'intimamente relacionados con la geometr\'ia ecuclidea del plano. Por ejemplo,
la suma de n\'umeros complejos se corresponde geom\'etricamente con la adici\'on de vectores en el plano mediante
la \emph{ley del paralelogramo}\index{Ley!del paralelogramo}. La interpretaci\'on geom\'etrica del producto se estudiar\'a m\'as adelante. Tambi\'en se puede definir una operaci\'on muy interesante, que se corresponde
geom\'etricamente con la reflexi\'on respecto del eje de abcisas; se trata de la \emph{conjugaci\'on compleja}\index{Operaciones!con complejos!conjugaci\'on compleja}. Si
$z=a+b\mathbf{i}$, su complejo conjugado es el n\'umero denotado $\overline{z}:=a-b\mathbf{i}$. 
\begin{center}
\begin{tikzpicture}
\draw[very thin,color=gray!20,step=0.5] (-2,-3.15) grid (3.15,3.15);
\draw[-] (-2,0) -- (3.15,0) node[right] {$\mathrm{Re}z$};
\draw[-] (0,-3.15) -- (0,3.15) node[above] {$\mathrm{Im}z$};
\draw[dashed,color=red] (1.5,0) -- (1.5,2) node[above=6pt,right] {$z=a+b\mathbf{i}=(a,b)$};
\node[red] at (1.7,-0.2) {$a$};
\draw[dashed,color=red] (0,2) node[left]{$b$} -- (1.5,2);
\draw[fill=red,color=red] (1.5,2) circle (0.05);
\draw[color=red] (0,0)--(1.5,2);
\draw[fill=red,color=red] (1.5,-2) circle (0.05) node[below=6pt,right] {$\overline{z}=a-b\mathbf{i}=(a,-b)$};
\draw[dashed,color=red] (0,-2) node[left]{$-b$} -- (1.5,-2);
\draw[color=red] (0,0)--(1.5,-2);
\draw[dashed,color=red] (1.5,0) -- (1.5,-2);
\end{tikzpicture}
\end{center}

El producto de n\'umeros complejos es conmutativo, asociativo, distributivo sobre la suma y, adem\'as, posee
un elemento neutro, el n\'umero $1=(1,0)=1+0\mathbf{i}$, como se comprueba inmediatamente. ?`Qu\'e ocurre con
el inverso? Se tiene que \emph{todo complejo distinto del $0=0+0\mathbf{i}$ tiene inverso respecto del producto}.\index{Operaciones!con complejos!inverso}
Es decir, dado cualquier complejo $z=a+b\mathbf{i}$, existe otro $z^{-1}=1/(a+b\mathbf{i})$ tal que $zz^{-1}=1=z^{-1}z$.
La construcci\'on de $z^{-1}$, que tambi\'en se denota $1/z$, es muy sencilla usando el conjugado de $z$:
$$
z^{-1}=\frac{1}{z}=\frac{\overline{z}}{z\overline{z}}=\frac{a-b\mathbf{i}}{a^2+b^2}\, .
$$
En efecto:
$$
zz^{-1}=z\frac{1}{z}=(a+b\mathbf{i})\frac{a-bi}{a^2+b^2}=\frac{a^2+b^2}{a^2+b^2}=1\,.
$$
As\'i, el inverso de $1-3\mathbf{i}$ es
$$
\frac{1}{1-3\mathbf{i}}=\frac{1+3\mathbf{i}}{1+3^2}=\frac{1}{10}+\frac{3}{10}\mathbf{i}\,.
$$
Con esto ya se puede calcular la divisi\'on entre dos n\'umeros complejos, dado que dividir $z$ entre $w$ no
es m\'as que multiplicar $z$ por el inverso de $w$. Por ejemplo, para dividir $2-2\mathbf{i}$ entre $1+3\mathbf{i}$,
proceder\'iamos como sigue:
\begin{align*}
\frac{2-2\mathbf{i}}{1+3\mathbf{i}}&=(2-2\mathbf{i})\frac{1-3\mathbf{i}}{1+9}\\
&=\frac{2-6+(-6-2)\mathbf{i}}{10}\\
&=-\frac{4}{10}-\frac{8}{10}\mathbf{i}\\
&=-\frac{2}{5}-\frac{4}{5}\mathbf{i}\,.
\end{align*}
En los c\'alculos precedentes nos ha aparecido repetidamentela cantidad $z\overline{z}$, que tiene una importante
interpretaci\'on geom\'etrica. Si $z=a+b\mathbf{i}$, de acuerdo con el teorema de Pit\'agoras se tiene que\index{Pit\'agoras!teorema}
$$
z\overline{z}=a^2+b^2
$$
es el cuadrado de la diagonal del tri\'angulo rect\'angulo con v\'ertices $a+b\mathbf{i}$, $a$ y $b\mathbf{i}$.
\begin{center}
\begin{tikzpicture}
\draw[very thin,color=gray!20,step=0.5] (-2,-2) grid (3.15,3.15);
\draw[-] (-2,0) -- (3.15,0) node[right] {$\mathrm{Re}z$};
\draw[-] (0,-2) -- (0,3.15) node[above] {$\mathrm{Im}z$};
\draw[dashed,color=red] (1.5,0) -- (1.5,2) node[above=6pt,right] {$z=a+b\mathbf{i}$};
\node[red] at (1.7,-0.2) {$a$};
\draw[dashed,color=red] (0,2) node[left]{$b$} -- (1.5,2);
\draw[color=red] (0,0)-- node[sloped,midway,above]{$|z|=\sqrt{a^2+b^2}$} (1.5,2);
\draw[color=blue] (0.5,0) arc (0:53.1:0.5) ; \node[color=blue] at (20:0.75)  {$\theta$};
\end{tikzpicture}
\end{center}
La cantidad $\sqrt{a^2+b^2}$ se denomina el \emph{m\'odulo}\index{Num@N\'umero!complejo!m\'odulo} del n\'umero complejo $z=a+b\mathbf{i}$,
y como acabamos de ver se corresponde con la longitud del vector que va del origen al punto $(a,b)$ en el plano
complejo. Se denota $|z|=\sqrt{a^2+b^2}$.

Existe otra forma de representar los n\'umeros complejos que es muy \'util para realizar otros c\'alculos
m\'as complicados, como las potencias y ra\'ices de un complejo. Se denomina la forma \emph{polar}.\index{Num@N\'umero!complejo!forma polar}
La idea consiste en considerar al n\'umero $z=a+b\mathbf{i}$ a trav\'es de su m\'odulo $|z|$ y el \'angulo $\theta$
que forma con el eje real. Antes de nada, observemos que todo n\'umero complejo $z$ se puede escribir como
$z=|z|u$, donde $u$ es un n\'umero complejo que est\'a sobre la circunferencia unidad\footnote{De hecho,
$u=z/|z|=(a+b\mathbf{i})(a^2+b^2)$.}. Entonces, podemos calcular el seno y el coseno de $\theta$:
$$
\sin\theta =\frac{b}{\sqrt{a^2+b^2}},
$$ 
y
$$
\cos\theta =\frac{a}{\sqrt{a^2+b^2}}.
$$
Por tanto:
$$
a=\mathrm{Re}z=\sqrt{a^2+b^2}\cos\theta =|z|\cos\theta,
$$
y
$$
b=\mathrm{Im}z=\sqrt{a^2+b^2}\sin\theta =|z|\sin\theta,
$$
es decir, se puede escribir
$$
z=a+b\mathbf{i}=|z|\cos\theta +|z|\mathbf{i}\sin\theta =|z|(\cos\theta +\mathbf{i}\sin\theta ).
$$
Al \'angulo $\theta$ se le denomina el \emph{argumento} de $z$\index{Num@N\'umero!complejo!argumento}. N\'otese que el argumento no est\'a un\'ivocamente
definido, pues \'angulos que difieren en $2\pi$ conducen al mismo n\'umero. Si no se especifica lo contrario, siempre
que se hable del argumento de un n\'umero  complejo se entender\'a aqu\'el \'angulo $\theta$ que cumple la relaci\'on
anterior y est\'a comprendido entre $-\pi$ y $\pi$ (esto se conoce como la \emph{rama principal} del argumento).

Fij\'emonos en que, si multiplicamos dos n\'umeros complejos $z$ y $w$ escritos en forma polar
respectivamente como
$z=|z|(\cos\theta +\mathbf{i}\sin\theta )$ y $w=|w|(\cos\alpha +\mathbf{i}\sin\alpha )$, resulta:
\begin{align*}
zw &=|z||w|(\cos\theta +\mathbf{i}\sin\theta )(\cos\alpha +\mathbf{i}\sin\alpha ) \\
	&=|z||w|(\cos\theta\cos\alpha -\sin\theta\sin\alpha +\mathbf{i}(\cos\theta\sin\alpha +\sin\theta\cos\alpha )) \\
	&=|z||w|(\cos(\theta+\alpha)+\mathbf{i}(\sin(\theta+\alpha)) ,
\end{align*}
es decir, \emph{los m\'odulos se multiplican y los argumentos se suman}.
Geom\'etricamente, este hecho permite interpretar la operaci\'on de multiplicaci\'on
por el n\'umero complejo $z=a+bi=|z|(\cos\theta +\mathbf{i}\sin\theta )$ sobre otro
complejo $w$ de la siguiente manera: se gira el vector que representa a $w$ un \'angulo
$\theta$ en sentido contrario a las agujas del reloj, y su m\'odulo $\vert w\vert$ se
alarga (o encoge) por un factor $\vert z\vert$. En particular, multiplicar por el
n\'umero complejo $\mathbf{i}$ (cuyo m\'odulo es la unidad y cuyo argumento es $\theta =\pi/2$),
equivale a realizar una rotaci\'on en sentido contrario a las agujas del reloj de
$90\degree$ (un cuarto de vuelta).

\begin{center}
\begin{tikzpicture}
\draw[very thin,color=gray!20,step=0.5] (-2,-2) grid (3.15,5.15);
\draw[-] (-2,0) -- (3.15,0) node[right] {$\mathrm{Re}z$};
\draw[-] (0,-2) -- (0,5.15) node[above] {$\mathrm{Im}z$};
\draw[color=red] (0,0)--node[sloped,midway,above]{$|z|$} (1.5,2);
\draw[color=red,fill=red] (1.5,2) circle (0.05);
\node[red] at (1.6,2.25) {$z$};
\draw[color=red] (0.5,0) arc (0:53.1:0.5) ; 
\node[color=red] at (30:0.75)  {$\theta$};
\draw[color=blue] (0,0)--(2,0.5);
\draw[color=blue,fill=blue] (2,0.5) circle (0.05);
\node[blue] at (2.25,0.6) {$w$};
\draw[color=blue] (0,0)--node[sloped,midway,above]{$|z|\,|w|$}(2,4.75);
\draw[color=blue,fill=blue] (2,4.75) circle (0.05);
\node[blue] at (2.1,5) {$z\cdot w$};
\draw[color=red] (0.85,0.23) arc (0:70.5:0.72) ;
\node[color=red] at (30:1.25)  {$\theta$};
\end{tikzpicture}
\end{center}

El comportamiento de los argumentos,
semejante al de la funci\'on exponencial, condujo a Euler a introducir la exponencial compleja\index{Euler, Leonhard!f\'ormula exponencial}
$$
e^{\mathbf{i}\alpha}:=\cos\alpha +\mathbf{i}\sin \alpha \,,
$$
de manera que un n\'umero complejo $z=a+b\mathbf{i}$ se escribir\'a como $z=|z|e^{\mathbf{i}\theta}$, donde
$\theta =\arctan (b/a)$ es su argumento\footnote{N\'otese que para el valor particular
$\theta =\pi$ se obtiene la f\'ormula de Euler que ya nos apareci\'o en la secci\'on \ref{explog}:
$e^{\mathbf{i}\pi}=-1$.},
y si se tiene otro complejo $w=|w|\, e^{\mathbf{i}\alpha}$, entonces el
producto se puede calcular como
$$
zw=|z||w|\, e^{\mathbf{i}\theta}e^{\mathbf{i}\alpha}=|z||w|\, e^{\mathbf{i}(\theta+\alpha)}\,.
$$
Por ejemplo, consideremos los n\'umeros $z=1-\mathbf{i}$ y $w=2+2\mathbf{i}$. El primero tiene m\'odulo
$|z|=\sqrt{1+1}=\sqrt{2}$ y el segundo $|w|=\sqrt{4+4}=2\sqrt{2}$. Los argumentos
respectivos son $\theta =\arctan (-1/1)=-\pi/4$ y $\alpha =\arctan (2/2)=\pi/4$. De esta forma,
$$
z=\sqrt{2}\, e^{-\mathbf{i}\pi/4},\quad w=2\sqrt{2}\, e^{\mathbf{i}\pi/4}.
$$
Su producto puede calcularse entonces como
$$
zw=2\sqrt{2}\sqrt{2}\, e^{\mathbf{i}(-\pi/4+\pi/4)}=4\,.
$$

Como hemos mencionado, el uso de la forma polar simplifica mucho el c\'alculo de potencias y ra\'ices
de n\'umeros complejos. Las potencias son particularmente simples cuando se escribe $z=a+ib$ en forma
polar $z=|z|e^{\mathbf{i}\theta}$, pues en ese caso
$$
z^n=\left( |z|\, e^{\mathbf{i}\theta}\right)^n=|z|^n\,  e^{\mathbf{i}n\theta} .
$$
Respecto a las ra\'ices, se tiene que si $z\neq 0$, entonces $z$ tiene exactamente $n$ ra\'ices $n-$\'esimas,\index{Rai@Ra\'iz!$n$-\'esima}
dadas por
\begin{equation}\label{comproot}
z_k =\sqrt[n]{|z|}\, e^{\mathbf{i}\left( \theta+2k\pi \right)/n}\, ,\,  k\in\{0,1,\ldots ,n-1\} .
\end{equation}
Como ejemplo, veamos cu\'ales son las ra\'ices cuadradas de $z=-1$. Comenzamos por escribir $-1$ en forma
polar:
$$
-1=|1|\, e^{\mathbf{i}\pi},
$$
y a continuaci\'on aplicamos la f\'ormula \eqref{comproot} con $n=2$, que nos da
$$
z_0=\sqrt{|1|}\, e^{\mathbf{i}(\pi +2\cdot 0\cdot \pi)/2}=e^{\mathbf{i}\pi/2}\mbox{ para }k=0,
$$
y
$$
z_1=\sqrt{|1|}\, e^{\mathbf{i}(\pi +2\cdot 1\cdot \pi)/2}=e^{\mathbf{i}3\pi/2}\mbox{ para }k=1.
$$
Por tanto, las ra\'ices cuadradas de $-1$ son $z_0=e^{\mathbf{i}\pi/2}=\mathbf{i}$ y
$z_1=e^{\mathbf{i}3\pi/2}=-\mathbf{i}$.

\section{Aplicaciones}
Cuando se conecta un electrodom\'estico\index{Electrodom\'estico} a la red, su circuito recibe un voltaje\index{Voltaje}
que var\'ia con el tiempo en forma senoidal, con una frecuencia $f=60\, Hz$ y una
amplitud de $120\, V$. Desde el punto de vista f\'isico, esto quiere decir que el voltaje a lo largo del cable cambia $60$ veces por segundo entre los valores $120\, V$ y 
$-120\, V$, lo cual implica que los electrones
en el cable oscilan a derecha e izquierda con esa frecuencia. Por este motivo, la
corriente dom\'estica se dice que es \emph{alterna}\index{Corriente!alterna} \'o de tipo AC (la que proporciona
una bater\'ia, por el contrario, es \emph{continua} o de tipo DC\index{Corriente!continua} y en ella no hay oscilaciones: los electrones corren
todos siempre en la misma direcci\'on, y por eso las bater\'ias se descargan al cabo de un 
cierto tiempo). En cualquier caso, aunque la posici\'on promedio de los \emph{electrones} 
no cambia, si hay una propagaci\'on del \emph{campo electromagn\'etico} a lo largo del 
cable, y 
ese campo electromagn\'etico es lo que alimenta al electrodom\'estico.
\begin{center}
\includegraphics[scale=0.33]{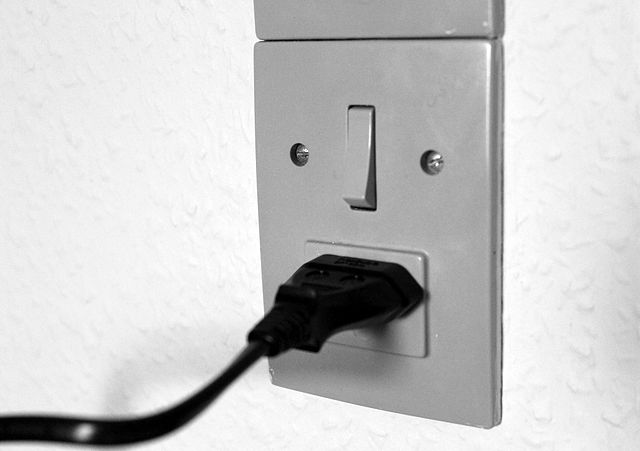} 
\end{center}

En t\'erminos matem\'aticos, esto significa que el voltaje est\'a
descrito por la funci\'on
\begin{equation}\label{senoidal}
V(t)=V_0\, \sin (wt+\phi)
\end{equation}
donde $w=2\pi\, 60\simeq 377$rad/s es la frecuencia angular (el n\'umero de ciclos por
segundo, que es $f=60$, multiplicado por los $2\pi$ radianes que hay en cada ciclo), 
$\phi$  es la fase\index{Fase} de la corriente y $V_0=120\, V$ es la amplitud\index{Amplitud} del 
voltaje\footnote{Realmente, el 
voltaje proporcionado por cualquier planta se expresa como el valor cuadr\'atico medio 
(RMS en ingl\'es). Como consecuencia, los picos en la amplitud no son de $120\, V$, sino 
de $120\sqrt{2}\simeq 170\, V$. Por simplicidad, aqu\'i omitiremos el factor $\sqrt{2}$.}.

Una vez que el campo electromagn\'etico comienza a alimentar un circuito AC,
la frecuencia $w$ de la corriente no var\'ia, as\'i que las magnitudes que hay
que considerar son $V_0$ y $\phi$. Un \emph{fasor}\index{Fasor} no es m\'as que un n\'umero
complejo\index{Num@N\'umero!complejo} que representa la funci\'on senoidal\index{Funci\'on!senoidal} \eqref{senoidal}. Concretamente,
a la funci\'on \eqref{senoidal} se le asocia el fasor $X(t)$ dado por\footnote{En
ingenier\'ia el\'ectrica se usa $i$ para representar la intensidad de la corriente,
por lo que resulta conveniente utilizar otra letra que denote la unidad imaginaria. Es 
com\'un tomar $j$ y as\'i lo haremos en esta secci\'on.}
\begin{equation}\label{fasor}
X=V_0\, e^{j\phi}\,,
\end{equation}
y viceversa. En particular, dado el fasor \eqref{fasor}, la funci\'on senoidal
\eqref{senoidal} se recupera multiplicando por $e^{jwt}$ y tomando la parte imaginaria
del resultado:
$$
\mathrm{Im}(Xe^{jwt})=\mathrm{Im}(V_0 e^{j\phi}e^{jwt})
=\mathrm{Im}(V_0 e^{j(\phi+wt)})=\mathrm{Im}(V_0 (\cos(\phi+wt)+j\sin(\phi+wt)))
=V_0 \sin(wt+\phi)\,.
$$
Conviene resaltar que en ingenier\'ia es frecuente utilizar la notaci\'on
$$
X=V_0\, e^{j\phi}=V_0 \angle \phi\,.
$$
La utilidad de los fasores se manifiesta en lo siguiente: cuando en un circuito se
superponen dos funciones senoidales $A\sin (wt+\phi)$ y $B\sin (wt+\psi)$, la
se\~nal resultante puede calcularse muy f\'acilmente a partir de la suma de los fasores
asociados $Ae^{j\phi}$ y $Be^{j\psi}$, pues si
$$
Ce^{j\varphi}=Ae^{j\phi}+Be^{j\psi}\,,
$$
entonces,
$$
A\sin (wt+\phi)+B\sin (wt+\psi)=\mathrm{Im}(Ae^{j(wt+\phi)}+Be^{j(wt+\psi)})
=\mathrm{Im}(e^{jwt}(Ae^{j\phi}+Be^{j\psi})
=\mathrm{Im}(e^{jwt}Ce^{j\varphi})\,;
$$
es decir, la se\~nal resultante de la superposici\'on es la correspondiente
al fasor suma de los fasores originales.

Una observaci\'on importante es que todo el an\'alisis se realiza de manera
an\'aloga si estamos tratando con ondas cosenoidales $V(t)=V_0\, \cos (wt+\phi)$. La
\'unica diferencia es que, en este caso, la se\~nal se recupera a partir del fasor
$X=V_0\, e^{j\phi}$ multiplicando por $e^{jwt}$ y tomando la parte \emph{real}
del resultado.

A la hora de operar con fasores, para el producto y el cociente obviamente es mejor utilizar la forma exponencial, mientras que para la suma y la resta es conveniente
pasar primero el fasor a la forma bin\'omica\index{Num@N\'umero!complejo!forma bin\'omica}, mediante los m\'etodos que ya conocemos. As\'i, dados $X=10\angle 60\degree$, $Y=5\angle 45\degree$, tendremos que
$$
X\cdot Y=50\angle 105\degree,\quad X/Y=2\angle 15\degree\,,
$$
y
\begin{align*}
X+Y=& 10\angle 60\degree +5\angle 45\degree \\
=& 10\,e^{j\pi/3}+5\,e^{j\pi/4} \\
=& 10(\cos (\pi/3)+j\sin (\pi/3))+5(\cos (\pi/4)+j\sin (\pi/4)) \\
=& 10(1/2+j\sqrt{3}/2)+5(\sqrt{2}/2+j\sqrt{2}/2)\\
=& 8{.}54+j12{.}2\\
=& 14{.}89\angle 55\degree\,.
\end{align*}

Veamos ahora c\'omo se aplican estos conceptos en el an\'alisis de circuitos
el\'ectricos. Los elementos b\'asicos de un circuito son los \emph{resistores} $R$\index{Resistor},
los \emph{capacitores} (o condensadores)\index{Capacitor}\index{Condensador} $C$ y los \emph{inductores}\index{Inductor} (o bobinas)\index{Bobina} $L$.
\begin{center}
\begin{tikzpicture}
\path (0,0) coordinate (ref_gnd);
\draw
  (ref_gnd) to[sV=$V_s$] ++(0,2)
            to[R=\(R\)] ++(3,0) 
            to[L=\(L\)] ++(3,0) 
            to[C=\(C\)] ++(0,-2) 
  -- (ref_gnd);
\end{tikzpicture}
\end{center}
Cada uno de ellos se caracteriza mediante un n\'umero complejo,
la \emph{impedancia}\index{Impedancia} $Z$ del elemento en 
cuesti\'on, cuya parte real se denomina \emph{resistencia}\index{Resistencia}
y cuya parte imaginaria se denomina \emph{reactancia}\index{Reactancia}. Los resistores no poseen reactancia (su impedancia es real) y tanto los 
capacitores como los inductores carecen de resistencia (su impedancia es imaginaria).
La siguiente tabla resume la asignaci\'on de impedancias a los elementos de un circuito
RLC alimentado por una se\~nal peri\'odica de frecuencia $w$:
\begin{center} 
\begin{tabular}{|c|c|c|c|}
\hline 
Elemento & Resistencia ($R$) & Reactancia ($X$) & Impedancia ($Z$)\\ 
\hline 
Resistor & $R$ & $0$ & $Z_R=R$ \\ 
\hline 
Inductor & $0$ & $wL$ & $Z_L=jwL$  \\ 
\hline 
Capacitor & $0$ & $-\frac{1}{wC}$ & $Z_C=-\frac{j}{wC}$ \\ 
\hline 
\end{tabular}
\end{center}
\medskip
La impedancia puede describirse tambi\'en mediante un fasor:
\begin{align}\label{impedancias}
Z_R=& R\angle 0\degree \nonumber\\[6pt]
Z_L=& wL\angle 90\degree \\ 
Z_C=& \frac{1}{wC}\angle -90\degree \nonumber\,. 
\end{align}

A la vista de estas expresiones, tenemos lo siguiente: por el circuito fluye una
corriente que podemos suponer (para simplificar) en la forma
\begin{equation}\label{intensidad}
i(t)=I_0\sin (wt)=I_0\angle 0\degree \,,
\end{equation}
es decir, su fase es $\phi=0$. Esta corriente es la misma en cada tramo del circuito.
Lo que cambia en cada tramo es la ca\'ida de voltaje $V$: entre los extremos del
resistor $R$, habr\'a una ca\'ida (o diferencia) de voltaje $V_R$, entre los del
inductor $L$, una ca\'ida $V_L$, y entre los del capacitor, una $V_C$. Para determinarlas,
se dispone de la Ley de Ohm\index{Ley!de Ohm}, que afirma que
$$
V=iZ
$$
de manera que, usando \eqref{impedancias} y \eqref{intensidad},
\begin{align}\label{voltajes}
V_R=& I_0R = I_0R\angle 0\degree \nonumber\\[6pt]
V_L=& jwI_0L = wI_0L\angle 90\degree \\ 
V_C=& -\frac{jI_0}{wC} =\frac{I_0}{wC}\angle -90\degree \nonumber\,. 
\end{align}
N\'otese que el voltaje a lo largo del resistor `est\'a en fase' con la corriente,
a lo largo del inductor est\'a `adelantado' a la corriente en $90\degree$ y a lo largo
del capacitor est\'a `retrasado' $-90\degree$. Precisamente, una de las ventajas de la 
notaci\'on de fasores es que hace evidente este tipo de relaciones (que, aunque en este
breve tratamiento no se aprecie, son de fundamental importancia en ingenier\'ia).\\
Ahora, podemos determinar el voltaje
de la fuente, $V_s$, aplicando otro principio b\'asico del an\'alisis de circuitos, la
Ley de Kirchoff\index{Ley!de Kirchoff}, que establece que en un circuito en serie:
$$
V_s=V_R+V_L+V_C\,.
$$
Dado que aparece una suma, conviene trabajar con las expresiones binomiales en
\eqref{voltajes}, as\'i que
$$
V_s=I_0R+j\left( wI_0L-\frac{I_0}{wC}\right)=
I_0\left( R+j\left( wL-\frac{1}{wC} \right)\right)\,.
$$
La se\~nal senoidal asociada a este fasor se puede recuperar expresando este resultado
en forma polar, multiplicando por $e^{jwt}$ y tomando la parte imaginaria:
\begin{align*}
V_s(t)=& \mathrm{Im}\left( I_0 e^{jwt}\sqrt{R^2+\left( wL-\frac{1}{wC}\right)^2}
e^{j\arctan \left( \frac{wL}{R}-\frac{1}{wCR}\right)} \right) \\
=& I_0\sqrt{R^2+\left( wL-\frac{1}{wC}\right)^2}\,\sin
\left( wt+\arctan \left( \frac{wL}{R}-\frac{1}{wCR}\right)\right)\,.
\end{align*}
Es interesante observar que tambi\'en se cumple una relaci\'on de la forma
$$
V_s=iZ\,,
$$
donde ahora $Z=Z_R+Z_L+Z_C$ es la \emph{impedancia del circuito}.

Se\~nalemos finalmente que en un circuito en paralelo el an\'alisis es similar,
solo que ahora la Ley de Kirchoff establece que
$$
I=I_R+I_L+I_C\,.
$$

\section{Ejercicios}
\begin{enumerate}  
\setcounter{enumi}{112}
\item 
Realizar las siguientes operaciones con n\'umeros complejos:
\begin{multicols}{3}
 \begin{enumerate}[(a)]
\item $(-5-4\mathbf{i}) + (3-6\mathbf{i})$
\item $(1,-7) + (2,5)$
\item $(\sqrt{3} -\sqrt{2}\mathbf{i}) - (\sqrt{2} - \sqrt{3}\mathbf{i})$
\item $-4 + (4+2\pi \mathbf{i})$
\item $(1,4) - (4,1)$
\item $(-3+2\mathbf{i})(1+2\mathbf{i})$
\item $(2,3)(1,5)$
\item $\overline{(-2+4\mathbf{i})}$
\item $\overline{(-\pi^2-\frac{1}{3}\mathbf{i})}$
\end{enumerate}
\end{multicols}

\item 
Si defini\'eramos el producto de complejos por componentes (como la suma), es decir,
$$
(a_1,b_1) (a_2,b_2) = (a_1a_2 , b_1b_2),
$$
mostrar que:
\begin{enumerate}[(a)]
\item S\'\i\ se cumplir\'ia la asociatividad
y habr\'ia un neutro (?`cu\'al?).
\item No todos los complejos tendr\'ian inverso (es decir, no podr\'iamos dividir).
\item No existir\'\i a un n\'umero que, elevado al cuadrado, resulte $-1$.
\end{enumerate}

\item 
Calcular los siguientes cocientes de n\'umeros complejos: 
\begin{multicols}{2}
 \begin{enumerate}[(a)]
\item $\dfrac{2+3i}{2-3i}$
\item $\dfrac{1-i}{2+5i}$
\item $\dfrac{(2,-4)}{(4,3)}$
\item $\dfrac{1}{i}$
\end{enumerate}
\end{multicols}

\item
Demostrar las siguientes propiedades de la conjugaci\'on de n\'umeros complejos:
\begin{multicols}{2}
 \begin{enumerate}[(a)]
\item $\overline{z_1+z_2} = \overline{z_1}+ \overline{z_2}$
\item $\overline{z_1z_2} = \overline{z_1} \overline{z_2}$
\item $\overline{\left (\frac{1}{z} \right )} = \frac{1}{\overline{z}}$, para $z \not= 0$
\item $\overline{\overline{z}} = z$
\end{enumerate}
\end{multicols}

\item 
Representar en el plano complejo los siguientes conjuntos de n\'umeros complejos: 
\begin{multicols}{2}
 \begin{enumerate}[(a)]
\item $2, -2, 2\mathbf{i}, -2\mathbf{i}$
\item $(2+2\mathbf{i}),(2-2\mathbf{i}), (-2+2\mathbf{i}), (-2-2\mathbf{i})$
\item $(1,5), (-2,3), (2\pi,-4), (-2\sqrt{2},-3\sqrt{2})$
\item Naturales
\item Enteros
\item Racionales
\item Reales
\item Imaginarios puros
\end{enumerate}
\end{multicols}

\item 
Escribir los siguientes n\'umeros complejos en forma polar:
\begin{multicols}{2}
 \begin{enumerate}[(a)]
\item $(1+2\mathbf{i}), (-2\sqrt{2}+2\sqrt{3}\mathbf{i}), (-3-\mathbf{i}), (\frac{\pi}{2}-\frac{\pi}{2} \mathbf{i})$ 
\item $3, -3, 3\mathbf{i}, -3\mathbf{i}$
\item $(3+3\mathbf{i}), (3-3\mathbf{i}), (-3+3\mathbf{i}), (-3-3\mathbf{i})$
\item $(1,\frac{3}{2}), (2,3), (4,6), (\sqrt{2},\frac{3\sqrt{2}}{2})$
\item $(3+5\mathbf{i}), (3-5\mathbf{i})$
\end{enumerate}
\end{multicols}

\item 
Escribir los siguientes n\'umeros complejos en forma bin\'omica:
\begin{multicols}{2}
 \begin{enumerate}[(a)]
\item $3e^{\mathbf{i}\pi/3}, 5e^{\mathbf{i}\pi/2}, e^{\mathbf{i}\pi}, \sqrt{6}e^{-\mathbf{i}\pi/4}$ 
\item $e^{\mathbf{i}\pi/2}, e^{\mathbf{i}\pi}, e^{\mathbf{i}3\pi/2}, e^{\mathbf{i}2\pi}$
\item $e^{\mathbf{i}0}, e^{\mathbf{i}2\pi/5}, e^{\mathbf{i}4\pi/5}, e^{\mathbf{i}6\pi/5}, e^{\mathbf{i}8\pi/5}, e^{\mathbf{i}10\pi/5}$
\item $e^{\mathbf{i}\pi/6}, 2e^{\mathbf{i}\pi/2}, 3e^{\mathbf{i}\pi/2}, \frac{1}{3}e^{\mathbf{i}\pi/2}, -2e^{\mathbf{i}\pi/2}$
\end{enumerate}
\end{multicols}

\item 
Realizar las siguientes operaciones en forma polar:
\begin{multicols}{4}
 \begin{enumerate}[(a)]
\item $3e^{\mathbf{i}2\pi/3} + e^{-\mathbf{i}\pi/4}$
\item $e^{\mathbf{i}\pi/3} 2e^{\mathbf{i}5\pi/6}$
\item $\overline{\sqrt{5}e^{\mathbf{i}3\pi/7}}$
\item $4e^{\mathbf{i}\pi/8} \overline{4e^{\mathbf{i}\pi/8}}$
\end{enumerate}
\end{multicols}


\item
Dibujar en el plano complejo dos n\'umeros complejos cualesquiera, llam\'emoslos $z_1$ y $z_2$ y, a partir de ellos, las representar las siguientes operaciones:
\begin{multicols}{4}
 \begin{enumerate}[(a)]
\item $-z_1$ y $-z_2$
\item $z_1 + z_2$
\item $z_1 z_2$
\item $\overline{z_1}$ y $\overline{z_2}$
\end{enumerate}
\end{multicols}

\item
Calcular el m\'odulo de los siguientes n\'umeros complejos:
\begin{multicols}{3}
 \begin{enumerate}[(a)]
\item $e^{\mathbf{i}\pi/4} + e^{-\mathbf{i}\pi/8}$
\item $\mathbf{i} e^{\mathbf{i}2\pi/7}$
\item $2 + 2e^{\mathbf{i}3\pi}$
\end{enumerate}
\end{multicols}

\item
Calcular las siguientes potencias:
\begin{multicols}{3}
 \begin{enumerate}[(a)]
\item $1^2, 1^3, 1^4, 1^5, \dots$
\item $\mathbf{i}^0, \mathbf{i}^1, \mathbf{i}^2, \mathbf{i}^3, \mathbf{i}^4, \mathbf{i}^5, \dots$
\item $\left (3e^{\mathbf{i}\pi/3} \right )^2$
\item $\left (\sqrt{5}e^{-\mathbf{i}\pi/3} \right )^3$
\item $(-2+\mathbf{i})^3$
\end{enumerate}
\end{multicols}

\item 
Demostrar la f\'ormula de A. De Moivre\index{De Moivre, Abraham} para cualquier $\alpha \in \mathbb{R}$ y cualquier $n \in \mathbb{N}$,
$$
(\cos \alpha + \mathbf{i} \sen \alpha)^n = \cos (n\alpha) + \mathbf{i} \sen (n\alpha). 
$$

\item 
Calcular y dibujar las ra\'{\i}ces cuadradas, c\'ubicas, cuartas, quintas, sextas y s\'eptimas de 1.

\item
Calcular y dibujar las siguientes ra\'{\i}ces:
\begin{multicols}{4}
 \begin{enumerate}[(a)]
\item $\left ( 9 e^{\mathbf{i} 2\pi/3}\right )^{1/3}$
\item $\left ( 64 e^{\mathbf{i} 2\pi/5}\right )^{1/2}$
\item $\sqrt[5]{\mathbf{i}}$
\item $(1+2\mathbf{i})^{1/4}$
\end{enumerate}
\end{multicols}

\item En este ejercicio comprobaremos que la f\'ormula \eqref{comproot} efectivamente caracteriza las
ra\'ices $n-$\'esimas de un n\'umero complejo.
\begin{enumerate}[(a)]
\item Comprueba que si $z_k$ designa a
$$
\sqrt[n]{|z|}\, e^{\mathbf{i}\left( \theta+2k\pi \right)}
$$
con $k\in\{0,1,\ldots ,n-1\}$, entonces $z^n_k=z$.
\item Comprueba que las ra\'ices $z_k$ son distintas dos a dos. Para ello, dados $0\leq k_1<k_2\leq n-1$,
estudia el cociente $z_{k_1}/z_{k_2}$.
\item Consideremos $w=|w|\, e^{\mathbf{i}\alpha}$ una ra\'iz $n-$\'esima cualquiera de $z$. Comprobar que
es una de las $z_k$ considerando que debe cumplirse $w^n=z$ (igualando m\'odulos y viendo que los
argumentos deben ser congruentes m\'odulo $n$).
\end{enumerate}

\item En el circuito de la figura el mult\'imetro mide la ca\'ida de
voltaje a lo largo del inductor $L$.
\begin{center}
\begin{circuitikz}
\draw[] (3,2) -- (3,3);
\draw[] (3,3) to[voltmeter] (5.5,3) to[] (5.5,2);
\path (0,0) coordinate (ref_gnd);
\draw
  (ref_gnd) to[sV=$120V\, @ 60Hz$] ++(0,2)
            to[R, l_=$12\Omega $] ++(3,0) 
            to[L, l_=$0{.}15H $] ++(3,0) 
            to[C=$100\mu F$] ++(0,-2) 
  -- (ref_gnd);
\end{circuitikz}
\end{center}
\noindent Determinar la impedancia total del circuito,
su corriente $I_0$ y la lectura del mult\'imetro.
\end{enumerate}

\section{Uso de la tecnolog\'ia de c\'omputo}

Maxima implementa la unidad compleja $\mathbf{i}$ en la forma de una constante\,
\texttt{\%i}\index{Maxima!constantes}. As\'i, podemos trabajar con los n\'umeros complejos en forma
rectangular como con cualquier otro n\'umero
\begin{maxcode}
\noindent
\begin{minipage}[t]{8ex}{\color{red}\bf
\begin{verbatim}
(%i1) 
\end{verbatim}
}
\end{minipage}
\begin{minipage}[t]{\textwidth}{\color{blue}
\begin{verbatim}
z:1-%i;
\end{verbatim}
}
\end{minipage}
\begin{math}\displaystyle
\parbox{8ex}{\color{labelcolor}(\%o1) }
1-i
\end{math}

\noindent
\begin{minipage}[t]{8ex}{\color{red}\bf
\begin{verbatim}
(%i2) 
\end{verbatim}
}
\end{minipage}
\begin{minipage}[t]{\textwidth}{\color{blue}
\begin{verbatim}
w:3+2*%i;
\end{verbatim}
}
\end{minipage}
\begin{math}\displaystyle
\parbox{8ex}{\color{labelcolor}(\%o2) }
2\,i+3
\end{math}

\noindent
\begin{minipage}[t]{8ex}{\color{red}\bf
\begin{verbatim}
(%i3) 
\end{verbatim}
}
\end{minipage}
\begin{minipage}[t]{\textwidth}{\color{blue}
\begin{verbatim}
conjugate(z);
\end{verbatim}
}
\end{minipage}
\begin{math}\displaystyle
\parbox{8ex}{\color{labelcolor}(\%o3) }
i+1
\end{math}

\noindent
\begin{minipage}[t]{8ex}{\color{red}\bf
\begin{verbatim}
(%i4) 
\end{verbatim}
}
\end{minipage}
\begin{minipage}[t]{\textwidth}{\color{blue}
\begin{verbatim}
conjugate(w);
\end{verbatim}
}
\end{minipage}
\begin{math}\displaystyle
\parbox{8ex}{\color{labelcolor}(\%o4) }
3-2\,i
\end{math}

\noindent
\begin{minipage}[t]{8ex}{\color{red}\bf
\begin{verbatim}
(%i5) 
\end{verbatim}
}
\end{minipage}
\begin{minipage}[t]{\textwidth}{\color{blue}
\begin{verbatim}
z+w;
\end{verbatim}
}
\end{minipage}
\begin{math}\displaystyle
\parbox{8ex}{\color{labelcolor}(\%o5) }
i+4
\end{math}

\noindent
\begin{minipage}[t]{8ex}{\color{red}\bf
\begin{verbatim}
(%i6) 
\end{verbatim}
}
\end{minipage}
\begin{minipage}[t]{\textwidth}{\color{blue}
\begin{verbatim}
z*w;
\end{verbatim}
}
\end{minipage}
\begin{math}\displaystyle
\parbox{8ex}{\color{labelcolor}(\%o6) }
\left( 1-i\right) \,\left( 2\,i+3\right) 
\end{math}
\end{maxcode}

Fij\'emonos en que Maxima realiza las operaciones simb\'olicamente (el producto
$zw$ se expresa exactamente como lo que es). Si queremos desarrollar (expandir)
el producto, debemos hacer
\begin{maxcode}
\noindent
\begin{minipage}[t]{8ex}{\color{red}\bf
\begin{verbatim}
(%i7) 
\end{verbatim}
}
\end{minipage}
\begin{minipage}[t]{\textwidth}{\color{blue}
\begin{verbatim}
expand(%);
\end{verbatim}
}
\end{minipage}
\begin{math}\displaystyle
\parbox{8ex}{\color{labelcolor}(\%o7) }
5-i
\end{math}
\end{maxcode}

Podemos extraer la parte real y la imaginaria de un complejo dado en forma rectangular:
\begin{maxcode}
\noindent
\begin{minipage}[t]{8ex}{\color{red}\bf
\begin{verbatim}
(%i8) 
\end{verbatim}
}
\end{minipage}
\begin{minipage}[t]{\textwidth}{\color{blue}
\begin{verbatim}
realpart(%o7);
\end{verbatim}
}
\end{minipage}
\begin{math}\displaystyle
\parbox{8ex}{\color{labelcolor}(\%o8) }
5
\end{math}

\noindent
\begin{minipage}[t]{8ex}{\color{red}\bf
\begin{verbatim}
(%i9) 
\end{verbatim}
}
\end{minipage}
\begin{minipage}[t]{\textwidth}{\color{blue}
\begin{verbatim}
imagpart(%o7);
\end{verbatim}
}
\end{minipage}
\begin{math}\displaystyle
\parbox{8ex}{\color{labelcolor}(\%o9) }
-1
\end{math}
\end{maxcode}

De la misma forma, podemos calcular otras magnitudes, como el m\'odulo y el argumento
(mediante los comandos \texttt{cabs} --de \textbf{c}omplex \textbf{abs}olute value
y \textbf{c}omplex \textbf{arg}ument, respectivamente)
\begin{maxcode}
\noindent
\begin{minipage}[t]{8ex}{\color{red}\bf
\begin{verbatim}
(%i10) 
\end{verbatim}
}
\end{minipage}
\begin{minipage}[t]{\textwidth}{\color{blue}
\begin{verbatim}
cabs(%o7);
\end{verbatim}
}
\end{minipage}
\begin{math}\displaystyle
\parbox{8ex}{\color{labelcolor}(\%o10) }
\sqrt{26}
\end{math}

\noindent
\begin{minipage}[t]{8ex}{\color{red}\bf
\begin{verbatim}
(%i11) 
\end{verbatim}
}
\end{minipage}
\begin{minipage}[t]{\textwidth}{\color{blue}
\begin{verbatim}
carg(%o7);
\end{verbatim}
}
\end{minipage}
\begin{math}\displaystyle
\parbox{8ex}{\color{labelcolor}(\%o11) }
-\mathrm{atan}\left( \frac{1}{5}\right) 
\end{math}
\end{maxcode}

El cociente se expresa como es habitual:
\begin{maxcode}
\noindent
\begin{minipage}[t]{8ex}{\color{red}\bf
\begin{verbatim}
(%i12) 
\end{verbatim}
}
\end{minipage}
\begin{minipage}[t]{\textwidth}{\color{blue}
\begin{verbatim}
z/w;
\end{verbatim}
}
\end{minipage}
\begin{math}\displaystyle
\parbox{8ex}{\color{labelcolor}(\%o12) }
\frac{1-i}{2\,i+3}
\end{math}
\end{maxcode}

Notemos que el resultado se expresa en forma simb\'olica. La expresi\'on
rectangular de $z/w$ se obtiene haciendo
\begin{maxcode}
\noindent
\begin{minipage}[t]{8ex}{\color{red}\bf
\begin{verbatim}
(%i13) 
\end{verbatim}
}
\end{minipage}
\begin{minipage}[t]{\textwidth}{\color{blue}
\begin{verbatim}
rectform(%);
\end{verbatim}
}
\end{minipage}
\begin{math}\displaystyle
\parbox{8ex}{\color{labelcolor}(\%o13) }
\frac{1}{13}-\frac{5\,i}{13}
\end{math}
\end{maxcode}

Comprobemos que el resultado es el que se obtiene multiplicando
$z$ por el inverso de $w$:
\begin{maxcode}
\noindent
\begin{minipage}[t]{8ex}{\color{red}\bf
\begin{verbatim}
(%i14) 
\end{verbatim}
}
\end{minipage}
\begin{minipage}[t]{\textwidth}{\color{blue}
\begin{verbatim}
rectform(conjugate(w)/(w*conjugate(w)));
\end{verbatim}
}
\end{minipage}
\begin{math}\displaystyle
\parbox{8ex}{\color{labelcolor}(\%o14) }
\frac{3}{13}-\frac{2\,i}{13}
\end{math}

\noindent
\begin{minipage}[t]{8ex}{\color{red}\bf
\begin{verbatim}
(%i15) 
\end{verbatim}
}
\end{minipage}
\begin{minipage}[t]{\textwidth}{\color{blue}
\begin{verbatim}
expand(z*%);
\end{verbatim}
}
\end{minipage}
\begin{math}\displaystyle
\parbox{8ex}{\color{labelcolor}(\%o15) }
\frac{1}{13}-\frac{5\,i}{13}
\end{math}
\end{maxcode}

Para pasar de la forma rectangular a la polar, disponemos del comando
\texttt{polarform}:
\begin{maxcode}
\noindent
\begin{minipage}[t]{8ex}{\color{red}\bf
\begin{verbatim}
(%i16) 
\end{verbatim}
}
\end{minipage}
\begin{minipage}[t]{\textwidth}{\color{blue}
\begin{verbatim}
polarform(z);
\end{verbatim}
}
\end{minipage}
\begin{math}\displaystyle
\parbox{8ex}{\color{labelcolor}(\%o16) }
\sqrt{2}\,{e}^{-\frac{i\,\pi }{4}}
\end{math}
\end{maxcode}

Los comandos \texttt{cabs} y \texttt{carg} funcionan con complejos en forma polar:
\begin{maxcode}
\noindent
\begin{minipage}[t]{8ex}{\color{red}\bf
\begin{verbatim}
(%i17) 
\end{verbatim}
}
\end{minipage}
\begin{minipage}[t]{\textwidth}{\color{blue}
\begin{verbatim}
cabs(%o16);
\end{verbatim}
}
\end{minipage}
\begin{math}\displaystyle
\parbox{8ex}{\color{labelcolor}(\%o17) }
\sqrt{2}
\end{math}

\noindent
\begin{minipage}[t]{8ex}{\color{red}\bf
\begin{verbatim}
(%i18) 
\end{verbatim}
}
\end{minipage}
\begin{minipage}[t]{\textwidth}{\color{blue}
\begin{verbatim}
carg(%o16);
\end{verbatim}
}
\end{minipage}
\begin{math}\displaystyle
\parbox{8ex}{\color{labelcolor}(\%o18) }
-\frac{\pi }{4}
\end{math}
\end{maxcode}

Por su parte, el comando \texttt{demoivre} pasa de la forma polar a la rectangular usando
la f\'ormula de Euler-DeMoivre $e^{\mathbf{i}\theta}=\cos\theta +\mathbf{i}\sin\theta$:
\begin{maxcode}
\noindent
\begin{minipage}[t]{8ex}{\color{red}\bf
\begin{verbatim}
(%i19) 
\end{verbatim}
}
\end{minipage}
\begin{minipage}[t]{\textwidth}{\color{blue}
\begin{verbatim}
demoivre(%o16);
\end{verbatim}
}
\end{minipage}
\begin{math}\displaystyle
\parbox{8ex}{\color{labelcolor}(\%o19) }
\sqrt{2}\,\left( \frac{1}{\sqrt{2}}-\frac{i}{\sqrt{2}}\right) 
\end{math}
\end{maxcode}

N\'otese que lo que resulta es el complejo original $z$, como debe ser:
\begin{maxcode}
\noindent
\begin{minipage}[t]{8ex}{\color{red}\bf
\begin{verbatim}
(%i20) 
\end{verbatim}
}
\end{minipage}
\begin{minipage}[t]{\textwidth}{\color{blue}
\begin{verbatim}
ratsimp(%);
\end{verbatim}
}
\end{minipage}
\begin{math}\displaystyle
\parbox{8ex}{\color{labelcolor}(\%o20) }
1-i
\end{math}
\end{maxcode}

El c\'alculo de las potencias no presenta dificultad alguna:
\begin{maxcode}
\noindent
\begin{minipage}[t]{8ex}{\color{red}\bf
\begin{verbatim}
(%i21) 
\end{verbatim}
}
\end{minipage}
\begin{minipage}[t]{\textwidth}{\color{blue}
\begin{verbatim}
expand(z^3);
\end{verbatim}
}
\end{minipage}
\begin{math}\displaystyle
\parbox{8ex}{\color{labelcolor}(\%o21) }
-2\,i-2
\end{math}

\noindent
\begin{minipage}[t]{8ex}{\color{red}\bf
\begin{verbatim}
(%i22) 
\end{verbatim}
}
\end{minipage}
\begin{minipage}[t]{\textwidth}{\color{blue}
\begin{verbatim}
polarform(%);
\end{verbatim}
}
\end{minipage}
\begin{math}\displaystyle
\parbox{8ex}{\color{labelcolor}(\%o22) }
{2}^{\frac{3}{2}}\,{e}^{-\frac{3\,i\,\pi }{4}}
\end{math}
\end{maxcode}

A la hora de calcular las ra\'ices $n-$\'esimas de un complejo dado en forma polar,
nos encontramos con un problema, pues Maxima no dispone de un comando para ello. Sin
embargo, no es dif\'icil construir una funci\'on que lo haga (v\'eanse las actividades
propuestas a continuaci\'on).

\section{Actividades de c\'omputo}

\begin{enumerate}[(A)]
\item Utilizando Maxima, comprobar que las potencias sucesivas de $z=|z|\, e^{\mathbf{i}\theta}$ est\'an situadas
sobre la curva parametrizada por $(|z|^t\cos (t\theta ),|z|^t\sin (t\theta))$ (una \emph{espiral de
Arqu\'imedes})\index{Arqu\'imedes!espiral}. Obs\'ervese que si $|z|>1$ la espiral se aleja del origen, mientras que si $|z|<1$ se acerca
a \'el (el caso $|z|=1$ conduce a la circunferencia unidad).
\item Construir una funci\'on en Maxima que tome como entrada un n\'umero complejo $z$ (en forma rectangular
o en forma polar) y un n\'umero natural $n$, y devuelva las ra\'ices $n-$\'esimas de $z$ en forma polar.\index{Rai@Ra\'iz!$n$-\'esima}
\item Utilizando la funci\'on creada en el \'item precedente, comprobar que las ra\'ices $n-$\'esimas de $z=1$ 
est\'an dispuestas sobre los v\'ertices de un pol\'igono regular de $n$ lados inscrito en la circunferencia unidad.

\end{enumerate}

\chapter{Polinomios y ra\'ices}
\vspace{2.5cm}\setcounter{footnote}{0}
\section{Nociones b\'asicas}

Los polinomios\index{Polinomio} son unas funciones esenciales en las Matem\'aticas\index{Funci\'on!polin\'omica}.
Como una muestra de su importancia,
mencionaremos que toda funci\'on suficientemente derivable, puede aproximarse por su polinomio de
Taylor al orden que se desee; esto es, puede sustitu\'irse una funci\'on $f(x)$, por complicada que
sea, por un polinomio de forma que la diferencia con la funci\'on original sea tan peque\~na como
queramos\footnote{El c\'alculo de este polinomio se estudia en los cursos de c\'alculo diferencial\index{C\'alculo diferencial}.}.
En la figura siguiente, se muestra la aproximaci\'on de la funci\'on $\sin (x)$ por un polinomio de
grado $3$, $P_3(x)=x-x^3/6$, alrededor del punto $x=0$.
\begin{center}
\begin{tikzpicture}[domain=-3:3]
\draw[] (-3.2,0) -- (3.2,0);
\draw[] (0,-1.25) -- (0,1.25);
\draw[color=red]   plot (\x,{sin(\x r)}) node[right] {$\sin (x)$};
\draw[color=orange]   plot (\x,{\x-\x^3/6}) node[right,below] {$P_3 (x)$};
\end{tikzpicture}
\end{center}

Otro hecho determinante para la importancia de los polinomios es que forman un conjunto de \emph{funciones}
que se comportan de manera muy parecida a como lo hacen los \emph{n\'umeros enteros}. Los polinomios se pueden sumar
y multiplicar entre s\'i y siguen resultando polinomios. Con respecto a la divisi\'on, no todo polinomio
tiene inverso para el producto (s\'olo lo tienen los polinomios constantes, de la misma forma que s\'olo el $1$ tiene inverso en
$\mathbb{Z}$) y, en general, el cociente de dos polinomios no es un polinomio, sino lo que se conoce 
como una  \emph{funci\'on racional}\index{Funci\'on!racional}
(de la misma forma que el cociente de dos enteros es un racional). Adem\'as,
muchas otras propiedades de los n\'umeros se siguen cumpliendo en el caso de los polinomios.\index{Num@N\'umero!polinomios}
En particular, se puede calcular el cociente y el resto de la divisi\'on de dos polinomios, as\'i como su
m\'inimo com\'un m\'ultiplo y su m\'aximo com\'un divisor con los mismos algoritmos (por ejemplo, el de Euclides).\index{Euclides!algoritmo}\index{Operaciones!con polinomios}
Por ejemplo, para dividir\index{Divisibilidad!polinomios} $5x^3-x^2+6$ entre $x-4$, 
proceder\'iamos como sigue:
\begin{center}
\begin{tabular}{lllll|l}
$5x^3$ & $-x^2$ & $+0$ & $+6$ &   & $x-4$  \\ \cline{6-1}
$\textcolor{red}{5x^3}$ & $\textcolor{red}{-20x^2}$ & &  &  & $ \textcolor{red}{5x^2}+\textcolor{blue}{19x}+\textcolor{orange}{76}$  \\ \cline{1-4}
 & $19x^2$ & $+0$ & $+6$ &  &   \\ 
 & $\textcolor{blue}{19x^2}$ & $\textcolor{blue}{-76x}$ &  &  &   \\ \cline{2-4}
 &  & $76x$ & $+6$ &  &   \\ 
 &  & $\textcolor{orange}{76x}$ & $\textcolor{orange}{-304}$ &  &   \\ \cline{3-4}
 &  &  & $+310$  &  &  \\ 
\end{tabular} 
\end{center}
\noindent de manera que
$$
5x^3-x^2+6=(5x^2+19x+76)(x-4)+310 \,.
$$

Por \'ultimo, tambi\'en se cumple la important\'isima propiedad de que todo polinomio se puede factorizar\index{Polinomio!factorizaci\'on} como
producto de unos polinomios m\'as sencillos, llamados irreducibles, y esta factorizaci\'on es \'unica, de manera
an\'aloga a como todo n\'umero entero se descompone en producto de n\'umeros primos\index{Factores primos}. Esta propiedad se resume diciendo
que el espacio de los polinomios con coeficientes reales, $\mathbb{R}[x]$, es un \emph{dominio de factorizaci\'on \'unica}.
Tambi\'en es un dominio de factorizaci\'on \'unica $\mathbb{Z}[x]$,\index{Zx@$\mathbb{Z}[x]$} el conjunto de todos los polinomios cuyos coeficientes
son n\'umeros enteros, ambos con la suma y producto de polinomios habituales.

Factorizar un polinomio es equivalente a encontrar sus ra\'ices\index{Rai@Ra\'iz!de un polinomio}. Una ra\'iz del polinomio 
$p(x)=a_nx^n+a_{n-1}x^{n-1}+\cdots +a_1x+a_0$ es un valor $\tilde{x}$ tal que anula el polinomio, esto es, $p(\tilde{x})=0$.
Si esto sucede, es $a_0=-(a_n\tilde{x}^n+\cdots +a_1\tilde{x})$ y el polinomio se puede factorizar en la forma
$$
p(x)=q(x)(x-\tilde{x})\, ,
$$
donde el grado del polinomio $q(x)$ es una unidad menor que el de $p(x)$. Si $\tilde{x}$ de nuevo es una ra\'iz
de $q(x)$, se podr\'a escribir $q(x)=r(x)(x-\tilde{x})$ y as\'i
$$
p(x)=r(x)(x-\tilde{x})^2\,.
$$
Se dice en ese caso que $\tilde{x}$ es una ra\'iz de \emph{multiplicidad} $2$, o ra\'iz doble. Las ra\'ices de multiplicidad
$k\geq 2$ se definen an\'alogamente. El 
\emph{Teorema fundamental del \'algebra}\index{Teorema!fundamental del \'Algebra} afirma que todo polinomio con
coeficientes reales de grado $n$, tiene $n$ ra\'ices complejas (alguna de las cu\'ales puede ser real, por supuesto),
contada cada una con su multiplicidad\index{Rai@Ra\'iz!de un polinomio!multiplicidad}. Por ejemplo, $x=3$
y $x=-1$ son ra\'ices del polinomio $x^3-5x^2+3x+9$, pero la primera es doble, pues se puede factorizar
\begin{equation}\label{ruffini}
x^3-5x^2+3x+9 =(x-3)^2(x+1)\,.
\end{equation} 
Otra forma de expresar el teorema, consiste en decir que si $x_1$, $x_2$,...,$x_r$ son las ra\'ices de un polinomio
$p(x)$, cada una con multiplicidad $m_1$, $m_2$,...,$m_r$ respectivamente, entonces,
$$
p(x)=(x-x_1)^{m_1}\cdots (x-x_r)^{m_r}\, .
$$
Observemos que si $z$ es una ra\'iz \emph{compleja} del polinomio con coeficientes \emph{reales} $p(x)$,
entonces su conjugado complejo
$\overline{z}$ \emph{tambi\'en} es una ra\'iz (?`por qu\'e?). En otras palabras: las ra\'ices complejas siempre
aparecen por pares, y si un factor $(x-z)$ interviene en la descomposici\'on de $p(x)$, tambi\'en lo hace
$(x-\overline{z})$, de manera que se tendr\'a
$$
p(x)=q(x)(x-z)(x-\overline{z})\,,
$$ 
para un cierto polinomio $q(x)$ cuyo grado es dos unidades menor que el de $p(x)$. Si escribimos $z=a+b\mathbf{i}$,
entonces $(x-z)(x-\overline{z})=(x-a)^2+b^2$ y por tanto cada ra\'iz compleja factoriza el polinomio en la forma
$$
p(x)=q(x)((x-a)^2+b^2)\,.
$$

En general no es sencillo encontrar ra\'ices de polinomios, pero la siguiente observaci\'on obvia 
es \'util a veces: dado un polinomio con coeficientes enteros
$p(x)=a_nx^n+a_{n-1}x^{n-1}+\cdots +a_1x+a_0$, si un polinomio lineal $x-b$ es uno de sus factores, entonces
$b$ debe ser un divisor del t\'ermino independiente $a_0$. As\'i, hubi\'eramos podido llegar a la descomposici\'on
\eqref{ruffini} observando que los divisores de $9$ son $\pm 1,\pm 3,\pm 9$ y estudiando si $x^3-5x^2+3x+9$ es
divisible de manera exacta por $(x\pm 1)$, $(x\pm 3)$, $(x\pm 9)$, respectivamente. \\
Como ejemplo, consideremos el polinomio $x^3+6x-20$. Como es de grado $3$, debe tener $3$ ra\'ices complejas, pero
como las que no son reales aparecen por pares, al menos una de las ra\'ices es real (en otro caso tendr\'iamos
\emph{cuatro} ra\'ices complejas, dos pares de conjugadas). Podemos entonces buscar un factor en la forma $(x-b)$, con
$b$ real que, adem\'as, debe estar entre los divisores de $20$: $\pm 1$, $\pm 2$, $\pm 4$, $\pm 5$. Se comprueba
r\'apidamente que $(x\pm 1)$ no dividen al polinomio, pero en el siguiente paso encontramos que
$$
x^3+6x-20 = (x-2)(x^2+2x+10)\,,
$$
de modo que $x=2$ es una ra\'iz. Las otras dos podr\'ian ser reales o un par de complejas conjugadas. Pero como
el factor cuadr\'atico es 
$$
x^2+2x+10=(x+1)^2+9\,,
$$
que tiene la forma $(x-a)^2+b^2$ para $a=-1$, $b=3$, conclu\'imos que
las restantes ra\'ices son $a\pm b\mathbf{i}=-1\pm 3\mathbf{i}$.

Veamos ahora otra situaci\'on en la que el uso de polinomios resulta esencial. Supongamos
que tenemos un conjunto de $N$ datos de la forma $(x_i,y_i)$, donde $i\in\{0,\ldots ,N-1\}$.
Por ejemplo, para $N=3$ podr\'iamos tener $(x_0,y_0)=(1,2{.}2)$, $(x_1,y_1)=(3{.}1,0{.}5)$ 
y $(x_2,y_2)=(4,-1)$. Si
queremos ajustarlos con una curva suave que pase por ellos, debemos buscar algo m\'as complejo 
que una recta (en general, los puntos \emph{no} estar\'an alineados), pero lo que sabemos
acerca de las rectas nos puede guiar para resolver este problema. 
En el cap\'itulo \ref{cap8} se vi\'o que dados dos puntos
$(x_0,y_0)$ y $(x_1,y_1)$ (es decir, $N=2$), la recta que pasa por ellos tiene la ecuaci\'on
$$
\frac{x-x_0}{x_1-x_0}=\frac{y-y_0}{y_1-y_0}\,,
$$
que se puede reescribir, despejando la $y$, como
$$
y=L(x)=\frac{x-x_1}{x_0-x_1}y_0 + \frac{x-x_0}{x_1-x_0}y_1\,.
$$
Deteng\'amonos un momento a analizar la estructura del polinomio $L(x)$. Claramente \'esta es
$$
L(x)=\phi_0(x)y_0 +\phi_1(x)y_1\,,
$$
donde $\phi_i (x)$, $i\in\{0,1\}$, es un polinomio en $x$ que vale $1$ para $x=x_i$ y se
anula en $x_j$ si $j\neq i$. Por ejemplo,
$$
\phi_1(x)=\frac{x-x_0}{x_1-x_0}\,,
$$
y se tiene $\phi_1(x_1)=1$, mientras que $\phi_1(x_0)=0$. La idea de Lagrange\index{Lagrange, Joseph Louis} consisti\'o
en generalizar estas propiedades para $N>2$, y as\'i, dados los $N$ datos
$(x_0,y_0)$,..., $(x_{N-1},y_{N-1})$, se define su \emph{polinomio interpolador}\index{Polinomio!interpolador} como
el polinomio de grado $N-1$ dado por
$$
y=L(x)=\phi_0 (x)y_0 +\cdots +\phi_{N-1}(x)y_{N-1}\,,
$$
donde
$$
\phi_i (x)=\frac{(x-x_0)(x-x_1)\cdots (x-x_{i-1})(x-x_{i+1})\cdots (x-x_{N-1})}
{(x_i-x_0)(x_i-x_1)\cdots (x_i-x_{i-1})(x_i-x_{i+1})\cdots (x_i-x_{N-1})}\,.
$$
N\'otese que, para $N=2$, el polinomio interpolador proporciona exactamente la 
recta entre dos puntos que ya conocemos.
Por construcci\'on, es obvio que para cada $j\in\{0,\ldots ,N-1\}$ es
$$
L(x_j)=y_j\,,
$$
que es la propiedad deseada. Para el ejemplo que hab\'iamos propuesto, con $N=3$, resulta
$$
L(x)=\frac{(x-3{.}1)(x-4)}{(1-3{.}1)(1-4)}2{.}2+
\frac{(x-1)(x-4)}{(3{.}1-1)(3{.}1-4)}0{.}5+
\frac{(x-1)(x-3{.}1)}{(4-1)(4-3{.}1)}(-1)
=-\frac{1}{105}(30x^2-38x-223)\,.
$$ 
\noindent La siguiente figura (obtenida con Maxima) ilustra los puntos en el plano junto con el polinomio interpolador. Obs\'ervese como, efectivamente, el polinomio pasa por los puntos dados,
proporcionando un ajuste suave para ellos.
\begin{center}
\includegraphics[scale=0.66]{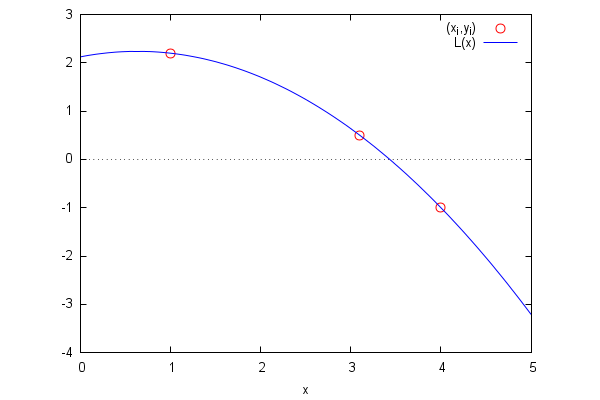} 
\end{center}

\section{Aplicaciones}

\subsection{Distancias de frenado}
El modelado m\'as simple del movimiento rectil\'ineo hace uso de funciones polinomiales.\index{Funci\'on!polin\'omica}
En este apartado, vamos a ver c\'omo calcular la distancia que necesita un autom\'ovil
para frenar totalmente, bajo distintas condiciones relacionadas con el conductor, el propio 
auto y la v\'ia por la que se transita. Utilizaremos para ello datos de la  
polic\'ia federal 
australiana (v\'ease la p\'agina web mencionada
en la actividad de c\'omputo \eqref{activA} de este cap\'itulo), que muestran distancias de frenado para piso seco y para piso mojado.
A partir de los datos para camino seco, experimentaremos ajustando con GeoGebra modelos 
polinomiales para la relaci\'on entre la velocidad del auto en 
el momento de frenar y la distancia recorrida hasta que se detiene totalmente. Analizaremos el caso de polinomios de primer y de segundo grado y compararemos los resultados.

Ajustar\index{Ajuste} un polinomio $p(x)$ a un conjunto de puntos 
$\{(x_1,y_1),(x_2,y_2), \ldots (x_n,y_n)\}$
consiste en determinar los coeficientes de $p(x)$ que hacen m\'inima
la suma de los cuadrados de los errores (es decir, las desviaciones del valor
del polinomio con respecto a los puntos que se tienen):
\[
  E = \sum_{i=1}^n (y_i -p(x_i))^2\,.
\]
El caso m\'as simple de ajuste lo proporciona un polinomio lineal $p(x) = ax$\index{Funci\'on!lineal}. Puede 
comprobarse como ejercicio que en este caso la suma de los cuadrados de los errores se
puede expresar como
\[
  E(a)=a^2\sum_{i=1}^n x_i^2-2a\sum_{i=1}^n x_i y_i+\sum_{i=1}^n y^2_i \doteq Aa^2+Ba+C
\]
(el s\'imbolo $\doteq$ quiere decir que lo que est\'a a su derecha est\'a definido por la ecuaci\'on), de manera que s\'olo depende de la variable $a$, el coeficiente del polinomio.
Completando el cuadrado, se ve que el polinomio $E(a)$ toma su valor
m\'inimo en el v\'ertice $a = -\frac{B}{2A}$. Para los datos particulares recogidos
en la siguiente tabla,
\begin{center}
  \begin{tabular}{c|c}
    velocidad (Km/h) & distancia (m) \\
    \hline 
	\rule{0pt}{3ex}
    50               & 35 \\
    55               & 40 \\
    60               & 45 \\
    65               & 50 \\
    70               & 55 \\
    75               & 65 \\
    80               & 70 \\
    \hline 
    \end{tabular}
\end{center}
\noindent una vez que se convierte la velocidad a $m/s$, se obtiene el valor $a = 2{.}881$.

\begin{center}
\includegraphics[scale=0.4]{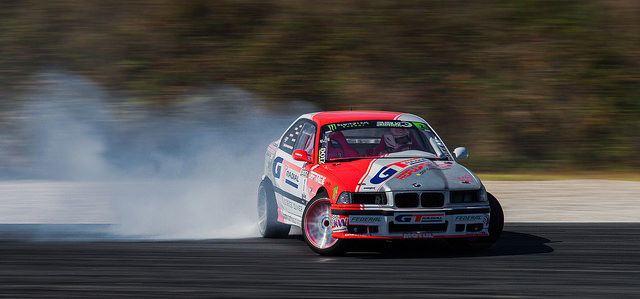} 
\end{center}

Para ajustar polinomios de grado $k$ se usa la misma idea\footnote{Sin embargo, el c\'alculo preciso de los coeficientes del polinomio de ajuste requiere t\'ecnicas de c\'alculo
diferencial, lo cual deber\'ia ser, una vez m\'as, una motivaci\'on para su estudio.}.
Aqu\'i veremos c\'omo obtener el polinomio de ajuste dentro de GeoGebra, usando las \'ordenes
\texttt{AjustePolin\'omico}\index{GeoGebra!ajuste polin\'omico} (que toma como argumentos una lista de puntos y el grado del polinomio)
y \texttt{SumaErroresCuadrados} (que toma como argumentos una lista de puntos y la funci\'on con
la que se quieren ajustar).

Como primer paso capturamos los datos en la vista de \texttt{Hoja de C\'alculo}\index{GeoGebra!hoja de c\'alculo}, convirtiendo la
velocidad de Km/h a m/s. Marcando los datos y picando con el bot\'on derecho
del rat\'on, se selecciona del men\'u contextual la opci\'on \texttt{Crea -> Lista de puntos}.
A continuaci\'on ajustaremos un polinomio de primer grado a nuestros datos.
Aunque el polinomio $p(x)$ obtenido aparentemente reproduce bien los datos,
su intersecci\'on con el eje de abcisas indica que un
auto viajando a $5{.}9$m/s (aproxim\'adamente $21$Km/h)
se detendr\'ia instant\'aneamente, lo cual es f\'isicamente imposible y, por tanto,
invalida este modelo.

\begin{center}
\includegraphics[scale=0.33]{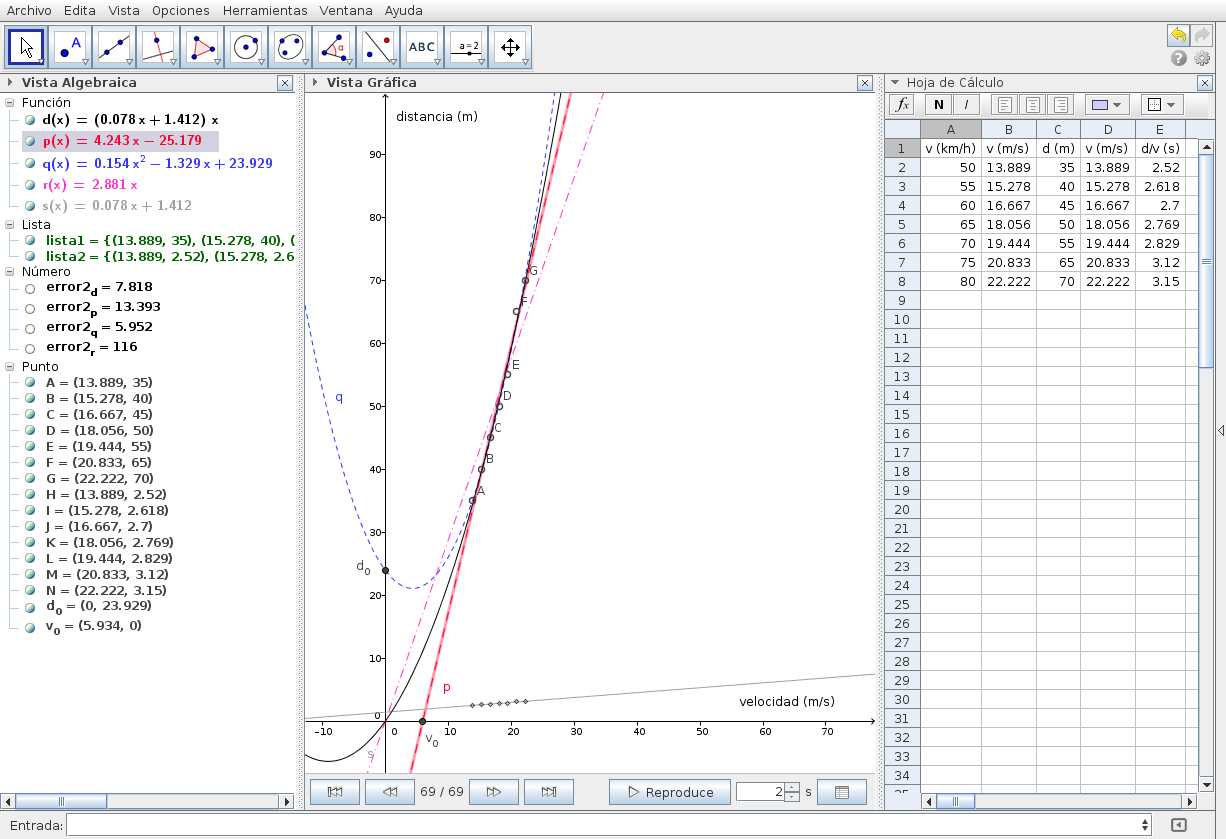} 
\end{center}

Probemos ahora a realizar el ajuste con un polinomio de segundo grado.
El polinomio $q(x)$ que se obtiene ajusta los datos con
menor error, sin embargo, su intersecci\'on con el eje $y$ indica que
un auto estacionado tendr\'ia
que recorrer $23{.}9$ m para detenerse (!) Es claro que el problema
reside en que debemos forzar a que los
polinomios pasen por el punto $(0,0)$, lo que obliga a que carezcan de
t\'ermino constante. Esta propiedad la cumple el
polinomio lineal $r(x) = 2.881x$ calculado previamente, pero ya hemos visto que
da un ajuste muy pobre a los datos (con un error de un orden de magnitud mayor
que los obtenidos). Aunque GeoGebra no tiene opciones para ajustar polinomios
sin t\'ermino constante, podemos utilizar el modelo
\[
  \frac{d}{v} = av + b, \quad v > 0
\]   
que es equivalente al modelo dado por la par\'abola\index{Par\'abola} $d = av^2 + bv$.
Adem\'as de proporcionar un error casi \'optimo, este \'ultimo modelo explica
la distancia de frenado como la
suma de la distancia que recorre el auto mientras que el conductor reacciona al evento, y la
distancia que recorre el auto una vez que se ha pisado el freno. El coeficiente 
$b = 1{.}412$ corresponde al tiempo promedio de reacci\'on de un conductor,
medido en segundos (v\'ease, por ejemplo,
\href{http://www.visualexpert.com/Resources/reactiontime.html}{\tt www.visualexpert.com}).
Por otra parte, el coeficiente $a = 0{.}078$ est\'a relacionado con el coeficiente
de fricci\'on $\mu$ entre las llantas y el suelo; expl\'icitamente $a = \frac{1}{2\mu g}$, 
con $g = 9.81$m/s$^2$ la aceleraci\'on de la gravedad. En nuestro ejemplo, despejando 
obtenemos un valor de $\mu = 0.653$, que concuerda bastante bien con los valores t\'ipicos 
para goma y asfalto seco (v\'ease \href{http://www.engineeringtoolbox.com/friction-coefficients-d\_778.html}{\tt www.engineeringtoolbox.com}).

\subsection{C\'odigos de Reed-Solomon}
Pasemos ahora a otra aplicaci\'on de los modelos polinomiales.
Como ya se ha visto, dado un conjunto de $n$ puntos distintos $\{(x_1,y_1),(x_2,y_2), \ldots (x_n,y_n)\}$, es posible ajustarlos mediante el polinomio interpolador de Lagrange, que
tiene la propiedad de pasar por todos ellos.
Para interpolar en GeoGebra $n$ puntos distintos se pueden usar las \'ordenes equivalentes
\texttt{AjustePolin\'omico}, que acabamos de ver, o \texttt{Polinomio}, que b\'asicamente llama a
\texttt{AjustePolin\'omico} con el segundo argumento igual a $n-1$, siendo $n$ la longitud de
la lista de puntos a interpolar.

La propiedad de los polinomios de poder pasar por una serie de puntos distintos tiene 
m\'ultiples aplicaciones, entre las que se encuentra corregir los errores que se producen
en la lectura de informaci\'on codificada (como fotograf\'ias digitalizadas, audio en un
CD o los c\'odigos QR\index{Co@C\'odigo!QR}). Uno de los m\'etodos m\'as usados para la correcci\'on de 
c\'odigos  corruptos es el de
Reed--Solomon, cuyo funcionamiento b\'asico describimos a continuaci\'on.\\
Supongamos que se quieren transmitir los cuatro n\'umeros $1{.}2,\ -3{.}2,\ -5{.}4,\ -1{.}1$
en este orden. Esto constituye un \emph{c\'odigo}, un mensaje que se desea hacer llegar a un 
receptor.\index{Co@C\'odigo!Reed-Solomon} 
Para asegurar la integridad del c\'odigo se a\~nade informaci\'on redundante a 
los datos originales, que en este caso ser\'a una secuencia de la misma longitud que la 
original\footnote{En general, el n\'umero de datos adicionales depende de la probabilidad
de error del sistema en cuesti\'on, mientras m\'as alta sea \'esta, m\'as informaci\'on redundante se requiere.}. Comenzamos por formar la serie de tiempo correspondiente
a nuestros datos originales, obteniendo la lista:
\begin{equation}\label{reedsol1}
S_e=[(1,1{.}2),(2,-3{.}2),(3,-5{.}4),(4,-1{.}1)]\,.
\end{equation}
A continuaci\'on, utilizando cualquiera de los dos comandos mencionados,
calculamos el polinomio interpolador,\index{Polinomio!interpolador}
\[
L_e(x) = 0{.}717x^3 - 3{.}2x^2 + 0{.}183x + 3{.}5\,.
\]
A partir de este polinomio, extrapolamos los valores $L_e(5),L_e(6), L_e(7),L_e(8)$, 
resultando la secuencia extendida:
$$
1{.}2,\ -3{.}2,\ -5{.}4,\ -1{.}1,\ 14,\ 44{.}2,\ 93{.}8,\ 167{.}1\,.
$$
Los ocho n\'umeros anteriores se transmiten en el orden en que aparecen en la lista
precedente. Supongamos que a lo m\'as puede
haber un 25\% de n\'umeros corruptos, y consideremos que la informaci\'on transmitida
est\'a dada por el c\'odigo
$$
C= [1{.}2,\ \color{red}{3{.}2},\ \color{black}{-5{.}4},\ -1{.}1,\ \color{red}{12{.}8},\ 
\color{black}{44.2},\ 93{.}8,\ 167{.}1] \,.
$$
En este ejemplo, pues, se da el m\'aximo n\'umero posible de errores,
en el segundo y el quinto lugar.
Como se est\'a duplicando el n\'umero de datos y el grado del polinomio utilizado
para calcular los $n$ \'ultimos es precisamente $n-1$, el receptor
puede determinar que en este caso (en que el n\'umero \emph{total} de datos trasmitidos
es $2n=8$) se utiliz\'o un polinomio c\'ubico para interpolar
los valores originales y extrapolar los datos adicionales. Con esta informaci\'on el
receptor puede determinar si ha habido o no errores en la transmisi\'on. Para ello,
forma la serie temporal correspondiente al c\'odigo recibido $C$, que en este caso ser\'ia
$$
S_r=
[(1,1{.}2),(2,3{.}2),(3,-5{.}4),(4,-1{.}1),(5,12{.}8),(6,44{.}2),(7,93{.}8),(8,167{.}1)]\,.
$$
Si no hubiese habido errores en la transmisi\'on, los cuatro primeros elementos
de esta lista deber\'ian ser los mismos que los que aparecen en la
serie temporal del emisor \eqref{reedsol1}. Sabemos que ha habido errores
porque el segundo elemento est\'a corrupto pero, ?`c\'omo puede saberlo el receptor
del c\'odigo $C$, si el no dispone de \eqref{reedsol1}? Lo que hace es calcular el
polinomio interpolador correspondiente a la lista $S_r$.
Si resulta un polinomio de
grado $3$, es que no hubo errores de transmisi\'on, pues en ese caso los \'ultimos
cuatro elementos de la lista son redundantes, no aportan informaci\'on nueva y el
polinomio que se obtiene es el mismo que resultar\'ia s\'olo de considerar los
cuatro primeros elementos del c\'odigo. Si al calcular el polinomio interpolador
de $S_r$ resulta que su grado es mayor a $3$, podemos asegurar que ha habido errores, la 
figura muestra el proceso en
GeoGebra.

\begin{center}
\includegraphics[scale=0.35]{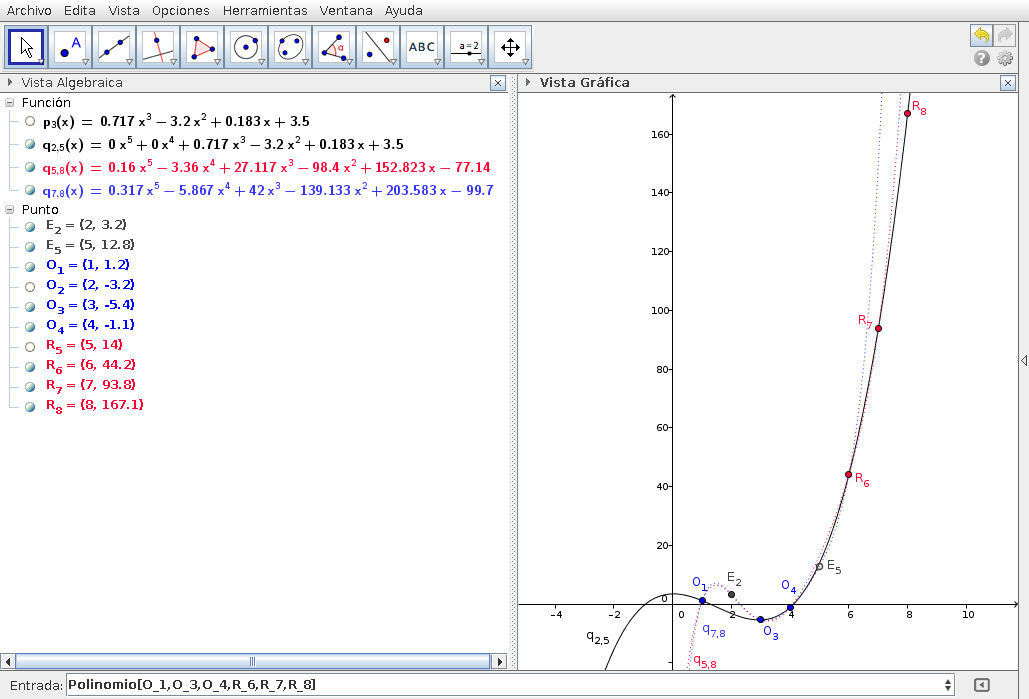} 
\end{center}

No s\'olo se puede determinar la existencia de errores en la transmisi\'on; el c\'odigo
de Reed--Solomon es \emph{corrector de errores}, es decir, puede determinar en qu\'e
posici\'on se ha producido el error y cu\'al es el valor correcto. Para ello, suponiendo
que al calcular el polinomio intepolante de $S_r$ su grado result\'o mayor a $3$, 
se omite uno a uno cada punto de la secuencia recibida y se repite el c\'alculo del
polinomio interpolador: si al omitir el $i-$\'esimo punto el polinomio obtenido,
digamos $q(x)$, es de grado 3, el $i-$\'esimo dato es el \'unico corrupto, y su valor corregido es $q(i)$. Si este procedimiento no conduce a un polinomio de grado $3$, es
que hay m\'as de un dato corrupto. Se procede entonces a eliminar los puntos dos a dos
repitiendo el c\'alculo hasta llegar a un polinomio interpolador de grado $3$.
En nuestro ejemplo, se omiten dos a dos
los puntos de la secuencia recibida: al omitir el segundo y el quinto punto
se obtiene el polinomio
interpolador c\'ubico $q_{2,5}(x) = 0{.}717x^3 - 3{.}2x^2 + 0{.}183x + 3{.}5$,
por lo que los \'unicos datos
corruptos son el segundo y el quinto, siendo su valores corregidos
$q_{2,5}(2) = -3{.}2$ y $q_{2,5}(5) = 14$, como cab\'ia esperar.

\section{Ejercicios}

\begin{enumerate}
\setcounter{enumi}{128}
\item
Realizar las siguientes sumas de polinomios:
\begin{multicols}{2}
\begin{enumerate}[(a)]
\item $(3+2x-5x^2+7x^3) + (-2+x-2x^2)$
\item $(x^5-x^4+x^3-x^2-1) - (-2-5x+4x^3)$
\item $(\sqrt{2}-x+\sqrt{3}x^2) + (\sqrt{3}-3x-2\sqrt{3}x^2)$
\item $(2-x)^4-(x^2-3)^3$
\item $(\sqrt[3]{4}-2x)^3 - (2x^5-3x^7+x^8)$
\item $(\pi + \pi ^2x +\pi ^3x^2+\pi ^4x^3) - (\pi -x)^6$
\end{enumerate}
\end{multicols}

\item
Realizar los siguientes productos de polinomios:
\begin{multicols}{2}
 \begin{enumerate}[(a)]
\item $(2+3x-4x^2)(-1-7x)$
\item $(3-2x+5x^2)(-1+3x^2-8x^4)$
\item $(x-a)(x-b)(x-c)$, $a$, $b$ y $c$ constantes.
\item $(x+\sqrt{2}x^3)(\sqrt{2}-3x^2)$
\item $(\pi x -e^2x^2 +\pi ^5x^4)(1-\sqrt{2}x)$
\end{enumerate}
\end{multicols}

\item
Calcular el cociente y el resto de las siguientes divisiones de polinomios:
\begin{multicols}{2}
 \begin{enumerate}[(a)]
\item $(2-3x+4x^2)/(2)$
\item $(3+3x-5x^2)/(1-6x)$
\item $(5-3x+x^2)/(4-x+2x^2)$
\item $(8-2x+3x^2)/(4+2x-5x^2+3x^3)$
\item $(\sqrt{3}-2x+\sqrt{2}x^2-4x^3)/(3+2x-x^2)$
\item $(2+\sqrt{5}x-\pi ^2x^2+\frac{7}{3}x^3)/(\sqrt{5}-x)$
\end{enumerate}
\end{multicols}

\item
Comprobar que, en las siguientes parejas de polinomios, uno es divisor de otro:
\begin{enumerate}[(a)]
\item $(x^2-4x+4)$ y $(2-x)$
\item $(x-3)$ y $(x^3-x^2-5x-3)$
\item $(x-\pi )$ y $(x^3 + \pi ^2x-\pi x^2-\pi ^3)$
\item $(2 + (2\sqrt{2}+\sqrt{6})x + (2\sqrt{3}-\sqrt{6})x^2 + (5\sqrt{2}-3)x^3 + 5\sqrt{3}x^4)$ y $(\sqrt{2} + \sqrt{3}x)$
\end{enumerate}

\item
Calcular el m\'aximo com\'un divisor de las siguientes parejas de polinomios:
\begin{multicols}{2}
 \begin{enumerate}[(a)]
\item $-2+3x+2x^2+2x^3+x^4$  y\\  $4+4x-7x^2-4x^3+2x^4+x^5$
\item $-1+x+x^3+x^4$ y\\ $x-x^2-2x^3+x^4+x^5$
\item $-\sqrt{2} + (1-\sqrt{2})x +x^2$ y $-2+x^2$
\item $\sqrt{8}x + (\sqrt{2}-2\pi )x^2 + (2-\pi )x^3 + x^4$ y\\ $2x+3x^2+x^3$
\item $1-x+2x^2+x^3$ y $-2+3x+x^2$
\end{enumerate}
\end{multicols}

\item
Calcular el m\'aximo com\'un divisor de los siguientes polinomios:
$$
1-x-x^2+x^3,\; 2x+2x^2+x^3+x^4 ,\; -1+4x+3x^2-x^3+x^4.
$$

\item
Calcular el m\'{\i}nimo com\'un m\'ultiplo de:
$$
-2-2x+x^2+x^3 ,\; 1+2x+x^2.
$$

\item
Calcular el valor de los siguientes polinomios en el punto que se indica:
\begin{multicols}{2}
 \begin{enumerate}[(a)]
\item $p_1(x)= x-3x^2+5x^{16}$, en $x=0$
\item $p_2(x)= -2-2x+x^2+x^3$, en $x=\sqrt{2}$
\item $p_3(x)= \pi ^3 + 3\pi ^2x + 3\pi x^2 + x^3$, en $x=1-\pi $
\item $p_4(x)= 3-5x+4x^3$, en $x=-2$
\end{enumerate}
\end{multicols}

\item
Demostrar que los polinomios en los que s\'olo aparecen potencias pares de la variable verifican:
$$
p(-x) = p(x).
$$
\noindent Recordemos que a las funciones que cumplen esta propiedad se las denomina \emph{funciones pares}. 
Hacer lo mismo para 
los polinomios en los que s\'olo aparecen potencias impares de la variable, que  verifican:
$$
p(-x) = -p(x),
$$
\noindent (es decir, son \emph{funciones impares}).

\item
Dibujar la gr\'afica de los siguientes monomios:
\begin{enumerate}[(a)]
\item $x^n$, donde $n=$0, 2, 4, 6, \dots
\item $x^n$, donde $n=$1, 3, 5, 7, \dots
\end{enumerate}

\item (Recu\'erdese el ejercicio \ref{derivadas} y que la derivada es \emph{lineal}). Calcular la derivada de los polinomios siguientes:
\begin{enumerate}[(a)]
\item $p(x) = 2-3x+5x^2-7x^8$
\item $q(x) = 4-16x^{10}$
\item $r(x) = -3b + (2-a)x - (4+5b)x^2$, con $a,b\in\mathbb{R}$ constantes.
\end{enumerate}

\item  Calcular todas las derivadas de los polinomios siguientes y evaluarlas en $x=0$:
\begin{enumerate}[(a)]
\item $p(x) = 4 + 3x - 2x^2 + 7x^3$
\item $q(x) = x^6$
\item $r(x) = e^{\pi/3} + 4e^{2\pi/5} x + 5e^{-\pi/2} x^5$
\end{enumerate}

\item
Demostrar la regla de derivaci\'on de un producto de tres polinomios:
\begin{displaymath}
(a(x)b(x)c(x))' = a'(x)b(x)c(x) + a(x)b'(x)c(x) + a(x)b(x)c'(x).
\end{displaymath}

\item
Si dos polinomios tienen la misma derivada, ?`Son esos dos polinomios iguales?

\item
?`Existe alg\'un polinomio que cumpla $p'(x) = p(x)$?

\item
Comprobar el teorema del resto en la divisi\'on de $2x^2-7x+5$ entre $x-3$.

\item
Hallar un polinomio de grado 3 tal que al dividirlo por $x-1$, $x-2$, $x+5$ d\'e, respectivamente, restos 2, 1, 8.

\item
Comprobar que $x^4-3x^3-3x^2+11x-6$ es divisible entre los polinomios lineales $x-1$, $x+2$ y $x-3$.
 
\item
Realizar las siguientes divisiones:
\begin{multicols}{2}
 \begin{enumerate}[(a)]
\item $(5x^4-3x^3+x^2-2x+6)/(x-1)$
\item $(x^4-2x^2+12)/(x+3)$
\item $(x^2-\frac{2}{3}x+1)/(x-\frac{1}{5})$
\item $((2-i)x^2-ix+(6+i))/(x-2i)$
\end{enumerate}
\end{multicols}

\item
Sabiendo que $-1$ y $2$ son ra\'{\i}ces de $x^5 - x^4 - 5x^3 + x^2 + 8x + 4$, hallar la multiplicidad de cada una.

\item
Hallar un polinomio cuyas ra\'{\i}ces sean las que se indican en cada caso:
\begin{enumerate}[(a)]
\item $x_1=1$, $x_2=2$ (doble) y $x_3=-2$
\item $x_1=2$, $x_2=\sqrt{5}$, $x_3=\pi$
\end{enumerate}

\item
Demostrar que calcular las ra\'{\i}ces de un polinomio de segundo grado es equivalente al problema de hallar dos n\'umeros cuya suma y producto sean conocidos.  
Determinar la relaci\'on entre la suma y producto conocidos y el polinomio de segundo grado correspondiente (los babilonios\index{Babilonios} resolv\'{\i}an as\'{\i} las ecuaciones de 
segundo grado hace cuatro mil a\~nos: transform\'andolas en el problema de dos n\'umeros con suma y producto conocidos).

\item
Hallar todas las ra\'{\i}ces (si es posible por m\'etodos exactos) de los siguientes polinomios: 
\begin{enumerate}[(a)]
\item $\pi x + 3$
\item $x^2+(1-\sqrt{2})x-\sqrt{2}$
\item $x^2-4x+5$
\item $3x^3+2x^2-x$
\item $x^3-2x-4$
\item $2x^3-x^2+3x-1$
\item $x^4 + 2(2-\pi)x^3 + (4-8\pi+\pi^2)x^2 + 4\pi(\pi-2)x + 4\pi^2$
\end{enumerate}

\item
Hallar las ra\'{\i}ces de $p(x)= 3x^4 + 15x^3 + 18x^2 - 12x -24$, indicando su multiplicidad, 
y comprobar que son ra\'{\i}ces del polinomio derivada, $p'(x)$, con multiplicidad una unidad menor.

\item
Sabiendo que un polinomio de coeficientes reales es de grado 5, y que $(2+i)$ y $(1-3i)$ son ra\'{\i}ces del mismo, ?`cu\'antas ra\'{\i}ces reales tiene?

\item
Hallar el polinomio m\'onico de menor grado de coeficientes reales que tenga a $3e^{i\pi/3}$ y $4e^{i2\pi/5}$ por ra\'{\i}ces.

\item
De las siguientes fracciones racionales, indicar cu\'ales son propias y cu\'ales impropias, y en  este \'ultimo caso, reducirlas a suma de polinomio m\'as fracci\'on propia:
\begin{multicols}{2} 
\begin{enumerate}[(a)]
\item $$\frac{2x^2-3x+4}{x^3-5x}$$
\item $$\frac{-5x^3+2x^2}{x^2-x+2}$$
\item $$\frac{3x^4-2x^3+4x-1}{5x^4+2x^3-3x^2+x+1}$$
\item $$\frac{x^2+3x-5}{2x^4}$$
\item $$\frac{(1+i)x^2+3ix-2}{x-2i}$$
\end{enumerate}
\end{multicols}

\item
Descomponer en suma de polinomio y fracciones simples, las siguientes fracciones racionales: 
\begin{multicols}{2}
\begin{enumerate}[(a)]
\item $$\frac{2x+3}{x^2+x-2}$$
\item $$\frac{x^2-2x+1}{x^3+3x^2-x-3}$$
\item $$\frac{3x-2}{2x-3}$$
\item $$\frac{-3x^2+4x-1}{2x^2-2x-6}$$
\item $$\frac{x^2+x+1}{x^3-x^2-x+1}$$
\item $$\frac{3x-2}{x^4-2x^3-11x^2+12x+36}$$
\item $$\frac{(2-i)x+i}{x^2-2ix+3}$$
\item $$\frac{2x^2-3i}{x^2-2ix-1}$$
\end{enumerate}
\end{multicols}

\end{enumerate}

\section{Uso de la tecnolog\'ia de c\'omputo}
Los polinomios se definen en Maxima como cualquier otra funci\'on:
\begin{maxcode}

\noindent
\begin{minipage}[t]{8ex}{\color{red}\bf
\begin{verbatim}
(%i1) 
\end{verbatim}
}
\end{minipage}
\begin{minipage}[t]{\textwidth}{\color{blue}
\begin{verbatim}
p(x):=x^4+2*x^3+2*x^2+3*x-2;
\end{verbatim}
}
\end{minipage}
\definecolor{labelcolor}{RGB}{100,0,0}
\begin{math}\displaystyle
\parbox{8ex}{\color{labelcolor}(\%o1) }
\mathrm{p}\left( x\right) :={x}^{4}+2\,{x}^{3}+2\,{x}^{2}+3\,x-2
\end{math}

\noindent
\begin{minipage}[t]{8ex}{\color{red}\bf
\begin{verbatim}
(%i2) 
\end{verbatim}
}
\end{minipage}
\begin{minipage}[t]{\textwidth}{\color{blue}
\begin{verbatim}
q(x):=x^5+2*x^4-4*x^3-7*x^2+4*x+4;
\end{verbatim}
}
\end{minipage}
\definecolor{labelcolor}{RGB}{100,0,0}
\begin{math}\displaystyle
\parbox{8ex}{\color{labelcolor}(\%o2) }
\mathrm{q}\left( x\right) :={x}^{5}+2\,{x}^{4}+\left( -4\right) \,{x}^{3}+\left( -7\right) \,{x}^{2}+4\,x+4
\end{math}
\end{maxcode}

Una vez definidos $p(x)$ y $q(x)$, disponemos de funciones para trabajar con sus
elementos b\'asicos. Por ejemplo, para calcular el coeficiente de $x^n$, tenemos
el comando \texttt{coeff(}\textit{expr}\texttt{,x,n)}:

\begin{maxcode}
\noindent
\begin{minipage}[t]{8ex}{\color{red}\bf
\begin{verbatim}
(%i3) 
\end{verbatim}
}
\end{minipage}
\begin{minipage}[t]{\textwidth}{\color{blue}
\begin{verbatim}
coeff(p(x),x,3);
\end{verbatim}
}
\end{minipage}
\definecolor{labelcolor}{RGB}{100,0,0}
\begin{math}\displaystyle
\parbox{8ex}{\color{labelcolor}(\%o3) }
2
\end{math}

\noindent
\begin{minipage}[t]{8ex}{\color{red}\bf
\begin{verbatim}
(%i4) 
\end{verbatim}
}
\end{minipage}
\begin{minipage}[t]{\textwidth}{\color{blue}
\begin{verbatim}
coeff(q(x),x,3);
\end{verbatim}
}
\end{minipage}
\definecolor{labelcolor}{RGB}{100,0,0}
\begin{math}\displaystyle
\parbox{8ex}{\color{labelcolor}(\%o4) }
-4
\end{math}
\end{maxcode}

El comando \texttt{hipow} porporciona el exponente m\'as alto al que
aparece elevada la inc\'ognita, es decir, el grado del polinomio:
\begin{maxcode}
\noindent
\begin{minipage}[t]{8ex}{\color{red}\bf
\begin{verbatim}
(%i5) 
\end{verbatim}
}
\end{minipage}
\begin{minipage}[t]{\textwidth}{\color{blue}
\begin{verbatim}
hipow(q(x),x);
\end{verbatim}
}
\end{minipage}
\definecolor{labelcolor}{RGB}{100,0,0}
\begin{math}\displaystyle
\parbox{8ex}{\color{labelcolor}(\%o5) }
5
\end{math}
\end{maxcode}
Para calcular ra\'ices de forma exacta, es \'util el comando \texttt{factor},
que calcula un polinomio como producto de factores cuando esto es posible,
\begin{maxcode}
\noindent
\begin{minipage}[t]{8ex}{\color{red}\bf
\begin{verbatim}
(%i6) 
\end{verbatim}
}
\end{minipage}
\begin{minipage}[t]{\textwidth}{\color{blue}
\begin{verbatim}
factor(p(x));
\end{verbatim}
}
\end{minipage}
\definecolor{labelcolor}{RGB}{100,0,0}
\begin{math}\displaystyle
\parbox{8ex}{\color{labelcolor}(\%o6) }
\left( x+2\right) \,\left( {x}^{3}+2\,x-1\right) 
\end{math}
\end{maxcode}

Si lo que se quiere es divir dos polinomios, se puede usar \texttt{divide}.
Toma como argumentos dos polinomios y la inc\'ognita, y devuelve una lista
donde el primer elemento es el cociente y el segundo el resto de la divisi\'on,

\begin{maxcode}
\noindent
\begin{minipage}[t]{8ex}{\color{red}\bf
\begin{verbatim}
(%i7) 
\end{verbatim}
}
\end{minipage}
\begin{minipage}[t]{\textwidth}{\color{blue}
\begin{verbatim}
divide(q(x),p(x),x);
\end{verbatim}
}
\end{minipage}
\definecolor{labelcolor}{RGB}{100,0,0}
\begin{math}\displaystyle
\parbox{8ex}{\color{labelcolor}(\%o7) }
[x,-6\,{x}^{3}-10\,{x}^{2}+6\,x+4]
\end{math}
\end{maxcode}

Podemos comprobar que el resultado es correcto multiplicando
el cociente por el divisor y sumando el resto,
\begin{maxcode}
\noindent
\begin{minipage}[t]{8ex}{\color{red}\bf
\begin{verbatim}
(%i8) 
\end{verbatim}
}
\end{minipage}
\begin{minipage}[t]{\textwidth}{\color{blue}
\begin{verbatim}
ratsimp(p(x)*x-6*x^3-10*x^2+6*x+4);
\end{verbatim}
}
\end{minipage}
\definecolor{labelcolor}{RGB}{100,0,0}
\begin{math}\displaystyle
\parbox{8ex}{\color{labelcolor}(\%o8) }
{x}^{5}+2\,{x}^{4}-4\,{x}^{3}-7\,{x}^{2}+4\,x+4
\end{math}
\end{maxcode}

El c\'alculo del m\'aximo com\'un divisor de dos polinomios se hace sin problemas:

\begin{maxcode}
\noindent
\begin{minipage}[t]{8ex}{\color{red}\bf
\begin{verbatim}
(%i9) 
\end{verbatim}
}
\end{minipage}
\begin{minipage}[t]{\textwidth}{\color{blue}
\begin{verbatim}
gcd(p(x),q(x));
\end{verbatim}
}
\end{minipage}
\definecolor{labelcolor}{RGB}{100,0,0}
\begin{math}\displaystyle
\parbox{8ex}{\color{labelcolor}(\%o9) }
x+2
\end{math}
\end{maxcode}

El c\'alculo del polinomio interpolador de Lagrange de una lista de puntos
como por ejemplo 
$[(1,1{.}2),(2,3{.}2),(3,-5{.}4),(4,-1{.}1),(5,12{.}8),(6,44{.}2),(7,93{.}8),(8,167{.}1)]$,
se puede realizar cargando el paquete \texttt{interpol}, que dispone de la orden
\texttt{lagrange} para ello.

\begin{maxcode}
\noindent
\begin{minipage}[t]{8ex}{\color{red}\bf
\begin{verbatim}
(%i10) 
\end{verbatim}
}
\end{minipage}
\begin{minipage}[t]{\textwidth}{\color{blue}
\begin{verbatim}
load(interpol)$
\end{verbatim}
}
\end{minipage}

\noindent
\begin{minipage}[t]{8ex}{\color{red}\bf
\begin{verbatim}
(%i11) 
\end{verbatim}
}
\end{minipage}
\begin{minipage}[t]{\textwidth}{\color{blue}
\begin{verbatim}
lista:[[1,1.2],[2,3.2],[3,-5.4],[4,-1.1],
[5,12.8],[6,44.2],[7,93.8],[8,167.1]]$
\end{verbatim}
}
\end{minipage}

\noindent
\begin{minipage}[t]{8ex}{\color{red}\bf
\begin{verbatim}
(%i12) 
\end{verbatim}
}
\end{minipage}
\begin{minipage}[t]{\textwidth}{\color{blue}
\begin{verbatim}
ratsimp(lagrange(lista));
\end{verbatim}
}
\end{minipage}
\definecolor{labelcolor}{RGB}{100,0,0}
\begin{math}\displaystyle
\parbox{8ex}{\color{labelcolor}(\%o12) }
\frac{31\,{x}^{7}-1009\,{x}^{6}+13513\,{x}^{5}-95995\,{x}^{4}+389074\,{x}^{3}-886516\,{x}^{2}+1020282\,x-437220}{1800}
\end{math}


\noindent
\begin{minipage}[t]{8ex}{\color{red}\bf
\begin{verbatim}
(%i12) 
\end{verbatim}
}
\end{minipage}
\begin{minipage}[t]{\textwidth}{\color{blue}
\begin{verbatim}
L(x):=''%$
\end{verbatim}
}
\end{minipage}
\end{maxcode}

Podemos comprobar que el polinomio de Lagrange tiene la propiedad de pasar por
los puntos, por ejemplo:

\begin{maxcode}
\noindent
\begin{minipage}[t]{8ex}{\color{red}\bf
\begin{verbatim}
(%i13) 
\end{verbatim}
}
\end{minipage}
\begin{minipage}[t]{\textwidth}{\color{blue}
\begin{verbatim}
L(5);
\end{verbatim}
}
\end{minipage}
\definecolor{labelcolor}{RGB}{100,0,0}
\begin{math}\displaystyle
\parbox{8ex}{\color{labelcolor}(\%o13) }
12{.}8
\end{math}
\end{maxcode}

Por \'ultimo, destaquemos que \texttt{factor} s\'olo trabaja con
factores cuyos coeficientes son n\'umeros enteros. Si se quiere
factorizar un polinomio con factores en los que los coeficientes sean de
la forma $x-(a+b\mathbf{i})$, se puede usar \texttt{gfactor}:

\begin{maxcode}
\noindent
\begin{minipage}[t]{8ex}{\color{red}\bf
\begin{verbatim}
(%i14) 
\end{verbatim}
}
\end{minipage}
\begin{minipage}[t]{\textwidth}{\color{blue}
\begin{verbatim}
gfactor(x^3-2*x-4);
\end{verbatim}
}
\end{minipage}
\definecolor{labelcolor}{RGB}{100,0,0}
\begin{math}\displaystyle
\parbox{8ex}{\color{labelcolor}(\%o14) }
\left( x-2\right) \,\left( x-i+1\right) \,\left( x+i+1\right) 
\end{math}
\end{maxcode}

\section{Actividades de c\'omputo}
\begin{enumerate}[(A)]

\item\label{activA} \sloppy Con base en los datos sobre las distancias de frenado para piso mojado mostrados
en la p\'agina 
\url{http://www.tmr.qld.gov.au/Safety/Driver-guide/Speeding/Stopping-distances/
Stopping-distances-on-wet-and-dry-roads.aspx},
realizar un ajuste en GeoGebra mediante un modelo cuadr\'atico que permita estimar en este caso el coeficiente de fricci\'on $\mu$ entre las llantas y el camino.
Comparar el resultado obtenido con los valores t\'ipicos
para hule y asfalto mojado mostrados en la p\'agina
\href{http://www.engineeringtoolbox.com/friction-coefficients-d\_778.html}{\tt http://www.engineeringtoolbox.com/friction-coefficients-d\_778.html}.

\item La siguiente tabla muestra la altitud de un transbordador espacial
en funci\'on del tiempo transcurrido desde su ingreso a la atm\'osfera terrestre.
   \begin{center}
   \begin{tabular}{c|c}
      tiempo (s) & altitud (m) \\
      \hline
      \rule{0pt}{3ex}
       $7$    & $76\,850$\\
      $11$     & $73\,775$ \\
      $15$      & $72\,500$ \\
      $19$       & $68\,525$ \\
      $23$        & $67\,750$ \\
      $29$         & $66\,650$ \\
      \hline
   \end{tabular}
   \end{center}
   \begin{enumerate}
      \item Ajustar los datos en GeoGebra mediante un modelo cuadr\'atico y utilizarlo
      para predecir cu\'ando llegar\'a el veh\'iculo al mar. ?`Es razonable
      esta predicci\'on?
      \item Realizar ahora un ajuste c\'ubico. ?`Es \'este modelo mejor que el anterior?
      \item Finalmente, realizar un ajuste lineal. ?`Cu\'al de los tres modelos
            tiene un error cuadratico menor? ?`Cu\'al de ellos refleja de mejor forma el
            comportamiento de los datos?
   \end{enumerate}
   
\item Dada la lista de puntos
$$
[(1,0{.}4),(2,-1{.}2),(3,-0{.}7), (4,1.4), (5,3.2), (6,-3)]\,,
$$
construir su polinomio interpolador de Lagrange usando Maxima y representarlo junto
con los datos, comprobando gr\'aficamente que el polinomio efectivamente los interpola.

\item Supongamos que se recibe el c\'odigo siguiente:
$$
[2,3,6,7,11,22,48,100,192,341]\,.
$$
Utilizando la teor\'ia de Reed-Solomon, determinar si se trata de un c\'odigo corrupto
y, en caso afirmativo, corregirlo y dar el c\'odigo original.
\end{enumerate}

\chapter{Sistemas de ecuaciones. Matrices}\label{finalcap}
\vspace{2.5cm}\setcounter{footnote}{0}
\section{Nociones b\'asicas}

Hasta ahora hemos analizado la forma de resolver ecuaciones con varios grados de dificultad,
desde hallar un \emph{n\'umero} $x$ tal que $ax=b$ (con $a$ y $b$ conocidos) hasta encontrar un
\emph{polinomio} $p(x)$ tal que $p(x)q(x)=r(x)$ (con $q(x)$ y $r(x)$ conocidos). Ahora, nos
vamos a ocupar de c\'omo resolver un \emph{sistema de ecuaciones}\index{Ecuaciones!sistema},
es\index{Sistema!de ecuaciones} decir, un problema de la forma
\begin{equation}\label{sistemas1}
\begin{cases}
a_{11}x+a_{12}y = b_1 \\
a_{21}x+a_{22}y= b_2\,,
\end{cases}
\end{equation}
donde $a_{11},a_{12},a_{21},a_{22},b_1,b_2$ son n\'umeros reales\footnote{En realidad, todo lo que vamos a ver sirve igual para el caso de n\'umeros complejos, pero en la exposici\'on supondremos que son reales por simplicidad.} y $x,y$ son las inc\'ognitas a determinar.

Hay much\'isimas situaciones en las que aparecen sistemas como \eqref{sistemas1}. Una de ellas
es el conocido como \emph{problema de la dieta}\index{Problema!dieta}\index{Dieta!problema}:
una persona sometida a dieta solo puede ingerir una cantidad limitada de calorias diarias y
\'estas se tienen que distribuir entre varios grupos de alimentos, como prote\'inas, frutas y verduras, l\'acteos, grasas y carbohidratos. T\'ipicamente, el nutri\'ologo plantear\'a unas
necesidades a cubrir, como por ejemplo: `diariamente se debe ingerir $4$ unidades del primer grupo (prote\'inas), $6$ del segundo (frutas y verduras), $3$ del tercero (l\'acteos), $2$ unidades del cuarto (grasas) y $4$ del quinto (carbohidratos)'. Al paciente se le proporcionar\'a una tabla de equivalencias, conteniendo una lista de alimentos que especifica en qu\'e grupo se encuentran y cu\'antas unidades contiene. Por ejemplo, la lista puede contener
`una taza de cereal de desayuno', aclarando que est\'a en el grupo de carbohidratos y que
equivale a $2$ unidades de los mismos. As\'i que el paciente podr\'ia tomar dos tazas de cereales al d\'ia o bien una sola taza y alg\'un otro alimento con contenido de dos unidades de carbohidratos.

Supongamos, para simplificar, que solo disponemos de cinco alimentos, cuyo contenido en
nutrientes de cada grupo es como sigue:
\begin{center}
\begin{table}[!ht]
\centering
\begin{tabular}{c|ccccccc}
\hline
& Prote\'ina & Verdura & L\'acteos & Grasa & Carbohidratos  \\
\hline
\rowcolor{LightCyan}
Alimento 1 & 0 & 2 & 0 & 0 & 1  \\
Alimento 2 & 0 & 2 & 1 & 1 & 0  \\
\rowcolor{LightCyan}
Alimento 3 & 2 & 0 & 2 & 1 & 0  \\
Alimento 4 & 4 & 0 & 0 & 3 & 1 \\
\rowcolor{LightCyan}
Alimento 5 & 1 & 1 & 2 & 2 & 1  \\
\hline
\end{tabular}
\end{table}
\end{center}

Queremos determinar la cantidad $x_1$ del primer alimento, $x_2$ del segundo, $x_3$ del
tercero, $x_4$ del cuarto y $x_5$ del quinto que debemos consumir para cumplir con la
recomendaci\'on del nutri\'ologo. 

\begin{center}
\includegraphics[scale=0.2]{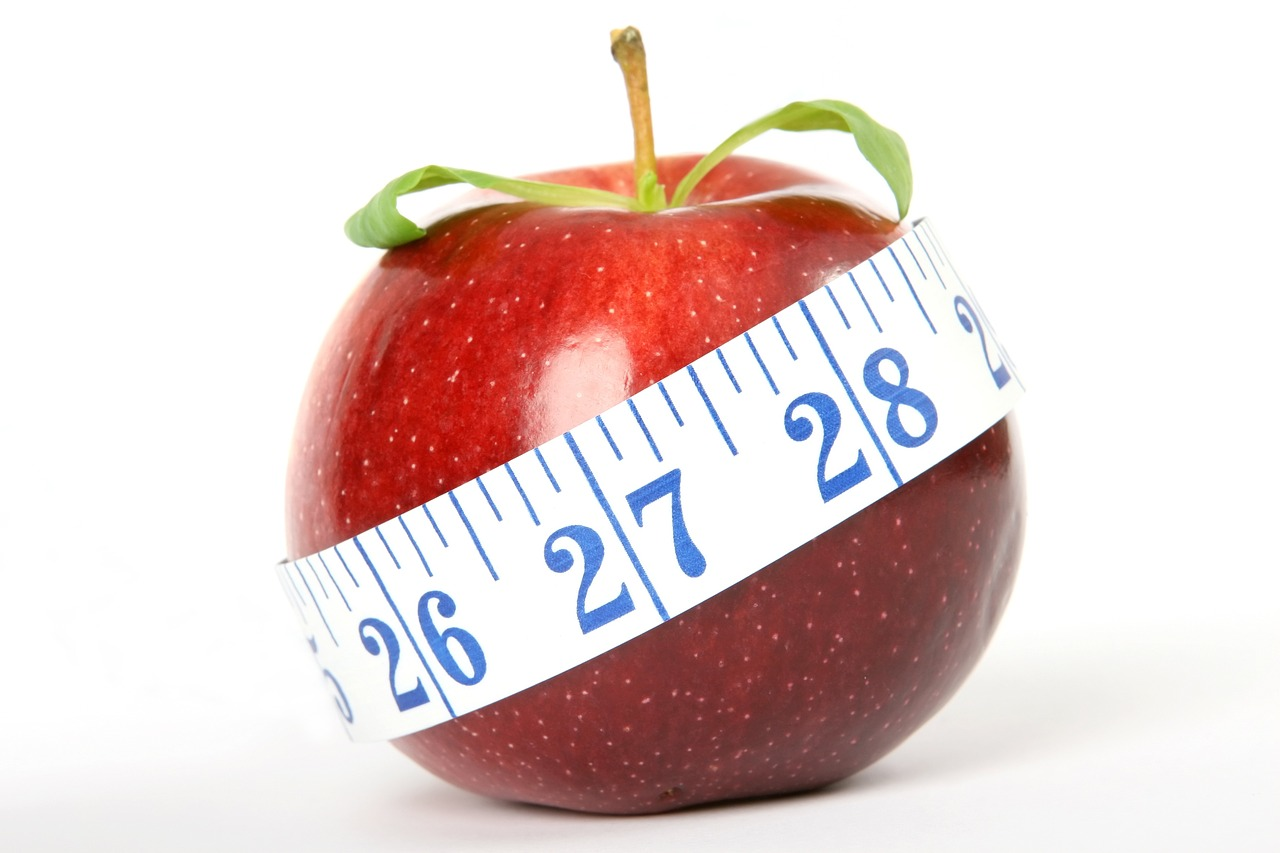} 
\end{center}

De la primera columna de la tabla precedente resulta
que el total de prote\'inas ingeridas ser\'a
$$
2x_3+4x_4+x_5
$$ 
y esa cantidad debe ser igual a $4$ (unidades de prote\'ina), de donde resulta la ecuaci\'on
$$
2x_3+4x_4+x_5 = 4\,.
$$

Procediendo de igual manera con las restantes columnas, llegamos al sistema de ecuaciones
\begin{equation}\label{sistemas1b}
\systeme{2x_3+4x_4+x_5=4,2x_1+2x_2+x_5=6,x_2+2x_3+2x_5=3,x_2+x_3+3x_4+2x_5=2,x_1+x_4+x_5=4}
\end{equation}

Fij\'emonos en c\'omo hemos colocado las inc\'ognitas $x_i$ (con $i\in\{1,2,3,4,5\}$),
alineadas verticalmente. Los `huecos' que quedan podemos suponer que est\'an ocupados
por expresiones de la forma $0x_i$, con el valor de $i$ que corresponda, es decir, el
sistema precedente es completamente equivalente a
\begin{equation}\label{sistemas2}
\systeme{0x_1+0x_2+2x_3+4x_4+x_5=4,
2x_1+2x_2+0x_3+0x_4+x_5=6,
0x_1+x_2+2x_3+0x_4+2x_5=3,
0x_1+x_2+x_3+3x_4+2x_5=2,
x_1+0x_2+0x_3+x_4+x_5=4}
\end{equation}

Al escribir el sistema de esta forma, enseguida nos damos cuenta de que podemos simplificar
la escritura especificando \'unicamente los coeficientes de las ecuaciones, dado que los signos
de suma y de igualdad ya sabemos d\'onde est\'an colocados. Tambi\'en necesitamos recordar
c\'omo estamos llamando a las inc\'ognitas y, teniendo esto en cuenta, podemos escribir
\eqref{sistemas2} en otra forma equivalente mucho m\'as compacta:
\begin{equation}\label{sistemas3}
\begin{pmatrix}
0 & 0 & 2 & 4 & 1 \\
2 & 2 & 0 & 0 & 1 \\
0 & 1 & 2 & 0 & 2 \\
0 & 1 & 1 & 3 & 2 \\
1 & 0 & 0 & 1 & 1
\end{pmatrix}
\begin{pmatrix}
x_1 \\ x_2 \\ x_3 \\ x_4 \\ x_5
\end{pmatrix}
=
\begin{pmatrix}
4 \\ 6 \\ 3 \\ 2 \\ 4
\end{pmatrix}\,.
\end{equation}

En \eqref{sistemas3} hemos introducido la notaci\'on de matrices para expresar el
sistema \eqref{sistemas2}: una \emph{matriz}\index{Matriz} no es m\'as que una tabla de
n\'umeros, llamados \emph{coeficientes}\index{Matriz!coeficientes} de la matriz, escrita 
entre par\'entesis. Las l\'ineas de n\'umeros horizontales se llaman
las \emph{filas}\index{Matriz!filas} o \emph{renglones}\index{Matriz!rengl\'on} de la matriz, 
mientras que las l\'ineas verticales son sus \emph{columnas}\index{Matriz!columna}. Si una 
matriz tiene $m$ filas y $n$ columnas, se dice que es de \emph{orden} $m\times n$
\index{Matriz!orden} y la matriz es \emph{cuadrada}\index{Matriz!cuadrada} si $m=n$. As\'i,
en \eqref{sistemas3} tenemos una matriz con cinco filas y cinco columnas (de orden
$5\times 5$), y dos matrices con cinco filas y una sola columna (de orden $5\times 1$). 
La matriz
$$
A=\begin{pmatrix}
0 & 0 & 2 & 4 & 1 \\
2 & 2 & 0 & 0 & 1 \\
0 & 1 & 2 & 0 & 2 \\
0 & 1 & 1 & 3 & 2 \\
1 & 0 & 0 & 1 & 1
\end{pmatrix}
$$
es la \emph{matriz de coeficientes}\index{Matriz!de coeficientes} del sistema equivalente
\eqref{sistemas2}, la matriz
$$
\mathbf{x}=\begin{pmatrix}
x_1 \\ x_2 \\ x_3 \\ x_4 \\ x_5
\end{pmatrix}
$$
es la matriz de \emph{inc\'ognitas}\index{Matriz!de inc\'ognitas} de \eqref{sistemas2} y, 
finalmente, la matriz
$$
\mathbf{b}=\begin{pmatrix}
4 \\ 6 \\ 3 \\ 2 \\ 4
\end{pmatrix}
$$
es la matriz de \emph{t\'erminos no homog\'enos}\index{Matriz!de t\'erminos no homog\'eneos} 
del sistema equivalente \eqref{sistemas2}, que matricialmente es
$$
A\mathbf{x}=\mathbf{b}\,.
$$

Resulta conveniente introducir la siguiente notaci\'on: para una matriz como $A$, el
coeficiente situado en la fila $i$ y la columna $j$ se denota $A_{ij}$. As\'i, se tiene
que $A_{23}=0$, $A_{44}=3$, etc.
Con esta notaci\'on, el sistema \eqref{sistemas1} se escribir\'a matricialmente como
$$
\begin{pmatrix}
a_{11} & a_{12} \\
a_{21} & a_{22}
\end{pmatrix}
\begin{pmatrix}
x \\ y
\end{pmatrix}
=
\begin{pmatrix}
b_1 \\ b_2 
\end{pmatrix}\,.
$$

Rec\'iprocamente, de una ecuaci\'on matricial se puede pasar a un sistema de ecuaciones.
Por ejemplo, dada la ecuaci\'on entre matrices
\begin{equation}\label{sistemas4}
\begin{pmatrix}
0 & -1 & 3 \\
1 & 2  & -1 \\
-2 & 3 & 1
\end{pmatrix}
\begin{pmatrix}
x \\ y \\z
\end{pmatrix}
=
\begin{pmatrix}
2 \\ -2 \\ 0 
\end{pmatrix}\,,
\end{equation}
vemos que su sistema equivalente de ecuaciones es
\begin{equation}\label{sistemas5}
\systeme{-y+3z=2,x+2y-z=-2,-2x+3y+z=0\,.}
\end{equation}

La pregunta ahora es: ?`c\'omo resolvemos un sistema tal como \eqref{sistemas2} o 
\eqref{sistemas5}? O, de otra manera, ?`c\'omo resolvemos \eqref{sistemas3} o 
\eqref{sistemas4}? Vamos a describir un m\'etodo muy eficiente que se remonta a 
Gauss\index{Gauss, Carl Friedrich}\footnote{Por supuesto, la esencia del m\'etodo en s\'i es mucho m\'as antigua y se remonta a la matem\'atica china (unos 200 a\~nos antes de nuestra era), tambi\'en fue descrito por Newton\index{Newton, Isaac}, pero en su forma actual y con
la notaci\'on que usaremos fue elaborado por Gauss.} y que,
por este motivo, se denomina \emph{eliminaci\'on gaussiana}\index{Eliminaci\'on gaussiana}.

Para fijar ideas, describiremos el m\'etodo con el sistema \eqref{sistemas5}. Comenzamos
por reordenar las ecuaciones de manera que en la primera el coeficiente
acompa\~nando a la primera inc\'ognita, que en este caso es $x$, sea distinto de cero.

No importa que ecuaci\'on se tome mientras se cumpla esta condici\'on. Como ejemplo,
nosotros tomaremos la tercera de \eqref{sistemas5} (el criterio seguido es que tiene
el mayor coeficiente de $x$ en valor absoluto) y de acuerdo con lo dicho, la 
intercambiamos con la primera.
De este modo, el sistema reordenado en el primer paso queda como
\begin{equation}\label{sistemas5b}
\systeme{-2x+3y+z=0,x+2y-z=-2,-y+3z=2\,.}
\end{equation}

El segundo paso consiste en utilizar la primera ecuaci\'on para eliminar el coeficiente de la 
$x$ en todas las dem\'as. En \eqref{sistemas5b}, la tercera ecuaci\'on ya tiene un $0$ como
coeficiente de la $x$, as\'i que no nos preocupamos por ella y pasamos directamente a la 
segunda ecuaci\'on, que tiene un $1$ acompa\~nando a la $x$. ?`C\'omo lo convertimos en
un $0$? Dividimos la primera ecuaci\'on por $-2$, para que el coeficiente de $x$ sea $1$.
Resulta:
$$
x-\frac{3}{2}y-\frac{1}{2}z=0\,.
$$
Ahora, restamos esta ecuaci\'on de la segunda en \eqref{sistemas5b}:
$$
\begin{array}{cl}
   &   x+2y-z=-2 \\[1pt]
\smash{\makebox[0pt][r]{\raisebox{8pt}{--}}}   &   x-\frac{3}{2}y-\frac{1}{2}z=0 \\\hline
   & \phantom{x+}\frac{7}{2}y-\frac{1}{2}z=-2\,.
\end{array}
$$

Finalizamos esta primera etapa reemplazando la segunda ecuaci\'on por esta que acabamos de obtener, en la que ya no aparece la $x$. El sistema equivalente es, pues,
$$
\systeme{-2x+3y+z=0,\frac{7}{2}y-\frac{1}{2}z=-2,-y+3z=2\,,}
$$
que podemos simplificar multiplicando la segunda ecuaci\'on por $2$ para eliminar
denominadores:
\begin{equation}\label{sistemas5c}
\systeme{-2x+3y+z=0,7y-z=-4,-y+3z=2\,.}
\end{equation}

En la siguiente etapa, vamos a eliminar los coeficientes de la variable $y$, as\'i que buscamos una fila (diferente de la primera, que ya no modificaremos en lo que resta) que
tenga el coeficiente de $y$ distinto de cero y la ponemos como segunda fila. En este caso,
no tenemos que hacer nada porque la segunda ecuaci\'on de \eqref{sistemas5c} ya cumple
esta condici\'on. La utilizaremos para hacer cero el coeficiente de $y$ en la tercera fila
y para ello basta con multiplicar la tercera por $7$ (resulta
$-7y+21z=14$) y sumarlas:
$$
\begin{array}{cl}

   &  \phantom{-} 7y-\phantom{21}z=-4  \\[1pt]
\smash{\makebox[0pt][r]{\raisebox{8pt}{+}}}   &  -7y+21z=14 \\\hline
   & \phantom{7y +1} 20z=10\,.
\end{array}
$$
Reemplazando la tercera ecuaci\'on por \'esta, nos queda
\begin{equation}\label{sistemas5d}
\systeme{-2x+3y+z=0,7y-z=-4,20z=10}
\end{equation}
un \emph{sistema triangular}\index{Sistema!de ecuaciones!triangular}
\index{Ecuaciones!sistema!triangular} que ya tiene soluci\'on 
inmediata mediante la llamada
\emph{sustituci\'on regresiva}\index{Sustituci\'on regresiva}: de la \'ultima ecuaci\'on
es inmediato que $z=1/2$, valor que se sustituye en la pen\'ultima para dar $7y-1/2=-4$, de
donde $y=-1/2$. Los valores obtenidos de $y$, $z$ se insertan en la primera ecuaci\'on
y queda una ecuaci\'on lineal para $x$, $-2x-3/2+1/2=0$, que se resuelve inmediatamente con
el resultado $x=-1/2$. La soluci\'on a \eqref{sistemas5} es
$$
\begin{pmatrix}
x \\ y \\ z
\end{pmatrix}
=
\begin{pmatrix}
-\frac{1}{2} \\ -\frac{1}{2} \\ \frac{1}{2}
\end{pmatrix}\,.
$$

Para sistemas con un mayor n\'umero de ecuaciones, el procedimiento que hemos descrito
se repite eliminando coeficientes de las sucesivas inc\'ognitas hasta que se deja en forma 
triangular, aplicando entonces la sustituci\'on regresiva.

Existe otro enfoque, m\'as te\'orico, para el estudio de los sistemas de ecuaciones y su 
relaci\'on con las matrices. Con la experiencia que hemos acumulado a lo largo del libro, 
tambi\'en podr\'iamos haber contestado a la pregunta de c\'omo se resuelve un sistema de 
ecuaciones diciendo `para resolver \eqref{sistemas3} o \eqref{sistemas4},
por ejemplo, basta con multiplicar ambos miembros de la ecuaci\'on por la inversa de
la matriz de coeficientes', de la misma manera que resolv\'iamos $ax=b$ multiplicando
ambos miembros por el inverso de $a$. 
Esta idea intuitiva est\'a bien, pero para dotarla de sentido primero hay que definir un 
producto entre matrices, de manera que podamos hablar de `matriz inversa'.\index{Operaciones!con matrices} 
Para justificar la forma en que multiplicaremos matrices, consideremos de nuevo la relaci\'on
entre un sistema de ecuaciones como
\begin{equation}\label{sistemas6}
\begin{cases}
b_{11}x+b_{12}y = \gamma_1 \\
b_{21}x+b_{22}y= \gamma_2\,,
\end{cases}
\end{equation}
y una ecuaci\'on matricial. En este caso, \eqref{sistemas6} es lo mismo que
\begin{equation}\label{sistemas7}
\begin{pmatrix}
b_{11} & b_{12} \\
b_{21} & b_{22}
\end{pmatrix}
\begin{pmatrix}
x \\ y
\end{pmatrix}
=
\begin{pmatrix}
\gamma_1 \\ \gamma_2
\end{pmatrix}\,.
\end{equation}
Si tenemos dos matrices como por ejemplo
$$
A=\begin{pmatrix}
a_{11} & a_{12} \\
a_{21} & a_{22}
\end{pmatrix}\quad\mbox{ y }\quad
B=\begin{pmatrix}
b_{11} & b_{12} \\
b_{21} & b_{22}
\end{pmatrix}\,,
$$
deduciremos la matriz que corresponde a su producto planteando la ecuaci\'on asociada a ese
producto. Para ello, multiplicamos \eqref{sistemas7} por la izquierda con $A$. Llamando
\begin{equation}\label{sistemas8}
\begin{pmatrix}
\beta_1 \\ \beta_2
\end{pmatrix}
=
\begin{pmatrix}
a_{11} & a_{12} \\
a_{21} & a_{22}
\end{pmatrix}
\begin{pmatrix}
\gamma_1 \\ \gamma_2
\end{pmatrix}
\end{equation}
resulta que
\begin{equation*}
\begin{pmatrix}
a_{11} & a_{12} \\
a_{21} & a_{22}
\end{pmatrix}
\begin{pmatrix}
b_{11} & b_{12} \\
b_{21} & b_{22}
\end{pmatrix}
\begin{pmatrix}
x \\ y
\end{pmatrix}
=
\begin{pmatrix}
\beta_1 \\ \beta_2
\end{pmatrix}
=
\begin{pmatrix}
a_{11} & a_{12} \\
a_{21} & a_{22}
\end{pmatrix}
\begin{pmatrix}
\gamma_1 \\ \gamma_2
\end{pmatrix}\,.
\end{equation*}
En particular, de la primera igualdad se deduce
\begin{equation}\label{sistemas9}
\begin{pmatrix}
a_{11} & a_{12} \\
a_{21} & a_{22}
\end{pmatrix}
\begin{pmatrix}
b_{11} & b_{12} \\
b_{21} & b_{22}
\end{pmatrix}
\begin{pmatrix}
x \\ y
\end{pmatrix}
=
\begin{pmatrix}
\beta_1 \\ \beta_2
\end{pmatrix}
\end{equation}

Ahora, suponiendo que el producto entre matrices que queremos definir es asociativo,
utilizamos la correspondencia entre \eqref{sistemas6} y \eqref{sistemas7} para escribir
$$
\begin{pmatrix}
a_{11} & a_{12} \\
a_{21} & a_{22}
\end{pmatrix}
\begin{pmatrix}
\gamma_1 \\ \gamma_2
\end{pmatrix}
=
\begin{pmatrix}
a_{11} & a_{12} \\
a_{21} & a_{22}
\end{pmatrix}
\begin{pmatrix}
b_{11}x+b_{12}y \\
b_{21}x+b_{22}y
\end{pmatrix}\,.
$$

El miembro izquierdo est\'a dado por \eqref{sistemas8}. Entonces, intercambiando los miembros derecho e izquierdo, hemos llegado a
$$
\begin{pmatrix}
a_{11} & a_{12} \\
a_{21} & a_{22}
\end{pmatrix}
\begin{pmatrix}
b_{11}x+b_{12}y \\
b_{21}x+b_{22}y
\end{pmatrix}
=
\begin{pmatrix}
\beta_1 \\ \beta_2
\end{pmatrix}\,.
$$

Esta ecuaci\'on matricial tiene un sistema de ecuaciones lineales equivalente, que no es otro
que
$$
\begin{cases}
a_{11}(b_{11}x+b_{12}y) +a_{12}(b_{21}x+b_{22}y) = \beta_1 \\
a_{21}(b_{11}x+b_{12}y) +a_{22}(b_{21}x+b_{22}y) = \beta_2\,.
\end{cases}
$$

Reordenamos las ecuaciones de este \'ultimo sistema para resaltar que nuestras inc\'ognitas realmente son $x$ e $y$:
$$
\begin{cases}
(a_{11}b_{11}+a_{12}b_{21})x + (a_{11}b_{12}+a_{12}b_{22})y = \beta_1 \\
(a_{21}b_{11}+a_{22}b_{21})x + (a_{21}b_{12}+a_{22}b_{22})y = \beta_2\,,
\end{cases}
$$
o bien, matricialmente,
\begin{equation}\label{sistemas10}
\begin{pmatrix}
a_{11}b_{11}+a_{12}b_{21} & a_{11}b_{12}+a_{12}b_{22} \\
a_{21}b_{11}+a_{22}b_{21} & a_{21}b_{12}+a_{22}b_{22}
\end{pmatrix}
\begin{pmatrix}
x \\ y
\end{pmatrix}
=
\begin{pmatrix}
\beta_1 \\ \beta_2
\end{pmatrix}
\end{equation}

Hemos llegado a la expresi\'on buscada para el \emph{producto de matrices}
\index{Matriz!producto}. Comparando \eqref{sistemas9} con \eqref{sistemas10}
obtenemos
\begin{equation}\label{sistemas11}
\begin{pmatrix}
a_{11} & a_{12} \\
a_{21} & a_{22}
\end{pmatrix}
\begin{pmatrix}
b_{11} & b_{12} \\
b_{21} & b_{22}
\end{pmatrix}
=
\begin{pmatrix}
a_{11}b_{11}+a_{12}b_{21} & a_{11}b_{12}+a_{12}b_{22} \\
a_{21}b_{11}+a_{22}b_{21} & a_{21}b_{12}+a_{22}b_{22}
\end{pmatrix}\,.
\end{equation}

N\'otese que el producto matricial \emph{no} es conmutativo: si se repite el proceso
precedente intercambiando $A$ con $B$, se llega a la expresi\'on
$$
\begin{pmatrix}
b_{11} & b_{12} \\
b_{21} & b_{22}
\end{pmatrix}
\begin{pmatrix}
a_{11} & a_{12} \\
a_{21} & a_{22}
\end{pmatrix}
=
\begin{pmatrix}
b_{11}a_{11}+b_{12}a_{21} & b_{11}a_{12}+b_{12}a_{22} \\
b_{21}a_{11}+b_{22}a_{21} & b_{21}a_{12}+b_{22}a_{22}
\end{pmatrix}\,,
$$
que \emph{no} coincide con \eqref{sistemas11}.

En general, es posible multiplicar una matriz $A$ de orden $m\times p$ por una matriz $B$ de 
orden $p\times n$ y el resultado es otra matriz $C=AB$ de orden $m\times n$. En ese caso, si 
$A,B$ son las matrices a multiplicar, los coeficientes $c_{ij}$ de la matriz producto $C=AB$ 
est\'an dados por
\begin{equation}\label{sistemas12}
c_{ij}=a_{i1}b_{1j}+\cdots +a_{ip}b_{pj}=\sum^p_{k=1}a_{ik}b_{kj}\,,
\end{equation}
con $i\in\{1,\ldots,m\}$ y $j\in\{1,\ldots,n\}$.

Las f\'ormulas para los coeficientes de la matriz producto pueden parecer muy complicadas (de 
hecho, lo son), pero se comprenden muy f\'acilmente mediante el 
\emph{esquema de Falk}\index{Falk, James!esquema} para su c\'alculo, que se basa en escribir 
el producto $AB=C$ en la forma
$$
\begin{tabular}{c|c}
 & B \\ 
\hline 
A & C \\ 
\end{tabular} 
$$

Con mayor detalle, cada elemento de la matriz producto $C$ se obtiene a partir de la fila
de $A$ a su izquierda, a la misma altura, y de la columna de $B$ situada sobre \'el, 
multiplicando elemento a elemento y sumando los resultados, como se indica en la siguiente 
figura:

\newcommand{\myunit}{1 cm}
\tikzset{
    node style sp/.style={draw,circle,minimum size=\myunit},
    node style ge/.style={circle,minimum size=\myunit},
    arrow style mul/.style={draw,sloped,midway,fill=white},
    arrow style plus/.style={midway,sloped,fill=white},
}
\begin{center}
\begin{tikzpicture}[>=latex]
\matrix (A) [matrix of math nodes,%
             nodes = {node style ge},%
             left delimiter  = (,%
             right delimiter = )] at (0,0)
{%
  a_{11} & a_{12} & \ldots & a_{1p}  \\
  |[node style sp]| a_{21}%
         & |[node style sp]| a_{22}%
                  & \ldots%
                           & |[node style sp]| {a_{2p}} \\
  \vdots & \vdots & \ddots & \vdots  \\
  a_{m1} & a_{m2} & \ldots & a_{mp}  \\
};

\matrix (B) [matrix of math nodes,%
             nodes = {node style ge},%
             left delimiter  = (,%
             right delimiter =)] at (6*\myunit,6*\myunit)
{%
  b_{11} & |[node style sp]| b_{12}%
                  & \ldots & b_{1n}  \\
  b_{21} & |[node style sp]| b_{22}%
                  & \ldots & b_{2n}  \\
  \vdots & \vdots & \ddots & \vdots  \\
  b_{p1} & |[node style sp]| b_{p2}%
                  & \ldots & b_{pn}  \\
};

\matrix (C) [matrix of math nodes,%
             nodes = {node style ge},%
             left delimiter  = (,%
             right delimiter = )] at (6*\myunit,0)
{%
  c_{11} & c_{12} & \ldots & c_{1n} \\
  c_{21} & |[node style sp,red]| c_{22}%
                  & \ldots & c_{2n} \\
  \vdots & \vdots & \ddots & \vdots \\
  c_{m1} & c_{m2} & \ldots & c_{mn} \\
};
\draw[blue] (A-2-1.north) -- (C-2-2.north);
\draw[blue] (A-2-1.south) -- (C-2-2.south);
\draw[blue] (B-1-2.west)  -- (C-2-2.west);
\draw[blue] (B-1-2.east)  -- (C-2-2.east);
\draw[<->,red](A-2-1) to[in=180,out=90]
    node[arrow style mul] (x) {$a_{21}\times b_{12}$} (B-1-2);
\draw[<->,red](A-2-2) to[in=180,out=90]
    node[arrow style mul] (y) {$a_{22}\times b_{22}$} (B-2-2);
\draw[<->,red](A-2-4) to[in=180,out=90]
    node[arrow style mul] (z) {$a_{2n}\times b_{p2}$} (B-4-2);
\draw[red,->] (x) to node[arrow style plus] {$+$} (y)%
                  to node[arrow style plus] {$+\raisebox{.5ex}{\ldots}+$} (z)%
                  to (C-2-2.north west);
\end{tikzpicture}
\end{center}

Como ejemplo, vamos a multiplicar las matrices
$$
A=\begin{pmatrix}
2 & 1 \\
-1 & 0 \\
4 & 3
\end{pmatrix}
\quad\mbox{ y }\quad
B=\begin{pmatrix}
-1 & 1 & -3/2\\
0 & 1/2 & 1
\end{pmatrix}\,.
$$
Como $A$ tiene orden $3\times 2$ y $B$ es de orden $2\times 3$, el resultado $C=AB$ ser\'a
una matriz de orden $3\times 3$. Utilizando el esquema de Falk o la f\'ormula
\eqref{sistemas12}, resulta
$$
C=AB=\begin{pmatrix}
2\cdot (-1)+1\cdot 0 & 2\cdot 1+1\cdot 1/2 & 2\cdot (-3/2)+1\cdot 1 \\
(-1)\cdot (-1)+0\cdot 0 & (-1)\cdot 1+0\cdot 1/2 & (-1)\cdot (-3/2)+0\cdot 1 \\
4\cdot (-1)+3\cdot 0 & 4\cdot 1+3\cdot (1/2)  &  4\cdot (-3/2)+3\cdot 1   
\end{pmatrix}
=
\begin{pmatrix}
-2 & 5/2 & -2 \\
1 & -1 & 3/2 \\
-4 & 11/2 & -3
\end{pmatrix}
$$

Aparentemente, esta forma de multiplicar matrices es demasiado complicada para resultar
de utilidad, pero nada m\'as lejos de la realidad: el producto matricial \eqref{sistemas12}
es una de las herramientas m\'as poderosas de la Matem\'atica actual y en la secci\'on
\ref{apl-mat} veremos algunas aplicaciones.

Una vez que se tiene la noci\'on de multiplicaci\'on de matrices, podemos buscar su
\emph{elemento neutro}\index{Matriz!producto!elemento neutro}: una matriz $I$ tal que
para cualquier otra matriz $A$ se tenga
$$
AI=A=IA\,.
$$

Si nos restringimos \'unicamente a matrices cuadradas de orden $n\times n$, es inmediato
comprobar que el elemento neutro es la \emph{matriz identidad}\index{Matriz!identidad}
$$
I=\begin{pmatrix}
1 &        & \\
  & \ddots & \\
  &        & 1
\end{pmatrix}\,.
$$

Armados con el concepto de matriz identidad, tambi\'en podemos buscar la \emph{matriz inversa}
de una matriz $A$\index{Matriz!inversa}: ser\'a otra matriz $A^{-1}$ tal que
\begin{equation}\label{sistemas13}
AA^{-1}=I=A^{-1}A\,.
\end{equation}

En lugar de dar una serie de f\'ormulas aun m\'as complejas que \eqref{sistemas12},
ilustraremos el uso de \eqref{sistemas13} con un ejemplo. Supongamos que queremos encontrar
la matriz inversa de
$$
A=\begin{pmatrix}
2 & -1 \\
0 & 3
\end{pmatrix}\,.
$$
Llamemos $b_{ij}$ a los coeficientes de la matriz inversa $A^{-1}$. Por definici\'on, 
estos coeficientes deben satisfacer \eqref{sistemas13}, esto es,
$$
\begin{pmatrix}
b_{11} & b_{12} \\
b_{21} & b_{22}
\end{pmatrix}
\begin{pmatrix}
2 & -1 \\
0 & 3
\end{pmatrix}
=
\begin{pmatrix}
1 & 0 \\
0 & 1
\end{pmatrix}\,.
$$
Efectuando el producto matricial en el miembro izquierdo, llegamos a
$$
\begin{pmatrix}
2b_{11} & -b_{11}+3b_{12} \\
2b_{21} & -b_{21}+3b_{22}
\end{pmatrix}
=
\begin{pmatrix}
1 & 0 \\
0 & 1
\end{pmatrix}\,.
$$
Esta ecuaci\'on matricial es equivalente al sistema lineal
$$
\begin{cases}
2b_{11}=1 \\
-b_{11}+3b_{12}=0 \\
2b_{21}=0 \\
-b_{21}+3b_{22}=1\,,
\end{cases}
$$
que puede resolverse f\'acilmente por sustituci\'on directa de la primera y la tercera ecuaci\'on en la segunda y cuarta, quedando solamente
$$
\begin{cases}
-2 +3b_{12}=0 \\
3b_{22} =1\,.
\end{cases}
$$

De la \'ultima ecuaci\'on resulta $b_{22}=1/3$ y de la primera $b_{12}=2/3$. Por tanto,
la matriz inversa de $A$ es
$$
A^{-1}=\begin{pmatrix}
1/2 & 1/6 \\
0 & 1/3
\end{pmatrix}\,.
$$

Una vez que tenemos la matriz inversa $A^{-1}$ es inmediato resolver cualquier sistema de ecuaciones que se exprese como $A\mathbf{x}=\mathbf{b}$, donde
$$
\mathbf{x}=\begin{pmatrix}
x \\ y
\end{pmatrix}\quad\mbox{ y }\quad
\mathbf{b}=\begin{pmatrix}
b_1 \\ b_2
\end{pmatrix}\,;
$$
la soluci\'on es, simplemente,
$$
\begin{pmatrix}
x \\ y
\end{pmatrix}
=
A^{-1}\begin{pmatrix}
b_1 \\ b_2
\end{pmatrix}\,.
$$

Para finalizar esta breve introducci\'on al \'algebra matricial, podemos preguntarnos qu\'e 
sucede con el producto que intuitivamente se le ocurrir\'ia definir a cualquiera, esto es, el 
\emph{producto componente a componente}\index{Matriz!producto!componente a componente} como en
$$
\begin{pmatrix}
a_{11} & a_{12} \\
a_{21} & a_{22}
\end{pmatrix}
\begin{pmatrix}
b_{11} & b_{12} \\
b_{21} & b_{22}
\end{pmatrix} =
\begin{pmatrix}
a_{11}b_{11} & a_{12}b_{12} \\
a_{21}b_{21} & a_{22}b_{22}
\end{pmatrix}\,.
$$

Este manera de multiplicar matrices tambi\'en se conoce como \emph{producto de Hadamard}
\index{Matriz!producto! de Hadamard}.
?`Tiene alguna aplicaci\'on este producto? La respuesta es que s\'i, como veremos en la
secci\'on \ref{imgproc}, dedicada a los principios del procesamiento de im\'agenes, donde 
usaremos el producto componente a componente para seleccionar partes de una imagen.

\section{Aplicaciones}\label{apl-mat}
\subsection{Matrices de incidencia en grafos}

Vamos a ver c\'omo modelar una \emph{red}\index{Red} mediante un \emph{grafo}\index{Grafo},
como ejemplo tomaremos la red de transporte a\'ereo de un peque\~no pa\'is con cinco
ciudades $A$, $B$, $C$, $D$, $E$, que poseen aeropuerto. Las ciudades se ubican en los
\emph{nodos}\index{Grafo!nodos} o \emph{v\'ertices}\index{Grafo!v\'ertices} del grafo,
y se traza una \emph{arista}\index{Grafo!arista} entre dos nodos si hay un vuelo que los
conecta, como en la figura.

\begin{center}
\begin{tikzpicture}
\Vertex[IdAsLabel,Math,opacity =.5]{A} 
\Vertex[IdAsLabel,Math,opacity =.5,x=2]{B}
\Vertex[IdAsLabel,Math,opacity =.5,y=2]{C}
\Vertex[IdAsLabel,Math,opacity =.5,y=-2]{D}
\Vertex[IdAsLabel,Math,opacity =.5,x=-2]{E}
\Edge[bend=45](A)(B)
\Edge[bend=-45](A)(B)
\Edge[bend=45](E)(A)
\Edge(E)(A)
\Edge[bend=-45](E)(A)
\Edge(E)(D)
\Edge[bend=-30](E)(D)
\Edge(A)(D)
\Edge[bend=-30](D)(C)
\Edge[bend=30](C)(B)
\Edge[bend=60](E)(B)

\end{tikzpicture}
\end{center}

A la vista de este grafo, podemos pensar que una de entre $A$ o $E$ es la capital del pa\'is,
porque son las que mayor n\'umero de conexiones tienen\footnote{Aunque no todas las ciudades 
se encuentran conectadas con ellas, como es el caso de $C$, quiz\'as porque se encuentran 
cerca por carretera, o porque existe una buena comunicaci\'on por otro medio alternativo, 
como el ferrocarril.}, mientras que $C$ parece la ciudad con peores comunicaciones.

La informaci\'on m\'as relevante de la red de comunicaciones a\'ereas de nuestro pa\'is se
puede representar por medio de la llamada 
\emph{matriz de incidencia}\index{Matriz!de incidencia} asociada a su grafo. Para
construirla, convenimos en asignarle (arbitrariamente) a cada ciudad un n\'umero.
Lo m\'as sencillo, desde luego, es usar la correspondencia $A\to 1$, $B\to 2$, $C\to 3$, $D\to 4$, $E\to 5$. A continuaci\'on, definimos la matriz $\Lambda$ de manera que su
coeficiente $(i,j)$ (esto es, ubicado en la intersecci\'on de la fila $i$ con la columna
$j$) sea el n\'umero de vuelos que conectan la ciudad $i$ con la $j$. En nuestro ejemplo,
tendremos
\begin{equation}\label{graph-1}
\Lambda =
\begin{pmatrix}
0 & 2 & 0 & 1 & 3 \\
2 & 0 & 1 & 0 & 1 \\
0 & 1 & 0 & 1 & 0 \\
1 & 0 & 1 & 0 & 2 \\
3 & 1 & 0 & 2 & 0
\end{pmatrix}\,.
\end{equation}

N\'otese que la matriz contiene ceros en la diagonal (no hay vuelos conectando una ciudad
con ella misma) y es \emph{sim\'etrica}\index{Matriz!sim\'etrica} (los vuelos que conectan
la ciudad $i$ con la $j$, son los mismos que conectan $j$ con $i$).

Supongamos ahora que una empresa multinacional interesada en invertir en nuestro hipot\'etico 
pa\'is env\'ia un representante para estudiar posibles negocios. El tiempo es oro, as\'i que 
la empresa decide que solo se visiten tres ciudades y que el representante haga un informe 
para ser presentado en el consejo de accionistas. Para programar el viaje, la empresa le
pide a una agencia de viajes que le informe de cu\'antas posibilidades hay (y del precio
de cada una). Por ejemplo, ?`cu\'antos viajes posibles existen entre $E$ y $C$ haciendo
exactamente una escala? Naturalmente, podr\'iamos hacer el c\'alculo manualmente,
enumerando todas las posibilidades y asegur\'andonos de que no exclu\'imos ninguna, y aunque
este procedimiento es factible en nuestro ejemplo, es f\'acil darse cuenta de que no es
posible hacerlo en el caso en que se tengan decenas de ciudades, por la cantidad de
tiempo y esfuerzo que requerir\'ia, as\'i que buscamos una t\'ecnica m\'as eficiente.

Resulta que la multiplicaci\'on de matrices \eqref{sistemas12} nos da la respuesta: si
calculamos $\Lambda^2 =\Lambda\,\Lambda$ encontraremos
\begin{equation}\label{graph-2}
\Lambda^2 =\begin{pmatrix}
14 & 3 & 3 & 6 & 4 \\
3 & 6 & 0 & 5 & 6 \\
3 & 0 & 2 & 0 & 3 \\
6 & 5 & 0 & 6 & 3 \\
4 & 6 & 3 & 3 & 14
\end{pmatrix}\,,
\end{equation}
y el coeficiente en la posici\'on $(i,j)$ (esto es, fila $i$ y columna $j$) nos da la
cantidad de vuelos entre $i$ y $j$ con una escala intermedia. Por ejemplo, 
ese n\'umero para las
ciudades $E$ y $C$ se encuentra en cualquiera de los coeficientes $\Lambda_{35}$ o
$\Lambda_{53}$ (pues recordemos que estas matrices son sim\'etricas) y es
$$
\Lambda_{53}=3\,.
$$

?`Por qu\'e funciona esta t\'ecnica? Si nos propusi\'eramos contar las conexiones entre
dos nodos $i$ y $j$ v\'ia un tercero $k$, tendr\'iamos que empezar por contar las conexiones
entre $i$ y $k$, luego contar\'iamos las que hay de $k$ a $j$, y finalmente las multiplicar\'iamos (porque para cada manera de ir de $i$ a $k$, podemos escoger entre todas las posibilidades de $k$ a $j$). O sea, el producto
$$
(\mbox{\texttt{\#}  de vuelos entre }i\mbox{ y }k)\times 
(\mbox{\texttt{\#} de vuelos entre }k\mbox{ y }j)
$$
proporciona todas las posibilidades de vuelos entre $i,j$ pasando por $k$. El n\'umero total
de posibilidades se obtiene sumando estos n\'umeros para cada $k$:
$$
\sum^5_{k=1} (\mbox{\texttt{\#} de vuelos entre }i\mbox{ y }k)\times 
(\mbox{\texttt{\#} de vuelos entre }k\mbox{ y }j)\,.
$$

Ahora bien, el n\'umero de vuelos entre $i$ y $j$ es precisamente $\Lambda_{ij}$. La
suma precedente se escribe entonces
$$
\sum_{k}\Lambda_{ik}\Lambda_{kj}\,,
$$
lo cual, comparando con \eqref{sistemas12}, es lo mismo que el coeficiente $(i,j)$ de
la matriz al cuadrado:
$$
\Lambda^2_{ij}=\sum_{k}\Lambda_{ik}\Lambda_{kj}\,.
$$

Una generalizaci\'on de este argumento nos lleva a observar que la entrada $(i,j)$ de
la matriz potencia $p-$\'esima $\Lambda^p$ proporciona el n\'umero de vuelos entre
los nodos $i$ y $j$ que tienen exactamente $p-1$ escalas intermedias. 
Si el lector intenta determinar este n\'umero directamente para un valor incluso tan 
peque\~no como $p=3$, enumerando las posibles rutas sobre el grafo, se convencer\'a
enseguida de la utilidad del producto matricial.

Por otra parte, es muy f\'acil hacer este tipo de c\'alculos en Maxima. 
Para construir una matriz
se usa el comando \texttt{matrix}\index{Matriz!\texttt{matrix}} y se dan como argumentos
las \emph{filas} de la matriz en orden, empezando por la primera (la superior):
\begin{maxcode}

\noindent
\begin{minipage}[t]{8ex}\color{red}\bf
(\% i1) 
\end{minipage}
\begin{minipage}[t]{\textwidth}\color{blue}
L:matrix([0,2,0,1,3],[2,0,1,0,1],[0,1,0,1,0],[1,0,1,0,2],[3,1,0,2,0]);
\end{minipage}
\[\displaystyle
\tag{L}
\begin{pmatrix}0 & 2 & 0 & 1 & 3\\
2 & 0 & 1 & 0 & 1\\
0 & 1 & 0 & 1 & 0\\
1 & 0 & 1 & 0 & 2\\
3 & 1 & 0 & 2 & 0\end{pmatrix}\mbox{}
\]

\noindent
\begin{minipage}[t]{8ex}\color{red}\bf
(\% i2)
\end{minipage}
\begin{minipage}[t]{\textwidth}\color{blue}
L.L;
\end{minipage}
\[\displaystyle
\tag{\% o2} 
\begin{pmatrix}14 & 3 & 3 & 6 & 4\\
3 & 6 & 0 & 5 & 6\\
3 & 0 & 2 & 0 & 3\\
6 & 5 & 0 & 6 & 3\\
4 & 6 & 3 & 3 & 14\end{pmatrix}\mbox{}
\]
\end{maxcode}

Observemos que el producto de matrices se denota por un punto ordinario `.', el
asterisco `$\ast$' se usa para el producto de Hadamard (componente a componente):
\begin{maxcode}

\noindent
\begin{minipage}[t]{8ex}\color{red}\bf
(\% i3)
\end{minipage}
\begin{minipage}[t]{\textwidth}\color{blue}
L*L;
\end{minipage}
\[\displaystyle
\tag{\% o3} 
\begin{pmatrix}0 & 4 & 0 & 1 & 9\\
4 & 0 & 1 & 0 & 1\\
0 & 1 & 0 & 1 & 0\\
1 & 0 & 1 & 0 & 4\\
9 & 1 & 0 & 4 & 0\end{pmatrix}\mbox{}
\]
\end{maxcode}

\subsection{Matrices en criptograf\'ia}\label{hill}

En 1929, el matem\'atico Lester S. Hill\index{Hill, Lester S.} public\'o un art\'iculo en la 
revista American Mathematical Monthly\footnote{Puede leerse gratuitamente en l\'inea en la
URL 
\url{http://links.jstor.org/sici?sici=0002-9890\%28192906\%2F07\%2936\%3A6\%3C306\%3ACIAAA\%3E2.0.CO\%3B2-J}.} 
en el que describ\'ia varios m\'etodos de encriptaci\'on. Uno de ellos, conocido hoy en 
d\'ia como \emph{cifrado de Hill}\index{Criptograf\'ia!cifrado de Hill}
ten\'ia la novedad de que las matem\'aticas empleadas ya no eran s\'olo la parte m\'as 
elemental de la teor\'ia de n\'umeros (aritm\'etica modular), sino que aparec\'ian objetos 
como matrices. 

Para introducir t\'ecnicas matem\'aticas en criptograf\'ia, resulta conveniente sustituir letras por n\'umeros de acuerdo con la siguiente tabla\label{tablacrypto} 
$$\begin{array}{c|c|c|c|c|c|c|c|c|c|c|c|c}
a & b & c & d & e & f & g & h & i & j & k & l & m \\\hline
0 & 1 & 2 & 3 & 4 & 5 & 6 & 7 & 8 & 9 & 10 & 11 & 12 

\end{array}
$$

$$\begin{array}{c|c|c|c|c|c|c|c|c|c|c|c|c}
n & o & p & q & r & s & t & u & v & w & x & y & z\\\hline
13& 14 & 15 & 16 & 17 & 18 & 19 & 20 & 21 & 22 & 23 & 24 & 25

\end{array}
$$

De esta manera, el alfabeto habitual $\mathcal{A}$ se sustituye por $\mathbb{Z}_{26}$.
El m\'etodo de encriptaci\'on que usemos comenzar\'a traduciendo el mensaje a encriptar
(o \emph{texto claro}\index{Texto claro}) en un conjunto de n\'umeros en
$\mathbb{Z}_{26}$, realizar\'a determinadas operaciones matem\'aticas
m\'odulo $26$ sobre ellos, de manera que el resultado siempre ser\'a un n\'umero
comprendido entre $0$ y $25$, y lo expresar\'a como una cadena de texto (el \emph{criptotexto}\index{Criptotexto} o texto encriptado) que, tradicionalmente, se
escribe en may\'usculas para distinguirlo del texto claro.

El criptosistema de Hill se basa en separar el mensaje \texttt{m} en bloques de texto de longitud $k$, y en cada uno de ellos aplicar una transformaci\'on mediante  una matriz $k\times k$, con operaciones m\'odulo $26$.
Por ejemplo, supongamos que queremos encriptar \texttt{m=hola}, usando una matriz
$2\times 2$ que funciona como clave:

$$
K=\begin{pmatrix}
3 & 2\\
5 & 3
\end{pmatrix}\,.
$$

Separamos el texto del mensaje en bloques de 2 caracteres: \texttt{'ho'} y \texttt{'la'}. 
Aplicamos la transformaci\'on indicada por $K$ al primer bloque,
$$K\begin{pmatrix}
h\\
o
\end{pmatrix}
=\begin{pmatrix}
3 & 2\\
5 & 3
\end{pmatrix}
\begin{pmatrix}
7\\
14
\end{pmatrix}
=\begin{pmatrix}
21+28\\
35+42
\end{pmatrix}
\equiv \begin{pmatrix}
23\\
25
\end{pmatrix}
 \mod\, 26 = \begin{pmatrix}
 X\\
 Z
 \end{pmatrix}\,,
 $$
y al segundo:
$$K\begin{pmatrix}
l\\
a
\end{pmatrix}
=\begin{pmatrix}
3 & 2\\
5 & 3
\end{pmatrix}
\begin{pmatrix}
11\\
0
\end{pmatrix}
=\begin{pmatrix}
33\\
55
\end{pmatrix}
\equiv \begin{pmatrix}
7\\
3
\end{pmatrix}
 \mod\, 26 = \begin{pmatrix}
 H\\
 D
 \end{pmatrix}\,.
 $$

Por tanto, el criptotexto $C$ que corresponde a \texttt{m=hola} mediante la matriz $K$ es $$
C=XZHD\,.
$$

La desencriptaci\'on se lleva a cabo realizando el mismo proceso descrito sobre el criptotexto, pero con la matriz inversa $K^{-1}$ m\'odulo $26$. ?`C\'omo se calcula?
Si revisamos el proceso seguido en la p\'agina \pageref{sistemas13} para calcular la
inversa de una matriz particular, nos convenceremos de que es posible hacer lo mismo
con una matriz arbitraria. Sin embargo, nosotros no procederemos as\'i, por el
contrario, utilizaremos las capacidades simb\'olicas de Maxima para obtener la inversa
(sin tener en cuenta el m\'odulo) como sigue:

\begin{maxcode}
\noindent
\begin{minipage}[t]{4em}\color{red}\bf
(\% i1) 
\end{minipage}
\begin{minipage}[t]{\textwidth}\color{blue}
M:matrix([a,b],[c,d]);

\end{minipage}
\[\displaystyle \tag{M}
\begin{pmatrix}a & b\\
c & d\end{pmatrix}\mbox{}
\]

\noindent
\begin{minipage}[t]{4em}\color{red}\bf
(\% i2)
\end{minipage}
\begin{minipage}[t]{\textwidth}\color{blue}
invert(M);

\end{minipage}
\[\displaystyle \tag{\% o2} 
\begin{pmatrix}\frac{d}{a d-b c} & -\frac{b}{a d-b c}\\
-\frac{c}{a d-b c} & \frac{a}{a d-b c}\end{pmatrix}\mbox{}
\]
\end{maxcode}

La cantidad $ad-bc$ es muy importante en el estudio del \'algebra
	matricial, se conoce como el \emph{determinante}\index{Matriz!determinante}
	de la matriz $A$, denotado $\det (A)$. A la vista de la expresi\'on precedente, 
	es claro que \'unicamente
	podremos calcular la inversa de $A$ si su determinante es distinto de cero
	mientras que, rec\'iprocamente, si $\det (A)=ad-bc\neq 0$, \texttt{(\%o2)} nos
	proporciona la matriz inversa $A^{-1}$. Ahora, podemos calcular la inversa
	de una matriz m\'odulo $26$ teniendo en cuenta que
	$a+b \mod\, n =a\mod\, n +b\mod\, n$ y $ab\mod\, n=(a\mod\, n)(b\mod\, n)$.
	Para ello, creamos la siguiente funci\'on (en la que, obs\'ervese, primero se
	comprueba si el determinante de $A$ es invertible m\'odulo el n\'umero $n$
	seleccionado, imprimiendo un mensaje de error si no lo es):
\begin{maxcode}
\noindent
\begin{minipage}[t]{4em}\color{red}\bf
(\% i3)
\end{minipage}
\begin{minipage}[t]{\textwidth}\color{blue}
minvmod(M,n):=block([dd:inv\_mod(determinant(M),n),A,B,C,D],

if is(not(dd)) then return("La matriz no es invertible"),

A:mod(M[2,2]*dd,n),

B:-mod(M[1,2]*dd,n),

C:-mod(M[2,1]*dd,n),

D:mod(M[1,1]*dd,n),

return(matrix([A,B],[C,D]))

)\$
\end{minipage}
\end{maxcode}

Ahora podemos calcular la inversa de la matriz $K$ m\'odulo $26$ sin m\'as que
ejecutar los siguientes comandos:
\begin{maxcode}
\noindent
\begin{minipage}[t]{4em}\color{red}\bf
(\% i4)
\end{minipage}
\begin{minipage}[t]{\textwidth}\color{blue}
K:matrix([3,2],[5,3]);
\end{minipage}
\[\displaystyle \tag{K}
\begin{pmatrix}3 & 2\\
5 & 3\end{pmatrix}\mbox{}
\]

\noindent
\begin{minipage}[t]{4em}\color{red}\bf
(\% i5)
\end{minipage}
\begin{minipage}[t]{\textwidth}\color{blue}
minvmod(K,26);
\end{minipage}
\[\displaystyle \tag{\% o5} 
\begin{pmatrix}23 & -24\\
-21 & 23\end{pmatrix}\mbox{}
\]
\end{maxcode}
Podemos ver que el resultado es correcto viendo que $K^{-1}.K=\mathrm{I}$:
\begin{maxcode}
\noindent
\begin{minipage}[t]{4em}\color{red}\bf
(\% i6)
\end{minipage}
\begin{minipage}[t]{\textwidth}\color{blue}
mod(\%.K,26);

\end{minipage}
\[\displaystyle \tag{\% o6} 
\begin{pmatrix}1 & 0\\
0 & 1\end{pmatrix}\mbox{}
\]
\end{maxcode}
y tambi\'en $K.K^{-1}=\mathrm{I}$:
\begin{maxcode}
\noindent
\begin{minipage}[t]{4em}\color{red}\bf
(\% i7)
\end{minipage}
\begin{minipage}[t]{\textwidth}\color{blue}
mod(K.\%th(2),26);

\end{minipage}
\[\displaystyle \tag{\% o7} 
\begin{pmatrix}1 & 0\\
0 & 1\end{pmatrix}\mbox{}
\]
\end{maxcode}

Como ejemplo, supongamos que recibimos el criptotexto $XZHD$, que sabemos ha sido encriptado con la matriz $K$. Como ya se mencion\'o, para desencriptarlo procedemos
como en la encriptaci\'on: dividiendo el mensaje en bloques de dos letras tendremos
$XZ$ y $HD$ que podemos traducir num\'ericamente en los vectores
$$
\begin{pmatrix}
23 \\ 25
\end{pmatrix}\quad\mbox{ y }\quad
\begin{pmatrix}
7 \\ 3
\end{pmatrix}\,.
$$
Multiplicando por la matriz inversa m\'odulo $26$ y reduciendo el resultado m\'odulo $26$
resulta
$$
\begin{pmatrix}
23 & -24 \\
-21 & 23
\end{pmatrix}
\begin{pmatrix}
23 \\ 25
\end{pmatrix}=
\begin{pmatrix}
-71 \\ 92
\end{pmatrix}\simeq
\begin{pmatrix}
7 \\ 14
\end{pmatrix}\mod\,26 =
\begin{pmatrix}
h \\ o
\end{pmatrix}
$$
y tambi\'en
$$
\begin{pmatrix}
23 & -24 \\
-21 & 23
\end{pmatrix}
\begin{pmatrix}
7 \\ 3
\end{pmatrix}=
\begin{pmatrix}
89 \\ -78
\end{pmatrix}\simeq
\begin{pmatrix}
11 \\ 0
\end{pmatrix}\mod\,26 =
\begin{pmatrix}
l \\ a
\end{pmatrix}
$$
con lo cual recuperamos el mensaje original: \texttt{hola}.

\section{Ejercicios}

\begin{enumerate}
\setcounter{enumi}{156}

\item Resolver los siguientes sistemas de ecuaciones.  
\begin{enumerate}[(a)]
\item\
$$
\systeme{x - 2y = -5,-3x + 4y =  6\,.}
$$

\item\
$$
\systeme{2x -  5y = 8,-4x + 10y = 6\,.}
$$

\item\
$$
\systeme{-5x + 5y = 10,4x - 4y = -8\,.}
$$
\end{enumerate}
Graficar las ecuaciones de cada uno de los incisos anteriores. ?`Qu\'e se
puede decir de las soluciones correspondientes?
\item Resolver los siguientes sistemas de ecuaciones.
\begin{enumerate}[(a)]
\item\
$$
\systeme{
 x - 2y + 3z  =  2,
5x + 8y - 2z  =  6,
2x + 6y - 4z  = -4\,.}
$$
\item\
$$
\begin{cases}
\phantom{(2-i)} x + (4+i)y  = -3 \\[4pt]
(2-i)x + \phantom{(221}3y  =  5i\,.
\end{cases}
$$
\end{enumerate}
\item Calcular la temperatura en los puntos interiores $T_1, T_2, T_3$ de
      la l\'amina  mostrada a continuaci\'on (la temperatura en cada punto
      interior es el promedio de la temperatura en los puntos adyacentes: norte,
      sur, este y oeste).
      \begin{center}
      \begin{tikzpicture}
  \fill [opacity=0.75,gray!25]
        (0,0) --node[below,black]{$100\degree$} (6,0) -- (6,2) --
        node[right,black]{$0\degree$} (6,4)-- (6,8) -- (0,2) -- cycle;      
  \fill [opacity=0] (0,0) circle [radius=\rr] coordinate (p0);     
  \fill [] (2,0) circle [radius=\rr] coordinate (p1);
  \fill [] (4,0) circle [radius=\rr] coordinate (p2);
  \fill [opacity=0] (6,0) circle [radius=\rr] coordinate (p3);         
  \fill [] (0,2) circle [radius=\rr] coordinate (p4)node[above left]{$60\degree$};
  \fill [] (2,2) circle [radius=\rr] coordinate (T2)node[below right]{$T_2$};
  \fill [] (4,2) circle [radius=\rr] coordinate (T3)node[below right]{$T_3$};
  \fill [] (6,2) circle [radius=\rr] coordinate (p5);
  \fill [] (2,4) circle [radius=\rr] coordinate (p6)node[above left]{$40\degree$};
  \fill [] (4,4) circle [radius=\rr] coordinate (T1)node[below right]{$T_1$};
  \fill [] (6,4) circle [radius=\rr] coordinate (p7); 
  \fill [] (4,6) circle [radius=\rr] coordinate (p8)node[above left]{$20\degree$}; 
  \fill [opacity=0] (6,8) circle [radius=\rr] coordinate (p9);
  \draw[very thick] (p0) -- (p3) -- (p9) -- (p4) -- cycle;
  \draw (p1) -- (T2) -- (p6);
  \draw (p2) -- (T3) -- (T1) -- (p8);
  \draw (p4) -- (T2) -- (T3) -- (p5);
  \draw (p6) -- (T1) -- (p7);         
\end{tikzpicture}
      \end{center}
\item Determinar la ecuaci\'on del polinomio $p(x) = a_2 x^2 + a_1 x + a_0$
      que pasa por los puntos $A(1,-2)$, $B(3,3)$ y $C(4,2)$.
\item Para una matrix dada $M$, se define su \emph{traspuesta}\index{Matriz!traspuesta}
$M^t$ como la obtenida intercambiando las filas y columnas de $M$ (es decir: la primera
fila de $M^t$ es la primera columna de $M$, la segunda
fila de $M^t$ es la segunda columna de $M$, etc.) Dadas las matrices
\[
A =  \begin{pmatrix}
             3 & -1 \\
            -2 &  4
           \end{pmatrix}
\,,\
B =  \begin{pmatrix}
             4 & 3 \\
             3 & 2
           \end{pmatrix}\,,
\]
y
$$
C = \begin{pmatrix}
            -4 &  2 \\
             0 & -1 \\
             3 &  6
           \end{pmatrix}\,,\
D = \begin{pmatrix}
             5 &  1 \\
            -2 & -3 \\
             1 &  4          
           \end{pmatrix}\,,
$$
se pide:
\begin{enumerate}[(a)]
\item Calcular $-2A + 3B,\ DA,\ AC^t,\ (A+B)^{-1}$ y $A^{-1}+B^{-1}$.
\item Calcular $AB$ y $BA$. ?`Se puede concluir algo respecto de la
      conmutatividad del producto de matrices?  
\item Recordemos que el determinante\index{Matriz!determinante} de la matriz
      \[
       M =  \begin{pmatrix}
                   a & b \\
                   c & d
             \end{pmatrix}
      \]
      se define como
      \[
      \det(M) = ad -bc\,.
      \]
      Calcular $\det(A),\ \det(B),\ \det(3A)$ y $\det(A+B)$.
      
      ?`Se puede concluir algo respecto del determinante de productos y
       sumas de matrices?
\item Comprobar que
      \begin{enumerate}
      \item $(C+D)^t = C^t + D^t$     
      \item $(CA)^t = A^t C^t$   
      \item $D(AB) = (DA)B)$
      \item $(AB)^{-1} = B^{-1}A^{-1}$
      \item $\det(AB) = \det(A)\det(B)$
      \item $\det(A^t) = \det(A)$
      \item $\det(A^{-1}) = 1/\det(A)$
      \end{enumerate}
\item Una matriz $M$ es \emph{sim\'etrica}\index{Matriz!sim\'etrica} si $M^t = M$.
      \begin{enumerate}
      \item ?`Cu\'ales de las matrices $A,\ B,\ C$ y $D$ son sim\'etricas?
      \item Supongamos que la matriz $M$ cumple que $M^t = -M$ (se dice en este caso
      que es \emph{antisim\'etrica}\index{Matriz!antisim\'etrica}). ?`Qu\'e
            coeficientes de $M$ son nulos?
      \end{enumerate}
\end{enumerate}
\item Dadas la matrices
      \[
      A = \begin{pmatrix}
                  2i & 3+2i \\
                   1 & -3+i
                 \end{pmatrix}\,,\
      B = \begin{pmatrix}
                  3+i & 2i \\
                  5-i & -4i
                 \end{pmatrix}\,\mbox{ y }\
      C = \begin{pmatrix}
                   4-i  & 5+2i & 6 \\
                  -2+2i & 1+i  & 5i
                 \end{pmatrix}\,,         
      \]
      calcular $(2A + iB)C$.
\item Sean
      \[
      A = \begin{pmatrix}
                  1 & 2 & 3 \\
                  2 & 1 & -3
                 \end{pmatrix}\,\
      B = \begin{pmatrix}
                  3 & 2 & -4 \\
                  5 & 4 & -3
                 \end{pmatrix}\,\mbox{ y }\
      C = \begin{pmatrix}
                   4 & 5 & 6 \\
                  -2 & 1 & 5
                 \end{pmatrix}\,.
\]
Resolver la ecuaci\'on matricial
\[
\frac{1}{3}(X + A) = 4(2X + B) - C\,.
\]
\item Encontrar una matriz $X$ de tama\~no $3 \times 3$ que satisfaga la
      ecuaci\'on propuesta en cada caso:
      \begin{enumerate}[(a)]
      \item $X^2 - 5X + 6I = 0$.
      \item $X^3 + 2X^2 + X + 2I = 0$.
      \end{enumerate}
\item Comprobar que la matriz
      \[
      A = \begin{pmatrix}
                  1 & -1 &  0 \\
                  0 &  1 & -1 \\
                  1 &  0 &  1
                 \end{pmatrix}
      \]
      satisface la ecuaci\'on $A^3 -3A^2 + 3A - 2I = 0$.
      Calcular $A^{-1}$ a partir de este hecho (esto es un ejemplo de un resultado
      general conocido como \emph{Teorema de Cayley-Hamilton}\index{Cayley-Hamilton!teorema}).
\item Una \emph{permutaci\'on}\index{Permutaci\'on} del conjunto
      $J_n \doteq \{1, 2, \ldots, n\}$ es una funci\'on biyectiva\index{Funci\'on!biyectiva}
      $\sigma: J_n \to J_n$. La \emph{matriz de permutaci\'on}\index{Matriz!de permutaci\'on}
      de $\sigma$ es la matriz $S = [s_{ij}]$ de tama\~no $n \times n$ tal que
      \[
      s_{ij} = \left\{ \begin{array}{ll}
                       1 & \mbox{si } \sigma(i) = j\\[6pt]
                       0 & \mbox{de otra forma} 
                       \end{array}
               \right.
      \]
      Por ejemplo, la permutaci\'on $\sigma: J_4 \to J_4$ dada por
      \begin{align*}
      \sigma(1) & = 1 \\
      \sigma(2) & = 3 \\
      \sigma(3) & = 4 \\
      \sigma(4) & = 2 \\
      \end{align*}
      tiene matriz de permutaci\'on
      \[
      S = \begin{pmatrix}
                  1 & 0 & 0 & 0 \\
                  0 & 0 & 1 & 0 \\
                  0 & 0 & 0 & 1 \\
                  0 & 1 & 0 & 0
           \end{pmatrix}\,.
      \]
      \begin{itemize}
      \item[(a)] Comprobar que $S^t$ es tambi\'en una matriz de
                 permutaci\'on de $J_4$. ?`C\'omo est\'an relacionadas las
                 permutaciones asociadas con $S$ y $S^t$?
      \item[(b)] Si la matriz de permutaci\'on $S$ es sim\'etrica, ?`qu\'e
                 puede decirse de la permutaci\'on $\sigma$ asociada con
                 $S$?
      \end{itemize}
\item Dado el grafo $\Gamma$ representado en la figura siguiente, se pide:

\begin{center}
\begin{tikzpicture}
\Vertex[IdAsLabel,Math,opacity =.5,y=2]{P1} 
\Vertex[IdAsLabel,Math,opacity =.5,y=1,x=2]{P2}
\Vertex[IdAsLabel,Math,opacity =.5,y=-1,x=2]{P3}
\Vertex[IdAsLabel,Math,opacity =.5,y=-2]{P4}
\Vertex[IdAsLabel,Math,opacity =.5,y=-1,x=-2]{P5}
\Vertex[IdAsLabel,Math,opacity =.5,y=1,x=-2]{P6}
\Edge(P1)(P3)
\Edge(P1)(P5)
\Edge(P2)(P4)
\Edge(P2)(P6)
\Edge(P3)(P5)
\Edge(P4)(P6)
\end{tikzpicture}
\end{center}

      \begin{enumerate}[(a)]
      \item Determinar la matriz de incidencia de $\Gamma$.
      \item ?`Cu\'antos caminos de longitud 3 existen para ir de
                 $P_1$ a $P_3$?
      \item Un \emph{ciclo}\index{Grafo!ciclo} en un grafo es un camino cerrado 
      			 que comienza y termina en un mismo v\'ertice. ?`Cu\'antos ciclos de
                 longitud 4 tiene $\Gamma$?
      \end{enumerate}
\item Mediante el m\'etodo de Hill, usando como clave la matriz
      \[
      K = \begin{pmatrix}
                  1 & 0 & 0 \\
                  1 & 6 & 1 \\
                  6 & 3 & 6
                 \end{pmatrix}\,,
      \]
      se pide:
      \begin{enumerate}[(a)]
      \item Encriptar el mensaje\footnote{El lector curioso quiz\'as agradezca saber que
      Luzbel\index{Luzbel} es el nombre de una banda de heavy metal muy popular en M\'exico en
      los an\~os 80, poseedora de una fiel base de fans. Y quiz\'as 
      tambi\'en le interese conocer el nombre de uno de sus miembros fundadores consultando
      la p\'agina de Wikipedia \url{https://es.wikipedia.org/wiki/Luzbel_(banda)}.} \texttt{luzbel}.
      \item Desencriptar el mensaje resultante del inciso precedente.
      \end{enumerate}

\end{enumerate}

\section{Uso de la tecnolog\'ia de c\'omputo}\label{imgproc}

En esta secci\'on veremos algunas aplicaciones muy b\'asicas de las matrices en el 
procesamiento de im\'agenes. La idea central es que una imagen no es m\'as que una matriz de
\emph{p\'ixeles}\index{P\'ixel}, de manera que cada posible operaci\'on algebraica con matrices
debe corresponder a una transformaci\'on de la imagen. Un p\'ixel es la unidad m\'inima de una 
imagen que puede representar la pantalla de una computadora, generalmente con la forma de un
cuadrado (o rect\'angulo) cuyo tama\~no depende de la \emph{resoluci\'on} con la que se est\'e
trabajando. Por ejemplo, una computadora con una pantalla cuya resoluci\'on es de
$1920\times 1280$, podr\'a representar precisamente $1920\cdot 1280 = 2\, 457\, 600$ p\'ixeles;
cualquier imagen se mostrar\'a en la pantalla como una composici\'on de esos $2$ millones y
medio de `cuadraditos', pero esa densidad es tan elevada que el ojo humano no los distinguir\'a 
y la imagen parecer\'a un continuo. Nosotros trabajaremos con im\'agenes formadas por un 
n\'umero peque\~no de p\'ixeles por motivos pr\'acticos (introduciremos la matriz que define
la imagen `a mano'), lo cual har\'a que las im\'agenes se ver\'an `angulosas' (o `pixeladas'), 
pero el lector no tendr\'a dificultad en convencerse de que cuanto mayor sea el n\'umero de 
p\'ixeles, m\'as `suavizada' se ver\'a la imagen.

En el siguiente comando, cargamos el paquete \texttt{draw} e instru\'imos a Maxima para que
omita mostrar algunos mensajes engorrosos (el comando es v\'alido en Linux o Mac OS, en un sistema Windows, debe sustituirse la cadena \texttt{/dev/null} por \texttt{NIL}).

\begin{maxcode}
\noindent
\begin{minipage}[t]{8ex}\color{red}\bf
(\% i1) 
\end{minipage}
\begin{minipage}[t]{\textwidth}\color{blue}
with\_stdout(``/dev/null'',load(draw))\$
\end{minipage}
\end{maxcode}

Comenzamos creando una imagen \emph{binaria}\index{Imagen!binaria} 
(cada p\'ixel solo puede ser blanco o negro)\footnote{Las
im\'agenes a color se forman de la misma manera, como matrices, pero en las cuales los
coeficientes pueden tomar valores entre $0$ y $2^b$, donde $b$ es el n\'umero de bits utilizados para representar el color. Por ejemplo, una imagen en color de $16$ bits, usa una matriz con coeficientes cuyos valores est\'an comprendidos entre $0$ y $2^{16}=65\,536$.}, representada por una matriz cuyas etiquetas $(i,j)$ representan la ubicaci\'on de los p\'ixeles, y cuyo valor $m_{ij}$ es $0$ si el p\'ixel es negro y $1$ si es blanco: 

\begin{maxcode}
\noindent
\begin{minipage}[t]{8ex}\color{red}\bf
(\% i2)
\end{minipage}
\begin{minipage}[t]{\textwidth}\color{blue}
felix:matrix(\newline
{\tt\small [1,1,1,1,1,1,0,1,1,1,1,1,1,1,1,1,1,1,1,1,1,1,1,1,0,1,1,1,1,1,1,1,1,1,1],
[1,1,1,1,1,0,0,0,1,1,1,1,1,1,1,1,1,1,1,1,1,1,1,0,0,0,1,1,1,1,1,1,1,1,1],
[1,1,1,1,1,0,0,0,0,1,1,1,1,1,1,1,1,1,1,1,1,1,1,0,0,0,1,1,1,1,1,1,1,1,1],
[1,1,1,1,1,0,0,0,0,0,1,1,1,1,1,1,1,1,1,1,1,1,0,0,0,0,1,1,1,1,1,1,1,1,1],
[1,1,1,1,1,0,0,0,0,0,0,1,1,1,1,1,1,1,1,1,1,0,0,0,0,0,0,1,1,1,1,1,1,1,1],
[1,1,1,1,1,0,0,0,0,0,0,0,1,0,0,0,0,0,0,0,0,0,0,0,0,0,0,0,0,0,1,1,1,1,1],
[1,1,1,1,1,0,0,0,0,0,0,0,0,0,0,0,0,0,0,0,0,0,0,0,0,0,0,0,0,0,1,1,1,1,1],
[1,1,1,1,1,0,0,0,0,0,0,0,0,0,1,1,1,1,1,1,1,1,0,0,0,0,1,1,0,0,1,1,1,1,1],
[1,1,1,1,1,0,0,0,0,0,0,0,0,1,1,1,1,1,1,1,1,1,1,0,0,0,1,1,1,0,0,1,1,1,1],
[1,1,1,1,1,0,0,0,0,0,0,0,1,1,1,1,1,1,1,1,1,1,1,1,0,0,1,1,1,1,0,0,1,1,1],
[1,1,1,1,0,0,0,0,0,0,0,0,1,1,1,1,1,1,1,1,1,1,1,1,0,0,1,1,1,1,1,0,1,1,1],
[1,1,1,1,0,0,0,0,0,0,0,1,1,1,1,1,1,1,1,1,1,1,1,1,1,0,1,1,1,1,1,1,0,1,1],
[1,1,1,1,0,0,0,0,0,0,0,1,1,1,1,1,1,1,1,1,1,1,1,1,1,0,1,1,1,0,0,1,0,0,1],
[1,1,1,1,0,0,0,0,0,0,1,1,1,1,1,1,1,1,1,0,0,0,1,1,1,0,1,1,1,0,0,0,1,0,1],
[1,1,1,0,0,0,0,0,0,0,1,1,1,1,1,1,1,1,0,0,0,0,1,1,1,0,1,1,1,0,0,0,1,0,0],
[1,1,1,0,0,0,0,0,0,0,1,1,1,1,1,1,1,1,0,0,0,0,1,1,1,0,0,1,1,0,0,0,1,0,0],
[1,1,0,0,0,0,0,0,0,0,1,1,1,1,1,1,1,1,0,0,0,0,1,1,1,0,0,1,1,0,0,0,1,0,0],
[1,0,0,0,0,0,0,0,0,0,1,1,1,1,1,1,1,1,0,0,0,0,1,1,1,0,0,1,1,1,0,0,1,1,0],
[1,0,0,0,0,0,0,0,0,0,0,1,1,1,1,1,1,1,1,0,0,0,1,1,1,0,0,0,1,1,1,1,1,1,0],
[0,0,0,0,0,0,0,0,0,0,0,1,1,1,1,1,1,1,1,1,1,1,1,1,0,0,1,0,1,1,1,1,1,1,0],
[1,0,0,0,0,0,0,0,0,0,0,0,1,1,1,1,1,1,1,1,1,1,1,1,0,1,1,0,0,0,0,0,1,0,0],
[1,1,0,0,0,0,0,1,1,1,0,0,0,1,1,1,1,1,1,1,1,1,1,0,0,1,0,0,0,0,0,0,0,1,0],
[1,1,1,0,0,0,1,1,1,1,1,1,0,0,0,1,1,1,1,1,1,1,0,0,1,1,0,0,0,0,0,0,1,0,0],
[1,1,0,0,0,0,1,1,1,0,0,0,1,1,0,0,0,0,0,0,0,0,0,1,1,1,0,0,0,0,0,0,0,0,0],
[1,0,0,0,0,0,1,1,1,0,0,0,1,1,1,1,0,0,0,0,1,1,1,1,1,1,1,0,0,0,0,0,0,0,0],
[1,1,1,1,0,0,1,1,1,1,1,0,0,1,1,1,1,1,1,1,1,1,1,1,1,1,1,1,1,1,1,0,1,0,1],
[1,1,1,1,1,0,0,1,1,1,1,1,1,0,0,0,1,1,1,1,1,1,1,1,1,1,1,1,1,0,0,1,0,0,1],
[1,1,1,1,1,0,0,1,1,1,1,1,1,0,0,0,0,0,0,0,1,1,1,1,0,0,0,0,0,0,1,1,0,1,1],
[1,1,1,1,1,1,0,0,1,1,1,1,1,1,0,0,0,0,0,0,0,0,0,0,0,0,0,0,0,0,1,0,1,1,1],
[1,1,1,1,1,1,1,0,0,1,1,1,1,1,0,0,0,1,1,1,0,0,0,0,0,0,0,0,0,1,0,0,1,1,1],
[1,1,1,1,1,1,1,1,0,0,1,1,1,1,1,0,1,1,1,1,1,1,0,1,1,1,0,0,1,0,0,1,1,1,1],
[1,1,1,1,1,1,1,1,1,0,0,1,1,1,1,1,0,0,1,1,1,1,1,1,1,0,0,1,0,0,1,1,1,1,1],
[1,1,1,1,1,1,1,1,1,1,0,0,0,1,1,1,1,0,0,0,0,0,0,0,0,0,0,0,0,1,1,1,1,1,1],
[1,1,1,1,1,1,1,1,1,1,1,1,0,0,0,0,0,1,1,1,1,1,0,0,0,0,0,1,1,1,1,1,1,1,1],
[1,1,1,1,1,1,1,1,1,1,1,1,1,1,0,0,0,0,0,0,0,0,0,0,0,1,1,1,1,1,1,1,1,1,1]\\
)\$}
\end{minipage}

\noindent
\begin{minipage}[t]{8ex}\color{red}\bf
(\% i3)
\end{minipage}
\begin{minipage}[t]{\textwidth}\color{blue}
wxplot\_size:[350,350]\$
\end{minipage}
\end{maxcode}
Para representar la matriz como una imagen, en Maxima tenemos el comando \texttt{image} 
(lo que hace este comando no es mas que colorear cada pixel con el n\'umero que indica el 
correspondiente coeficiente de la matriz). Notemos que la imagen de p\'ixeles se 
reetiqueta para comenzar en $(0,0)$ (mientras que en una matriz, el primer coeficiente
es $(1,1)$): 

\begin{maxcode}
\noindent
\begin{minipage}[t]{8ex}\color{red}\bf
(\% i4)
\end{minipage}
\begin{minipage}[t]{\textwidth}\color{blue}
wxdraw2d( 
palette=gray,colorbox=false,proportional\_axes=xy,\\
image(felix,0,0,34,34)
);
\end{minipage}
\[\displaystyle
\tag{\% t4} 
\includegraphics[scale=0.5]{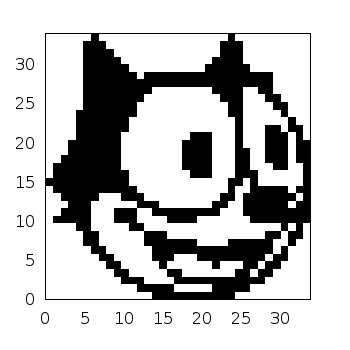}\mbox{}\]
\end{maxcode}

Las transformaciones algebraicas sobre matrices se corresponden con 
transformaciones gr\'aficas sobre la imagen. As\'i, intercambiar el orden de las 
columnas (de derecha a izquierda) invierte la orientaci\'on de la imagen\index{Imagen!orientaci\'on}.  Esto se hace mediante el comando 
\texttt{genmatrix} de Maxima. Se indica la
funci\'on \texttt{u[i,j]} que da la expresi\'on de los coeficientes en t\'erminos de
$(i,j)$, la expresi\'on $35-(j-1)$ aparece porque la etiqueta de los p\'ixeles inicia en $0$. Para la dimensi\'on solo se especifica $(35,35)$ (Maxima
entiende que el \'indice de filas y columnas inicia en $1$). 
\begin{maxcode}
\noindent
\begin{minipage}[t]{8ex}\color{red}\bf
(\% i5)
\end{minipage}
\begin{minipage}[t]{\textwidth}\color{blue}
f[i,j]:=felix[i,35-(j-1)]\$
\end{minipage}

\noindent
\begin{minipage}[t]{8ex}\color{red}\bf
(\% i6)
\end{minipage}
\begin{minipage}[t]{\textwidth}\color{blue}
flipfelix:genmatrix(f,35,35)\$
\end{minipage}

\noindent
\begin{minipage}[t]{8ex}\color{red}\bf
(\% i7)
\end{minipage}
\begin{minipage}[t]{\textwidth}\color{blue}
wxdraw2d( 
palette=gray,colorbox=false,proportional\_axes=xy,\\
image(flipfelix,0,0,34,34)
);
\end{minipage}
\[\displaystyle
\tag{\% t7} 
\includegraphics[scale=0.5]{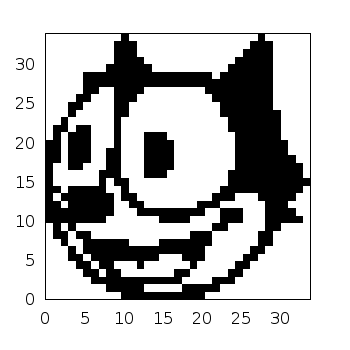}\mbox{}\]
\end{maxcode}

Trasponer\index{Matriz!traspuesta} una matriz (intercambiar filas por columnas) tiene el efecto de reflejar la imagen respecto del eje $Oy$ y luego hacer un giro de 90\degree . 
Maxima dispone de \texttt{transpose} para esto: 
\begin{maxcode}
\noindent
\begin{minipage}[t]{8ex}\color{red}\bf
(\% i8)
\end{minipage}
\begin{minipage}[t]{\textwidth}\color{blue}
wxdraw2d(palette=gray,colorbox=false,proportional\_axes=xy,\\
image(transpose(felix),0,0,34,34));
\end{minipage}
\[\displaystyle
\tag{\% t8} 
\includegraphics[scale=0.5]{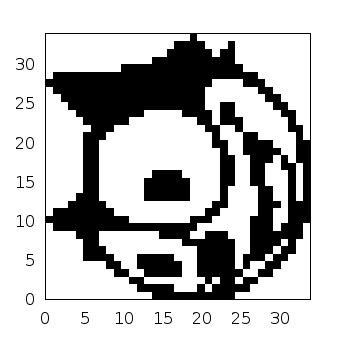}\mbox{}\]
\end{maxcode}

El \emph{negativo}\index{Imagen!negativo} de una imagen binaria resulta
de intercambiar blanco ($1$) y negro ($0$). Para ello, basta con sumar $1$ m\'odulo $2$ 
a cada coeficiente de la matriz original:

\begin{maxcode}
\noindent
\begin{minipage}[t]{8ex}\color{red}\bf
(\% i9)
\end{minipage}
\begin{minipage}[t]{\textwidth}\color{blue}
u[i,j]:=1\$
\end{minipage}

\noindent
\begin{minipage}[t]{8ex}\color{red}\bf
(\% i10)
\end{minipage}
\begin{minipage}[t]{\textwidth}\color{blue}
unos:genmatrix(u,35,35)\$
\end{minipage}

\noindent
\begin{minipage}[t]{8ex}\color{red}\bf
(\% i11)
\end{minipage}
\begin{minipage}[t]{\textwidth}\color{blue}
felixneg:mod(felix + unos,2)\$
\end{minipage}

\noindent
\begin{minipage}[t]{8ex}\color{red}\bf
(\% i12)
\end{minipage}
\begin{minipage}[t]{\textwidth}\color{blue}
wxdraw2d( 
palette=gray,colorbox=false,proportional\_axes=xy,\\
image(felixneg,0,0,34,34)
);
\end{minipage}
\[\displaystyle
\tag{\% t12} 
\includegraphics[scale=0.5]{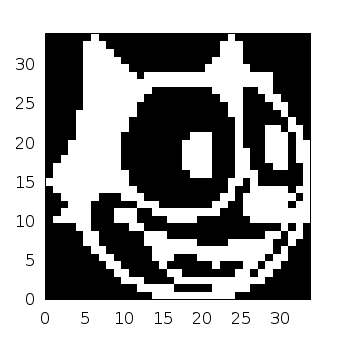}\mbox{}\]
\end{maxcode}

Tambi\'en podemos crear una ventana para recortar partes de la imagen. Una \emph{ventana}\index{Imagen!ventana} no es mas que una imagen de la misma dimensi\'on que la original, con todos sus p\'ixeles negros excepto los de una region rectangular en la que son blancos. Al sumar esta ventana (m\'odulo $2$) con la imagen que desemaos transformar, todo se vuelve negro a excepci\'on de la region rectangular, que permanece sin cambios. Por ejemplo, para nuestra imagen de F\'elix, cuyas dimensiones son $35\times 35$, podemos crear una ventana que seleccione la regi\'on que va desde $(3,3)$ como esquina superior izquierda, hasta $(12,8)$ como esquina inferior derecha. Para lograrlo, primero generamos un fondo negro: 

\begin{maxcode}
\noindent
\begin{minipage}[t]{8ex}\color{red}\bf
(\% i13)
\end{minipage}
\begin{minipage}[t]{\textwidth}\color{blue}
Z:zeromatrix(35,35)\$
\end{minipage}
\end{maxcode}
Y a continuaci\'on la ventana blanca: cambiamos los pixeles que nos interesean de $0$ (negro) a 
$1$ (blanco): 

\begin{maxcode}
\noindent
\begin{minipage}[t]{8ex}\color{red}\bf
(\% i14)
\end{minipage}
\begin{minipage}[t]{\textwidth}\color{blue}
for i:3 thru 12 do (for j:3 thru 30 do (Z[i,j]:1))\$
\end{minipage}

\noindent
\begin{minipage}[t]{8ex}\color{red}\bf
(\% i15)
\end{minipage}
\begin{minipage}[t]{\textwidth}\color{blue}
ventana:Z\$
\end{minipage}
\end{maxcode}

Por \'ultimo, aplicamos el filtro. Notese que \emph{no} usamos el producto matricial (que en
Maxima se denota por .), sino el producto componente a 
componente\index{Matriz!producto!componente a componente} (que se denota por *) 
m\'odulo 2: 

\begin{maxcode}
\noindent
\begin{minipage}[t]{8ex}\color{red}\bf
(\% i16)
\end{minipage}
\begin{minipage}[t]{\textwidth}\color{blue}
felixfiltro:felix * ventana\$
\end{minipage}

\noindent
\begin{minipage}[t]{8ex}\color{red}\bf
(\% i17)
\end{minipage}
\begin{minipage}[t]{\textwidth}\color{blue}
wxdraw2d( 
palette=gray,colorbox=false,proportional\_axes=xy,\\
image(felixfiltro,0,0,34,34)
);
\end{minipage}
\[\displaystyle
\tag{\% t17} 
\includegraphics[scale=0.5]{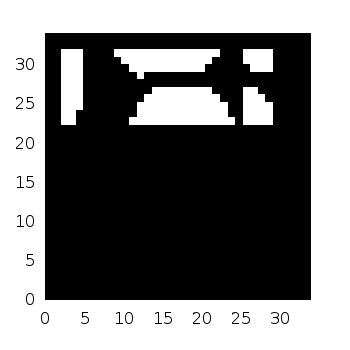}\mbox{}\]
\end{maxcode}

Vamos a crear otra imagen, esta vez con una \emph{escala de grises}\index{Imagen!escala de grises}. Esto quiere decir que usaremos $8$ bits para representar cada coeficiente de la matriz que, por tanto, podr\'a tomar cualquier valor entero entre $0$ (negro) y
$2^8=256$ (blanco).

\begin{maxcode}
\noindent
\begin{minipage}[t]{8ex}\color{red}\bf
(\% i18)
\end{minipage}
\begin{minipage}[t]{\textwidth}\color{blue}
kitty:matrix(\newline
    [1,1,1,1,1,1,1,1,1,1,1,1,1,1,1,1,1,1,1,1,1,1,1,1,1,1,1,1,1,1,1,1,1,1,1],
    [1,1,1,1,1,1,1,1,1,1,1,1,1,1,1,1,1,1,1,1,1,1,1,1,1,1,1,1,1,1,1,1,1,1,1],
    [1,1,1,1,1,1,1,1,1,1,1,1,1,1,1,1,1,1,1,1,1,1,1,1,1,1,1,1,1,1,1,1,1,1,1],
    [1,1,1,1,1,1,1,1,1,1,1,1,1,1,1,1,1,1,1,1,1,1,1,1,1,1,1,1,1,1,1,1,1,1,1],
    [1,1,1,1,1,1,1,1,1,1,1,1,1,1,1,1,1,1,1,1,1,1,1,1,1,1,1,1,1,1,1,1,1,1,1],
    [1,1,1,1,1,1,1,1,1,1,1,1,1,1,1,1,1,1,1,1,1,1,1,1,1,1,1,1,1,1,1,1,1,1,1],
    [1,1,1,1,1,1,1,1,1,1,1,1,1,1,1,1,1,1,1,1,1,1,1,1,1,1,1,1,1,1,1,1,1,1,1],
    [1,1,1,1,1,1,1,1,1,1,1,1,1,1,1,1,1,1,1,1,1,1,1,1,1,1,1,1,1,1,1,1,1,1,1],
    [1,1,1,1,1,1,1,1,1,1,1,1,1,1,1,1,1,1,1,1,1,1,1,1,1,1,1,1,1,1,1,1,1,1,1],
    [1,1,1,1,1,1,1,1,1,1,1,1,1,1,1,1,1,1,1,0,0,1,1,1,1,1,1,1,1,1,1,1,1,1,1],
    [1,1,1,1,1,1,1,1,1,1,1,0,0,0,1,1,1,1,0,0.5,0.5,0,1,0,0,0,1,1,1,1,1,1,1,1,1],
    [1,1,1,1,1,1,1,1,1,1,0,1,1,1,0,0,0,0,0.5,0.5,0.5,0.5,0,1,1,1,0,1,1,1,1,1,1,1,1],
    [1,1,1,1,1,1,1,1,1,1,0,1,1,1,1,1,1,0,0.5,0.5,0.5,0,0,0,1,1,0,1,1,1,1,1,1,1,1],
    [1,1,1,1,1,1,1,1,1,1,0,1,1,1,1,1,1,0,0.5,0.5,0,0.5,0.5,0.5,0,0,0,1,1,1,1,1,1,1,1],
    [1,1,1,1,1,1,1,1,1,1,1,0,1,1,1,1,1,1,0,0,0,0.5,0.5,0.5,0,0.5,0.5,0,1,1,1,1,1,1,1],
    [1,1,1,1,1,1,1,1,1,1,0,1,1,1,1,1,1,1,1,1,1,0,0,0,0.5,0.5,0.5,0,1,1,1,1,1,1,1],
    [1,1,1,1,1,1,1,1,1,1,0,1,1,1,1,1,1,1,1,1,1,1,1,0,0.5,0.5,0,1,1,1,1,1,1,1,1],
    [1,1,1,1,1,1,1,1,1,1,0,1,1,1,1,1,1,1,1,1,1,1,1,1,0,0,1,1,1,1,1,1,1,1,1],
    [1,1,1,1,1,1,1,1,1,0,1,1,1,1,0,1,1,1,1,1,1,0,1,1,1,1,0,1,1,1,1,1,1,1,1],
    [1,1,1,1,1,1,1,0,0,0,0,1,1,1,0,1,1,1,1,1,1,0,1,1,1,0,0,0,0,1,1,1,1,1,1],
    [1,1,1,1,1,1,1,1,1,1,0,0,1,1,1,1,1,0,0,1,1,1,1,1,0,0,1,1,1,1,1,1,1,1,1],
    [1,1,1,1,1,1,1,1,0,0,0,1,1,1,1,1,1,1,1,1,1,1,1,1,1,0,0,0,1,1,1,1,1,1,1],
    [1,1,1,1,1,1,1,1,1,1,1,0,1,1,1,1,1,1,1,1,1,1,1,1,0,1,1,1,1,1,1,1,1,1,1],
    [1,1,1,1,1,1,1,1,1,1,1,1,0,0,1,1,1,1,1,1,1,1,0,0,1,1,1,1,1,1,1,1,1,1,1],
    [1,1,1,1,1,1,1,1,1,1,1,1,1,1,0,0,0,0,0,0,0,0,1,1,1,1,1,1,1,1,1,1,1,1,1],
    [1,1,1,1,1,1,1,1,1,1,1,1,1,1,1,1,1,1,1,1,1,1,1,1,1,1,1,1,1,1,1,1,1,1,1],
    [1,1,1,1,1,1,1,1,1,1,1,1,1,1,1,1,1,1,1,1,1,1,1,1,1,1,1,1,1,1,1,1,1,1,1],
    [1,1,1,1,1,1,1,1,1,1,1,1,1,1,1,1,1,1,1,1,1,1,1,1,1,1,1,1,1,1,1,1,1,1,1],
    [1,1,1,1,1,1,1,1,1,1,1,1,1,1,1,1,1,1,1,1,1,1,1,1,1,1,1,1,1,1,1,1,1,1,1],
    [1,1,1,1,1,1,1,1,1,1,1,1,1,1,1,1,1,1,1,1,1,1,1,1,1,1,1,1,1,1,1,1,1,1,1],
    [1,1,1,1,1,1,1,1,1,1,1,1,1,1,1,1,1,1,1,1,1,1,1,1,1,1,1,1,1,1,1,1,1,1,1],
    [1,1,1,1,1,1,1,1,1,1,1,1,1,1,1,1,1,1,1,1,1,1,1,1,1,1,1,1,1,1,1,1,1,1,1],
    [1,1,1,1,1,1,1,1,1,1,1,1,1,1,1,1,1,1,1,1,1,1,1,1,1,1,1,1,1,1,1,1,1,1,1],
    [1,1,1,1,1,1,1,1,1,1,1,1,1,1,1,1,1,1,1,1,1,1,1,1,1,1,1,1,1,1,1,1,1,1,1],
    [1,1,1,1,1,1,1,1,1,1,1,1,1,1,1,1,1,1,1,1,1,1,1,1,1,1,1,1,1,1,1,1,1,1,1])\$
\end{minipage}
\end{maxcode}

De nuevo utilizamos el comando \texttt{image} para colorear los p\'ixeles de acuerdo con los
valores de los coeficientes en la matriz. N\'otese que las dimensiones de la matriz
\texttt{kitty} son las mismas que las de \texttt{felix}.

\begin{maxcode}
\noindent
\begin{minipage}[t]{8ex}\color{red}\bf
(\% i19)
\end{minipage}
\begin{minipage}[t]{\textwidth}\color{blue}
wxdraw2d( 
palette=gray,colorbox=false,proportional\_axes=xy,\\
image(kitty,0,0,34,34)
);
\end{minipage}
\[\displaystyle
\tag{\% t19} 
\includegraphics[scale=0.5]{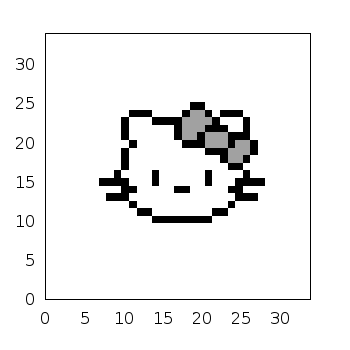}\mbox{}\]
\end{maxcode}

Como ambas im\'agenes tienen el mismo tama\~no ($35\times 35$) podemos 
hacer un \emph{fundido}\index{Imagen!fundido} entre ellas, 
sumando las im\'agenes multiplicadas por un factor de manera que 
en el instante inicial tengamos a F\'elix\index{F\'elix} y en el final a Kitty\index{Kitty}; 
esto se llama en Matem\'aticas 
una \emph{combinaci\'on lineal convexa}\index{Combinaci\'on!lineal convexa}\footnote{Haciendo 
click con el bot\'on derecho del rat\'on puede iniciarse y detenerse la animaci\'on en 
wxMaxima.}: 

\begin{maxcode}
\noindent
\begin{minipage}[t]{8ex}\color{red}\bf
(\% i20)
\end{minipage}
\begin{minipage}[t]{\textwidth}\color{blue}
wxanimate\_draw(t,20,\\
         image((1-t/20)*felix + (t/20)*kitty,0,0,34,34),\\
         palette=gray,colorbox=false,proportional\_axes=xy),\\
wxanimate\_framerate=5;
\end{minipage}
\[\displaystyle
\tag{\% t20} 
\includegraphics[scale=0.55]{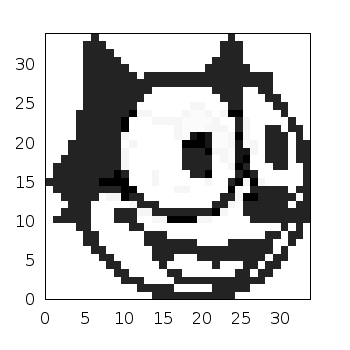}
\mbox{}\]
\end{maxcode}

Sin duda alguna el lector reconocer\'a numerosos ejemplos de videoclips y pel\'iculas donde se hace un uso intensivo de este efecto de animaci\'on. Otros muchos son posibles
utlizando operaciones matem\'aticas sobre matrices m\'as complicadas, motivo por el cual
no los analizaremos aqu\'i, pero esperamos que esta breve introducci\'on sirva como motivaci\'on para su estudio en cursos m\'as avanzados.

\section{Actividades de c\'omputo}
\begin{enumerate}[(A)]

\item Escribir un programa en Maxima que implemente el cifrado de Hill descrito en la 
subsecci\'on \ref{hill}. El programa
debe tomar como argumento una matriz $K$ de tama\~no $2\times 2$ y una cadena de
texto \texttt{s} (en min\'usculas, sin espacios ni guiones), y debe devolver el 
criptotexto como una cadena en may\'usculas.

\item Una de las ventajas del m\'etodo de encriptaci\'on de Hill es que es muy resistente al 
\emph{an\'alisis de frecuencias}\index{Criptograf\'ia!an\'alisis de frecuencias}. 
Para comprender qu\'e significa esto, consideremos primero un m\'etodo de encriptaci\'on 
sencillo, como el de \emph{cifrado af\'in}\index{Criptograf\'ia!cifrado af\'in}.
Supongamos que el mensaje a transmitir es \texttt{envienrefuerzos}. Como siempre, el proceso
de encriptaci\'on inicia transformando el texto claro en una cadena de n\'umeros con la
ayuda de la tabla en la p\'agina \pageref{tablacrypto}:
$$
4,13,21,8,4,13,17,4,5,20,4,17,25,14,18\,.
$$
A continuaci\'on, sobre cada n\'umero se aplica una transformaci\'on af\'in (de ah\'i el nombre) de la forma $f(x)=ax+b$ y el resultado se reduce m\'odulo 2$26$, es decir, la transformaci\'on completa ser\'a
$$
\Phi(x)=f(x)\mod\,26 =(ax+b)\mod\,26\,.
$$
Por ejemplo, la letra \texttt{e} corresponde al n\'umero $4$. Si elegimos $a=1$, $b=3$
(par\'ametros que corresponden al cifrado de Julio C\'esar\index{Criptograf\'ia!cifrado de
Julio C\'esar}\index{C\'esar, Julio}), este n\'umero se transforma en
$\Phi(4)=4+3\mod\,26=7\mod\,26=7$, que corresponde a la letra \texttt{H}. De hecho, el 
cifrado de Julio C\'esar tambi\'en se conoce como \emph{cifrado de desplazamiento}\index{Criptograf\'ia!cifrado de desplazamiento}
porque corresponde a desplazar las letras tres lugares a su derecha en el 
alfabeto\footnote{Nos dice Suetonio en su \emph{Vida de los C\'esares}, p\'arrafo 56:
\begin{quote}
\emph{Si ten\'ia que decir algo confidencial, lo escrib\'ia cifrado, esto es, cambiando el orden de las letras del alfabeto para que ni una palabra pudiera entenderse. Si alguien quiere decodificarlo y entender su significado, debe sustituir la cuarta letra del alfabeto, es decir, la D, por la A, y as\'i con las dem\'as}.
\end{quote}} La
siguiente funci\'on en Maxima toma un mensaje dado como una cadena de texto (que en
Maxima se indica escribi\'endolo entre comillas), junto con los par\'ametros $a$ y $b$,
y devuelve el mensaje encriptado:
\begin{maxcode}
\noindent
\begin{minipage}[t]{4em}\color{red}\bf
(\% i1)
\end{minipage}
\begin{minipage}[t]{\textwidth}\color{blue}
cifrado\_afin(s,a,b):=block([tmp1,tmp2,tmp3],

tmp1:map(cint,charlist(s))-97,

tmp2:mod(a*tmp1+b,26),

tmp3:map(ascii,tmp2+65),

simplode(tmp3))\$
\end{minipage}
\noindent%
\end{maxcode}
De esta manera, la encriptaci\'on de \texttt{notenemosreservadeagua} con los par\'ametros de 
Julio C\'esar es
\begin{maxcode}
\noindent
\begin{minipage}[t]{8ex}{\color{red}\bf
\begin{verbatim}
(%i2) 
\end{verbatim}
}
\end{minipage}
\begin{minipage}[t]{\textwidth}{\color{blue}
\begin{verbatim}
cifrado_afin("notenemosreservadeagua",1,3);
\end{verbatim}
}
\end{minipage}
\definecolor{labelcolor}{RGB}{100,0,0}
\begin{math}\displaystyle
\parbox{8ex}{\color{labelcolor}(\%o2) }
QRWHQHPRVUHVHUYDGHDJXD
\end{math}
\end{maxcode}

Si queremos desencriptar un mensaje que ha sido encriptado utilizando esta t\'ecnica pero
no conocemos los par\'ametros que se han usado, podemos 
partir de la observaci\'on de que la frecuencia relativa de las letras en el mensaje original
se mantiene en el mensaje encriptado. En otras palabras,
la letra que m\'as aparece en el mensaje encriptado ser\'a la letra transformada mediante
$\Phi(x)$ de la letra m\'as frecuente en el mensaje original. 
A continuaci\'on, consultamos una 
tabla de las letras m\'as frecuentes en el idioma espa\~nol, como la siguiente\footnote{Obtenida de Wikipedia: \freqletras}: 
\begin{center}
\begin{tabular}{c|c}
Letra & Frecuencia \\ 
\hline 
e & $13{.}68$\% \\ 
a & $12{.}53$\% \\ 
o & $8{.}68$\% \\ 
s & $7{.}98$\% \\ 
r & $6{.}87$\% \\ 
n & $6{.}71$\% \\ 
\hline 
\end{tabular} 
\end{center}
e identificamos de qu\'e letra se trata. 

Por ejemplo,
en el mensaje cifrado \emph{QRWHQHPRVUHVHUYDGHDJXD} la letra m\'as frecuente es la $H$,
que aparece $5$ veces, y le sigue la $D$, que aparece tres veces. Este recuento se puede hacer directamente en Maxima con el comando \texttt{discrete\_freq} del paquete
\texttt{descriptive}:

\begin{maxcode}
\noindent
\begin{minipage}[t]{4em}\color{red}\bf
(\% i3)
\end{minipage}
\begin{minipage}[t]{\textwidth}\color{blue}
load(descriptive);
\end{minipage}

\noindent
\begin{minipage}[t]{4em}\color{red}\bf
(\% i4)
\end{minipage}
\begin{minipage}[t]{\textwidth}\color{blue}
discrete\_freq(charlist("QRWHQHPRVUHVHUYDGHDJXD"));\\
\end{minipage}

\definecolor{labelcolor}{RGB}{100,0,0}
\begin{math}\displaystyle
\parbox{8ex}{\color{labelcolor}(\%o4) } 
[[D,G,H,J,P,Q,R,U,V,W,X,Y],[3,1,5,1,1,2,2,2,2,1,1,1]]\mbox{}
\end{math}
\end{maxcode}

Esto nos hace sospechar que la \emph{H} es la imagen de la letra m\'as frecuente en espa\~nol,
que es la $e$, mientras que $D$ es la imagen de la segunda letra m\'as frecuente, la $a$. Es decir, sospechamos que $\Phi (e)=H$ y $\Phi (a)=D$ o bien, en t\'erminos num\'ericos,
$$
\systeme{a 4 +b=7,a 0 +b=3\,.}
$$
La soluci\'on a este sistema es, como cab\'ia esperar, $a=1$, $b=3$. 

La t\'ecnica funciona con este cifrado porque en el criptotexto se detectan claramente dos letras con mayor frecuencia relativa de aparici\'on, que son las que ponemos en correspondencia con las dos letras m\'as frecuentes en el idioma espa\~nol (a mucha distancia
de la tercera). Podemos ver esto en Maxima gr\'aficamente con el siguiente comando (donde
los argumentos de \texttt{bars} son listas de tres elementos conteniendo la posici\'on
de la barra sobre el eje $x$, su altura y su anchura, respectivamente):
\begin{maxcode}

\noindent
\begin{minipage}[t]{4em}\color{red}\bf
(\% i5)
\end{minipage}
\begin{minipage}[t]{\textwidth}\color{blue}
wxdraw2d(fill\_color = red,fill\_density = 0.6,

  bars([1,3,0.2],[2,1,0.2],[3,5,0.2],[4,1,0.2],[5,1,0.2],[6,2,0.2],

       [7,2,0.2],[8,2,0.2],[9,2,0.2],[10,1,0.2],[11,1,0.2],[12,1,0.2]),

  xaxis = true,

  xtics = {[\texttt{"}D\texttt{"},1],[\texttt{"}G\texttt{"},2],[\texttt{"}H\texttt{"},3],[\texttt{"}J\texttt{"},4],[\texttt{"}P\texttt{"},5],[\texttt{"}Q\texttt{"},6],

           [\texttt{"}R\texttt{"},7],[\texttt{"}U\texttt{"},8],[\texttt{"}V\texttt{"},9],[\texttt{"}W\texttt{"},10],[\texttt{"}X\texttt{"},11],[\texttt{"}Y\texttt{"},12]}

)\$

\end{minipage}
\[\displaystyle \tag{\% t5} 
\includegraphics[scale=0.5]{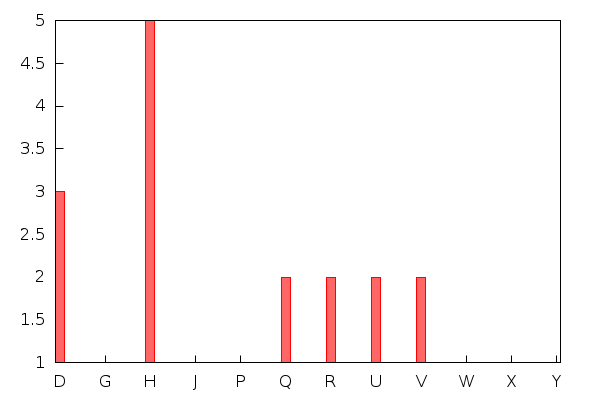} 
\]
\end{maxcode}

Como actividad, se pide realizar este an\'alisis encriptando
el mismo mensaje pero esta vez
usando el cifrado de Hill con la matriz $K$ de la subsecci\'on \ref{hill}.
Representar las frecuencias de las letras y comparar el resultado con el
obtenido en el cifrado af\'in.

%

\item En esta actividad estudiaremos las llamadas
\emph{transformaciones lineales}\index{Transformaci\'on!lineal} en el plano. 
Dada una matriz $A = (a_{ij})$ de
tama\~no $2 \times 2$, podemos definir una funci\'on $T_A:\R^2 \to \R^2$ como
\[
T_A(x,y) = \begin{pmatrix}
                   a_{11} & a_{12} \\
                   a_{21} & a_{22} 
                  \end{pmatrix}
          \begin{pmatrix}
                   x \\
                   y 
                  \end{pmatrix}
\]
Utilicemos GeoGebra para ver que sucede al aplicar la transformaci\'on
$T_A$ a los puntos del cuadrado unitario
\[
C = \{(x,y) \in \R^2:0 \le x,y \le 1\}\,.
\]
\begin{center}
\scalebox{0.35}{\includegraphics{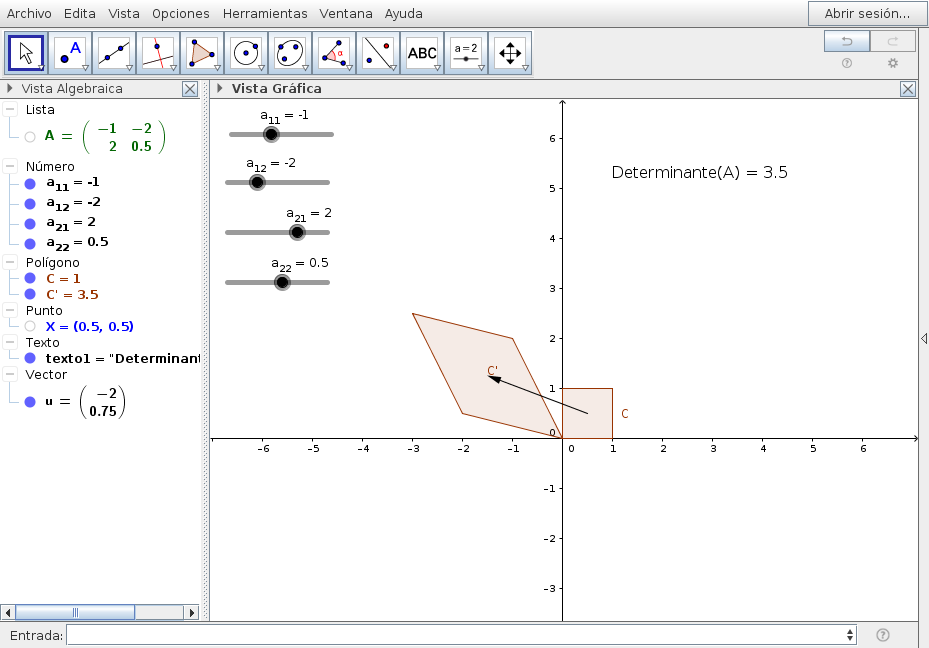}}
\end{center}
Los pasos a seguir son los siguientes:
\begin{enumerate}
\item Con la herramienta \texttt{Deslizador} creamos los controles
\texttt{a\_{11},a\_{12},a\_{21},a\_{22}},
que se usar\'an como elementos de la matriz $A$.
\item Creamos la matriz $A$, con la orden \texttt{A = {{a\_{11}, a\_{12}},{a\_{21}, a\_{22}}}}.
\item Creamos el centroide del cuadrado unitario, con la orden
\texttt{X=(1/2,1/2)}.
\item Creamos el vector que va del punto $X$ a su imagen $AX$, con la orden \texttt{Vector[X, A*X]}.
\item Creamos el cuadrado unitario $C$, con la orden \texttt{C = Poligono[(0,0),(0,1),4]}.
\item Aplicamos la transformaci\'on $T_A$ a los puntos del cuadrado unitario,
      mediante la orden \texttt{AplicaMatriz[A,C]}.
\item Calculamos el determinante de la matriz $A$, con la orden \texttt{Determinante[A]}.
\end{enumerate}
Modificando los coeficientes de la matriz $A$ con los deslizadores
creados se ve que el \'area del paralelogramo resultante
\[
C' = \left\{T_A(x,y):(x,y) \in C\right\}
\]
siempre es igual al valor absoluto del determinante de $A$.

Como actividad de c\'omputo se pide repetir el experimento anterior,
aplicando la transformaci\'on $T_A$ al circulo de \'area unitaria 
\[
{\cal C} = \left\{(x,y) \in \R^2:x^2 + y^2 \le \sqrt{\frac{1}{\pi}}\right\}\,.
\]
?`Puede el lector dar una explicaci\'on del resultado observado?

\end{enumerate}
%

\chapter*{Ap\'endice: Maxima y GeoGebra}\addcontentsline{toc}{chapter}{Ap\'endice: Maxima y GeoGebra}
\vspace*{2.5cm}

Ambos programas se pueden instalar, en los sistemas operativos que los soportan, como cualquier otro, 
por lo que no nos detendremos en el proceso de instalaci\'on y supondremos que el usuario lector los 
tiene disponibles en su computadora. Se encuentran en \url{http://maxima.sourceforge.net/} y 
\url{https://www.geogebra.org/}, respectivamente.

Maxima es un sistema de \'algebra computacional (o CAS\index{CAS} por sus siglas en ingl\'es, \textbf{C}omputer \textbf{A}lgebra \textbf{S}ystem) que originalmente se ejecuta en una consola 
(\emph{shell}, en ingl\'es), esto es, en l\'inea de comandos. Sin embargo existen varias interfaces gr\'aficas disponibles, como Imaxima, TeXmacs, Kayali o el popular sistema SAGE. Nosotros usamos wxMaxima\index{Maxima!wxMaxima}\index{wxMaxima} (disponible en \url{http://andrejv.github.io/wxmaxima/}), que tiene una apariencia que recuerda vagamente a la de otros CAS como Maple\textregistered\  o Mathematica\textregistered\ , de modo que el usuario pueda
trasladar f\'acilmente su experiencia con cualquiera de ellos a los dem\'as.

Al arrancar una sesi\'on en wxMaxima, por defecto nos encontraremos con una hoja de trabajo en blanco
(y, quiz\'as, con una ventana de historial a la derecha). Picando con el rat\'on en cualquier
punto de la hoja, podemos comenzar a introducir \'ordenes, como \texttt{2+2}. Para ejecutarlas, 
por defecto hay que pulsar simult\'aneamente las teclas \texttt{Shift + Enter}, que se encuentran 
una sobre la otra en la parte derecha del teclado (este comportamiento est\'a tomado de 
Mathematica\textregistered\ ). 
Una vez hecho esto, veremos que Maxima a\~nade autom\'aticamente un punto y coma al final del comando 
que hab\'iamos escrito (pasando a \texttt{2+2;}) y devuelve el resultado, \texttt{4}.
\begin{center}
\includegraphics[scale=0.375]{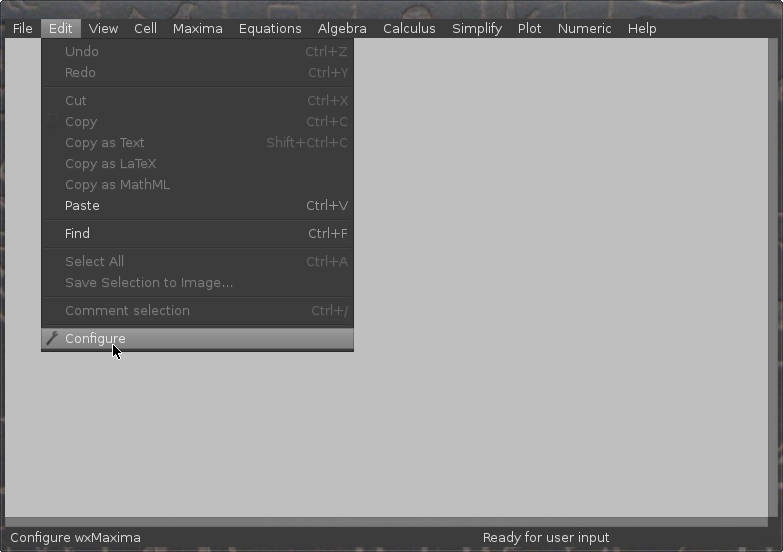} 
\end{center}
La interfaz wxMaxima se puede configurar a gusto del usuario de muchas maneras. Para ello, se accede
al men\'u \texttt{Editar -> Configurar}, como en la figura precedente. Una vez all\'i encontraremos
varias opciones para modificar la apariencia de la hoja de trabajo, o del comportamiento de la interfaz (n\'otese que las capturas de pantalla que acompa\~nan a este ap\'endice corresponden a una versi\'on
de wxMaxima configurada de forma muy distinta a la que se presenta por defecto, se trata de una configuraci\'on en ingl\'es y dise\~nada para pasar mucho tiempo frente a la pantalla). Por ejemplo, es habitual desear que la ejecuci\'on de los comandos se realice sin m\'as que
pulsar \texttt{Enter}, y hay una opci\'on disponible a tal efecto.
\begin{center}
\includegraphics[scale=0.375]{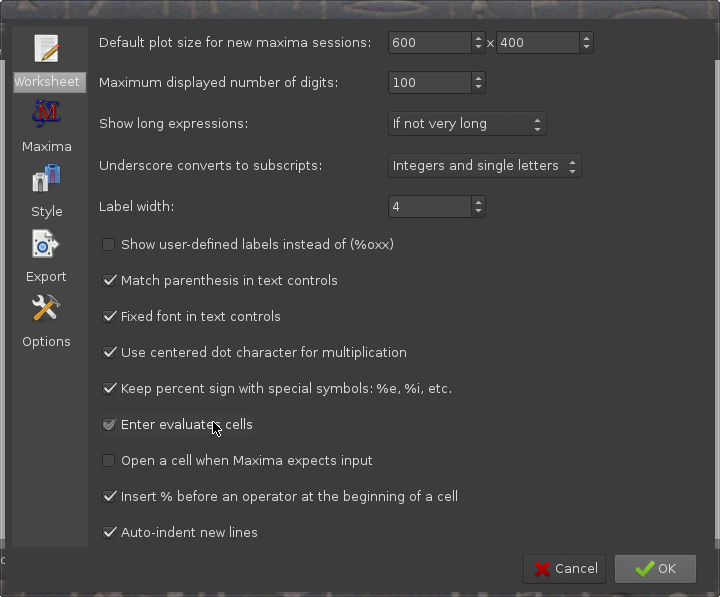} 
\end{center}
El entorno de una hoja de trabajo en wxMaxima se organiza alrededor del concepto de \emph{celda}.\index{Maxima!wxMaxima!celda}\index{wxMaxima!celda} Una celda no es m\'as que el bloque formado por una instrucci\'on
y el resultado de su ejecuci\'on. Las celdas se pueden seleccionar, borrar, plegar y desplegar
usando las
barras verticales que aparecen a su izquierda. Existen tambi\'en celdas de texto, muy \'utiles
para incluir comentarios en la hoja de trabajo\footnote{Para los \emph{connoisseurs}: es incluso posible utilizar el formato de \LaTeX\,, para ello basta con comenzar el texto de la celda con los caracteres \texttt{TeX:} en una sola l\'inea.}. Se accede a los diferentes tipos de celda a trav\'es del
men\'u.
\begin{center}
\includegraphics[scale=0.4]{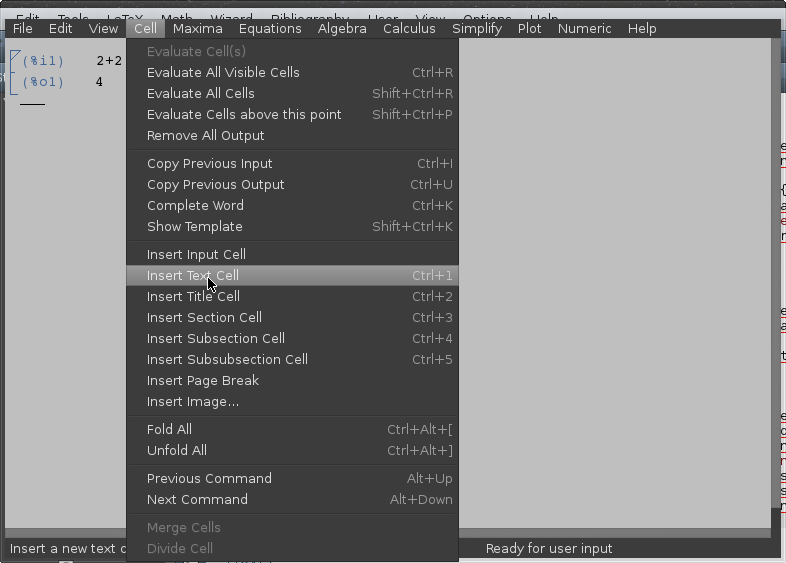} 
\end{center}
Como un ejemplo sencillo, en la siguiente figura se ha creado un documento con una celda de t\'itulo,
una secci\'on, un p\'arrafo de texto y un par de celdas con c\'alculos.
\begin{center}
\includegraphics[scale=0.425]{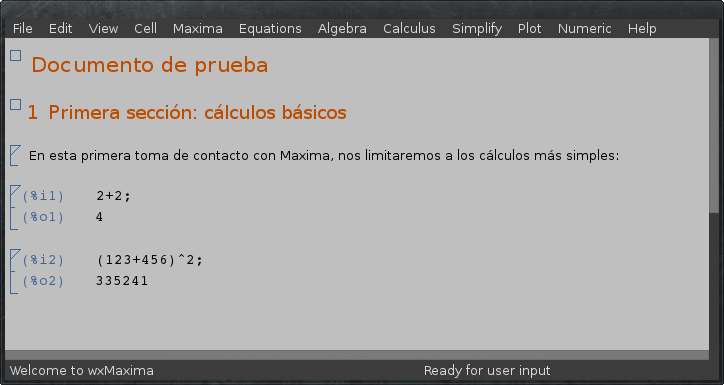} 
\end{center}
Una vez hayamos acabado nuestro trabajo, podemos guardar el documento como es habitual, yendo al
men\'u y siguiendo la ruta \texttt{Archivo -> Guardar como}. La interfaz wxMaxima ofrece dos posibilidades para
guardar el documento, una con extensi\'on \texttt{.wxm} y otra con extensi\'on \texttt{.xwxm}. La
diferencia entre ellos simplemente es que el segundo guarda tambi\'en los gr\'aficos y otros objetos
extendidos que pudiera haber en el documento de trabajo (esta opci\'on, por tanto, es la m\'as
recomendable para los no expertos). Otro aspecto interesante es la posibilidad de \emph{exportar}
nuestro documento como una p\'agina web, en formato \texttt{HTML}, lo cual resulta extremadamente
\'util a la hora de elaborar materiales destinados a uso en clases que pueden compartirse por
intranet o internet\footnote{Tambi\'en existe la posibilidad de exportar el documento en formato
\TeX, que posteriormente puede compilarse para proporcionar una versi\'on en \texttt{pdf}.}.

Por \'ultimo, mencionaremos que Maxima dispone de un sistema de ayuda muy completo, aunque las
explicaciones frecuentemente tienen un car\'acter marcadamente t\'ecnico (que puede aprovecharse
para aprender muchos detalles sobre teor\'ia de c\'omputo en general).
\begin{center}
\includegraphics[scale=0.33]{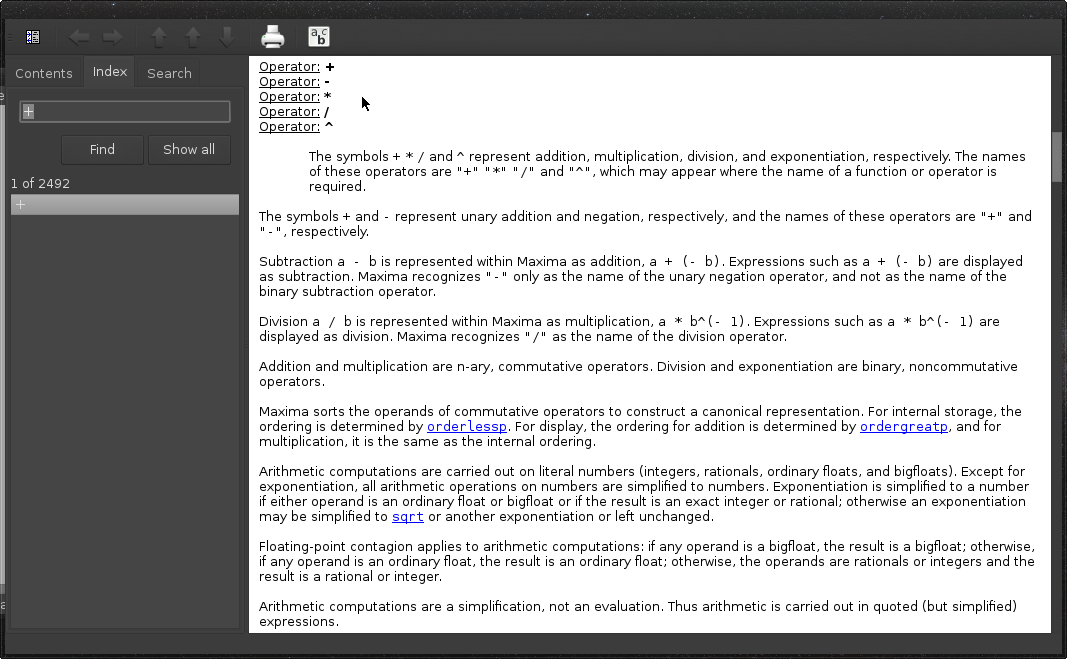} 
\end{center}

Por su parte, GeoGebra es un entorno de geometr\'ia din\'amica (o DGE por sus siglas en
ingl\'es, \textbf{D}ynamic \textbf{G}eometry \textbf{E}nvironment), similar a otros como 
Cabri\textregistered\ , 
Cinderella\textregistered\  o The Geometer's Sketchpad\textregistered\ , con el a\~nadido de tratarse, 
como Maxima, de software libre.
\begin{figure}[H]
\begin{center}
\scalebox{0.3}{\includegraphics{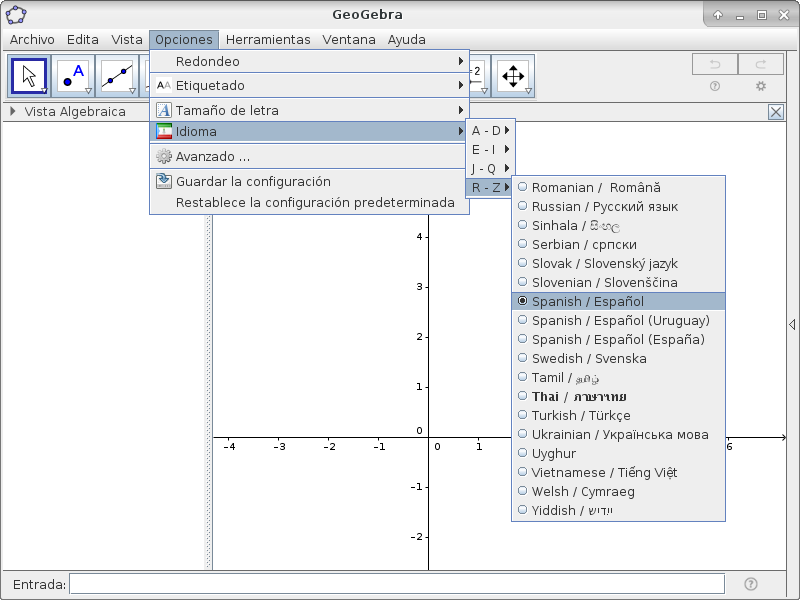}}
\end{center}
\end{figure}
Suponiendo que se ha configurado el idioma como \texttt{Espa\~nol} (como se muestra
arriba), al arrancar una sesi\'on de GeoGebra se presenta en el margen derecho de
la ventana principal una ventana desplegable para seleccionar una apariencia
predeterminada\footnote{Las \texttt{Apariencias} pueden variar con la versi\'on
de GeoGebra usada, que en estas notas es la 5.0.}: \texttt{\'Algebra},
\texttt{Geometr\'ia}, \texttt{Hoja de C\'alculo}, \texttt{CAS},
\texttt{Gr\'aficos 3D} \'o \texttt{Probabilidad}.
Cada vista se compone de una barra de cajones de herramientas, una o m\'as
vistas y una o m\'as ventanas de entrada, predeterminadas seg\'un la vista elegida. 
Por ejemplo, la apariencia por defecto, la de \texttt{Algebra}, consta de una
barra con 12 cajones de herramientas, la \texttt{Vista Algebraica} y la
\texttt{Vista Gr\'afica}, y una ventana de \texttt{Entrada}. Adem\'as de estos
elementos siempre est\'a presente una barra de men\'us desplegables:
\texttt{Archivo}, \texttt{Edita}, \texttt{Vista}, \texttt{Opciones},
\texttt{Herramientas}, \texttt{Ventana} y \texttt{Ayuda}.   
\begin{figure}[H]
\begin{center}
\scalebox{0.3}{\includegraphics{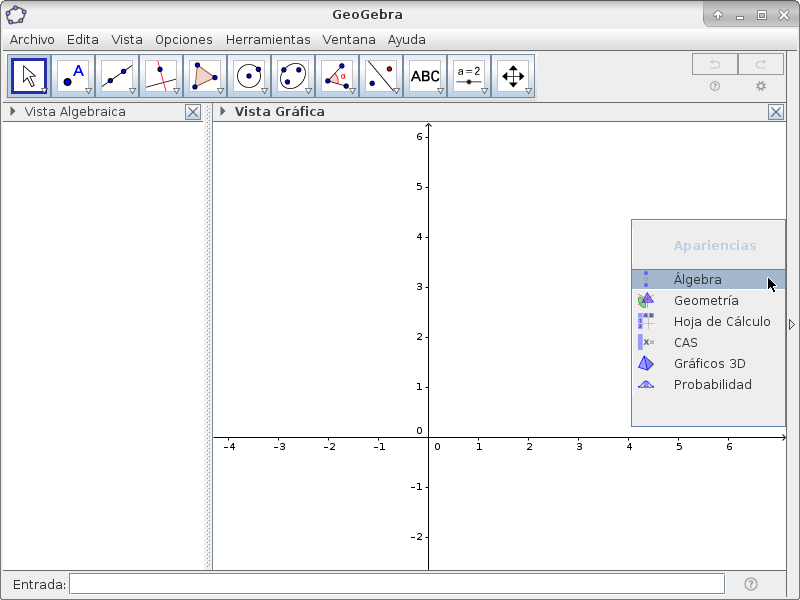}}
\end{center}
\end{figure}

La forma m\'as com\'un de crear objetos en GeoGebra es me\-diante acciones
apropiadas realizadas con el rat\'on sobre la \texttt{Vista Gr\'afica}.
Por ejemplo, para crear un punto en el plano, seleccionamos la herramienta
\texttt{Punto} del segundo caj\'on de la barra de herramientas y picamos
con el bot\'on izquierdo del rat\'on en el lugar de la
\texttt{Vista Gr\'afica} en el que queremos crear el punto. Notemos que en
la \texttt{Vista Algebraica} aparece la representaci\'on en coordenadas
cartesianas del punto creado $A(3,2)$.
\begin{figure}[h]
\begin{center}
\scalebox{0.36}{\includegraphics{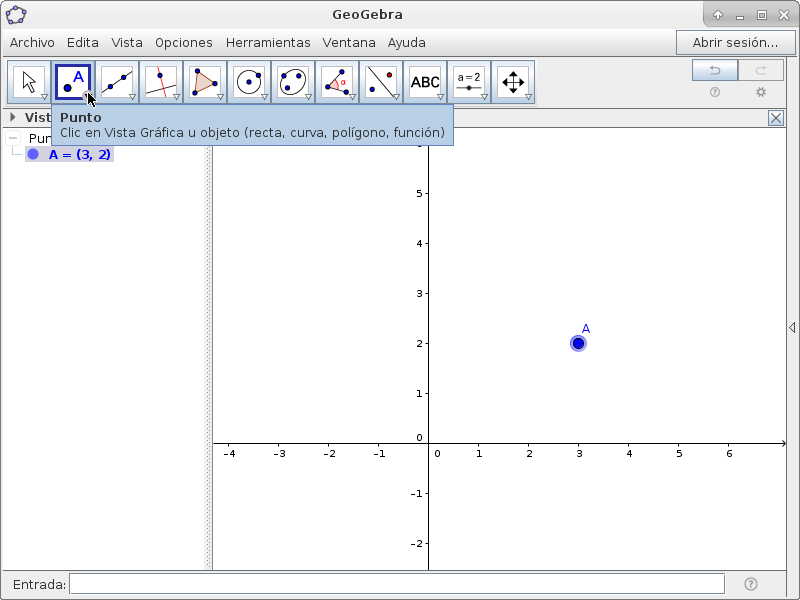}}
\end{center}
\end{figure}

Tambi\'en es posible crear objetos tecleando la orden apropiada en la
ventana de \texttt{Entrada}, situada en la barra inferior. 
As\'i se ha creado el punto $B(-1,1)$ en la figura siguiente.
\begin{figure}[h]
\begin{center}
\scalebox{0.36}{\includegraphics{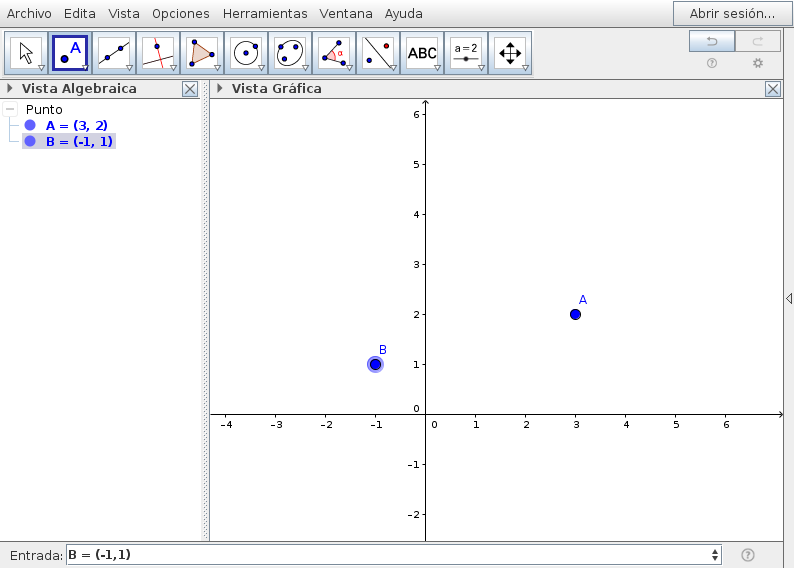}}
\end{center}
\end{figure}

Nuestra \'ultima construcci\'on es el punto medio $C$ entre $A$ y $B$, el cual
determinamos usando la herramienta \texttt{Medio o Centro}.
\begin{figure}[h]
\begin{center}
\scalebox{0.33}{\includegraphics{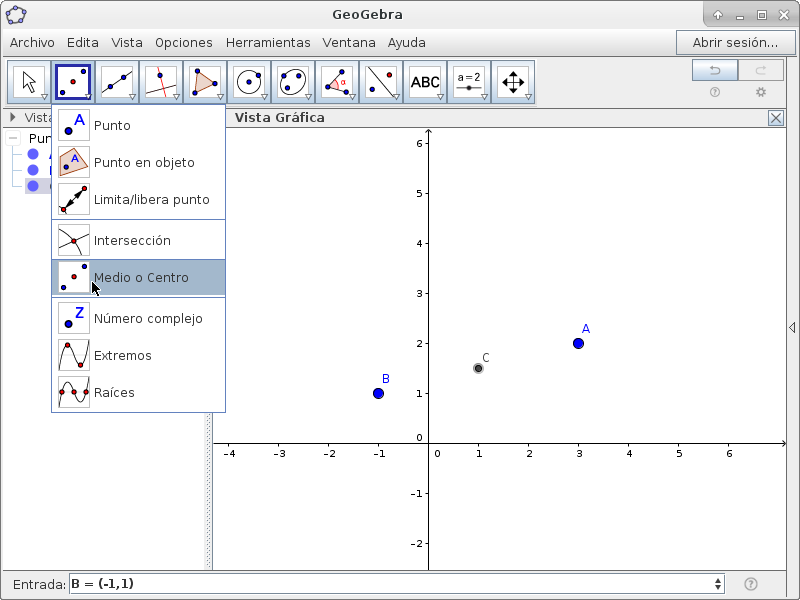}}
\end{center}
\end{figure}

Al picar con el bot\'on derecho del rat\'on sobre un objeto se despliega
un men\'u contextual que permite realizar diversas acciones; por ejemplo,
cambiar el nombre del objeto en cuesti\'on, como se muestra abajo.
\begin{figure}[h]
\begin{center}
\scalebox{0.36}{\includegraphics{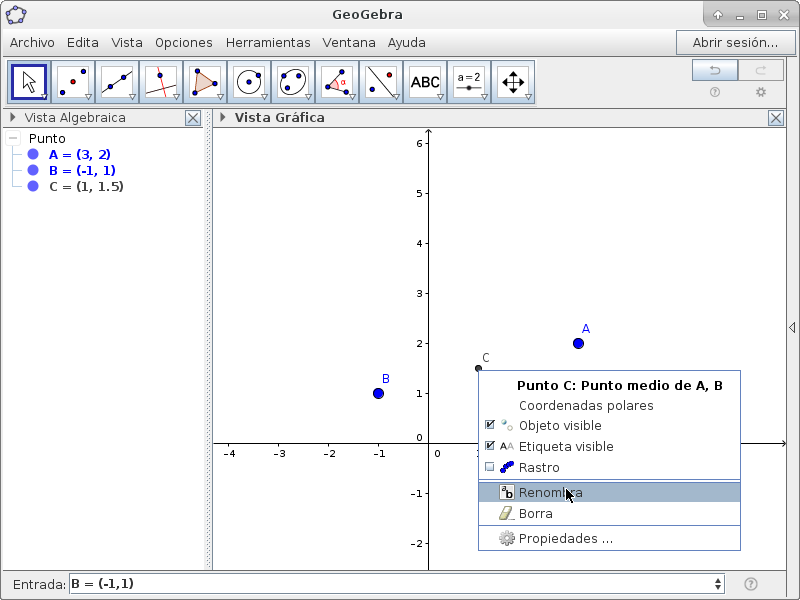}}
\end{center}
\end{figure}

Las construcciones se actualizan de forma din\'amica al modificar la
posici\'on de los objetos que las componen utilizando, por ejemplo, la
herramienta \texttt{Elige y Mueve} del primer caj\'on de herramientas, o
redifini\'endolos picando dos veces sobre ellos.  
Finalmente, el submen\'u \texttt{Exporta} del men\'u \texttt{Archivo} permite
crear una hoja din\'amica en formato \texttt{HTML} con la construcci\'on hecha, o bien
guardar la \texttt{Vista Gr\'afica} como un archivo en formato \texttt{png} o \texttt{eps},
entre otras posibilidades.
\begin{figure}[H]
\begin{center}
\scalebox{0.3}{\includegraphics{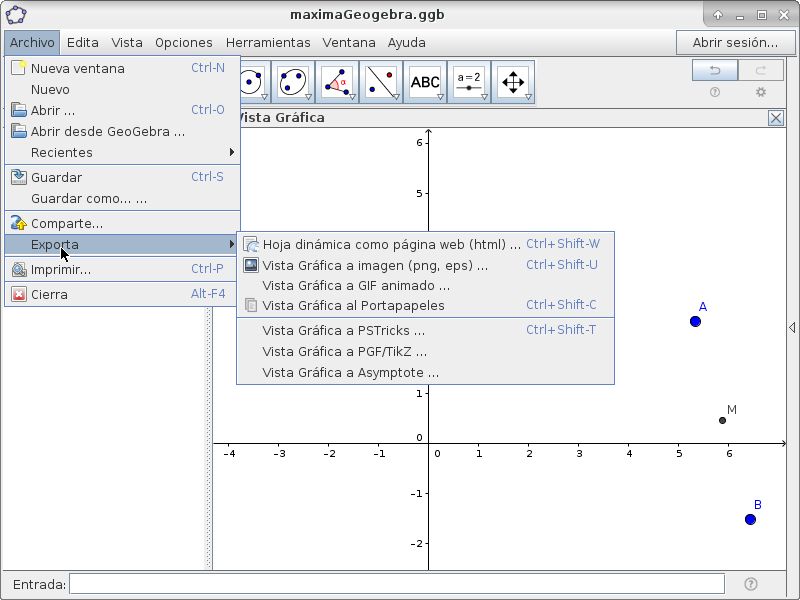}}
\end{center}
\end{figure}

Para finalizar, mencionaremos la existencia del repositorio anteriormente conocido como
`GeoGebraTube' y actualmente disponible en \url{https://www.geogebra.org/materials/}. Se
trata de una p\'agina web desde la que se pueden descargar miles de applets de GeoGebra
y a la que cualquiera puede subir los suyos. Aunque la calidad de los contenidos es muy
heterog\'enea, es muy recomendable echar un vistazo de vez en cuando y desde aqu\'i
animamos al lector interesado a contribuir con sus propias creaciones.


\newpage
\cleardoublepage\phantomsection
\addcontentsline{toc}{chapter}{\indexname}
\titlespacing{\chapter}{2.5cm}{-\baselineskip}{2.5cm}
\printindex

\newpage
\cleardoublepage\phantomsection
\addcontentsline{toc}{chapter}{Procedencia de las im\'agenes}
\thispagestyle{empty}
\vspace*{1cm}
\begin{center}
\textbf{Procedencia de las im\'agenes}
\end{center}

\bigskip

\begin{enumerate}\sloppy
\item P\'agina 7 (el Hombre de Vitrubio): Luc Viatour/www.Lucnix.be, en \url{https://en.wikipedia.org/wiki/Vitruvian_Man#/media/File:Da_Vinci_Vitruve_Luc_Viatour.jpg} (Imagen de dominio p\'ublico).
\item P\'agina 8 (Sof\'ia Loren): Mauricio Navarrete Contreras, en \url{https://www.flickr.com/photos/127592731@N07/15788618971}
(Imagen bajo licencia Creative Commons Attribution 2.0 Generic). Modificada por los
autores para superimponer la m\'ascara de Marquardt. Tambi\'en disponible en
\url{http://idbelleza.com/2016/05/26/iconos-de-la-belleza-sophia-loren/}.
\item P\'agina 17 (reloj): wilhei, en  \url{https://pixabay.com/en/clock-time-pointer-time-of-watches-705672/} (Imagen de dominio p\'ublico).
\item P\'agina 20 (libros): Fotograf\'ia tomada por uno de los autores (JAV).
\item P\'agina 22 (ej\'ercito de terracota): Ingo.staudacher, en \url{https://commons.wikimedia.org/wiki/File:Terracotta_army_5256.jpg} (Imagen bajo licencia Creative Commons Attribution-Share Alike 3.0 Unported).
\item P\'agina 31 (cardenales espiando): Henri Adolphe Laissement, en
\url{https://en.wikipedia.org/wiki/Eavesdropping#/media/File:
Henri_Adolphe_Laissement_Kardin%C3%A4le_im_Vorzimmer_1895.jpg} (Imagen de dominio p\'ublico).
\item P\'agina 33 (sat\'elite): NASA, en 
\url{https://commons.wikimedia.org/wiki/File:GPS_Satellite_NASA_art-iif.jpg} (Imagen de dominio p\'ublico). 
\item P\'agina 35 (camellos): amira\_a, en \url{https://commons.wikimedia.org/wiki/File:Camels.convoy%288008893480%29.jpg} (Imagen bajo licencia Creative Commons Attribution 2.0 Generic).
\item P\'agina 36 (nueces): Melchoir, en 
\url{https://commons.wikimedia.org/wiki/File:Mixed_nuts_small_white1.jpg} (Imagen bajo licencia Creative Commons Attribution-Share Alike 3.0 Unported).
\item P\'agina 41 (domin\'o): \url{http://www.publicdomainpictures.net/view-image.php?image=34179&picture=la-caida-de-domino-principio} (Imagen de dominio p\'ublico).
\item P\'agina 45 (loter\'ia): Hermann, en \url{https://pixabay.com/en/lotto-lottery-ticket-bill-profit-484782/} (Imagen de dominio p\'ublico).
\item P\'agina 51 (buffet): Spencer195, en 
\url{https://commons.wikimedia.org/wiki/File:Chinese_buffet2.jpg} (Imagen bajo licencia Creative Commons Attribution-Share Alike 3.0 Unported).\clearpage\thispagestyle{empty}   
\item P\'agina 60 (salto): Guillaume Baviere, en \url{https://commons.wikimedia.org/wiki/File:Eloyse_Lesueur_G%C3%B6teborg_2013.jpg#/media/File:Eloyse_Lesueur_G%C3%B6teborg_2013.jpg} (Imagen bajo licencia Creative Commons Attribution 2.0 Generic).
\item P\'agina 69 (tablero ajedrez): McGeddon, en 
\url{https://commons.wikimedia.org/wiki/File:Wheat_and_chessboard_problem.jpg} (Imagen bajo
licencia Creative Commons Attribution-Share Alike 4.0 International).
\item P\'agina 78 (bacteria E.Coli): Agricultural Research Service (USA), en
\url{https://es.wikipedia.org/wiki/Archivo:E_coli_at_10000x,_original.jpg} (Imagen de dominio p\'ublico).
\item P\'agina 80 (sism\'ografo): Yamaguchi, en 
\url{https://commons.wikimedia.org/wiki/File:Kinemetrics_seismograph.jpg}
(Imagen bajo licencia Creative Commons Attribution-Share Alike 3.0 Unported).
\item P\'agina 81 (conejos): Worm That Turned, en
\url{https://commons.wikimedia.org/wiki/File:Wild_Rabbits_at_Edinburgh_Zoo.jpg}
(Imagen bajo licencia Creative Commons Attribution-Share Alike 3.0 Unported).
\item P\'agina 111 (antena parab\'olica): Richard Bartz, en
\url{https://en.wikipedia.org/wiki/Parabolic_antenna#/media/File:Erdfunkstelle_Raisting_2.jpg}
(Imagen bajo licencia Creative Commons Attribution-Share Alike 2.5 Generic).
\item P\'agina 112 (chimenea): BluebearsLani, en
\url{https://www.flickr.com/photos/7292946@N08/15459161906} (Imagen bajo licencia Creative Commons Attribution-NoDerivs 2.0 Generic).
\item P\'agina 113 (techo): rockford\_chimney, en 
\url{http://www.instructables.com/file/FBWUB5EIH7Z9NPM/?size=MEDIUM}
(Imagen bajo licencia Creative Commons Attribution-NonCommercial-ShareAlike 2.5 Generic).
\item P\'agina 114 (dibujo elipse): Dino, en
\url{https://en.wikipedia.org/wiki/File:Drawing_an_ellipse_via_two_tacks_a_loop_and_a_pen.jpg}
(Imagen bajo licencia Creative Commons Attribution-Share Alike 3.0 Unported).
\item P\'agina 130 (enchufe): Tomasz Sienicki, en
\url{https://commons.wikimedia.org/wiki/File:Gniazdo_el_z_przelacznikiem_%28ubt%29.JPG}
(Imagen bajo licencia Creative Commons Attribution 3.0 Unported). 
\item P\'agina 144 (frenazo): www.hacienda-la-colora.com, en \url{https://www.flickr.com/photos/xavier33300/15052804730}
(Imagen bajo licencia Creative Commons Attribution 2.0 Generic).
\item P\'agina 160 (manzana): moreharmony, en \url{https://pixabay.com/en/appetite-apple-calories-catering-1239057/}
(Imagen bajo licencia Creative Commons CC0 1.0 Universal).
\end{enumerate}


\end{document}